\documentclass[11pt]{amsart}
\usepackage{amsmath}
\usepackage{amscd}
\usepackage{pb-diagram}
\usepackage{amssymb}
\usepackage{graphicx}
\usepackage{subcaption}
\usepackage[table]{xcolor}
\usepackage[bookmarks=true,bookmarksnumbered=true]{hyperref}

\usepackage{amssymb}

\def\refeq#1{\if\workingver y(\ref{#1})-[[#1]]\else(\ref{#1})\fi}
\def\refth#1{\if\workingver y\ref{#1}-[[#1]]\else\ref{#1}\fi}
\def\mylabel#1{\if\workingver y\label{#1}{\bf\ \ [[#1]]\ \ }\else\label{#1}\fi}
\def\mybibitem#1{\if\workingver y\bibitem{#1}{\bf\ \ [[#1]]\ \
}\else\bibitem{#1}\fi}

\def\articletheorems{
\newtheorem{thm}{Theorem}[section]
\newtheorem{lem}[thm]{Lemma}

\newtheorem{defn}[thm]{Definition}
\newtheorem{cor}[thm]{Corollary}
\newtheorem{prop}[thm]{Proposition}
\newtheorem{propdef}[thm]{Proposition and Definition}
\newtheorem{thmdef}[thm]{Theorem and Definition}
\newtheorem{ex}[thm]{Example}
\newtheorem{algo}{Algorithm}[section] 
\newtheorem{alg}[thm]{Algorithm}  

}

\newcommand{\mvmap}{\multimap}
\newcommand{\mto}{\multimap}
\newcommand{\pto}{\nrightarrow}

\renewcommand{\emptyset}{\varnothing}
\renewcommand{\rho}{\varrho}
\renewcommand{\phi}{\varphi}
\renewcommand{\epsilon}{\varepsilon}

\def\cA{\text{$\mathcal A$}}

\def\cC{\text{$\mathcal C$}}
\def\cD{\text{$\mathcal D$}}
\def\cE{\text{$\mathcal E$}}
\def\cF{\text{$\mathcal F$}}

\def\cL{\text{$\mathcal L$}}
\def\cM{\text{$\mathcal M$}}
\def\cN{\text{$\mathcal N$}}

\def\cP{\text{$\mathcal P$}}

\def\cT{\text{$\mathcal T$}}
\def\cU{\text{$\mathcal U$}}
\def\cV{\text{$\mathcal V$}}
\def\cW{\text{$\mathcal W$}}
\def\cX{\text{$\mathcal X$}}

\def\bA{\text{$\mathbf A$}}
\def\bB{\text{$\mathbf B$}}
\def\bC{\text{$\mathbf C$}}
\def\bD{\text{$\mathbf D$}}
\def\bE{\text{$\mathbf E$}}
\def\bF{\text{$\mathbf F$}}
\def\bG{\text{$\mathbf G$}}

\newcommand{\id}{\operatorname{id}}
\newcommand{\cl}{\operatorname{cl}}

\newcommand{\inte}{\operatorname{int}}

\newcommand{\dom}{\operatorname{dom}}

\newcommand{\card}{\operatorname{card}}

\newcommand{\im}{\operatorname{im}}

\newcommand{\conv}{\protect\mbox{\rm conv\,}}
\renewcommand{\emptyset}{\varnothing}

\newcommand{\Con}{\operatorname{Con}}
\newcommand{\Inv}{\operatorname{Inv}}
\newcommand{\Fix}{\operatorname{Fix}}

\def\proof{{\bf Proof:\ }}

\def\begeq#1{\begin{equation}\mylabel{#1}}
\def\endeq{\end{equation}}

\def\mathobj#1{\mbox{$#1$}}

\def\NN{\mathobj{\mathbb{N}}}
\def\PP{\mathobj{\mathbb{P}}}

\def\RR{\mathobj{\mathbb{R}}}

\def\ZZ{\mathobj{\mathbb{Z}}}

\def\rep#1{\mbox{$\langle#1\rangle$}}

\def\scalprod#1{\langle #1 \rangle}
\def\spn#1{\langle #1 \rangle}

\def\implies{\;\Rightarrow\;}
\def\iff{\;\Leftrightarrow\;}

\def\setof#1{\mbox{$\{\,#1\,\}$}}

\newcommand{\bdy}{{\bf \partial}}

\def\0#1{\hbox{\kern25pt}$ #1 $\\}
\def\1#1{\hbox{\kern40pt}$ #1 $\\}
\def\2#1{\hbox{\kern55pt}$ #1 $\\}
\def\3#1{\hbox{\kern70pt}$ #1 $\\}

\newcounter{li}

\def\begalg#1{\begin{algo}\mylabel{#1}\normalshape:\small\baselineskip 10pt\\}
\def\endalg{\end{algo}}

\def\Figures(include=#1,cat=#2){
  \renewcommand{\textfraction}{.20}
  \renewcommand{\topfraction}{.80}
  \renewcommand{\bottomfraction}{.80}
  \renewcommand{\floatpagefraction}{.80}
  \newcount\figcount
  \figcount=0
  \let\includefigures=#1
  \def\figcat{#2}
}

\def\FigureFromFile[#1][#2](#3)#4
{
  \begin{figure}[htbp]
     \global\advance\figcount by 1
     \if\includefigures y\special{anisoscale #1.wmf, \the\hsize #2}\fi
     \vspace{#2}
     \caption{#4}
     \mylabel{#3}
   \end{figure}
}

\def\FigureFromFileTwoD[#1][#2,#3](#4)#5
{
  \begin{figure}[htbp]
     \global\advance\figcount by 1
     \if\includefigures y\special{anisoscale #1.wmf, #2 #3}\fi
     \vspace{#2}
     \caption{#5}
     \mylabel{#4}
   \end{figure}
}

\def\FigureF<#1>[#2](#3)#4
{
  \begin{figure}[htbp]
     \global\advance\figcount by 1
     \if\includefigures y\special{anisoscale \figcat/fig\number\figcount.wmf,
       \the\hsize #2}
     \fi
     \if\includefigures p
       \leavevmode
       \epsfxsize=\hsize
       \epsffile{#1}
     \fi
     \if\includefigures y
          \vspace{#2}
     \fi
     \caption{#4}
     \mylabel{#3}
   \end{figure}
}

\def\Figure[#1](#2)#3
{
  \begin{figure}[htbp]
     \global\advance\figcount by 1
     \if\includefigures y\special{anisoscale \figcat/fig\number\figcount.wmf,
       \the\hsize #1}
     \fi
     \if\includefigures p
       \leavevmode
       \epsfxsize=\hsize
       \epsffile{fig\number\figcount.eps}
     \fi
     \if\includefigures y
          \vspace{#1}
     \fi
     \caption{#3}
     \mylabel{#2}
   \end{figure}
}

\let\visiblecomments y 

\graphicspath{
  {./images/}
}

\newcommand{\currentDate}{\today}
\def\PP{\mathobj{\mathbb{P}}}

\newcommand{\bdop}{\partial}

\def\Down{\mbox{$\operatorname{Down}$}}

\def\sing{\mbox{$\operatorname{sing}$}}
\def\Clsd{\mbox{$\operatorname{Clsd}$}}
\def\uim{\mbox{$\operatorname{uim}$}}
\def\AP{\mbox{$\operatorname{AP}$}}

\newcommand{\mouth}{\operatorname{mo}}
\newcommand{\mo}{\operatorname{mo}}     
\def\adhl{\prec}

\def\vclass#1{[#1]_{\cV}}

\definecolor{lgray}{rgb}{0.9, 0.9, 0.9}
\def\cg{}

\newcommand{\lep}[2][]{#2^{\sqsubset#1}}
\renewcommand{\rep}[2][]{#2^{\sqsupset#1}}

\def\category#1{\text{\sc #1}}
\newcommand{\PfCC}{\category{PfCc}}
\newcommand{\EgPfCC}{\category{EgPfCc}}
\newcommand{\PgCC}{\category{PgCc}}
\newcommand{\GMod}{\category{GMod}}
\newcommand{\FMod}{\category{FMod}}
\newcommand{\CCR}{\category{Cc}}
\newcommand{\ChCC}{\category{ChCc}}
\newcommand{\ChPfCC}{\category{ChPfCc}}
\newcommand{\Ch}{\category{Ch}}
\newcommand{\DSet}{\category{DSet}}
\newcommand{\DPSet}{\category{DPSet}}

\renewcommand{\conv}{\protect\mbox{\rm conv}}

\def\boundaryless{boundaryless}
\def\reduced{reduced}
\def\Psub{P^\sharp}
\def\PPsub{\PP^\sharp}
\def\leqsub{\leq^\sharp}

\def\xa{\mathbf{p}}
\def\xb{\mathbf{q}}
\def\xc{\mathbf{r}}
\def\xd{\mathbf{s}}
\def\xe{\mathbf{t}}
\def\xf{\mathbf{u}}

\def\bA{\mathbf{A}}
\def\bB{\mathbf{B}}
\def\bC{\mathbf{C}}
\def\bD{\mathbf{D}}
\def\bE{\mathbf{E}}
\def\bF{\mathbf{F}}

\def\ba{\mathbf{a}}
\def\bb{\mathbf{b}}
\def\bc{\mathbf{c}}
\def\balpha{{\boldsymbol\alpha}}

\def\phir{\phi_|}
\def\psir{\psi_|}

\newcommand{\exend}{\hspace*{\fill}$\Diamond$}

\articletheorems

\begin{document}

\author{Marian Mrozek}
\address{Marian Mrozek, Division of Computational Mathematics,
  Faculty of Mathematics and Computer Science,
  Jagiellonian University, ul.~St. \L{}ojasiewicza 6, 30-348~Krak\'ow, Poland
}
\email{Marian.Mrozek@uj.edu.pl}
\author{Thomas Wanner}
\address{Thomas Wanner, Department of Mathematical Sciences,
George Mason University, Fairfax, VA 22030, USA
} \email{twanner@gmu.edu}
\date{today}
\thanks{Research of  M.M.\ was partially supported by
  the Polish National Science Center under Ma\-estro Grant No. 2014/14/A/ST1/00453
  and Opus Grant No. 2019/35/B/ST1/00874.
  T.W.\ was partially supported by NSF grant DMS-1407087 and by the Simons Foundation
  under Award~581334.
}

\title[Connection matrices in combinatorial topological dynamics]
{Connection matrices in combinatorial topological dynamics}

\date{Version compiled on \currentDate}

\begin{abstract}
Connection matrices are one of the central tools in Conley's approach to the study
of dynamical systems, as they provide information on the existence of connecting
orbits in Morse decompositions. They may be considered as a generalization of the 
Morse complex boundary operator in Morse theory. 
Their computability has recently been addressed by Harker, Mischaikow, 
and Spendlove~\cite{HMS2021} in the context of lattice filtered chain complexes. 
In the current paper, we extend the newly introduced Conley theory for combinatorial
vector and multivector fields on Lefschetz complexes~\cite{LKMW2020}
by transferring the concept of connection matrix to this setting. 
This is accomplished by connection matrices for arbitrary poset filtered chain complexes, 
as well as an associated equivalence, which allows for changes in the underlying  posets.
We show that for the special case of gradient combinatorial 
vector fields in the sense of Forman~\cite{Fo98a}, connection matrices are necessarily unique.
Thus, the classical results of Reineck~\cite{reineck:90a, reineck:95a} have a natural
analogue in the combinatorial setting.
\end{abstract}

\maketitle


\newpage

\tableofcontents

\newpage


\section{Introduction}
\label{sec:intro}

Classical Morse theory concerns a compact smooth manifold together with a
smooth real-valued function with non-degenerate critical points. 
An important part of the theory introduces the Morse complex
which is a chain complex whose $i$th chain group is a free abelian group 
spanned by critical points of Morse index $i$ and whose boundary homomorphism 
is defined by counting the (oriented) flow lines between critical points
in the gradient flow induced by the Morse function (see~\cite[Section 4.2]{Knudson2015}).
One of the fundamental results of classical Morse theory states that the homology
of the manifold is isomorphic to the homology of the Morse complex.

The stationary points of the gradient flow in Morse theory, 
which are precisely the critical points of the Morse function,
provide the simplest example of an isolated invariant set, a key concept
of Conley theory~\cite{Conley1978}. For every isolated invariant set there
is a homology module associated with it. It is called the homology Conley index.
A (minimal) Morse decomposition is a decomposition of space into a partially
ordered collection of isolated invariant sets, called Morse sets, such that
every recurrent trajectory (in particular every stationary or periodic trajectory)
is located in a Morse set and every non-recurrent trajectory is a heteroclinic
connection between Morse sets from a higher Morse set to a lower Morse set in
the poset structure of the Morse decomposition. The collection of stationary
points of the gradient flow of a Morse function provides the simplest example
of a Morse decomposition in which the Morse sets are just the stationary points
and the Conley index of a stationary point coincides with the homology of a
pointed $k$-dimensional sphere with $k$ equal to the Morse index of the point. 

Conley theory in its simplest form may be viewed as a twofold generalization
of Morse theory. On the one hand, it substantially weakens the general
assumptions by replacing the smooth manifold by a compact metric space
and the gradient flow of the Morse function by an arbitrary (semi)flow.
On the other hand, it replaces the collection of critical points of the
Morse function by the more general Morse decomposition in which the
counterpart of the Morse complex takes the form of the direct sum of
the Conley indexes of all Morse sets. The homology of this generalized
complex, as in the Morse theory, is isomorphic to the homology of the space.
The boundary operator in this setting is called the connection matrix of the
Morse decomposition. This definition of connection matrix was introduced
by Franzosa in~\cite{Fr1989}. His definition, based on homology braids,
is technically complicated, in part because the generalized setting
captures the situations of bifurcations when, unlike for Morse theory,
the connection matrix need not be uniquely determined by the flow.
In a later paper Robbin and Salamon~\cite{RoSa1992} simplify slightly
the definition by replacing homology braids with filtered chain complexes
which helps with separating dynamics from algebra. This separation is
even more visible in the recent algorithmic approach to connection matrices 
by Harker, Mischaikow, and Spendlove~\cite{HMS2021}. Actually, the
separation of dynamics and algebra allows the authors in~\cite{HMS2021,
RoSa1992} to set up the definition of a connection matrix of a Morse
decomposition along the following pipeline, where we omit a number
of technical details:
\begin{itemize}                                
   \item[(i)]   Consider a Morse decomposition $\cM:=\{M_p\}_{p\in P}$
                of the phase space indexed by a poset $P$.
   \item[(ii)]  For each down set $I\subset P$ consider the associated
                attractor~$M_I$ consisting of points on trajectories whose 
                limits sets are in $\bigcup_{p\in I}M_p$ and construct an
                attracting neighborhood~$N_I$ for the attractor~$M_I$ in
                such a way that $I\mapsto N_I$ is a lattice homomorphism.
   \item[(iii)] Under some smoothness assumptions the family $\{N_I\}$ 
                of step (ii) induces a $P$-filtered chain complex.
   \item[(iv)]   A connection matrix of $\cM$ is an algebraic object
                associated with a poset filtered chain complex, in
                particular with the $P$-filtered chain complex of step~(iii).
\end{itemize}
Conley theory, in particular via connection matrices, is a very useful tool
in the qualitative study of dynamical systems. However, to apply it one
requires a well-defined dynamical system on a compact metric space. This
is not the case when the dynamical system is exposed only via a finite
set of samples as in the case of time series collected from observations
or experiments. The study of dynamical systems known only from samples
becomes an important part of the rapidly growing field of data science. 
In this context a generalization of Morse theory presented by Robin
Forman~\cite{Fo98a} turns out to be very fruitful. In this generalization
the smooth manifold is replaced by a finite CW complex and the gradient
vector field of the Morse function by the concept of a combinatorial
vector field. These structures may be easily constructed from data
and analyzed by means of the combinatorial, also called discrete,
Morse theory by Forman. 

Recently, the concepts of isolated invariant set and Conley index
have been  carried over to this combinatorial setting~\cite{batko:etal:20a,
KaMrWa2016, LKMW2020, Mr2017, mrozek:wanner:21a}. In the present paper we
extend these ideas by constructing connection matrices of a Morse
decomposition of a combinatorial multivector field. To achieve this,
we have to modify the connection matrix pipeline discussed above, because
in the combinatorial case, and even in the case of a multiflow, step~(ii)
in the pipeline cannot be completed in general. To overcome this difficulty
we enlarge the original poset to guarantee the existence of the necessary
lattice of attracting neighborhoods. The added elements are then removed by
introducing a certain equivalence in the category of poset filtered chain complexes
with a varying poset structure. A side benefit of this approach is the
shortening of the pipeline by omitting step~(ii).
Such a shortening is possible, because,  under the enlarged poset, 
the filtered chain complex of step (iii) may be obtained
directly from the Morse decomposition via an associated partition of the phase space indexed by
the extended poset. The proposed theory of connection matrices for
combinatorial multivector fields lets us prove the main result of
the paper: The uniqueness of connection matrices for gradient
combinatorial Forman vector fields.

Although the main motivation for our paper is the adaptation
of connection matrix theory to the combinatorial setting
of multivector fields, we believe that the approach presented here
has broader potential. Formally speaking, connection matrices are
purely algebraic objects and are presented this way from the very
beginning. Nevertheless, in the early papers they are strongly tied
to dynamical considerations. In particular, the important concept of
uniqueness and non-uniqueness in these papers is addressed only via the
underlying dynamics. As we already mentioned, the decoupling of algebra
and dynamics started with Robbin and Salamon~\cite{RoSa1992}, and is
even stronger in Harker, Mischaikow, and Spendlove~\cite{HMS2021}.
However, to the best of our knowledge, so far there has been no purely
algebraic definition of uniqueness. Harker, Mischaikow, and
Spendlove~\cite{HMS2021} prove that any two connection matrices of
the same filtered chain complex are conjugate via a filtered isomorphism.
Clearly, this is not the uniqueness concept used in the context of dynamics.
In our approach we propose a stronger algebraic equivalence of connection
matrices which allows for filtered chain complexes with no unique connection
matrix. We use this definition in the already mentioned main result of our
paper on the uniqueness of connection matrices for gradient combinatorial
vector fields. We believe that the detachment of connection matrix theory
from dynamics is worth the effort, because it may bring applications in new
fields. The potential of applications in topological data analysis is
already discussed in~\cite{HMS2021}, and the ties between connection
matrices and persistent homology are also visible in~\cite{dey:etal:p22}.
Moreover, combinatorial vector fields have been used successfully in the 
study of algebraic and combinatorial problems~\cite{JW2009, Ko08, Sk2006}.

Finally, combinatorial multivector fields make it possible to construct examples of
a variety of complex dynamical phenomena in a straightforward way. It is therefore
our hope that the results of this paper make topological methods in dynamics, and
in particular the concept of connection matrices, more accessible to a broader
mathematical audience.

The remainder of this paper is organized as follows. In Section~\ref{sec:main} we
provide an informal overview over our main definitions and results. By providing
numerous examples, we hope that this section will serve as a guide for the more technical
later sections of the paper. In addition, this section can serve as a quick and informal 
introduction to the concept of connection matrices. After presenting necessary preliminary
definitions and facts in Section~\ref{sec:prelim}, the main technical part of the paper
starts with a discussion of poset graded and poset filtered chain complexes in
Section~\ref{sec:posetfilteredcc}. Based on these definitions, the algebraic
connection matrix can then be introduced in Section~\ref{sec:algconnmatrix}.
The remaining three sections of the paper apply the earlier algebraic constructions
to multivector fields on Lefschetz complexes. While the connection matrix in this
setting is considered in Section~\ref{sec:lefconnmatrix}, the dynamics of combinatorial
multivector fields is the subject of Section~\ref{sec:dcmvf}. Finally, in
Section~\ref{sec:cm-gcvf} we show that if one considers a gradient Forman vector
field on a Lefschetz complex, then the associated connection matrix is necessarily
uniquely determined.


\section{Main results}
\label{sec:main}

In this section we informally summarize our definition of connection matrices in the combinatorial setting,
together with the necessary background material, examples,  and main results of the paper.
Precise definitions, statements, and proofs will be given in the sequel. We emphasize that throughout the paper
we assume field coefficients in all considered modules, in particular in chain complexes and homology modules.

\subsection{Combinatorial dynamical systems and multivector fields}
By a {\em combinatorial dynamical system} we mean a  multivalued map $F:X\mto X$ defined 
on a finite topological space $X$. Alternatively, it may be viewed as a finite directed graph whose set 
of vertices is the topological space $X$, and with $F$ interpreted as the map sending a vertex to the collection 
of its neighbors connected via an outgoing directed edge. We call it the {\em $F$-digraph}.
In most applications, this set $X$ of vertices of the $F$-digraph 
is a collection or certain subcollection of cells of a finite cellular complex,
for instance a simplicial complex, with its
topology induced by the associated face poset via the Alexandrov Theorem~\cite{Al1937}.
This topology is very different from the metric topology of the geometric realization 
of the cellular complex in terms of its separation properties, but the same 
in terms of homotopy and homology groups via McCord's Theorem~\cite{MC1966}. 
Hence, as far as algebraic topological invariants are concerned, these topologies may
be used interchangeably, with the finite topology having the advantage of better
explaining some peculiarities of combinatorial dynamics. 

As in classical multivalued dynamics, given a combinatorial dynamical system $F$,
we define a {\em solution} of $F$ to be a map $\gamma:I\to X$
defined on a subset $I\subset\ZZ$ of integers and such that $\gamma(i+1)\in F(\gamma(i))$ for $i,i+1\in I$.
The solution $\gamma$ is called {\em full} if $I=\ZZ$, and it is a {\em path} if $I$ is the intersection of $\ZZ$ with a compact 
interval in $\RR$. In the directed graph interpretation, a solution 
may be viewed as a finite or infinite directed walk through the graph. 
Finally, a subset $A\subset X$ is called {\em invariant}  if for every point $a\in A$ 
there exists a full solution $\gamma:\ZZ\to A$ through the point $a$, that is,
a solution satisfying the identity $\gamma(0)=a$.  

In this paper we are interested in a class of combinatorial dynamical systems 
induced by combinatorial vector and multivector fields. We recall that 
a {\em combinatorial multivector field}~$\cV$ is a partition of a finite topological space~$X$ 
into multivectors, where a {\em multivector} is a {\em locally closed set}, that is, a set which
is the difference of two closed sets (see~\cite{En1989}). 
A multivector is called a {\em vector}, if its cardinality is one or two. If  a combinatorial multivector field contains only vectors, then
we call it a {\em combinatorial vector field}, a concept introduced already by Forman~\cite{Fo98a}.  
In applications, the partition indicates the resolution below which there is not enough evidence to 
analyse the dynamics either due to insufficient amount of data or because of lacking computational power.

The results of this paper require the combinatorial multivector field to be defined 
on a special finite topological space, namely a Lefschetz complex~\cite{Le1942}, called by Lefschetz and in~\cite{HMS2021} simply a complex.
In short, a {\em Lefschetz complex} is just a basis  $X$ of a finitely generated free chain complex $(C,d)$ with emphasis 
on the basis. This means that $X$ is the Lefschetz complex and $(C,d)$, denoted $C(X)$,  is the chain complex associated with the Lefschetz complex.
In terms of applications, a typical example of a Lefschetz complex is the set of cells of a cellular complex
or simplices of a simplicial complex. These simplices or cells constitute a natural basis for the associated chain complex.
Also, every locally closed subset of a Lefschetz complex is a Lefschetz complex, a feature which, in particular,
facilitates constructing concise examples.
By the homology of a Lefschetz complex $X$ we mean the homology of the associated chain complex $C(X)$. 
We view a Lefschetz complex as a finite topological space via the Alexandrov Theorem. For this, we need a well-defined face relation, which is given by the transitive closure of the facet relation, 
where an element $x\in X$ is called a {\em facet} of $y\in X$ if~$x$ appears with non-zero coefficient in the boundary~$d(y)$.

\begin{figure}
  \begin{center}
    \includegraphics[width=0.39\textwidth]{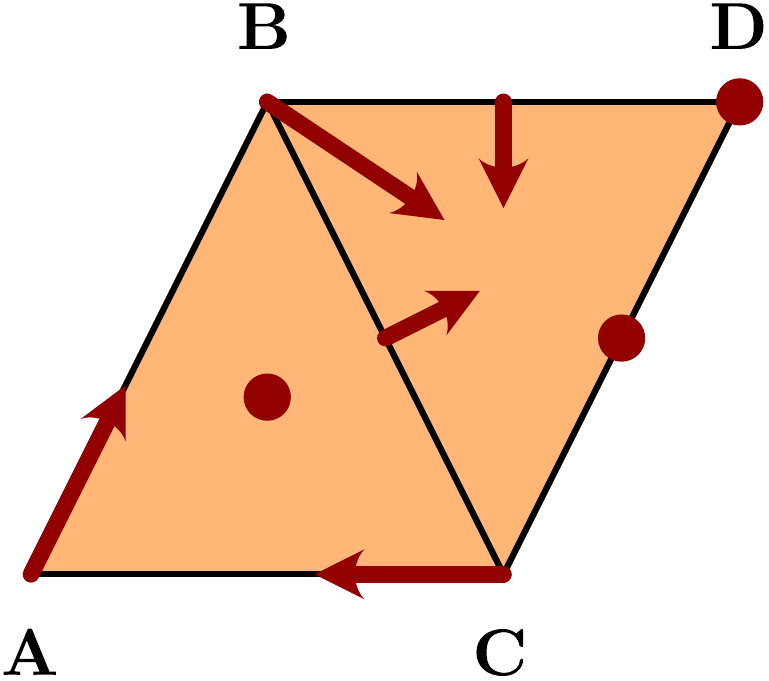}\quad
    \includegraphics[width=0.55\textwidth]{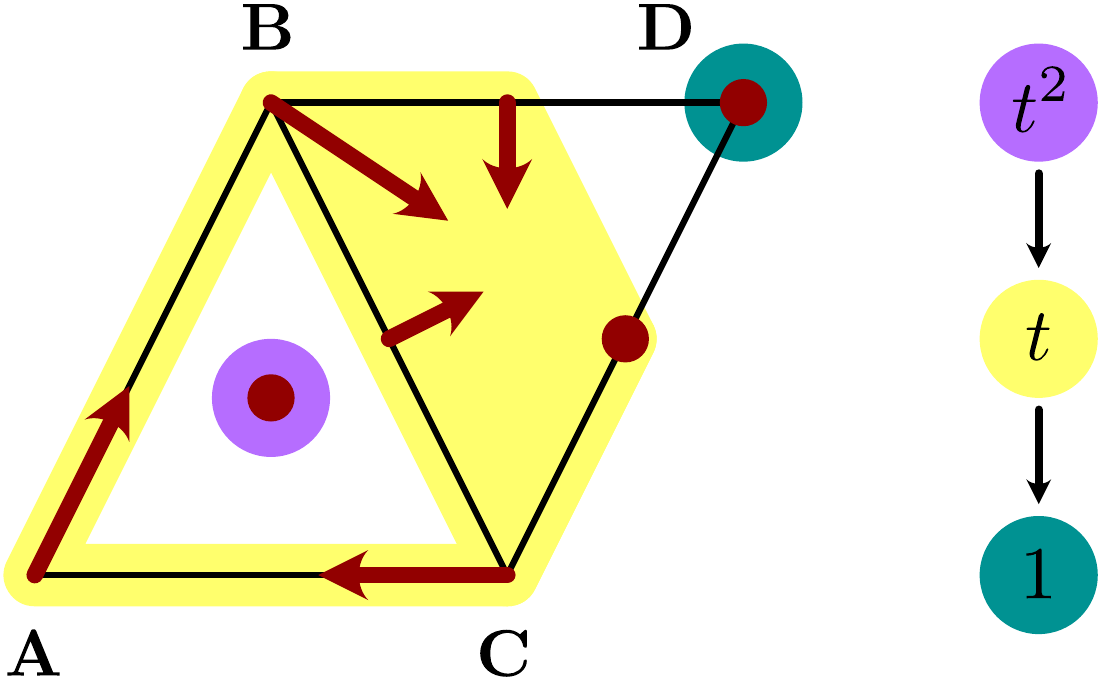}
  \end{center}
  \caption{{\em A first multivector field.} The left panel shows a combinatorial multivector
           field on a simplicial complex which consists of two triangles, five edges,
           and four vertices. The multivector field has three critical cells, marked as
           red dots, two vectors on the bottom and left-most edges, as well as one
           multivector, which consists of the right triangle~$\bB\bC\bD$, its two
           edges~$\bB\bC$ and~$\bB\bD$, and the vertex~$\bB$. In the middle 
           panel, we indicate the three Morse sets of the associated combinatorial
           dynamical system~$F_\cV$ in purple, yellow, and green, while the Conley-Morse
           graph of this Morse decomposition is shown in the panel on the right.}
  \label{fig:multivectorfield}
\end{figure}
\begin{ex}[{\em A first multivector field}]
\label{ex:multivectorfield}
{\em
Figure~\ref{fig:multivectorfield}(left) presents an example of a combinatorial multivector field 
on a Lefschetz complex $X$ which is just a simplicial complex consisting of two triangles
$\bA\bB\bC$, $\bB\bC\bD$, five edges $\bA\bB$, $\bA\bC$, $\bB\bC$, $\bB\bD$, $\bC\bD$, and
four vertices $\bA$,  $\bB$, $\bC$, $\bD$. Hence,~$X$ is a finite topological space consisting of
eleven elements. The combinatorial multivector field in Figure~\ref{fig:multivectorfield} consists
of three singletons $\{\bA\bB\bC\}$, $\{\bC\bD\}$, and~$\{\bD\}$ indicated in the figure by a red
dot, two doubletons $\{\bA,\bA\bB\}$ and~$\{\bC,\bA\bC\}$ marked with a red arrow, and one
multivector $\{\bB,\bB\bC,\bB\bD,\bB\bC\bD\}$ consisting of four simplices and marked with red arrows joining 
each simplex in the multivector with its top-dimensional coface in the same multivector. Why this method of drawing
multivectors is natural will become clear when we discuss the combinatorial dynamical system induced by the combinatorial
multivector field.
\exend}
\end{ex} 
 
Every combinatorial multivector field $\cV$ induces a combinatorial dynamical 
system $F_\cV: X\mto X$ given for $x\in X$ by 
\begin{equation}
\label{eq:F-cV}
F_\cV(x):=\cl x\cup [x]_\cV \; .
\end{equation}
In this definition, $\cl x$ stands for the closure of $x$ in $X$, that is, it is given by the collection of all faces of~$x$. Moreover, $[x]_\cV$ denotes the unique multivector~$V$ in the partition $\cV$ with $x\in V$.
The formula for $F_\cV$ is related to the resolution interpretation of a combinatorial multivector field mentioned earlier.
Namely, inside a multivector we do not exclude any movement, treating a multivector as kind of a black box.
This justifies the presence of~$[x]_\cV$ in~\eqref{eq:F-cV}.
And, since the dynamics of a combinatorial multivector field models a flow, any movement between multivectors 
is possible only through their common boundary. This justifies the presence of the closure~$\cl x$ in  \eqref{eq:F-cV}.

\begin{ex}[{\em A first multivector field}, continued]
\label{ex:multivectorfield-2}
{\em
Even for a simple multivector field such as the one shown in
Figure~\ref{fig:multivectorfield}(left) the induced combinatorial
dynamical system is quite large. Interpreted as a directed graph, it
has eleven vertices (one for each simplex) and~$42$ directed edges.
Hence, instead of drawing such a digraph, we interpret
Figure~\ref{fig:multivectorfield}(left) as a digraph with
vertices in the centers of mass of simplices (not marked) and only a minimum 
of directed edges (arrows) which cannot be deduced from \eqref{eq:F-cV}. In
fact, the only arrows which do depend on~$\cV$ are arrows of the form $x\to y$
where $y\in [x]_\cV\setminus \cl x$. Moreover, since vertices inside the same
multivector form a clique in the directed graph, it suffices to mark only some
of the arrows joining them and, as we already mentioned, we mark arrows joining
each simplex in the multivector with its top-dimensional cofaces in the same
multivector. Interpreting Figure~\ref{fig:multivectorfield}(left) this way,
it is not difficult to infer the solutions of $F_\cV$ from the figure. For
instance, we have a solution $\gamma$ defined for $n\in\ZZ$ by 
\[
\gamma(n):=\begin{cases}
             \bA\bB\bC & n<0,\\
             \bB\bC & n=0,\\
             \bB\bC\bD & n=1,\\
             \bB\bD & n=2,\\
             \bD & n>2.
           \end{cases}
\]
In the sequel we will write solutions in a compact form 
\[
\gamma=\overleftarrow{\bA\bB\bC}\cdot\bB\bC\cdot\bB\bC\bD\cdot \bB\bD\cdot \overrightarrow{\bD}
\]
were the left and right arrows on top of the left and right ends, respectively,
indicate that the same pattern is repeated up to negative and positive infinity.
\exend}
\end{ex} 

An undesired consequence of the otherwise natural formula \eqref{eq:F-cV} 
is that every $x\in X$ is a fixed point of $F_\cV$, that is
$x\in F_\cV(x)$. This immediately implies that
every subset of $X$ is invariant for combinatorial dynamical systems induced by combinatorial 
multivector fields. Clearly, the interest is in invariant sets whose invariance goes beyond the fact 
that all points are fixed points of $F_\cV$. In order to properly describe such invariant sets we need the following 
observations and definitions.

According to \eqref{eq:F-cV}, a solution of the combinatorial dynamical system may leave a multivector $V\in\cV$
only via a point in $\mouth V:=\cl V\setminus V$, a set which we call the {\em mouth} of $V$.
Indeed, if $\sigma(i)\in V$ and $\sigma(i+1)\not \in V$, then one necessarily has $\sigma(i+1)\in\mouth V$.
This means that the mouth of~$V$ acts as the exit set for~$V$. The mouth is closed, because $V$
is locally closed. This allows us to interpret the closure~$\cl V$ as a small isolating block for
the multivector~$V$, and the relative homology~$H(\cl V,\mouth V)$ as its {\em Conley index}.
We call a multivector {\em critical} if $H(\cl V,\mouth V)\neq 0$, and {\em regular} otherwise.
Motivated by the fact that in the classical setting a zero Conley index of an isolating block does
not guarantee the existence of a full solution inside the block, 
we say that a full solution $\gamma$ is {\em essential}, if for every regular multivector~$V$ the preimage
$\gamma^{-1}(V)$ does not contain an infinite interval of integers. In other words, every essential
solution must leave a regular multivector~$V$ in forward and backward time before it can enter again.

An invariant set $A$ is called an {\em essential invariant set}, 
if every $a\in A$ admits an essential solution through~$a$ which is completely contained in~$A$. 
We are interested in {\em isolated invariant sets}, which are defined as essential invariant sets admitting 
a closed superset $N\supset S$  such that  $ F(S)\subset N $, and such that 
every path in $N$ with endpoints in $S$ is itself contained in $S$.
We call such a set $N$ an {\em isolating set}.
Isolating sets are the combinatorial counterparts of  isolating neighborhoods in the classical setting.
Nevertheless, we prefer the name {\em isolating set}, because they do not need to be neighborhoods in general.
An isolating set is called an {\em attracting set} if every path in~$X$ with left endpoint in~$N$ is contained in~$N$,
and the associated isolated invariant set~$S$ is then called an {\em attractor}. Similarly one defines a {\em repelling set} and a {\em repeller}.
We note that an essential invariant set is isolated if and only if it is locally closed
and {\em $\cV$-compatible}, that is, it equals the union of all multivectors contained in
it, see for example~\cite[Proposition 4.10, 4.12, 4.13]{LKMW2020}. Also, an essential invariant
set is an attractor (or repeller) if and only if it is $\cV$-compatible and closed (or open),
as was shown in~\cite[Theorem 6.2]{LKMW2020}.

\begin{ex}[{\em A first multivector field}, continued]
\label{ex:multivectorfield-3}
{\em
  It is easy to verify that every singleton in a multivector field is critical, and every doubleton is
  necessarily regular. In contrast, for a general multivector of cardinality greater than two it is not
  possible to automatically determine whether it is critical or regular just based on its cardinality.
  For example, in Figure~\ref{fig:multivectorfield}(left) the multivector
  $V:=\{\bB,\bB\bC,\bB\bD,\bB\bC\bD\}$ is regular, 
  since its closure may be homotopied to its mouth, i.e., one has $H(\cl V,\mouth V)= 0$.
  Yet, it is not difficult to see that detaching the vertex~$\bB$ from~$V$ and instead
  attaching the edge~$\bC\bD$ to~$V$ preserves the cardinality of the multivector, but makes~$V$
  critical. An example of a solution in Figure~\ref{fig:multivectorfield}(left) 
  which is not essential is
\[
\overleftarrow{\bC\cdot\bA\bC}\cdot\bA\cdot\bA\bB\cdot\overrightarrow{\bB\cdot\bB\bC\bD}
\]
and of a solution which is essential is
\[
\overleftarrow{\bC\cdot\bA\bC\cdot\bA\cdot\bA\bB\cdot\bB\cdot\bB\bC\bD}\cdot\overrightarrow{\bD}.
\]
Furthermore, if we consider $A:=\{\bA,\bB,\bC,\bA\bB,\bA\bC,\bB\bC\}$
in Figure~\ref{fig:multivectorfield}(left), then the solution
\begin{equation} \label{ex:multivectorfield-3a}
  \overleftrightarrow{\bA\cdot\bA\bB\cdot\bB\cdot\bB\bC\cdot\bC\cdot\bA\bC}
\end{equation}
is a periodic essential solution which passes through all elements of~$A$. This demonstrates
that the set~$A$ is an essential invariant set. However, it is not an isolated invariant set, because
it is not $\cV$-compatible. We may modify $A$ by taking $A' := A \cup \{\bB\bD,\bB\bC\bD\}$.
Then the modified set $A'$ is $\cV$-compatible, but still not an isolated invariant set, 
since it is not locally closed. In fact, the smallest isolated invariant set containing~$A$
is $A'':=A\cup \{\bB\bD,\bB\bC\bD,\bC\bD\}$, which contains all simplices of~$X$ except
for the triangle~$\bA\bB\bC$ and the vertex~$\bD$.
\exend}
\end{ex}

\subsection{Conley index and Morse decompositions}  
Given an isolated invariant set~$S$ of a combinatorial multivector field~$\cV$ on a Lefschetz complex~$X$, 
we define its {\em Conley index} as the relative homology $H(\cl S,\mo S)$. This definition is motivated
by the fact that the topological pair~$(\cl S,\mo S)$ is one of possibly many index pairs for~$S$, see for
example~\cite[Definition~5.1 and Proposition~5.3]{LKMW2020} for more details. Since in this paper
we are only interested in homology with field coefficients, and since such homology is uniquely 
determined by the associated Betti numbers, it is convenient to identify the homology and, in particular,
the Conley index with the associated {\em Poincar\'e polynomial} whose coefficients are the consecutive
Betti numbers. In the case of the Conley index we refer to this polynomial as the {\em Conley polynomial}
of the isolated invariant set.

\begin{ex}[{\em A first multivector field}, continued]
\label{ex:multivectorfield-4}
{\em
  In the situation of the combinatorial multivector field from Figure~\ref{fig:multivectorfield},
  the multivector given by the singleton $S:=\{\bA\bB\bC\}$ is an isolated invariant set. One can
  easily see that the relative homology~$H(\cl S,\mouth S)$ is the homology of a pointed $2$-sphere.
  Therefore, the Conley polynomial of~$\{\bA\bB\bC\}$ is given by~$t^2$. Similarly, it was shown
  in Example~\ref{ex:multivectorfield-3} that the smallest isolated invariant set containing the
  periodic solution~(\ref{ex:multivectorfield-3a}) is given by the set~$A''$ of all simplices of~$X$,
  except for the triangle~$\bA\bB\bC$ and the vertex~$\bD$. We leave it for the reader to verify that
  the Conley index~$H(\cl A'', \mo A'')$ then has the Conley polynomial~$t$. Finally,
  the Conley index of the singleton~$\{\bD\}$ is the homology of a pointed $0$-sphere,
  and therefore the associated Conley polynomial is~$1$.
\exend}
\end{ex}

A family $\{M_r\}_{r\in P}$ of mutually disjoint, non-empty, isolated invariant sets
indexed by a poset $P$ is called a {\em Morse decomposition} of~$X$ if for every
essential solution~$\gamma$ either all values of~$\gamma$ are contained in the same
set~$M_r$, or there exist indices $q>p$ in~$P$ and $t_p,t_q\in\ZZ$ such that
$\gamma(t)\in M_q$ for $t\leq t_q$ and $\gamma(t)\in M_p$ for $t\geq t_p$.
In the latter case, the solution~$\gamma$ is called a {\em connection} from~$M_q$
to~$M_p$. Furthermore, the sets~$M_r$ are called the {\em Morse sets} of the Morse
decomposition. 

One can check that the collection of strongly connected components of the $F$-digraph
is a Morse decomposition. In fact, this Morse decomposition is always the finest
possible Morse decomposition, see for example~\cite[Theorem~4.1]{DJKKLM2017}. In this
respect, the combinatorial case differs from the classical case where the finest
Morse decomposition may not exist.

It is customary to condense the information about a Morse decomposition in the form
of its {\em Conley-Morse graph}. This graph consists of the Hasse diagram of the
poset~$P$ whose vertices, representing the individual Morse sets~$M_p$, are
labeled with the respective Conley polynomials.

\begin{ex}[{\em A first multivector field}, continued]
\label{ex:multivectorfield-5}
{\em
Consider $P:=\{\xa,\xb,\xc\}$  linearly ordered  by $\xa<\xb<\xc$,
and again the combinatorial multivector field presented in
Figure~\ref{fig:multivectorfield}(left). It is not difficult to
verify that the family~$\cM$ consisting of $M_\xc=\{\bA\bB\bC\}$,
$M_\xb=\{\bA,\bB,\bC,\bA\bB,\bA\bC,\bB\bC,\bB\bD,\bC\bD,\bB\bC\bD\}$,
and~$M_\xa=\{\bD\}$ is a Morse decomposition. In
Figure~\ref{fig:multivectorfield}(middle), the three Morse sets
are indicated with three different colors. In our example all three
Morse sets have one-dimensional Conley indices, respectively, in dimension
zero for~$M_\xa$, one for~$M_\xb$, and two for~$M_\xc$. Hence, the
Conley polynomials are given by~$1$ for~$M_\xa$, $t$ for~$M_\xb$,
and~$t^2$ for~$M_\xc$. The associated Conley-Morse graph is
visualized in Figure~\ref{fig:multivectorfield}(right).

Although we have not yet given the definition of the connection matrix
for a Morse decomposition, we can characterize some of its features
on the basis of this example. As we will see later, a connection matrix
is just a matrix representation of an abstract boundary homomorphism
acting on the direct sum of Conley indexes of Morse sets. The entries
in the matrix are homomorphisms between the individual Conley indexes.
As in the classical Morse theory, the homology of the resulting chain
complex coincides with the homology of the underlying Lefschetz complex. Assuming
homology coefficients in the field~$\ZZ_2$, the only possible homomorphisms
between Conley indexes of Morse sets in our example are either zero or
an isomorphism. Denote the Conley index of the Morse set~$M_p$ by~$C_p$.
Marking the isomorphisms by~$1$ and leaving empty fields for the zero
homomorphism, the connection matrix for our example turns out to be
\[
\begin{array}{c||c|c|c|}
    &  C_\xa &  C_\xb &    C_\xc \\
  \hline
  \hline
  C_\xa &    &  0 &     \\
  \hline                                                                                    
  C_\xb &    &    &  1  \\
  \hline
 C_\xc &    &    &     \\
  \hline
\end{array} \; .
\]
Hence, all entries of the connection matrix are zero, except for the
homomorphism from~$C_\xc$ to~$C_\xb$. In our simple example, each
Conley index is non-trivial in a different grade, and this in turn 
implies that~$C_\xa$, $C_\xb$, and~$C_\xc$ are also chain subgroups
in grades zero, one, and two, respectively. Therefore, the chain complex
may be written in the compact form
\[
   \begin{diagram}
  \dgARROWLENGTH 1.6em
    \node{\cdots}
    \node{0}
    \arrow{w,t}{0}
    \node{C_\xa}
    \arrow{w,t}{0}
    \node{C_\xb}
    \arrow{w,t}{0}
    \node{C_\xc}
    \arrow{w,t}{1}
    \node{0}
    \arrow{w,t}{0}
    \node{\cdots}
    \arrow{w,t}{0}
  \end{diagram}.
\]
One can easily see that its homology is trivial except in grade
zero --- and this provides precisely the homology of the Lefschetz
complex in the considered example. In fact, the connection matrix
acts as an algebraic version of the Conley-Morse graph. Its unique
non-zero entry reflects the heteroclinic connection from the Morse
set~$M_\xc$ to the set~$M_\xb$. Notice that we also have two different
heteroclinic connections from~$M_\xb$ to~$M_\xa$, but they algebraically
annihilate each other and thus lead to the corresponding entry being
zero in the connection matrix.
\exend}
\end{ex}

\begin{figure}
  \begin{center}
    \includegraphics[width=0.46\textwidth]{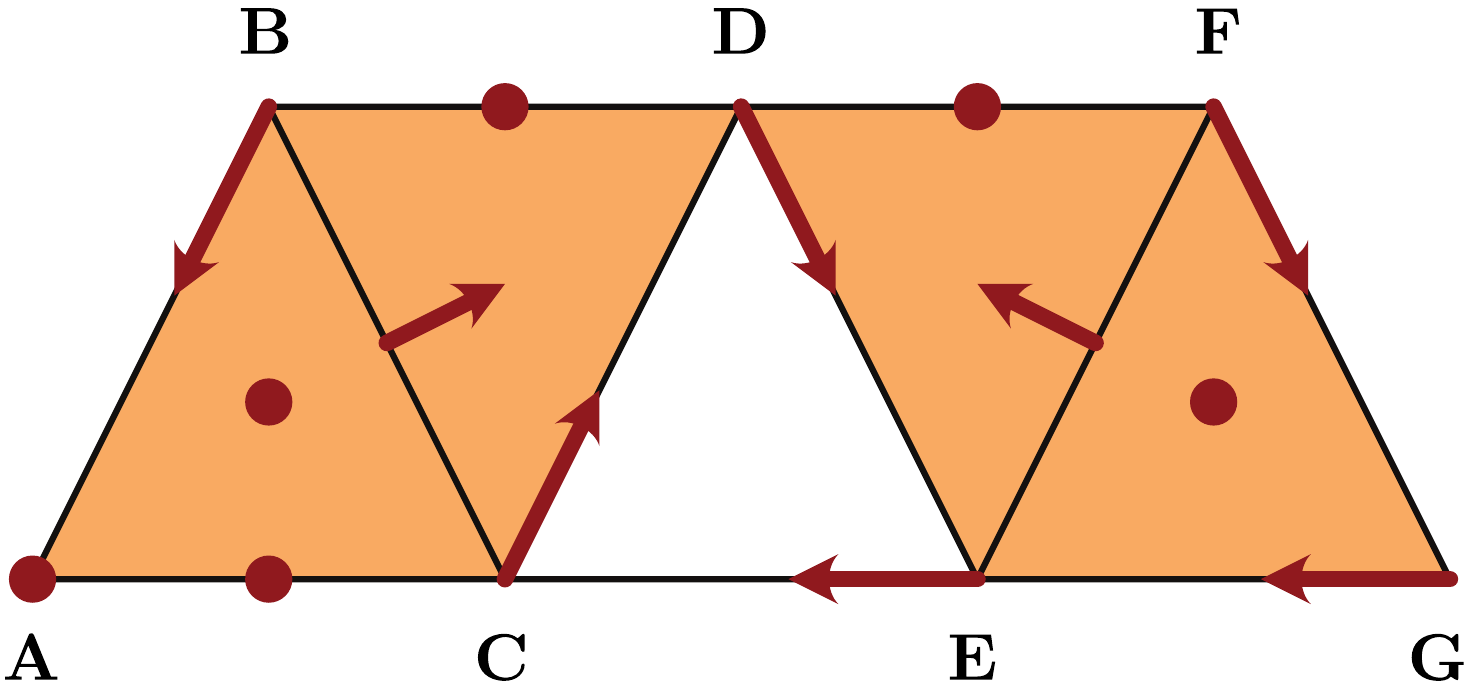}\qquad
    \includegraphics[width=0.46\textwidth]{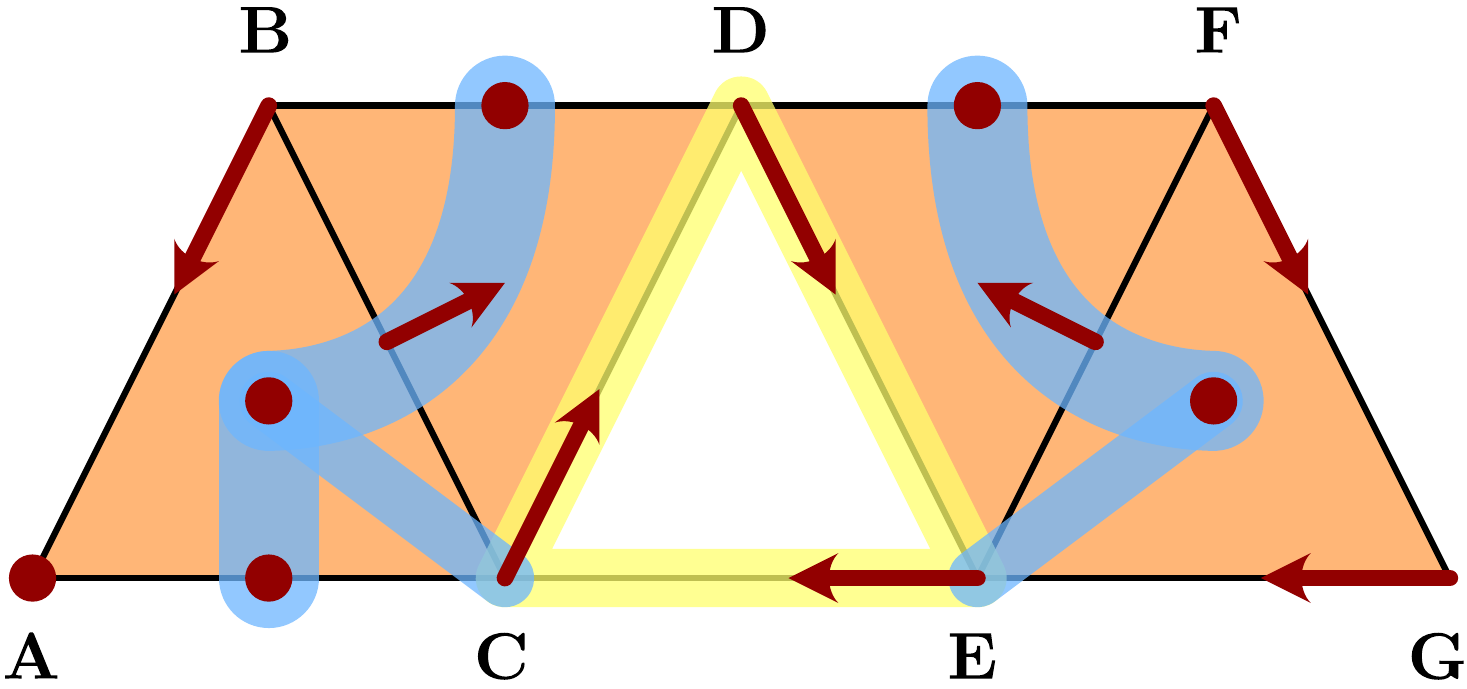} \\[4ex]
    \includegraphics[width=0.85\textwidth]{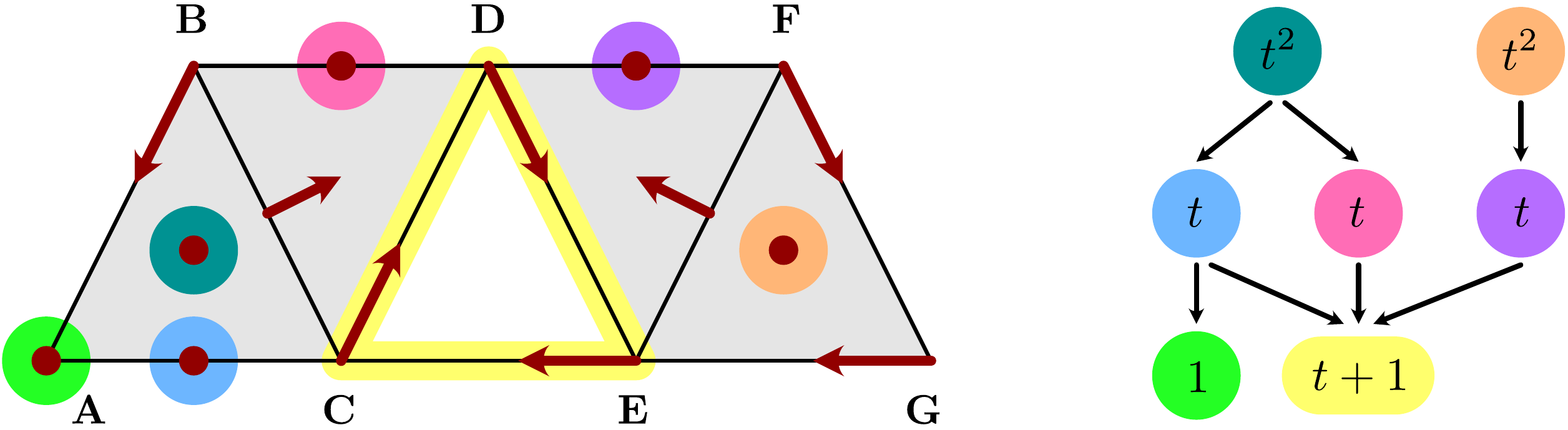}
  \end{center}
  \caption{{\em A Forman vector field with periodic orbit.}
           The top left panel depicts a Forman vector field on a simplicial
           complex, which consists of six critical cells and eight vectors
           given by doubletons. Notice that the edges~$\bC\bD$, $\bD\bE$,
           and~$\bE\bC$ give rise to a periodic solution. In the top right
           panel, five connections originating in the two critical cells in
           dimension two are shown. While three of these connect to critical
           cells of dimension one, two connect to the periodic orbit. The two 
           panels on the bottom illustrate the Morse sets of the
           combinatorial vector field in different colors, together
           with the associated Conley-Morse graph.
           }
  \label{fig:periodicex0}
\end{figure}

\begin{ex}[{\em A Forman vector field with periodic orbit}]
\label{ex:formanperiodic-1}
{\em
An example with a more elaborate Morse decomposition can be found in
Figure~\ref{fig:periodicex0}. In this case, the underlying Lefschetz
complex is the simplicial complex indicated in the top left panel of
the figure. The multivector field~$\cV$ is actually a Forman combinatorial
vector field, which consists of six critical cells, together with eight
vectors given by doubletons. Of the critical cells, two have Morse index~$2$,
three have Morse index~$1$, while only one has index~$0$ --- where the Morse
index of a critical cell is defined as the degree of its Conley polynomial.
We would like to point out that in the case of a critical cell, this polynomial
is always a monomial. In addition to the critical cells, the edges~$\bC\bD$,
$\bD\bE$, and~$\bE\bC$ give rise to a periodic solution. Together, these
seven sets constitute the Morse sets of the minimal Morse decomposition of the vector field~$\cV$.
They are indicated in different colors in the panel on the lower left. 
In order to get the full structure of the Morse decomposition, one needs to also determine the partial
order between these Morse sets. For this, the top right panel shows five
connections originating in the two critical cells in dimension two. While
three of these connect to critical cells of dimension one, two connect to
the periodic orbit. Similarly, one can easily determine which index~$1$
critical cells connect to the periodic orbit or the vertex~$\bA$. From this,
one can readily determine the Conley-Morse graph shown in the lower right 
panel.
\exend}
\end{ex}

Intuitively, one can immediately see the connecting orbit structure 
of the above example, and therefore also its associated Morse decomposition.
In fact, it was shown in~\cite{mrozek:wanner:21a} that for any Forman
vector field on a simplicial complex~$X$ one can explicitly construct
a classical semiflow on~$X$ which exhibits the same Morse sets and
Conley-Morse graph. Conversely, it was shown in~\cite{mrozek:etal:22a}
that suitable phase space subdivisions in a classical dynamical system
combined with ideas from multivector fields can be used to rigorously
establish the existence of classical periodic solutions from the 
existence of combinatorial counterparts. These results show that the
above-mentioned intuition is more than just a coincidence.
\begin{figure}
  \begin{center}
    \includegraphics[width=0.23\textwidth]{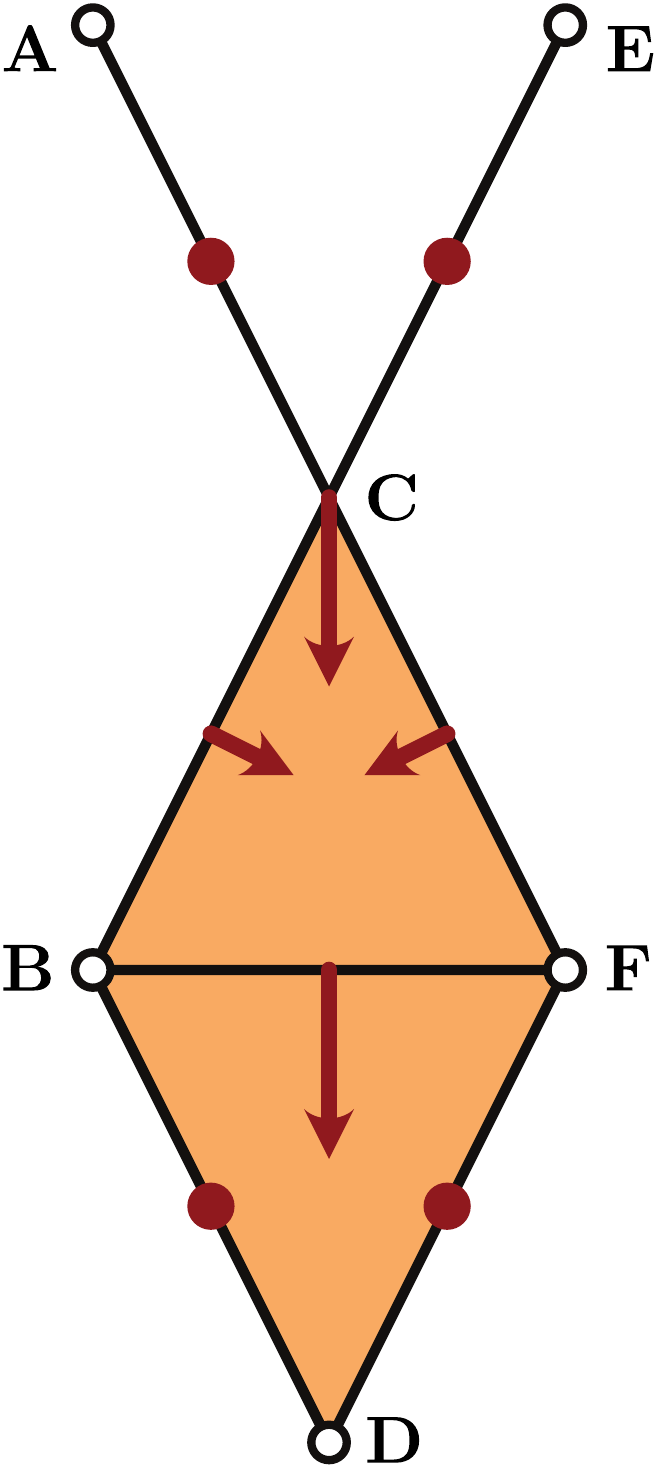}\quad\quad\quad
    \includegraphics[width=0.65\textwidth]{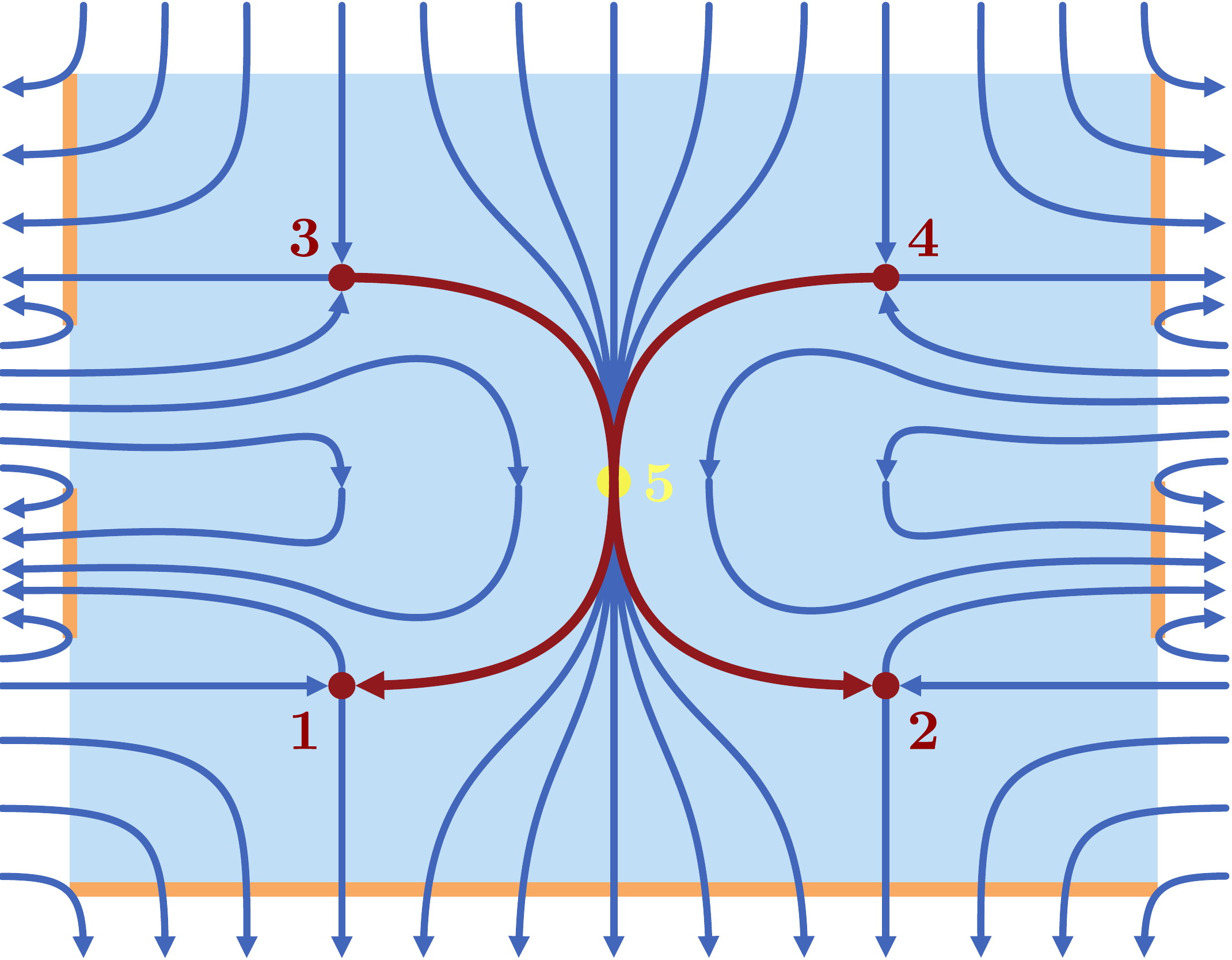}
  \end{center}
  \caption{{\em A multiflow without lattice of attractors.}
           The panel on the left shows a multivector field on a Lefschetz
           complex, which consists of the depicted simplicial complex, but
           without the vertices~$\bA$, $\bB$, $\bD$, $\bE$, and~$\bF$, shown
           as white circles. The multivector field consists of four singletons
           and one doubleton, as well as the multivector $\{ \bC, \bB\bC,
           \bC\bF, \bB\bC\bF \}$. The dynamics of this example can be 
           represented as the multiflow shown on the right, where the only
           point of forward or backward nonuniqueness is the yellow
           point~$\mathbf{5}$. For neither example one can construct a
           lattice of attracting sets.
           }
  \label{fig:multiflowex}
\end{figure}

\begin{ex}[{\em A multiflow without lattice of attractors}]
\label{ex:no-lat-hom}
{\em
As our next example, consider the combinatorial multivector field
presented in Figure~\ref{fig:multiflowex}(left). It is defined on
the Lefschetz complex~$X$ which is obtained by removing the closed
subset $\{\bA,\bB,\bD,\bE,\bF\}$ from the simplicial complex consisting
of the triangles $\bB\bC\bF$, $\bB\bD\bF$, the edges $\bA\bC$, $\bC\bE$,
and all their faces. The multivector field consists of four singleton
edges and the doubleton $\{ \bB\bF, \bB\bD\bF \}$, as well as the
multivector $\{ \bC, \bB\bC, \bC\bF, \bB\bC\bF \}$. The singletons
of the four critical edges $\bA\bC$, $\bB\bD$, $\bC\bE$, $\bD\bF$
are four isolated invariant sets. Furthermore, if we define the
poset $P=\{\xa,\xb,\xc,\xd\}$ via the Hasse diagram
\begin{equation}
\label{ex:hasseP1}
   \begin{diagram}
  \dgARROWLENGTH 1.6em
    \node{\xc}
    \arrow{s,-}
    \arrow{se,-}
    \node{\xd}
    \arrow{s,-}
    \arrow{sw,-}\\
    \node{\xa}
    \node{\xb}
  \end{diagram}
\end{equation}
and the family $\{M_p\}_{p\in P}$ with $M_\xa:=\{\bB\bD\}$, $M_\xb:=\{\bD\bF\}$,
$M_\xc:=\{\bA\bC\}$, as well as $M_\xd:=\{\bC\bE\}$, then we obtain a well-defined
Morse decomposition. Observe that all inequalities in the poset~$P$ matter, because
there are connections from~$M_\xc$ to both~$M_\xa$ and~$M_\xb$, as well as from~$M_\xd$
to both~$M_\xa$ and~$M_\xb$.
For instance,  
\[
\gamma=\overleftarrow{\bA\bC}\cdot\bC\cdot\bB\bC\bF\cdot\bB\bF\cdot \bB\bD\bF\cdot \overrightarrow{\bB\bD}
\]
is an essential solution and a connection from~$M_\xc$ to~$M_\xa$. 
A multiflow counterpart of this example is presented in Figure~\ref{fig:multiflowex}(right).
Its four saddle points are the stationary points marked in the figure 
as~$\mathbf{1}$, $\mathbf{2}$, $\mathbf{3}$, and~$\mathbf{4}$. The only
point of non-uniqueness of solutions is the yellow point marked~$\mathbf{5}$.
Similarly to the combinatorial case, there are connections through this point
from the saddle~$\mathbf{3}$ to both saddles~$\mathbf{1}$ and~$\mathbf{2}$,
as well as from saddle~$\mathbf{4}$ to both saddles~$\mathbf{1}$ and~$\mathbf{2}$.
\exend}
\end{ex} 

After this sequence of examples, we return to our general discussion of
Morse decompositions. Consider the family of {\em down sets} in the poset~$P$,
that is, the subsets $I\subset P$ such that with every $p\in I$ also all
elements below~$p$ are contained in~$I$. We denote this family as~$\Down(P)$. 
It is easy to see that~$\Down(P)$ is a {\em lattice of sets}, which means that
it is closed under union and intersection. We now associate with every
$I\in \Down(P)$ the set~$M_I$ consisting of all right endpoints of a path
with left endpoint in the union~$\bigcup_{r\in I}M_r$. One can check that
the so-defined set~$M_I$ is an attractor. However, the family of all such
attractors is not a lattice in general, because the intersection of two
attractors does not need to be an attractor itself.

In the classical setting of flows or semiflows it is possible to overcome
this difficulty by constructing an attracting neighborhood~$N_I$ for each
attractor~$M_I$ in such a way that the family $\setof{N_I\mid I\in\Down(P)}$
is again a lattice, and such that the map $I\mapsto N_I$ is a lattice
homomorphism, i.e., it preserves unions and intersections. In addition,
the attracting neighborhoods~$N_I$ have to be constructed in such a way
that~$M_I$ is the largest invariant set in~$N_I$. In other words, the
attracting set~$N_I$ determines the attractor~$M_I$. In the classical setting,
this lattice $\{ N_I \mid I\in\Down(P) \}$ is then used to proceed with the
construction of the connection matrix. However, as the following example
indicates, such a lattice and associated lattice homomorphism may not exist
in the case of a multiflow and, similarly, in the case of a combinatorial
multivector field which is inherently multivalued. 

\begin{ex}[{\em A multiflow without lattice of attractors}, continued]
\label{ex:no-lat-hom2}
{\em
Consider again the multivector field and associated multiflow shown in
Figure~\ref{fig:multiflowex}, which were already discussed in
Example~\ref{ex:no-lat-hom}. In both of these cases, the Morse
decomposition is indexed by the same poset $P=\{\xa,\xb,\xc,\xd\}$
with Hasse diagram \eqref{ex:hasseP1}. One can easily see that the
lattice of down sets of this poset is given by $\{\emptyset, \{\xa\},
\{\xb\}, \{\xa,\xb\}, \{\xa,\xb,\xc\}, \{\xa,\xb,\xd\}, P \}$.
First observe that in the combinatorial case $A:=X\setminus\{\bA\bC,\bC\bE\}$ 
does not admit an essential solution in $A$ through $\bC$, $\bB\bC$, $\bC\bF$,
and~$\bB\bC\bF$. Therefore, $A$ is not an isolated invariant set and, in
consequence, not an attractor. Notice, however, that
$M_{\{\xa,\xb,\xc\}}\cap M_{\{\xa,\xb,\xd\}}=A$. Hence, $\setof{M_I\mid I\in\Down(P)}$
is not a lattice in this case. We claim that the above-mentioned work-around for
flows with the lattice of attracting neighborhoods~$N_I$ does not work in our
combinatorial example. More precisely, there exists no lattice homomorphism
which sends each down set~$I$ to an attracting set~$N_I$ for the corresponding
attractor, and such that~$M_I$ is the largest invariant subset of~$N_I$.
To show this, assume to the contrary that $I \mapsto N_I$ is such a
homomorphism. Then one immediately obtains 
\[
  \bC\in M_{\{\xa,\xb,\xc\}} \cap M_{\{\xa,\xb,\xd\}}\subset
  N_{\{\xa,\xb,\xc\}} \cap N_{\{\xa,\xb,\xd\}} = N_{\{\xa,\xb\}} =
  N_{\{\xa\}}  \cup N_{\{\xb\}}.
\]
Suppose now that we have $\bC\in N_{\{\xa\}}$. Since
$\bC\cdot\bB\bC\bF\cdot\bB\bF\cdot\bB\bD\bF\cdot\bD\bF$ is a path
from~$\bC$ to~$\bD\bF$, and since~$N_{\{\xa\}}$ is an attracting
set, this implies the inclusion $\bD\bF\in N_{\{\xa\}}$ ---
and therefore $\{\bB\bD,\bD\bF\}$ is an invariant subset of~$N_{\{\xa\}}$.
Yet, this is not possible, since we assumed that~$M_{\{\xa\}}=\{\bB\bD\}$
is the largest invariant subset of~$N_{\{\xa\}}$. Analogously one can
rule out $\bC\in N_{\{\xb\}}$. Altogether this indeed proves that such
a lattice homomorphism does not exist. The argument for the multiflow
is similar.
\exend}
\end{ex}       

\subsection{Acyclic partitions of Lefschetz complexes}

Before we explain how the problem with lattices of attracting sets
indicated in the previous section can be addressed, we will first
consider the special case when the family of attractors
$\setof{M_I\mid I\in\Down(P)}$ is indeed a lattice and $I\mapsto M_I$
is a lattice homomorphism. One can show that in the combinatorial setting
an attractor is always an attracting set of itself. Thus, in this case
we also have a lattice of attracting sets --- and this puts us in the
setting when the approach of~\cite{HMS2021, RoSa1992} works. Actually,
in view of the closedness of attracting sets, we may consider an arbitrary
sublattice~$\cL$ of the lattice of closed sets in a Lefschetz complex~$X$.
Given a set $L\in\cL$, consider now the union of all proper subsets of~$L$
in~$\cL$, and denote it by~$L^\star$. If the inequality $L^\star\neq L$
holds, then we call~$L$ {\em join-irreducible}. Clearly, the set difference
$L^\circ:=L\setminus L^\star$ is locally closed, as a difference of two
closed sets. Moreover, one can in fact prove that the family
$\AP(\cL):=\setof{ L^\circ\mid \text{$L$ is join-irreducible}}$ is
a partition of~$X$, which in addition is {\em acyclic} in the sense
that the transitive closure of the relation $L_1^\circ\preceq L_2^\circ$
defined by $L_1^\circ\cap \cl L_2^\circ\neq\emptyset$ makes~$\AP(\cL)$
into a poset. Moreover, this resulting poset turns out to be isomorphic
to the poset of join-irreducibles of~$\cL$ ordered by inclusion. Actually,
the presented way of passing from a lattice of sets to an acyclic partition
can be viewed as a version of the celebrated Birkhoff
Theorem~\cite{Bi1937} for abstract lattices, yet for the
special case of lattices of sets.

\begin{ex}[{\em A first multivector field}, continued]
\label{ex:multivectorfield-6}
{\em
For the Morse decomposition~$\cM$ discussed in
Example~\ref{ex:multivectorfield-5} the lattice of down sets is
given by the collection $\{\emptyset,\{\xa\},\{\xa,\xb\},P\}$,
and one can easily see that the associated family of attractors
$\cA:=\{\emptyset,\{\bD\},X\setminus\{\bA\bB\bC\},X\}$ forms a
lattice. Moreover, each attractor except the empty set is
join-irreducible. Thus, in this case the corresponding acyclic
partition~$\AP(\cA)$ of~$X$ coincides with the family~$\cM$                                          
of Morse sets, and the mapping $P\ni p\mapsto M_p\in \AP(\cA)$
is indeed an order isomorphism.
\exend
}
\end{ex}

As we explain in Section~\ref{sec:alg-conn-matr} and in detail  
in Section~\ref{sec:c-matr-acyc-part}, an acyclic partition of a Lefschetz 
complex may be used to define a connection matrix. 
Unfortunately, as we indicated in the previous section, the family of
attractors $\setof{M_I\mid I\in\Down(P)}$ may not be a lattice in the
combinatorial setting. Therefore, in order to obtain the required acyclic
partition we have to modify the presented approach.  First observe that 
given an arbitrary family~$\cU$ of subsets of~$X$, there is a smallest
lattice of subsets of~$X$ containing~$\cU$. We refer to it as the
{\em lattice extension} of~$\cU$ and denote it by~$\cU'$. One can
easily verify that if~$\cU$ is a finite family of closed sets, then so
is its lattice extension. Since, clearly, the union and intersection of
attracting sets is again an attracting set, the lattice extension~$\cA'$
of the family of attractors given by $\cA:=\setof{M_I\mid I\in\Down(P)}$
is a lattice of attracting sets. Hence, we can in fact construct a
lattice, but it is no longer indexed by the down sets of the original poset~$P$.
Nevertheless, one can show that the acyclic partition~$\cD:=\AP(\cA')$,
considered as a poset, is order isomorphic to an extension~$\hat{P}$ of
the original poset~$P$. This enables us to write the family~$\cD$ in the
form $\cD=\{D_p\}_{p\in\hat{P}}$ where $\hat{P}\ni p\mapsto D_p\in\cD$ is
the mentioned order isomorphism. Note that every element $D_p$ for
$p\in\hat{P}$, as a locally closed subset of the Lefschetz complex~$X$,
is itself a Lefschetz complex. Moreover, one can deduce from the
definition of Morse decomposition that $H(D_p)=0$ for all 
$p\in\hat{P}\setminus P$. This observation will then let us get rid
of the elements in the set difference~$\hat{P}\setminus P$ via an
equivalence relation introduced in the final, algebraic step of
our construction. Moreover, it turns out that instead of searching for the lattice
extension $\cA'$ and the associated partition, we can just take a refined acyclic 
partition $\cD$ consisting of all Morse sets together with all the multivectors not 
contained in a Morse set. 
It is not difficult to see that such multivectors
must be regular. Therefore, although this may lead to an acyclic partition 
indexed by a larger poset, the partition elements are explicitly given
and we may get rid of the elements in ~$\hat{P}$ coming from the regular
multivectors outside Morse sets by the the same algebraic equivalence.
In fact, this straightforward construction replaces the laborious and generally not possible 
step (ii) of the connection matrix pipeline discussed in the introduction. 

\begin{ex}[{\em A multiflow without lattice of attractors}, continued]
\label{ex:no-lat-hom3}
{\em
As we pointed out in Example~\ref{ex:no-lat-hom2} the intersection
$A=M_{\{\xa,\xb,\xc\}}\cap M_{\{\xa,\xb,\xd\}}$ is not an attractor
which implies that the family of attractors is not a lattice.
Actually, in this example $A$ is the only set in the lattice
extension~$\cA'$ of the family of attractors $\cA$ which is not
an attractor. The set~$A$ is join-irreducible in this lattice, 
therefore it provides 
\[
A^\circ=\{\bC,\bB\bC,\bC\bF,\bB\bF,\bB\bC\bF,\bB\bD\bF\}
\]
to the associated acyclic partition $\AP(\cA')$. However, $A^\circ$
splits as the disjoint union of two multivectors
$V_1:=\{\bC,\bB\bC,\bC\bF,\bB\bC\bF\}$ and $V_2:=\{\bB\bF,\bB\bD\bF\}$.
Therefore, we may take as the acyclic partition family
\[
\cE:=\{\,\{\bA\bC\},\{\bC\bE\},\{\bB\bD\},\{\bD\bF\},V_1,V_2\,\} , 
\]
which consists of the four singleton Morse sets and the two regular
multivectors~$V_1$ and~$V_2$.
To write it as an indexed family $\{E_p\}_{p\in \hat{P}}$, we take the
extended poset $\hat{P}$ as the union $P\cup\{\xe,\xf\}$ 
with partial order given by the Hasse diagram
\[
   \begin{diagram}
  \dgARROWLENGTH 1.0em
    \node{\xc}
    \arrow{se,-}
    \node[2]{\xd}
    \arrow{sw,-}\\
    \node[2]{\xe}
    \arrow{s,-}\\
    \node[2]{\xf}
    \arrow{se,-}
    \arrow{sw,-}\\
    \node{\xa}
    \node[2]{\xb}
  \end{diagram}.
\]
Moreover, we define~$E_p$ as the Morse set~$M_p$ for $p\in P$, and
additionally let~$E_\xe:=V_1$ and~$E_\xf:=V_2$.
\exend
}
\end{ex}       

\subsection{Algebraic connection matrices} 
\label{sec:alg-conn-matr}
The acyclic partition~$\{D_p\}_{p\in \hat{P}}$ of the Lefschetz complex~$X$
discussed in the previous section lets us decompose the associated chain
complex~$C(X)$ as a direct sum $C(X)=\bigoplus_{p\in \hat{P}}C(D_p)$.
It turns out that the order isomorphism $p\mapsto D_p$ gives this direct
sum the special structure of a poset filtered chain complex, which in
turn leads directly to a purely algebraic concept of connection matrix.
This will be explained in more detail in the following.
 
We first explain the case when $\hat{P}=P$ is fixed. In terms of
applications to dynamics this corresponds to the situation when we have
a lattice of attracting sets indexed by down sets of the original poset
in the Morse decomposition. Recall that we assume field coefficients for
all considered modules, in particular chain complexes and homology modules.
Given a chain complex $(C,d)$ together with a direct sum decomposition 
\begin{equation}
\label{eq:C-bigoplus-Cp}
C=\bigoplus_{p\in P} C_p
\end{equation}
and a down set $I\in\Down(P)$ we introduce the abbreviation $C_I:=\bigoplus_{p\in I} C_p$.
We say that the boundary homomorphism $d$ is {\em filtered} if $d(C_I)\subset C_I$
is satisfied for every $I\in\Down(P)$. If this is the case, then the chain complex~$(C,d)$
together with the decomposition \eqref{eq:C-bigoplus-Cp} is called a $P$-{\em filtered
chain complex}. Similarly, any homomorphism, and in particular, any chain map $ h:C\to C'$
between two $P$-filtered chain complexes, is called {\em filtered} if 
\begin{equation}
\label{eq:varphi-C-I}
 h(C_I)\subset C'_I
 \quad\text{ for every }\quad
 I\in\Down(P).
\end{equation}
We would like to point out that every acyclic partition~$\cE$ of a Lefschetz
complex~$X$ makes~$C(X)$, the chain complex of~$X$, a filtered chain complex
via the decomposition
\[
   C(X)=\bigoplus_{E\in\cE}C(E).
\]
The decomposition is well-defined, because each $E\in\cE$, as a locally closed
subset of~$X$, is itself a Lefschetz complex, and the fact that the boundary
homomorphism is filtered may be concluded from the assumption that the
partition is acyclic, see also Proposition~\ref{prop:Lefschetz-filtered-chain-complex}.

Now, two $P$-filtered chain complexes are {\em filtered chain homotopic} if
there exist filtered chain maps $h:C\to C'$ and $h':C'\to C$ such that the
composition~$h' h$ is filtered chain homotopic to~$\id_C$, and $hh'$ is
filtered chain homotopic to~$\id_{C'}$. In this context, a filtered chain
homotopy is a chain homotopy which is itself filtered as a homomorphism.
Robbin and Salamon~\cite{RoSa1992} prove that every $P$-filtered chain
complex is $P$-filtered chain homotopic to a {\em reduced} filtered
chain complex, that is, a filtered chain complex~$(C',d')$ such that
$d'_{pp}=0$ for all $p\in P$. Moreover, Harker, Mischaikow, and
Spendlove~\cite{HMS2021} prove that if~$(C,d)$ is also filtered chain
homotopic to another reduced filtered complex~$(C'',d'')$, then the
filtered chain complexes~$(C',d')$ and~$(C'',d'')$ are in fact filtered
isomorphic. Hence, up to filtered chain isomorphism every filtered chain
complex has exactly one reduced representative. By definition, this
representative is called its {\em Conley complex}, and the matrix of the
boundary operator of the Conley complex is referred to as its
{\em connection matrix}.

\begin{ex}[{\em An algebraic example}]
\label{ex:alg-conn-matrix}
{\em
   Consider free module $C=R\spn{X}$ spanned by the set of symbols 
\[
   X:=\{\bA,\bB,\ba,\bb,\bc,\balpha\}
\]
and with partition $X=X_0\cup X_1\cup X_2$ given by 
\[
   X_0:=\{\bA,\bB\},\quad X_1:=\{\ba,\bb,\bc\}, \quad X_2:=\{\balpha\}.
\]
This partition lets us treat the free module~$C$ as a $\ZZ$-graded module
with gradation $C=\bigoplus_{i=0}^2 R\spn{X_i}$. In order to make~$C$ a free chain
complex we assume that~$R$ is the field~$\ZZ_2$ and consider a homomorphism
$d:C\to C$ defined on the basis~$X$ by the matrix
\[
\begin{array}{c||cccc|c|c|}
       d  &  \bA &  \bB &    \ba &  \bc &   \bb &    \balpha      \\
  \hline
  \hline
  \bA     &      &      &      1 &    1 &       1 &               \\
  \bB     &      &      &      1 &    1 &       1 &               \\
  \ba     &      &      &        &      &         &           1   \\
  \bc     &      &      &        &      &         &               \\
  \hline
  \bb     &      &      &        &      &         &           1   \\
  \hline                                                     
  \balpha &      &      &        &      &         &               \\
  \hline
\end{array}\;.
\]
It is not difficult to check that $(C,d)$ is indeed a free chain complex.
Furthermore, it can be made into a poset filtered chain complex by
considering the poset $P:=\{\xa,\xb,\xc\}$ linearly ordered by
$\xa<\xb<\xc$ and setting
\[
C_\xa:=R\spn{\{ \bA,\bB,\ba,\bc\}}, \quad C_\xb:=R\spn{\{ \bb\}}, \quad C_\xc:=R\spn{\{ \balpha\}}.
\]
Again, one can easily check that the gradation $C=\bigoplus_{p\in P}C_p$
turns~$(C,d)$ into a $P$-filtered chain complex. Since we have $d_{\xa\xa}\neq 0$,
this complex~$C$ is not reduced. 

One can also consider a submodule~$\bar{C}$ of~$C$ which is obtained by
removing the generators~$\bB$ and~$\bc$, i.e., we have
\[
\bar{C}_\xa:=R\spn{\{ \bA,\ba\}}, \quad
\bar{C}_\xb:=R\spn{\{ \bb\}}, \quad
\bar{C}_\xc:=R\spn{\{ \balpha\}},
\]
together with the boundary homomorphism $\bar{d}:\bar{C}\to \bar{C}$
given by the matrix
\[
\begin{array}{c||cc|c|c|}
  \bar{d} &  \bA &  \ba &    \bb &    \balpha   \\
  \hline
  \hline
  \bA     &      &      &        &              \\
  \ba     &      &      &        &          1   \\
  \hline
  \bb     &      &      &        &          1   \\
  \hline                                                   
  \balpha &      &      &        &              \\
  \hline
\end{array}\;.
\]
We leave it to the reader to verify that this indeed defines
a $P$-filtered complex~$(\bar{C},\bar{d})$. Furthermore,
this chain complex is reduced and, in fact, is a Conley complex of~$(C,d)$
which makes the matrix of~$\bar{d}$ a connection matrix of~$(C,d)$.
To see this in more detail, consider the two specific homomorphisms
$h:C\to\bar{C}$ and $g:\bar{C}\to C$ given by the matrices
\[
\begin{array}{c||cccc|c|c|}
        h &  \bA &  \bB &    \ba &  \bc &   \bb &    \balpha   \\
  \hline
  \hline
  \bA     &    1 &    1 &        &      &         &               \\
  \ba     &      &      &      1 &      &         &               \\
  \hline
  \bb     &      &      &        &      &      1  &               \\
  \hline                                                     
  \balpha &      &      &        &      &         &         1     \\
  \hline
\end{array}
\qquad\text{ and }\qquad
\begin{array}{c||cc|c|c|}
       g  &  \bA &   \ba &   \bb &    \balpha   \\
  \hline
  \hline
  \bA     &    1 &      &        &              \\
  \bB     &      &      &        &              \\
  \ba     &      &    1 &        &              \\
  \bc     &      &    1 &      1 &              \\
  \hline
  \bb     &      &      &      1 &              \\
  \hline                                                     
  \balpha &      &      &        &          1   \\
  \hline
\end{array}\;.
\]
One can immediately verify that both~$g$ and~$h$ are indeed $P$-filtered chain maps,
and that the identity $h \circ g = \id_{\bar{C}}$ is satisfied. Moreover, the composition~$g \circ h$
can be computed as the matrix 
\[
\begin{array}{c||cccc|c|c|}
g \circ h &  \bA &  \bB &    \ba &  \bc &   \bb &    \balpha   \\
  \hline
  \hline
  \bA     &    1 &    1 &        &      &         &               \\
  \bB     &      &      &        &      &         &               \\
  \ba     &      &      &      1 &      &         &               \\
  \bc     &      &      &      1 &      &       1 &               \\
  \hline
  \bb     &      &      &        &      &       1 &               \\
  \hline                                                     
  \balpha &      &      &        &      &         &            1  \\
  \hline
\end{array}\;.
\]
Therefore, this composition~$g \circ h$ is filtered chain homotopic
to~$\id_C$ via the $P$-filtered chain homotopy $\gamma:C\to C$
which sends all generators to zero, except~$\bB$ which is sent to~$\bc$.
In fact, this chain homotopy is an example of an elementary reduction
via the pair of generators~$(\bB,\bc)$ discussed in~\cite[Section~4.3]{KaMiMr2004},
see also \cite{KMS1998}. We presented this example purely algebraically, but
in fact its algebra provides the connection matrix for the Morse decomposition
of a combinatorial multivector field which we will discuss in more detail
in Example~\ref{ex:nonuniqueper-1} below.
\exend
}
\end{ex}

As we mentioned already earlier, one has to modify this definition to accommodate the
situation when we cannot keep the original poset~$P$ as in
Examples~\ref{ex:no-lat-hom},~\ref{ex:no-lat-hom2}, and~\ref{ex:no-lat-hom3}.
Under a changing poset the
condition~\eqref{eq:varphi-C-I} in the definition of a filtered
homomorphism needs to be modified. For this, let~$P$ and~$P'$ be two
posets and consider a $P$-filtered chain complex~$(C,d)$, as well as
a $P'$-filtered chain complex~$(C',d')$. In order to speak about a
filtered homomorphism in this setting we need a way to relate down
sets in the two posets~$P$ and~$P'$. Note that if $\alpha:P'\to P$
is order preserving, then we have $\alpha^{-1}(I) \in \Down(P')$
for every $I\in\Down(P)$. Therefore, we can define a morphism 
from a $P$-filtered chain complex $(C,d)$ to a  $P'$-filtered
chain complex~$(C',d')$ as a pair~($\alpha,h)$ where $\alpha:P'\to P$ 
is order preserving and $h:C\to C'$ is a chain map satisfying 
\begin{equation}
\label{eq:varphi-C-I-prime}
 h(C_I) \subset C'_{\alpha^{-1}(I)}
 \quad\text{ for every }\quad
 I\in\Down(P).
\end{equation}
A similar modification allows us to also extend the definition of
filtered chain homotopy to the setting of varying posets.

We still need some further modifications to guarantee that in
this algebraic step we can eliminate the elements added to the
poset to obtain a lattice of attracting sets. For this, recall
that subcomplexes~$C_p$ for the added values of~$p$ are Lefschetz complexes
with zero homology or, equivalently, they are chain homotopic to zero.
Therefore, if~$\bar{C}$ denotes the Conley complex of~$C$ computed
under a fixed poset~$P$, the chain groups~$\bar{C}_p$ for an added
element~$p$ are zero. In consequence, the respective rows and columns
in the connection matrix are empty, because the only basis of the zero
group is empty. To formalize the removal of these empty rows and columns
we do four things:
\begin{itemize}
\item We add a distinguished subset $P_\star\subset P$ and require that
all chain groups~$C_p$ for $p\in P\setminus P_\star$ are chain homotopic
to zero.
\item We extend the definition of a reduced chain complex by additionally
requesting that $P_\star = P$ in a reduced complex.
\item We allow partial order preserving maps $\alpha:P'\pto P$ to relate
down sets in the posets~$P$ and~$P'$, but we require that~$\alpha$ is
defined at least for all $p' \in P'_\star$. The partial map~$\alpha$ does
not need to be defined for elements $p'\in P'\setminus P'_\star$, since a 
term~$\bar{C}'_p=0$ in the Conley complex contributes nothing to the whole
Conley complex.
\item Finally, since the preimage $\alpha^{-1}(I)$ for $I\in\Down(P)$ no
longer needs to be a down set under a partially defined order preserving
map, we also have to replace~\eqref{eq:varphi-C-I-prime} by the condition
\begin{equation}
\label{eq:varphi-C-I-prime-leq}
  h(C_I) \subset C'_{\alpha^{-1}(I)^{\leq}}
  \quad\text{ for every }\quad
  I\in\Down(P),
\end{equation}
where~$A^\leq$ denotes the smallest down set containing~$A$, for $A\subset P'$.
We refer to chain maps satisfying \eqref{eq:varphi-C-I-prime-leq} as
{\em $\alpha$-filtered}.\footnote{In fact, when we formally define this
concept later, we instead use the definition~\eqref{eq:filt-hom} for
technical reasons. Its equivalence with~\eqref{eq:varphi-C-I-prime-leq}
is established in Proposition~\ref{prop:filt-hom}.}
\end{itemize}

Summarizing, the above extensions lead to a well-defined category~$\PfCC$
whose objects are of the form~$(P,C,d)$, where~$P$ is a poset with a
distinguished subset~$P_\star$ and~$(C,d)$ is a $P$-filtered chain complex.
In addition, morphisms from~$(P,C,d)$ to~$(P',C',d')$ are of the form~$(\alpha,h)$,
where $\alpha:P'\pto P$ is a partial order preserving map whose domain
contains~$P'_\star$ and which satisfies $\alpha(P'_\star) \subset P_\star$,
and $h: C \to C'$ is a chain map satisfying~\eqref{eq:varphi-C-I-prime-leq}.
As it turns out, the main results of~\cite{HMS2021, RoSa1992} can be extended
to this new category. More precisely, we prove later in the paper the following
result.
\begin{thm}[Existence of the Conley complex and connection matrix]
\label{thm:PfCC}
Consider the category~$\PfCC$ of poset filtered chain complexes with varying
posets introduced above. Then the following hold:
\begin{itemize}
\item[(i)]  Every object in~$\PfCC$ is filtered chain homotopic to a reduced object. 
\item[(ii)] If two reduced objects in~$\PfCC$ are filtered chain homotopic, then
they automatically are isomorphic in~$\PfCC$.
\end{itemize}
In other words, up to isomorphism every object in~$\PfCC$ has exactly one reduced
representative, which is called its {\em Conley complex}. The matrix of the
boundary operator of the Conley complex is called {\em connection matrix}.
\end{thm}
The above result provides both the Conley complex and the connection matrix also
in the case when there is no lattice homomorphism from down sets in the poset of
a Morse decomposition to attracting sets in this Morse decomposition, as discussed
earlier. A bonus, which comes as a side effect, is that we can take a shortcut
in the connection matrix pipeline by skipping the construction of the lattice
of attracting neighborhoods and passing immediately from the Morse decomposition
in step~(i) to a filtered chain complex in step (iii) via
an acyclic partition of the phase space associated with the Morse decomposition
under the extended poset.
\begin{ex}[{\em A multiflow without lattice of attractors}, continued]
\label{ex:no-lat-hom4}
{\em
Consider the acyclic partition $\cE=\{E_p\}_{p\in\hat{P}}$ introduced in
Example~\ref{ex:no-lat-hom3}, which corresponds to the Morse decomposition
of the multivector field on the Lefschetz complex $X$ in Figure~\ref{fig:multiflowex}(left).
The boundary homomorphism $d$ of the associated chain complex $C(X)$ of~$X$ 
has the matrix 
{\small
\[
\begin{array}{c||c|c|cc|cccc|c|c||}
       d    &  \bB\bD &  \bD\bF &    \bB\bF & \bB\bD\bF &  \bC & \bB\bC & \bC\bF & \bB\bC\bF &   \bA\bC & \bC\bE   \\
\hline
\hline
  \bB\bD    &         &         &           &         1 &      &        &        &           &          &          \\ 
\hline
  \bD\bF    &         &         &           &         1 &      &        &        &           &          &          \\ 
\hline
  \bB\bF    &         &         &           &         1 &      &        &        &         1 &          &          \\ 
  \bB\bD\bF &         &         &           &           &      &        &        &           &          &          \\ 
\hline
  \bC       &         &         &           &           &      &      1 &      1 &           &        1 &        1 \\ 
  \bB\bC    &         &         &           &           &      &        &        &         1 &          &          \\ 
  \bC\bF    &         &         &           &           &      &        &        &         1 &          &          \\ 
  \bB\bC\bF &         &         &           &           &      &        &        &           &          &          \\ 
\hline
  \bA\bC    &         &         &           &           &      &        &        &           &          &          \\ 
\hline
  \bC\bE    &         &         &           &           &      &        &        &           &          &          \\ 
\hline
\end{array}\;.
\]
}
One can immediately verify that $(C(X),d)$ gives a $\hat{P}$-filtered chain complex.
For $p\in \hat{P}\setminus P=\{\xe,\xf\}$ the induced Lefschetz complex $X_p$ is chain
homotopic to zero. Therefore, we take $\hat{P}_*:=P$ as the distinguished subset
in~$\hat{P}$. This way we obtain an object  $(\hat{P},C(X),d)$ of $\PfCC$.
We note that $d$ is not reduced, because $d_{\xe\xe}\neq 0$, $d_{\xf\xf}\neq 0$
and, additionally, $C(X_p)$ is chain homotopic to zero for $p \in\hat{P}\setminus P$.

Consider now the free module $R\spn{\bB\bD,\bD\bF,\bA\bC,\bC\bE}$ and turn it into a
chain complex~$\bar{C}$ by assuming the boundary map to be zero. Clearly, the definitions 
$\bar{C}_\xa:=\spn{\bB\bD}$, $\bar{C}_\xb:=\spn{\bD\bF}$, $\bar{C}_\xc:=\spn{\bA\bC}$,
and $\bar{C}_\xd:=\spn{\bC\bE}$ render~$\bar{C}$ a $P$-filtered chain complex. After
finally letting $\bar{P}:=\bar{P}_\star:=P$, we now claim that $(\bar{P},\bar{C},0)$
is a Conley complex of the Morse decomposition associated with the multivector field  
on the Lefschetz complex $X$ in Figure~\ref{fig:multiflowex}(left). In order to see
this, let $\alpha: \bar{P} \hookrightarrow \hat{P}$ denote the inclusion map and let
$\beta:\hat{P}\pto \bar{P}$ be the partial map which, as a relation, is the inverse
of~$\alpha$. Consider the homomorphisms $h: C(X) \to \bar{C}$ and $g: \bar{C}\to C(X)$
given by the matrices
{\small
\[
\begin{array}{c||c|c|cc|cccc|c|c||}
       h    &  \bB\bD &  \bD\bF &    \bB\bF & \bB\bD\bF &  \bC & \bB\bC & \bC\bF & \bB\bC\bF &   \bA\bC & \bC\bE   \\
\hline
\hline
  \bB\bD    &       1 &         &         1 &           &      &        &      1 &           &          &          \\ 
\hline
  \bD\bF    &         &       1 &         1 &           &      &        &      1 &           &          &          \\ 
\hline
  \bA\bC    &         &         &           &           &      &        &        &           &        1 &          \\ 
\hline
  \bC\bE    &         &         &           &           &      &        &        &           &          &      1   \\ 
\hline
\end{array}
\]
}
and
{\small
\[
\begin{array}{c||c|c|c|c||}
       g    &  \bB\bD &  \bD\bF &   \bA\bC & \bC\bE   \\
\hline
\hline
  \bB\bD    &       1 &         &          &          \\ 
\hline
  \bD\bF    &         &       1 &          &          \\ 
\hline
  \bB\bF    &         &         &          &          \\ 
  \bB\bD\bF &         &         &          &          \\ 
\hline
  \bC       &         &         &          &          \\ 
  \bB\bC    &         &         &        1 &       1  \\ 
  \bC\bF    &         &         &          &          \\ 
  \bB\bC\bF &         &         &          &          \\ 
\hline
  \bA\bC    &         &         &        1 &          \\ 
\hline
  \bC\bE    &         &         &          &       1  \\ 
\hline
\end{array}\;.
\] 
}
One can immediately verify that $h$ is an $\alpha$-filtered chain map,
and that~$g$ is a $\beta$-filtered chain map. Hence, both $(\alpha,h)$
and $(\beta,g)$ are morphisms in $\PfCC$. It turns out that $(\alpha,h)
\circ (\beta,g) = \id_{(\bar{P},\bar{C},0)}$, and that $(\beta,g)\circ
(\alpha,h)$ is filtered chain homotopic to $\id_{(\hat{P},C(X),d)}$ via
a $\hat{P}$-filtered homotopy $\gamma:C(X)\to C(X)$ defined by the matrix
{\small
\[
\begin{array}{c||c|c|cc|cccc|c|c||}
  \gamma    &  \bB\bD &  \bD\bF &    \bB\bF & \bB\bD\bF &  \bC & \bB\bC & \bC\bF & \bB\bC\bF &   \bA\bC & \bC\bE   \\
\hline
\hline
  \bB\bD    &         &         &           &           &      &        &        &           &          &          \\ 
\hline
  \bD\bF    &         &         &           &           &      &        &        &           &          &          \\ 
\hline
  \bB\bF    &         &         &           &           &      &        &        &           &          &          \\ 
  \bB\bD\bF &         &         &         1 &           &      &        &      1 &           &          &          \\ 
\hline
  \bC       &         &         &           &           &      &        &        &           &          &          \\ 
  \bB\bC    &         &         &           &           &    1 &        &        &           &          &          \\ 
  \bC\bF    &         &         &           &           &      &        &        &           &          &          \\ 
  \bB\bC\bF &         &         &           &           &      &        &      1 &           &          &          \\ 
\hline
  \bA\bC    &         &         &           &           &      &        &        &           &          &          \\ 
\hline
  \bC\bE    &         &         &           &           &      &        &        &           &          &          \\ 
\hline
\end{array}\;.
\]
}%
This confirms that $(\bar{P},\bar{C},0)$ is a Conley complex of the Morse
decomposition associated with the multivector field on the Lefschetz
complex $X$ in Figure~\ref{fig:multiflowex}(left).
\exend
}
\end{ex}

\begin{figure}
  \begin{center}
    \includegraphics[width=0.99\textwidth]{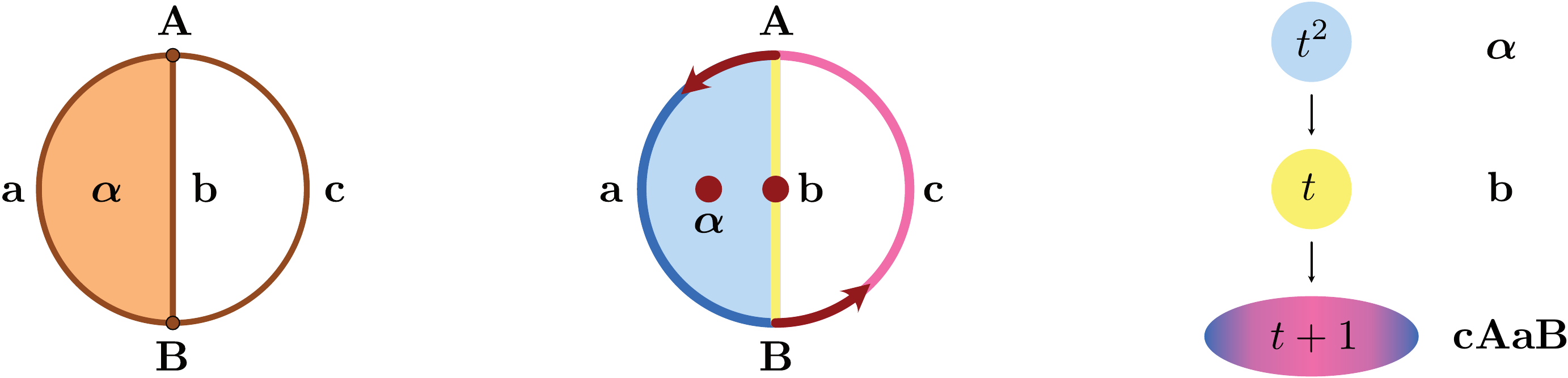}
  \end{center}
  \caption{{\em Small Lefschetz complex with periodic orbit}.
           The left panel shows a small Lefschetz complex which consists
           of a $2$-cell~$\mathbf{\balpha}$, three $1$-cells~$\mathbf{a}$,
           $\mathbf{b}$, $\mathbf{c}$, and two $0$-cells~$\mathbf{A}$,
           $\mathbf{B}$. On this complex, we study the combinatorial 
           vector field shown in the middle panel, which consists of
           two singletons and two doubletons. The associated Conley-Morse
           graph is shown on the right, with Morse sets given by the
           critical cells~$\{ \balpha \}$ and~$\{ \mathbf{b} \}$, as
           well as the periodic orbit~$\{ \mathbf{c}, \mathbf{A},
           \mathbf{a}, \mathbf{B} \}$.  
           }
  \label{fig:nonuniqueper1}
\end{figure}
\begin{figure}
  \begin{center}
    \includegraphics[width=0.99\textwidth]{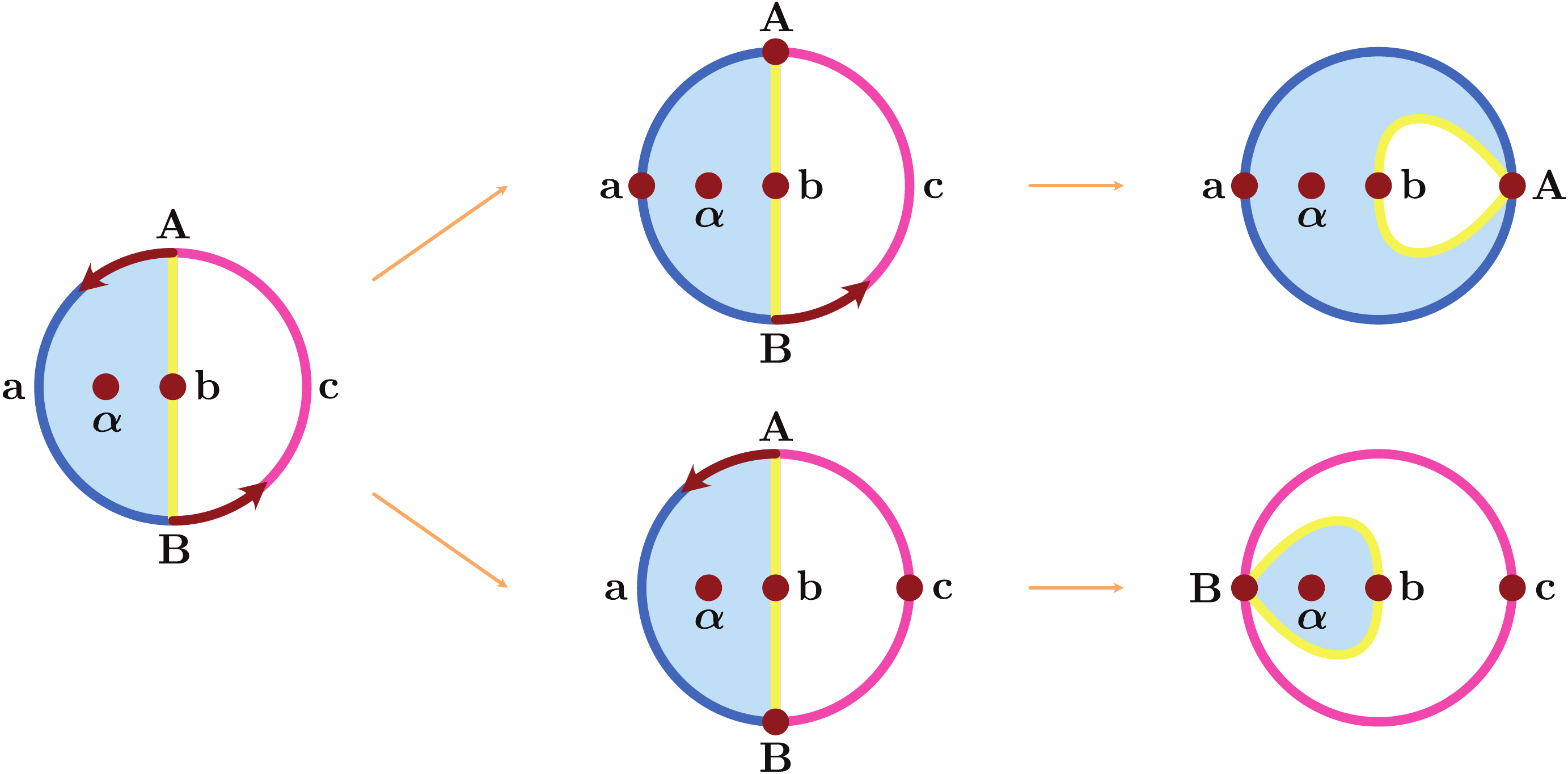}
  \end{center}
  \caption{{\em Small Lefschetz complex with periodic orbit}.
           For the combinatorial vector field shown in
           Figure~\ref{fig:nonuniqueper1}, the above diagram sketches
           the computation of two different connection matrices, based 
           on two different elementary chain complex reductions.
           The top arrow sequence uses the reduction pair~$(B,c)$,
           while the bottom sequence is for~$(A,a)$. In both cases,
           the final chain complex is reduced and boundaryless.
           }
  \label{fig:nonuniqueper2}
\end{figure}

We have already presented one example of a connection matrix in
the context of Example~\ref{ex:multivectorfield-5}. In that case,
it could easily be derived from the fact that the homology of the
Conley complex is isomorphic to the homology of the underlying
Lefschetz complex. A more complicated example will further 
illustrate Theorem~\ref{thm:PfCC}.

\begin{ex}[{\em Small Lefschetz complex with periodic orbit}]
\label{ex:nonuniqueper-1}
{\em
As a more elaborate example, yet one that still can be discussed
directly in detail, consider the Lefschetz complex~$X$ shown in the
left panel of Figure~\ref{fig:nonuniqueper1}. This complex consists
of a semi-circle shaped two-dimensional cell~$\balpha$, whose
boundary consists of the vertices~$\bA$ and~$\bB$, together with
the two $1$-cells~$\mathbf{a}$ and~$\mathbf{b}$. In addition, the
complex contains a third $1$-cell~$\mathbf{c}$ which joins the two
vertices. On~$X$ we consider the combinatorial vector field~$\cV$
sketched in the middle panel of Figure~\ref{fig:nonuniqueper1},
and which consists of two singletons and two doubletons. The
associated Conley-Morse graph is shown on the right, with Morse
sets given by the critical cells~$\{ \balpha \}$ and~$\{ \mathbf{b} \}$,
as well as the periodic orbit~$\{ \mathbf{c}, \mathbf{A}, \mathbf{a},
\mathbf{B} \}$.
\begin{table}
  \begin{center}
  \begin{tabular}{c||cc|c|c|}
      & $\mathbf{A}$ & $\mathbf{a}$ & $\mathbf{b}$ & $\balpha$ \\ \hline\hline
      $\mathbf{A}$ &   & 0 & 0 &   \\ 
      $\mathbf{a}$ &   &   &   & 1 \\ \hline
      $\mathbf{b}$ &   &   &   & 1 \\ \hline
      $\balpha$    &   &   &   &   \\ \hline
  \end{tabular}
  \hspace*{2cm}
  \begin{tabular}{c||cc|c|c|}
      & $\mathbf{B}$ & $\mathbf{c}$ & $\mathbf{b}$ & $\balpha$ \\ \hline\hline
      $\mathbf{B}$ &   & 0 & 0 &   \\ 
      $\mathbf{c}$ &   &   &   & 0 \\ \hline
      $\mathbf{b}$ &   &   &   & 1 \\ \hline
      $\balpha$    &   &   &   &   \\ \hline
  \end{tabular} \\[2ex]
  \end{center}
  \caption{{\em Small Lefschetz complex with periodic orbit}.
           Two connection matrices for the combinatorial vector
           field from Figure~\ref{fig:nonuniqueper1}. The matrix
           on the left is the result of the reduction process 
           indicated via the top arrows in Figure~\ref{fig:nonuniqueper2},
           while the bottom arrows lead to the matrix shown on the right.
           }
  \label{table:nonuniqueper}
\end{table}

Due to the small size of the Lefschetz complex~$X$, one can immediately
reduce the associated chain complex via elementary reduction pairs.
For example, if one uses the reduction pair~$(B,c)$, one obtains the
sequence of Lefschetz complexes shown in the top arrow sequence
in Figure~\ref{fig:nonuniqueper2}, and this leads to the connection
matrix shown in Table~\ref{table:nonuniqueper}(left). We discussed algebraic details 
of the computation of this connection matrix in Example~\ref{ex:alg-conn-matrix}.
In contrast,
if one uses the reduction pair~$(A,a)$, then the sequence of Lefschetz
complexes is as in the bottom arrow sequence of Figure~\ref{fig:nonuniqueper2},
and one obtains the connection matrix in Table~\ref{table:nonuniqueper}(right).
The algebraic details of this computation are similar to those presented in 
Example~\ref{ex:alg-conn-matrix}.
\exend
}
\end{ex}
%

\subsection{Uniqueness of connection matrices}
Theorem~\ref{thm:PfCC}(ii) implies that connection matrices are uniquely
determined up to an isomorphism in~$\PfCC$. While this statement is clearly
correct, it does not convey the complete story. This statement only means
that any two connection matrices of a given filtered chain complex are similar
via the matrix of a filtered chain map. However, there is also a stronger
equivalence relation between connection matrices, which is important in its
own right.

To explain this, recall that every Conley complex of a given filtered chain
complex, considered as an object of~$\PfCC$, may be considered together with
the filtered chain maps establishing the filtered chain homotopy. The
composition of these filtered chain maps between the two Conley complexes of
a given filtered chain complex is a filtered chain map. Nevertheless, this
composition may or may not be a {\em graded chain map\/}, which we define as a
chain map $h:C_p\to C'_p$ such that $h(C_p) \subset C^\prime_{\alpha^{-1}(p)}$.
This leads to a stronger equivalence relation between connection matrices
which lets us speak about uniqueness and nonuniqueness of connection matrices. 
In particular, one of the main results of the paper is the following theorem,
whose precise formulation can be found in Theorem~\ref{thm:cm-for-forman-gradient}.
\begin{thm}[Unique connection matrix for gradient vector fields]
\label{thm:uniqueness}
   Assume that~$\cV$ is a gradient combinatorial vector field on a regular
   Lefschetz complex~$X$. Then the Morse decomposition consisting of all the critical
   cells of~$\cV$ has precisely one connection matrix. It coincides with the
   matrix of the boundary operator of the associated Conley complex.
\end{thm}
In addition to being unique, the connection matrix of a gradient combinatorial
vector field also allows us to establish connections between Morse sets, as the
following result shows. Its precise formulation is the subject of
Theorem~\ref{thm:co-for-forman-gradient}, and this result also explains how the
connections can be found explicitly.
\begin{thm}[Existence of connecting orbits]
\label{thm:connectionex}
   Assume that~$\cV$ is a gradient combinatorial vector field on a regular
   Lefschetz complex~$X$, and let~$M_p$ and~$M_q$ with $p < q$ denote two
   critical cells of~$\cV$ such that the entry in the connection matrix
   corresponding to these cells is non-zero. Then there exists a connection
   between~$M_q$ and~$M_p$.
\end{thm}
In fact, the above result under slightly stronger assumptions 
remains true for general multivector fields,
see Theorem~\ref{thm:h-conn}. We illustrate these two results in the
following example.
\begin{figure}
  \begin{center}
    \includegraphics[width=0.45\textwidth]{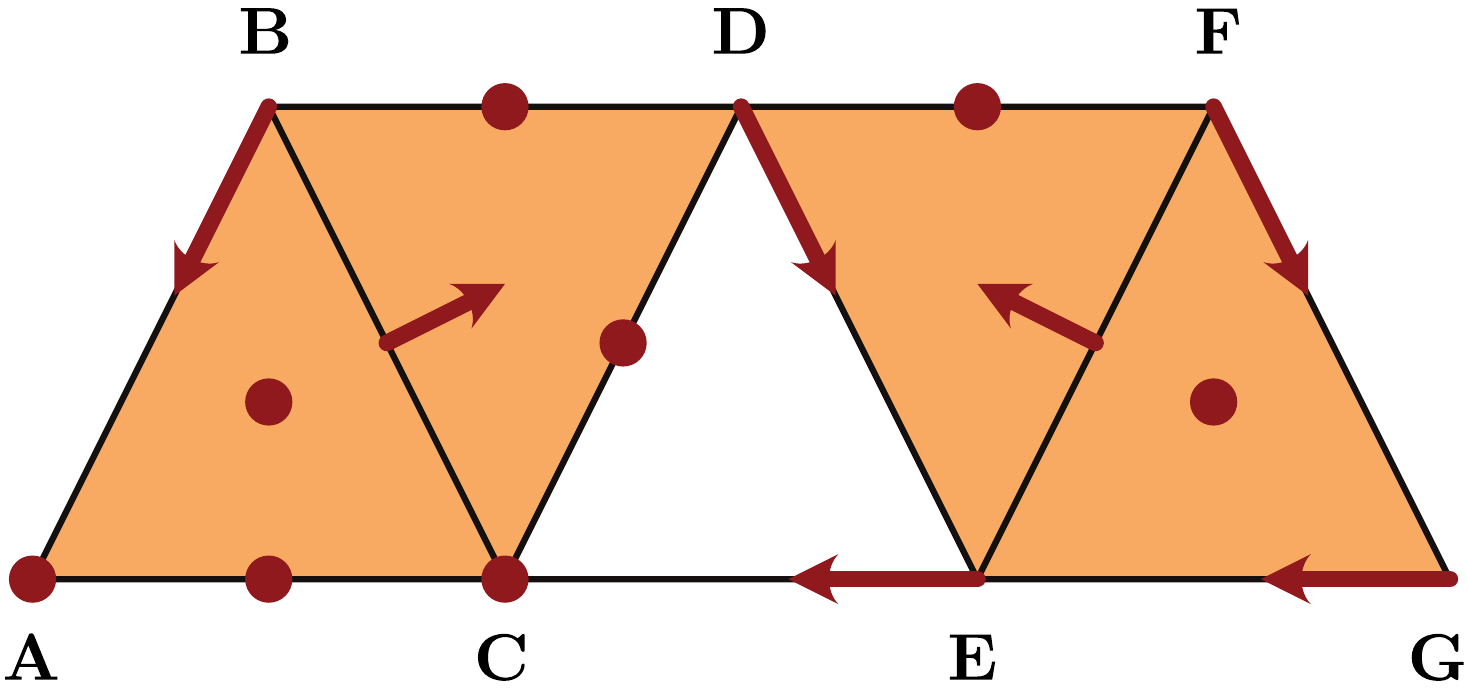}\quad\quad
    \includegraphics[width=0.45\textwidth]{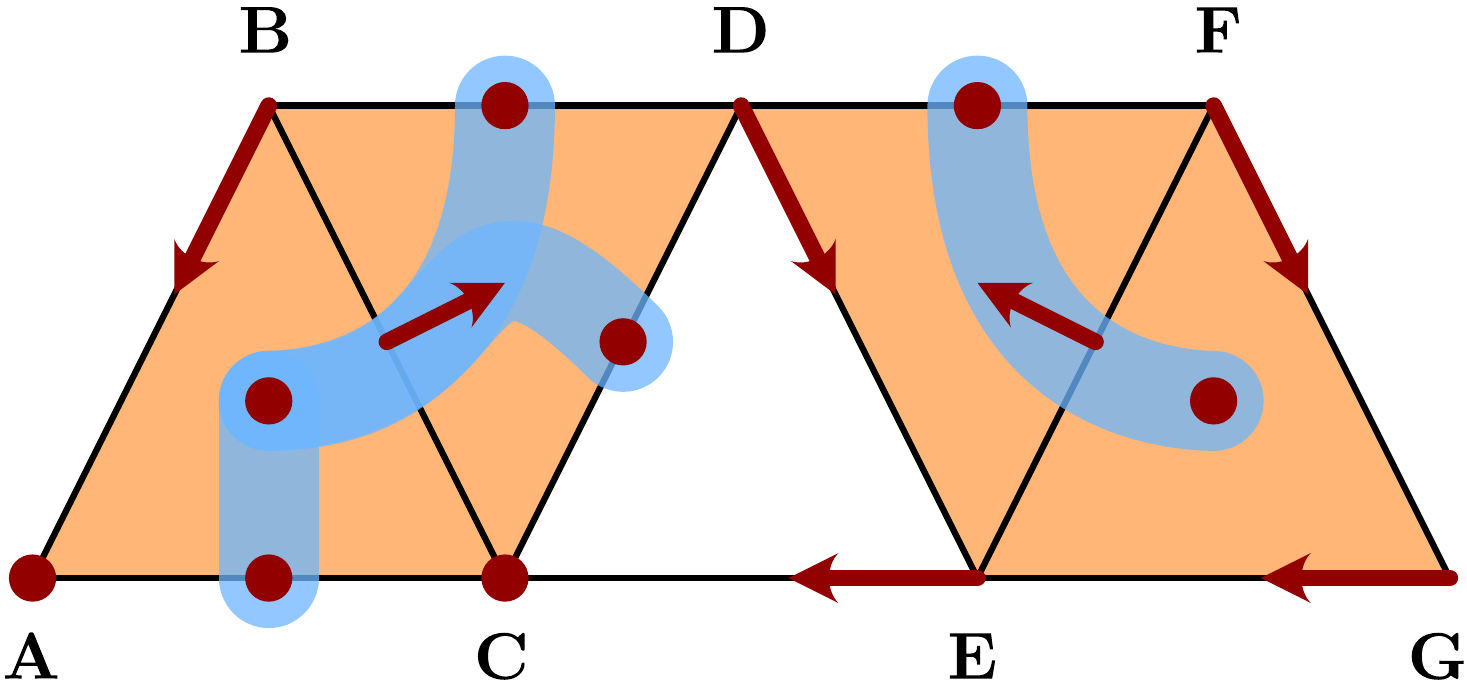} \\[3ex]
    \includegraphics[width=0.45\textwidth]{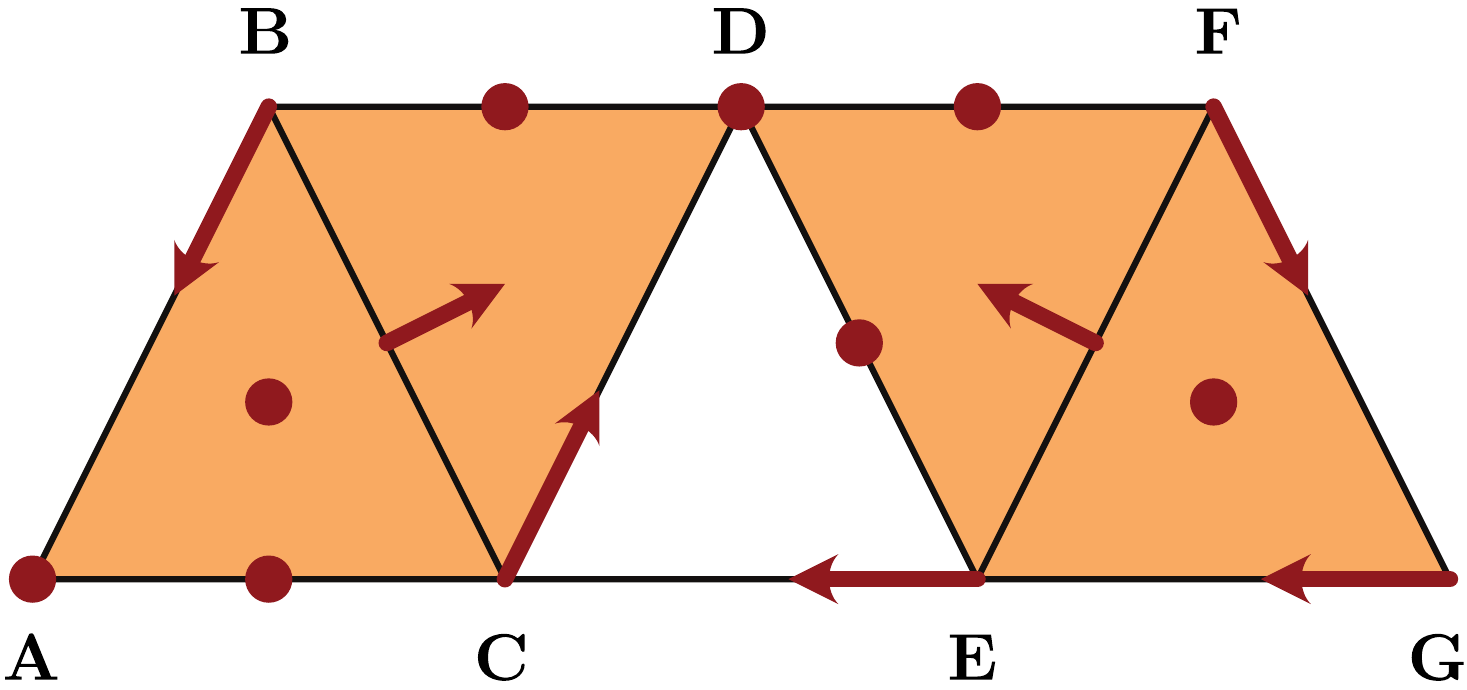}\quad\quad
    \includegraphics[width=0.45\textwidth]{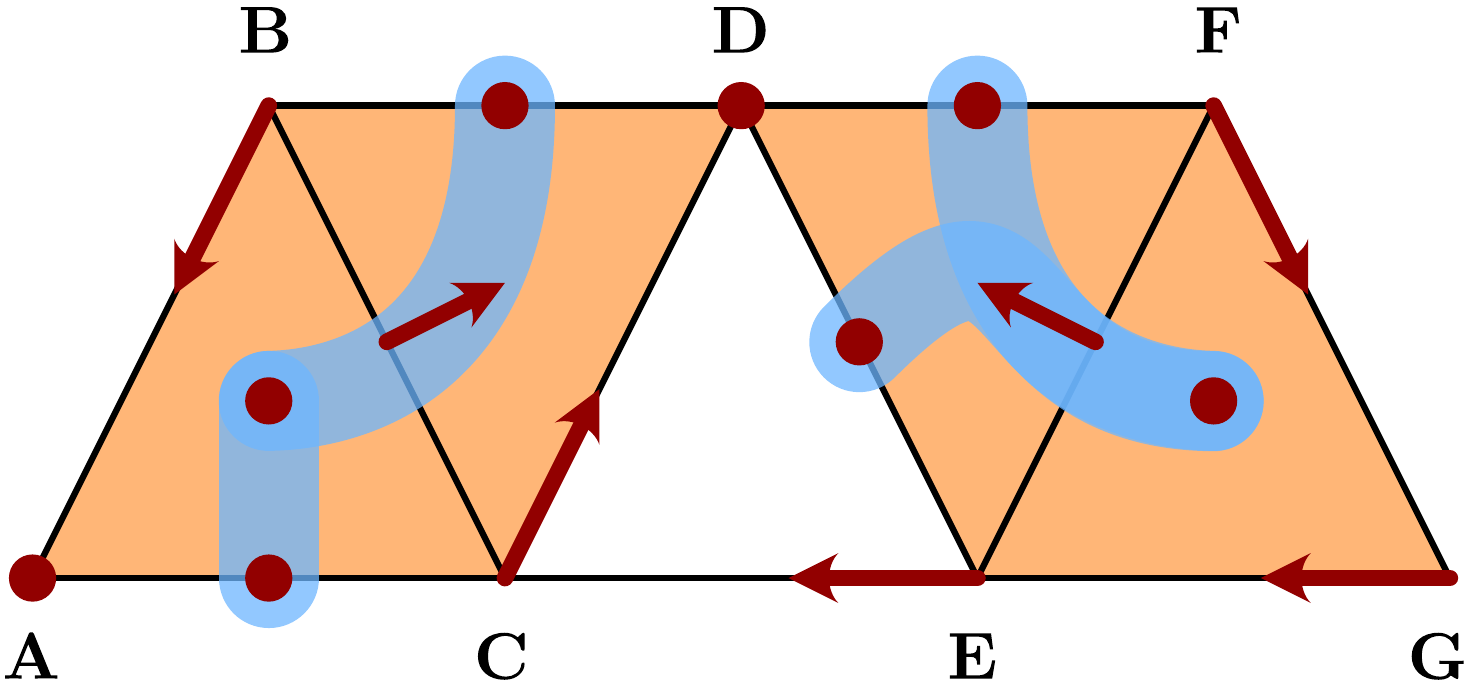} \\[3ex]
    \includegraphics[width=0.45\textwidth]{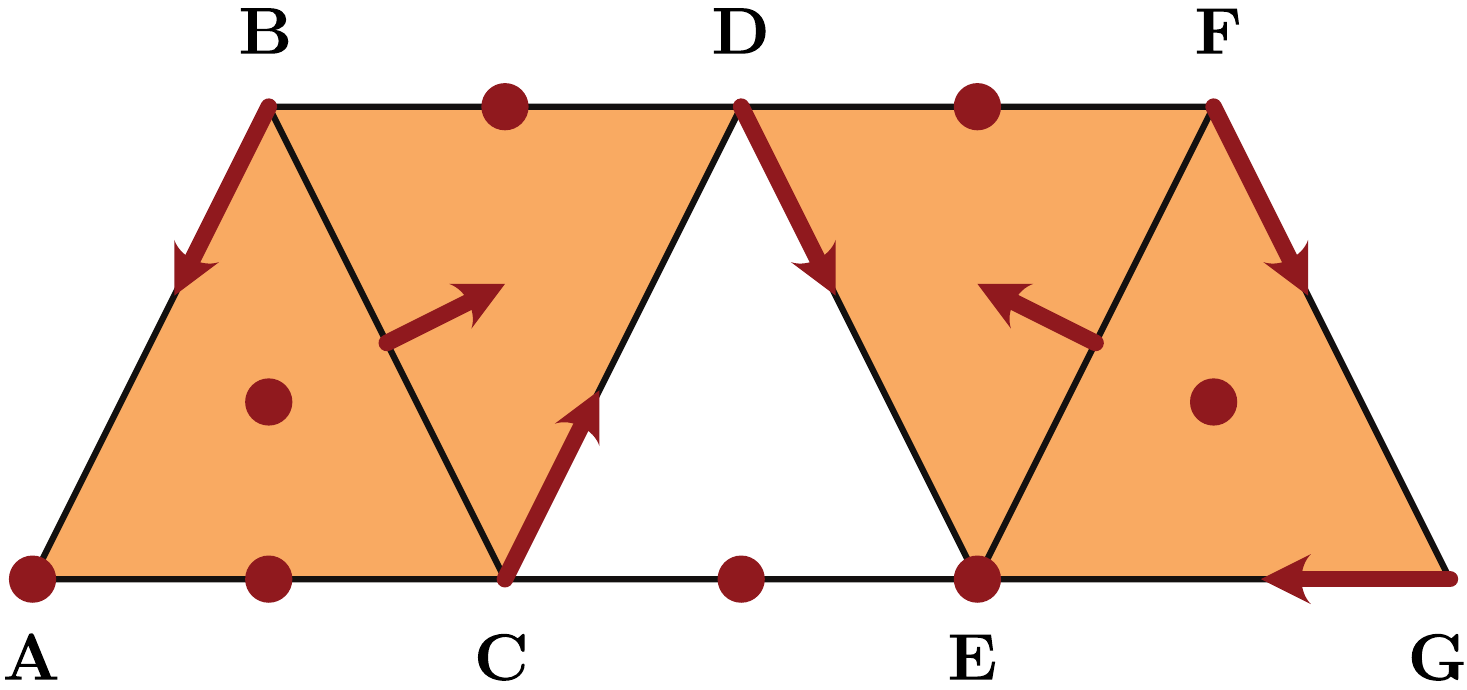}\quad\quad
    \includegraphics[width=0.45\textwidth]{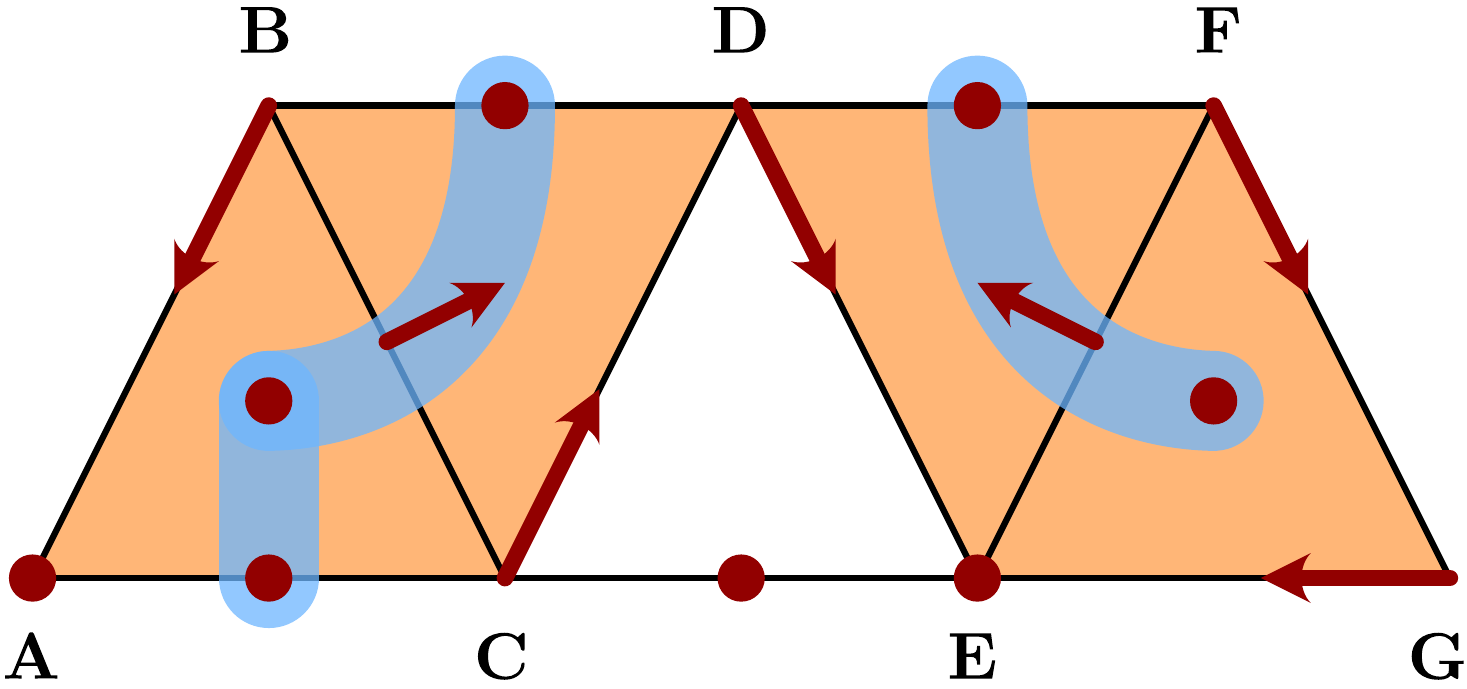}
  \end{center}
  \caption{{\em Three Forman gradient vector fields}.
           The above panels show three different combinatorial Forman
           gradient vector fields. All of them are related to the
           combinatorial vector field from Figure~\ref{fig:periodicex0},
           yet they destroy its periodic orbit which traverses the
           vertices~$\bC$, $\bD$, and~$\bE$. This is achieved by
           replacing precisely one vector in the periodic orbit
           by two critical cells --- one edge and one vertex. Since the
           periodic orbit consists of three vectors, there are three
           different ways of breaking it, leading to the vector
           fields in the left column. In the right column, we indicate
           for each case which connections exist between the index~$2$
           and index~$1$ critical cells.
           }
  \label{fig:periodicex123}
\end{figure}
\begin{figure}
  \begin{center}
    \includegraphics[width=0.95\textwidth]{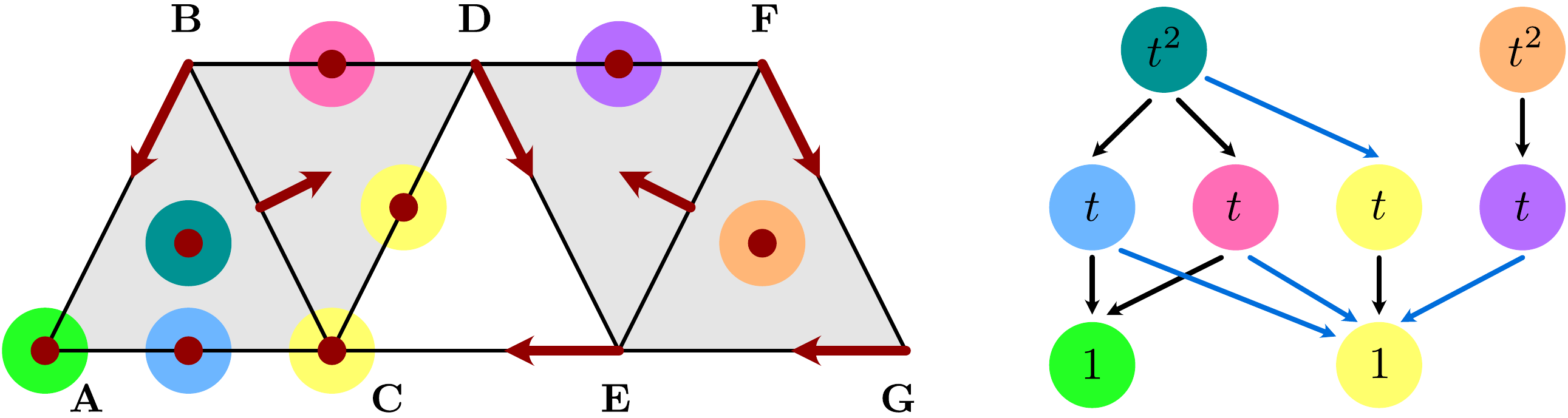} \\[5ex]
    \includegraphics[width=0.95\textwidth]{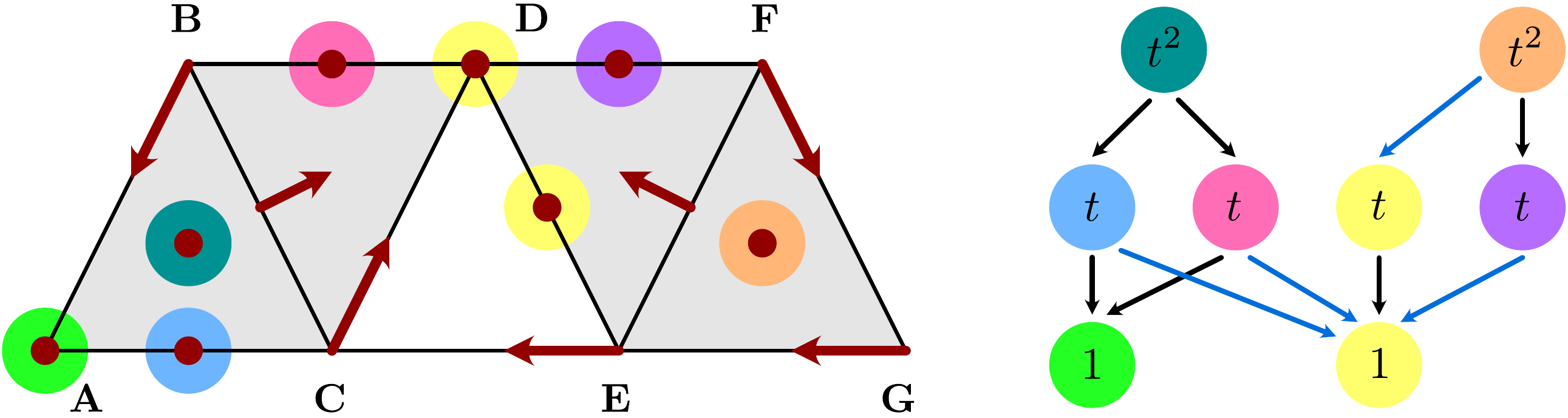} \\[5ex]
    \includegraphics[width=0.95\textwidth]{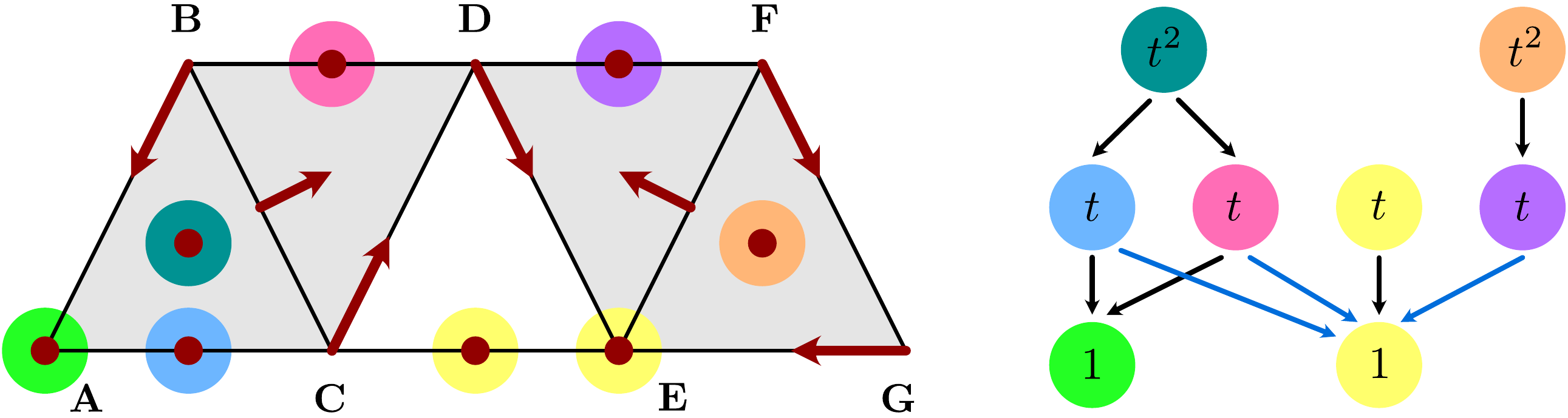}
  \end{center}
  \caption{{\em Three Forman gradient vector fields}.
           For the combinatorial vector fields shown in
           Figure~\ref{fig:periodicex123}, the above panels depict
           their associated Morse decompositions. The Morse sets
           are indicated by different colors in the left column,
           while the Conley-Morse graphs are shown on the right.
           In each case, the two Morse sets obtained via breaking
           the periodic orbit in Figure~\ref{fig:periodicex0}
           are shown in yellow.
           }
  \label{fig:morsedecompex123}
\end{figure}
\begin{table}
  \begin{center}
  \begin{tabular}{c||cccc|cc|}
      & $\mathit{CD}$ & $\mathbf{AC}$ & $\mathbf{BD}$ & $\mathbf{DF}$ &
        $\mathbf{ABC}$ & $\mathbf{EFG}$ \\ \hline\hline
      $\mathbf{A}$  & 0 & 1 & 1 & 0 &   &   \\
      $\mathit{C}$  & 0 & 1 & 1 & 0 &   &   \\ \hline
      $\mathit{CD}$ &   &   &   &   & 1 & 0 \\
      $\mathbf{AC}$ &   &   &   &   & 1 & 0 \\
      $\mathbf{BD}$ &   &   &   &   & 1 & 0 \\
      $\mathbf{DF}$ &   &   &   &   & 0 & 1 \\ \hline
  \end{tabular} \\[3ex]
  \begin{tabular}{c||cccc|cc|}
      & $\mathit{DE}$ & $\mathbf{AC}$ & $\mathbf{BD}$ & $\mathbf{DF}$ &
        $\mathbf{ABC}$ & $\mathbf{EFG}$ \\ \hline\hline
      $\mathbf{A}$  & 0 & 1 & 1 & 0 &   &   \\
      $\mathit{D}$  & 0 & 1 & 1 & 0 &   &   \\ \hline
      $\mathit{DE}$ &   &   &   &   & 0 & 1 \\
      $\mathbf{AC}$ &   &   &   &   & 1 & 0 \\
      $\mathbf{BD}$ &   &   &   &   & 1 & 0 \\
      $\mathbf{DF}$ &   &   &   &   & 0 & 1 \\ \hline
  \end{tabular} \\[3ex]
  \begin{tabular}{c||cccc|cc|}
      & $\mathit{CE}$ & $\mathbf{AC}$ & $\mathbf{BD}$ & $\mathbf{DF}$ &
        $\mathbf{ABC}$ & $\mathbf{EFG}$ \\ \hline\hline
      $\mathbf{A}$  & 0 & 1 & 1 & 0 &   &   \\
      $\mathit{E}$  & 0 & 1 & 1 & 0 &   &   \\ \hline
      $\mathit{CE}$ &   &   &   &   & 0 & 0 \\
      $\mathbf{AC}$ &   &   &   &   & 1 & 0 \\
      $\mathbf{BD}$ &   &   &   &   & 1 & 0 \\
      $\mathbf{DF}$ &   &   &   &   & 0 & 1 \\ \hline
  \end{tabular} \\[1ex]
  \end{center}
  \caption{{\em Three Forman gradient vector fields}.
           For the three combinatorial vector fields shown in
           Figure~\ref{fig:periodicex123} and associated Morse 
           decompositions in Figure~\ref{fig:morsedecompex123},
           the three unique connection matrices are listed above.
           Only the two right-most entries in the third
           row are changing, and they correspond to the connections
           from the index~$2$ critical cells to the index~$1$
           critical cell on the destroyed periodic orbit. Note
           that the nonzero entries in the lower right
           $4 \times 2$-submatrices correspond exactly to the
           connections indicated in the right column of
           Figure~\ref{fig:periodicex123}.
           }
  \label{table:periodicex123}
\end{table}

\begin{ex}[{\em Three Forman gradient vector fields}]
\label{ex:formangradient-1}
{\em
Consider the three combinatorial Forman vector fields shown in the
left column of Figure~\ref{fig:periodicex123}. All of them are
defined on the same simplicial complex, and they consist of eight
critical cells and seven doubletons. In fact, one can easily see 
that all of these vector fields are related to the combinatorial
vector field from Figure~\ref{fig:periodicex0}, which has a periodic
orbit traversing the vertices~$\bC$, $\bD$, and~$\bE$. While the
latter vector field is therefore not gradient, it is possible to
destroy its recurrent behavior by replacing precisely one vector
in the periodic orbit by two critical cells --- one edge and one
vertex. Since the periodic orbit consists of three vectors, this 
leads to the three different gradient vector fields shown in the
left column of Figure~\ref{fig:periodicex123}. Their corresponding
Morse decompositions are depicted in Figure~\ref{fig:morsedecompex123}.

Using our results from Section~\ref{sec:cm-gcvf}, one can easily
determine the unique connection matrix of a combinatorial gradient
vector field. This basically amounts to finding suitable bases
of chains in each dimension and creating the matrix of the boundary
operator with respect to these bases, see the detailed description
in Theorem~\ref{thm:cm-for-forman-gradient} and Proposition~\ref{prop:barX}.
For the gradient combinatorial vector field in the last row of
Figure~\ref{fig:periodicex123} this is illustrated in detail in
Example~\ref{ex:formangradient-4}, and this leads to the last 
connection matrix in Table~\ref{table:periodicex123}. For the first
two gradient vector fields, their connection matrices are given by
the first two matrices in Table~\ref{table:periodicex123},
respectively, and they can be determined using the basis elements
described in Example~\ref{ex:formangradient-3}. Notice that the
nonzero entries in all three of these matrices between critical
cells of index~$2$ and critical cells of index~$1$ yield precisely
the connecting orbits shown in the right column of
Figure~\ref{fig:periodicex123}. Similarly, nonzero entries 
between critical cells of index~$1$ and critical cells of index~$0$
do correspond to connections in the combinatorial vector fields.
Note, however, that the entry in the connection matrix records
also the multiplicity of such connections. This is the reason why
the connection matrices record the connections between~$\bB\bD$
or~$\bA\bC$ and the vertex~$\bA$, as well as with the additional
critical vertex in the support of the destroyed periodic orbit.
In contrast, connections originating at either~$\bD\bF$ or the
critical cell of dimension one on the former periodic orbit do
not give rise to nonzero entries in the connection matrix, since
in these cases the two different connections end at the same
critical vertex.
\exend
}
\end{ex}

\begin{ex}[{\em A Forman vector field with periodic orbit}, continued]
\label{ex:formanperiodic-2}
{\em
Based on the previous example, we can now revisit the combinatorial
vector field from Figure~\ref{fig:periodicex0}, which has already
been discussed in Example~\ref{ex:formanperiodic-1}. Notice that this
vector field has a periodic orbit which traverses the vertices~$\bC$,
$\bD$, and~$\bE$. The existence of this periodic orbit makes the system
non-gradient, and therefore we cannot apply Theorem~\ref{thm:uniqueness}
to obtain the existence of a unique connection matrix.

Nevertheless, Theorem~\ref{thm:uniqueness} guarantees the existence
of a connection matrix in this situation as well. In order to explain 
how such a matrix can be determined using the previous example, note that
the Conley polynomial of the periodic orbit is given by~$t + 1$, i.e., 
the associated homology groups are one-dimensional in dimensions~$0$
and~$1$. In Example~\ref{ex:formangradient-1} we saw that by replacing
a vector on the periodic orbit by one critical cell each of dimension~$0$,
one obtains a gradient system with a unique connection matrix. In fact,
the sum of the homology groups of these newly added critical cells is
isomorphic to the Conley index of the periodic orbit, and therefore one
would expect that all three connection matrices from Table~\ref{table:periodicex123}
should be connection matrices for the vector field in Figure~\ref{fig:periodicex0}.
As we will see later, this is indeed the case, and thus in all three
cases the nontrivial entries correspond to connections between the
respective Morse sets. For a more detailed explanation we refer the reader
to Example~\ref{ex:formanperiodic-3}.
\exend
}
\end{ex}


\section{Preliminaries}
\label{sec:prelim}

\subsection{Sets, maps, relations, and partial orders}
\label{sec:sets-and-maps}
We denote  the sets of reals, integers, positive integers, and non-negative
integers by~$\RR$, $\ZZ$, $\NN$, and~$\NN_0$, respectively. Given a
set~$X$, we write~$\card X$ for the number of elements of~$X$ and
we denote the family of all subsets of~$X$ by~$\cP(X)$.
A {\em partition} of a set~$X$ is a family $\cE\subset\cP(X)$
of mutually disjoint  non-empty subsets of~$X$ such that $\bigcup\cE=X$.
For a subfamily $\cE'\subset \cE$ of a partition~$\cE$ of~$X$
we will use the compact notation~$|\cE'|$ to denote the union
of sets in~$\cE'$, that is, we define $|\cE'|=\bigcup \cE'$.

A {\em $K$-indexed family} of subsets of a given set $X$ is an
injective map of the form $K\ni k\mapsto X_k\in\cP(X)$.
We denote such a family by $(X_k)_{k\in K}$.
Note that every family $\cE\subset\cP(X)$ may be considered as
an $\cE$-indexed family via the map $\cE\ni E\mapsto E\in\cP(X)$.
In this case we say that $\cE$ is a {\em self-indexing} family.
By a {\em $K$-gradation} of a set~$X$ we mean a $K$-indexed family~$(X_k)_{k\in K}$
such that~$\{X_k\}_{k\in K}$ is a partition of~$X$. Note, in particular,
that every surjective map $f:X\to K$ induces the $K$-gradation
$(f^{-1}(k))_{k\in K}$ of $X$. We refer to this $K$-gradation as the
{\em $f$-gradation of $X$}. Given a subset $A\subset X$ and a
gradation~$(X_k)_{k\in K}$ of $X$, the {\em induced gradation}
is the gradation $(A \cap X_k)_{k\in K_A}$ of $A$ where
$K_A:=\setof{k\in K\mid A\cap X_k\neq\emptyset}$.

We write $f:X\pto Y$ for a {\em partial map} from $X$ to $Y$, that is, a map
defined on a subset $\dom f\subset X$, called the {\em domain} of $f$,
and such that the set of values of $f$, denoted by
$\im f$, is contained in $Y$. Partial maps are composed in the obvious way:
If $f:X\pto Y$ and $g:Y\pto Z$ are partial maps, then the domain of their
composition~$g \circ f : X \pto Z$ is given by $\dom(g \circ f) = f^{-1}(\dom g)$,
and on this domain we have $(g\circ f)(x) = g(f(x))$.

In the sequel, we work with the {\em category of finite sets with a
distinguished subset}, which is denoted by~$\DSet$ and defined as follows.
The objects of $\DSet$ are pairs $(X,X_\star)$, where $X$ is a finite set
and $X_\star\subset X$ is a distinguished subset. The morphisms from
$(X,X_\star)$ to  $(Y,Y_\star)$ in $\DSet$ are {\em subset preserving
partial maps}, that is, partial maps $f:X\pto Y$ such that the inclusions
$X_\star\subset \dom f$ and $f(X_\star)\subset Y_\star$ are satisfied.
One easily verifies that~$\DSet$ with the composition of morphisms defined
as the composition of partial maps and the identity morphism defined as the
identity map is indeed a category.

Given a set $X$ and a binary relation $R\subset X\times X$, we write~$xRy$
meaning $(x,y)\in R$. The {\em inverse} of~$R$, denoted $R^{-1}$, is the relation
$R^{-1}:=\setof{(y,x)\mid xRy}$. By the {\em transitive closure} of~$R$ we mean
the relation $\bar{R}\subset X\times X$ given by $x\bar{R}y$ if there exists
a sequence $x=x_0,x_1,\ldots, x_n=y$ such that $n\geq 1$ and  $x_{i-1}Rx_i$
for $i=1,2,\ldots, n$. 

A multivalued map $F:X\mvmap Y$ is a map $F:X\to\cP(Y)$.
For $A\subset X$ we define the {\em image of $A$} by
$
    F(A):=\bigcup \setof{F(x) \mid x\in A}
$
and for $B\subset Y$ we define the {\em preimage of $B$} by
$
    F^{-1}(B):=\setof{x\in X\mid F(x)\cap B\neq\emptyset}.
$

Given a relation $R$, we associate with it a multivalued map $F_R:X\mvmap X$
via the definition $F_R(x):=R(x)$, where
$
   R(x):=\setof{y\in X\mid xRy}
$
is the {\em image of $x\in X$ in~$R$}.
Obviously, $R\mapsto F_R$  is a one-to-one correspondence between
binary relations in $X$ and multivalued maps from $X$ to $X$.
Often, it will be convenient to interpret the relation $R$ as a directed graph
whose set of vertices is $X$ and a directed arrow goes from $x$ to $y$ whenever $xRy$.
The three concepts relation, multivalued map, and directed graph are
equivalent on the formal level and the distinction is used only to emphasize
different points of view. However, in this paper it will be convenient to use
all these concepts interchangeably.

Recall that a {\em partial order} in a set $P$
is a reflexive, antisymmetric and transitive relation in $P$.
If not stated otherwise we denote a partial order by~$\leq$ and
its inverse by~$\geq$. We also write $<$ and $>$ for the associated
strict partial orders, that is, the relations~$\leq$ and~$\geq$ excluding equality.
We also recall that a {\em poset} is a pair  $\PP=(P,\leq)$
where $P$ is a set and $\leq$ is a partial order in $P$.
Whenever the the partial order  is clear from the context,
we refer to $P$ as the poset.

Given a preordered set $P$ we say that a $q\in P$ {\em covers} a $p\in P$
if $p<q$ and there is no $r\in P$ such that
$p<r<q$. We say that $p$ is a {\em predecessor} of~$q$ if~$q$ covers~$p$.
Recall that a subset $A\subset P$ is {\em convex} if the conditions
$x,z\in A$ and $x\leq y\leq z$ imply $y\in A$.
Clearly, the intersection of a family of convex sets is a convex set. 
This lets us define the {\em convex hull} of $J\subset P$, 
denoted $\conv_P(J)$, as the intersection of the family of all convex sets in $P$
containing $J$. It is straightforward to observe that
\begin{equation}
\label{eq:conv-hull}
\conv_P(J)=\setof{p\in P\mid \exists p_-,p_+\in J\;\mbox{ with }\; p_-\leq p \leq p_+}.
\end{equation}

A set $A\subset P$ is a {\em down set} or a {\em lower set} if the conditions
$z\in A$ and $y\leq z$ yield $y\in A$.
Dually, the subset $A\subset P$ is an {\em upper set}  if 
$z\in A$ and $z\leq y$ yield $y\in A$.

We denote the family of all down sets in~$P$ by $\Down(P)$.
For $A\subset P$ we further write
$A^{\leq}:=\setof{x\in P\mid \exists_{a\in A}\;x\leq a}$ and
$A^{<}:=A^{\leq}\setminus A$.

\begin{prop}
\label{prop:interval}
   For any subset $I \subset P$ the set~$I^{\leq}$ is a down set.
   In addition, if~$I$ is convex, then~$I^{<}$ is a down set as well.
\end{prop}
\proof
  The verification that  $I^{\leq}$ is a down set is straightforward.
To see that $I^{<}$ is a down set take an $x\in I^{<}$. Then we have $x\not\in I$ and
$x<z$ for some $z\in I$. Let $y\leq x$. Then  $y\in I^{\leq}$.
Since $I$ is convex, we cannot have the inclusion $y\in I$. It follows that $y\in I^{<}$.
\qed

\medskip
Let $P$ and $P'$ be posets. We say that a partial map $f:P\pto P'$ is {\em order preserving}
if the conditions $x,y\in\dom f$ and $x\leq y$ imply $f(x)\leq f(y)$.
We say that $f$ is an {\em order isomorphism} if $f$ is an order preserving bijection
such that $f^{-1}$ is also order preserving.

We define the category $\DPSet$ of posets with a distinguished subset as follows.
Its objects are pairs $(P,P_\star)$ where $P$ is a finite poset and $P_\star\subset P$ is
a distinguished subset. A morphism from $(P,P_\star)$ to $(P',P'_\star)$ in $\DPSet$
is an order preserving partial map $f: P\pto P'$ such that $P_\star\subset\dom f$ and $f(P_\star)\subset P'_\star$.
We call such a morphism $f$ {\em strict} if $\dom f=P_*$.
One easily verifies that $\DPSet$ with the composition of morphisms defined as the composition of partial maps
and the identity morphism defined as the identity map
is indeed a category. Moreover, since every set may be considered as a poset partially ordered by the identity
and, clearly, every partial map preserves identities, we may consider $\DSet$ as a subcategory of $\DPSet$.
In this case, every partial map is automatically order preserving.

To simplify notation, in the sequel we denote objects of $\DSet$ and $\DPSet$ with a single capital letter and
the distinguished subset by the same letter with subscript $\star$.

\subsection{Topological spaces}
\label{sec:top}
For our terminology concerning topological
spaces we refer the reader to~\cite{En1989,munkres:00a}.
Recall that a {\em topology} on a set $X$ is a family $\cT$ of subsets
of $X$ which is closed under finite intersections and
arbitrary unions, and which satisfies $\emptyset,X\in\cT$.
A {\em topological space} is a pair $(X,\cT)$ where $\cT$ is a topology on~$X$.
We often refer to $X$ as a topological space assuming that the topology~$\cT$
on~$X$ is clear from context.
The sets in $\cT$ are called {\em open}.
The {\em interior} of~$A$, denoted $\inte A$, is the union of all open subsets of~$A$.
A subset $A\subseteq X$ is {\em closed} if $X\setminus A$ is open.
We denote by $\Clsd(X)$ the family of closed subsets of $X$.
The {\em closure} of $A$, denoted $\cl_\cT A$, or by $\cl A$ if $\cT$ is clear from context,
is the intersection of all closed supersets of $A$.

Given two topological spaces $(X,\cT)$ and $(X',\cT')$ we say that
$f:X\to X'$ is {\em continuous} if  $U\in\cT'$ implies $f^{-1}(U)\in\cT$.
If we want to emphasize the involved topologies we say that
$f: (X,\cT) \to (X',\cT')$ is continuous.

Recall that a subset~$A$ of a topological space~$X$
is called {\em locally closed},
if every point $x\in A$ admits a neighborhood $U$ in $X$ such that
$A\cap U$ is closed in $U$.
Locally closed sets as well as $\mouth A:=\cl A\setminus A$, which we refer to as the {\em mouth} of $A$,
are important in the sequel.

\begin{prop}[see Problem 2.7.1 in~\cite{En1989}]
\label{prop:lcl}
Assume that~$A$ is a subset of a topological space $X$.
Then the following conditions are equivalent.
\begin{itemize}
   \item[(i)] $A$ is locally closed,
   \item[(ii)] the mouth of $A$ is closed in $X$,
   \item[(iii)] $A$ is a difference of two closed subsets of $X$,
   \item[(iv)] $A$ is an intersection of an open set in $X$ and a closed set in $X$.
\end{itemize}
\end{prop}
\proof
   It is obvious that (ii) implies (iii), (iii) implies (iv), and that~(iv) implies~(i).
Thus, it suffices  to prove that if (ii) fails, then so does (i).
Hence, assume that the mouth $\mouth A=\cl A\setminus A$
is not closed in $X$. Then there exists a point
$x\in  \cl(\cl A\setminus A)\setminus(\cl A\setminus A)$.
It follows that $x\in A$. Let $U$ be an arbitrary neighborhood of~$x$ in~$X$.
We will prove that $A\cap U$ is not closed in~$U$.
Since $x\in   \cl(\cl A\setminus A)$, we can find a point $y\in (\cl A\setminus A)\cap U$.
Let $V$ be an open neighborhood of~$y$ in~$X$.
Then also $V\cap U$ is a neighborhood of~$y$ in~$X$.
Since $y\in \cl A\setminus A \subset \cl A$, we see that $V\cap U\cap A\neq\emptyset$.
Since all open neighborhoods of $y$ in $U$ are of the form $V\cap U$ where $V$ is open in $X$,
we see that $y\in\cl_U\left(U\cap A\right)$.
However, $y\not\in A$, hence, $y\not\in U\cap A$.
It follows that the intersection~$U\cap A$ is not closed in $U$.
This holds for every neighborhood~$U$ of~$x$ in~$X$.
Therefore, $A$ is not locally closed.
\qed

\medskip
The topology $\cT$ is $T_2$
or {\em Hausdorff} if for any two different points $x,y\in X$ there exist disjoint sets
$U,V\in\cT$ such that $x\in U$ and $y\in V$.
It is $T_0$ or {\em Kolmogorov} if for any two different points $x,y\in X$ there exists a $U\in\cT$
such that $U\cap\{x,y\}$ is a singleton.

Given $Y\subset X$, the family $\cT_Y:=\setof{U\cap Y\mid U\in \cT}$ is a topology on~$Y$ called
the {\em induced topology}.
The topological space $(Y,\cT_Y)$ is called a {\em subspace} of $(X,\cT)$.
For $A\subset Y$ we write $\cl_Y A:=\cl_{\cT_Y}A$.

Finally, we say that a topological space~$X$ is a {\em finite topological
space} if the underlying set~$X$ is finite. Finite topological spaces differ
from general topological spaces, because the only Hausdorff topology
on a finite topological space~$X$ is the discrete topology consisting
of all subsets of~$X$. Moreover, the celebrated Alexandrov Theorem~\cite{Al1937}
states that every finite, $T_0$ topological space can equivalently be considered as
a  poset, if we set $x\leq_\cT y$ whenever $x \in \cl_\cT y$.
The results of this paper apply only to a very special, but crucial for applications,
class of finite topological spaces, namely Lefschetz complexes (see Section~\ref{sec:lefschetz}).
However, the view through the lens of the non-standard features of finite topological spaces
facilitates a better understanding of the peculiarities of combinatorial topological dynamics.
We refer the interested reader to~\cite{barmak:11a,barmak:etal:20a,LKMW2020} for more information.

\subsection{Modules and their gradations}
\label{sec:modules}

For the terminology used in the paper concerning modules we refer the reader to~\cite{lang:02a}.
Here, we briefly recall that a {\em module} over a ring $R$ is a triple $(M,+,\cdot)$ such that
$(M,+)$ is an abelian group and $\cdot: R\times M\to M$ is a scalar multiplication which is distributive
with respect to the additions in~$M$ and~$R$, and which satisfies both $a\cdot(b\cdot x)=(ab)\cdot x$
and $1_R\cdot x=x$ for $x\in M$ and $a,b\in R$.
We shorten the notation $a\cdot x$ to $ax$.
We often say that $M$ is a module assuming that the operations~$+$
and~$\cdot$ are clear from context. Recall that a {\em submodule} of~$M$
is a subset $N\subset M$ such that~$N$ with the operations~$+$ and~$\cdot$
restricted to~$N$ is itself a module.
The associated {\em quotient module} is denoted $M/N$.
For submodules $N_1,N_2,\ldots,N_k$ of a given module $M$,
their {\em algebraic sum}
\[
N_1+N_2+\cdots N_k:=\setof{x_1+x_2+\ldots + x_k\in M\mid x_i\in N_i}
\]
is easily seen to be a submodule of $M$. We call it a {\em direct sum} and denote it by
the symbol $N_1\oplus N_2\oplus\ldots\oplus N_k$,
if for each $x\in N_1+N_2+\ldots + N_k$ the representation
 $x=x_1+x_2+\ldots+x_k$ with $x_i\in N_i$ is unique.

The following
straightforward result will be useful for establishing the existence of
connection matrices later on.

\begin{prop}
\label{prop:XAB}
Assume that~$A$ and~$B$ are submodules of a module $X$ with~$X=A+B$.
If~$B'$ is a submodule of~$B$ such that $B = (A\cap B) \oplus B'$,
then we have $X=A\oplus B'$.
\qed
\end{prop}

A subset  $Z\subset M$ {\em generates} $M$ if for every element $x\in M\setminus\{0\}$ there exist elements $x_1,x_2,
\ldots, x_n\in Z$ and $a_1,a_2,\ldots, a_n\in R$ such that $x=\sum_{i=1}^na_ix_i$.
A subset  $Z\subset M$ is {\em linearly independent} if for any choice of
$x_1,x_2,\ldots, x_n\in Z$ and  $a_1,a_2,\ldots, a_n\in R$ the equality $\sum_{i=1}^na_ix_i=0$
implies $a_1=\ldots=a_n=0$. A linearly independent subset $B\subset M$ which generates $M$ is called
a {\em basis} of $M$. A module may not have a basis. If it does have a basis, the module is called {\em free}.
We note that $\emptyset$ is the unique basis of the trivial module
$\{0\}$.

Assume that $M$ is a free module and that $B\subset M$ is a fixed basis.
Then we have the associated bilinear map $\scalprod{\cdot,\cdot}_B:
M\times M\to R$, called {\em scalar product}, and defined on basis elements $b,b'\in B$ by
$\scalprod{b,b'}_B:=0$ for $b \neq b'$, as well as $\scalprod{b,b'}_B:=1$ for
$b=b'$.
For an element $x\in M$ we define its {\em support} with respect to $B$
as $|x|_B:=\setof{b\in B\mid \scalprod{x,b}_B\neq 0}$. One easily verifies that
for all $x,y\in M$
\begin{equation}
\label{eq:support-union}
   |x+y|_B\subset |x|_B\cup |y|_B.
\end{equation}
Assume that~$X$ is a finite set. Then the collection $R\spn{X}:=\setof{f:X\to R}$
of all functions forms a module
with respect to pointwise addition and scalar multiplication.
We note that $R\spn{\emptyset}$ is the zero module.
For every $x\in X$ we have a function $\bar{x}:M\to R$ which sends $x$ to $1$ and any other element of $X$
to $0$. One easily verifies that  $\bar{X}:=\setof{\bar{x}\mid x\in X}$ is a basis of $R\spn{X}$.
Therefore, $R\spn{X}$ is a free module, and we call  $R\spn{X}$ the {\em free module spanned by $X$}.
In the sequel we identify~$\bar{x}$ with~$x$.

Let $M$ and $M'$ be two modules over $R$. A map $h:M\to M'$  is called a {\em module homomorphism}
if we have $h(a_1x_1+a_2x_2)=a_1h(x_1)+a_2h(x_2)$ for arbitrary $x_1,x_2\in M$ and $a_1,a_2\in R$.
The {\em kernel} of $h$, denoted by $\ker h$, is the submodule $\setof{x\in M\mid h(x)=0}$,
while the {\em image} of $h$, denoted by $\im h$, is the submodule $\setof{y\in M'\mid \exists_{x\in M} h(x)=y}$.

We now turn our attention to gradations of modules.
Assume that~$R$ is a fixed ring and $M$ is a module over $R$. Let~$K$ be an
arbitrary set. Then a {\em $K$-gradation} of~$M$ is a collection of
submodules $\{M_k\}_{k\in K}$ of~$M$ indexed by~$K$, and such that
$M=\bigoplus_{k\in K}M_k$ is the direct sum of the submodules~$M_k$.
This definition requires that $K\neq\emptyset$,
but it is convenient to assume that the direct sum
over the empty set is always the zero module.
Thus, the only module admitting the $\emptyset$-gradation
is the zero module.
By a {\em $K$-graded module} over~$R$ we mean a module~$M$
over~$R$ together with a fixed, implicitly given, $K$-gradation.
Note that the zero module admits not only the $\emptyset$-gradation,
but also a unique $K$-gradation for each non-empty
set $K$. Thus, for each set $K$ the zero module is a $K$-graded module.
We denote it by $0_K$.

If $B$ is a basis of a free module $M$ and $\cA$ is a partition of $B$, then 
$M$ is $\cA$-graded with the gradation
\[
    M=\bigoplus_{A\in\cA} R\cdot A,
\]
where the expression on the right-hand side is defined as
\[
    R \cdot A = \left\{ \sum_{i=1}^n r_i a_i \; \mid \;
                r_i \in R, \; a_i \in A, \; i = 1,\ldots,n, \; n \in \NN \right\}.
\]
Note that in the case that~$A$ is finite, one can identify~$R \cdot A$
with~$R\spn{A}$. As a special case, consider the situation when the partition
of~$B$ consists only of singletons. 
Then the gradation becomes
\[
    M=\bigoplus_{b\in B} R\cdot b.
\]
We refer to this gradation as the {\em $B$-basis gradation} of $M$.

We say that a submodule $M'$ of $M$ is a {\em $K$-graded submodule} if $M'$ is also
$K$-graded with the decomposition
\begin{equation*}
   M'=\bigoplus_{k\in K}M'_k
\end{equation*}
and such that $M'_k\subset M_k$. For  $I\subset K$ we have an {\em $I$-graded submodule}
\begin{equation}
\label{eq:M-I}
    M_I:=\bigoplus_{i\in I}M_i.
\end{equation}
The homomorphisms
\[
    \iota_I: M_I\ni x \mapsto x\in M
\]
and
\[
    \pi_I: M=\bigoplus_{k\in K}M_k\ni \sum_{k\in K}x_k \mapsto \sum_{i\in I}x_i\in M_I
\]
are called the associated {\em canonical inclusion\/} and {\em canonical projection\/}.
In the case $I=\setof{i}$ we abbreviate the above notation to~$\iota_i$ and~$\pi_i$.

Assume $K'$ is another set and $M'$ is a $K'$-graded module.
For a module homomorphism $f:M\to M'$ and subsets $J\subset K$ and $I\subset K'$ we have
an induced homomorphism $f_{IJ}:M_J\to M'_I$ defined as the composition
\[
  f_{IJ}:=\pi_I \circ f \circ \iota_J,
\]
where~$\pi_I : M' \to M_I'$ and~$\iota_J : M_J \to M$.
Again, if $I=\{i\}$ and $J=\{j\}$ we abbreviate the notation to $f_{ij}$.
One easily verifies the following proposition.
\begin{propdef}
\label{prop:matrix-rep}
Assume $K$ and $K'$ are finite. Then
\begin{equation}
\label{eq:matrix-rep}
   f=\sum_{i\in K'} \iota_i \sum_{j\in K}f_{ij} \pi_j.
\end{equation}
In particular, $f$ is uniquely determined by
the matrix $[f_{ij}]_{i\in K',j\in K}$ of homomorphisms $f_{ij}:M_j\to M'_i$.
We refer to this matrix as the {\em $(K,K')$-matrix of the homomorphism~$f$}.
\qed
\end{propdef}

Notice that in the case when the $K$- and $K'$-gradations are basis
gradations given by bases $B=\{b_1,b_2,\ldots, b_n\}$ in~$M$ 
and $B'=\{b'_1,b'_2,\ldots, b'_m\}$ in~$M'$, respectively, then the homomorphisms  $f_{ij}$
take the form
\[
  f_{ij} \; : \; t\cdot b_j\mapsto t \scalprod{f(b_j),b_i}\cdot b_i,
\]
which means that the $(K,K')$-matrix of $f$ may be identified in this case 
with the matrix of coefficients $\scalprod{f(b_j),b_i}$.

As in the case of classical matrices one easily verifies the following proposition.
\begin{prop}
\label{prop:composition}
Consider modules $M$, $M'$, and $M''$ which are graded
by~$K$, $K'$, and $K''$, respectively.
Let $f:M\to M'$ and $f':M'\to M''$ be module homomorphisms.
Then the $(K,K'')$-matrix of $f$
consists of the homomorphisms
\begin{equation}
\label{eq:composition}
  (f'f)_{ik}=\sum_{j\in K'} f'_{ij}f_{jk}
\end{equation}
for $i\in K''$ and $k\in K$.
\qed
\end{prop}

Given a $\ZZ$-graded module $M$, a homomorphism $f:M\to M$
is of {\em degree} $k\in\ZZ$ if $f_{ij}\neq 0$ implies $i-j=k$ for $i,j\in\ZZ$.
If $f$ is of degree $k$ we write $f_j$ meaning $f_{ij}$ with $i=j+k$.

\subsection{Chain complexes and their homology}
\label{sec:chain-complexes}

Recall that a  {\em chain complex} is a pair~$(C,d)$
consisting of a $\ZZ$-graded $R$-module~$C = (C_k)_{k \in \ZZ}$ and a $\ZZ$-graded
homomorphism $d:C\to C$ which satisfies $d^2=0$, and which has degree $-1$, i.e.,
which satisfies $d(C_k) \subset C_{k-1}$ for all $k \in \ZZ$.
The homomorphism~$d$ is called the {\em boundary homomorphism} of the chain complex~$(C,d)$.
Let $(C',d')$ be another chain complex.
A {\em chain map} $\varphi:(C,d)\to (C',d')$ is a module
homomorphism $\varphi:C\to C'$ of degree zero, satisfying $\varphi d=d'\varphi$.
The following proposition is straightforward.
\begin{prop}
\label{prop:inverse-chain-map}
If $\varphi:(C,d)\to (C',d')$ is a chain map and $\varphi$ is an isomorphism of $\ZZ$-graded modules
then $\varphi^{-1}$ is also a chain map.
\qed
\end{prop}
We denote by $\CCR$ the category whose object are chain complexes and whose
morphisms are chain maps. One easily verifies that this is indeed a category.

A subset $C'\subset C$ is a {\em chain subcomplex} of $C$ if $C'$ is a $\ZZ$-graded submodule of $C$
such that $d(C')\subset C'$. Recall that given a subcomplex $C'\subset C$ we have a well-defined
{\em quotient complex} $(C/C',d')$ where $d':C/C'\to C/C$ is the boundary homomorphism induced by $d$.

We now turn our attention to a fundamental equivalence relation on chain maps.
A pair of chain maps $\varphi,\varphi':(C,d)\to (C',d')$ are called
{\em chain homotopic} if there exists a {\em chain homotopy} joining $\varphi$
and $\varphi'$, i.e., a module homomorphism $S:C\to C'$ of degree~$+1$
such that $\varphi'-\varphi=d' S + S d$. The existence of a chain homotopy
between two chain maps is easily seen to be an equivalence relation
in the set of chain maps~$\CCR((C,d),(C',d'))$. Given a chain map
$\varphi\in \CCR((C,d),(C',d'))$ we denote by $[\varphi]$
the equivalence class of $\varphi$ with respect to this equivalence relation.
We define the {\em homotopy category}~$\ChCC$ of chain complexes by taking
chain complexes as objects, equivalence classes of morphisms in~$\CCR$
as morphisms in~$\ChCC$, and the formula
\begin{equation}
\label{eq:eq-comp}
   [\psi]\circ[\varphi]:=[\psi\varphi]
\end{equation}
for $\psi\in \CCR((C',d'),(C'',d''))$ as the definition of composition of
morphisms in $\ChCC$. Note that then the equivalence classes of identities
in~$\CCR$ are the identities in~$\ChCC$.

\begin{prop}
\label{prop:ChCC-category}
The category $\ChCC$ is well-defined.
\end{prop}
\proof
The only non-obvious part of the argument is the verification
that the composition given by \eqref{eq:eq-comp} is well-defined.
Thus, we need to prove that if  $\varphi,\varphi':(C,d)\to (C',d')$ and $\psi,\psi':(C',d')\to (C'',d'')$
are chain homotopic, then $\psi\varphi$ and $\psi'\varphi'$ are chain homotopic.
Let $S:C\to C'$ and $S':C'\to C''$ be
chain homotopies between $\varphi$, $\varphi'$
and $\psi$, $\psi'$, respectively.
Moreover, consider the map $S'':=\psi' S + S'\varphi$.
Then~$S''$ is clearly a degree~$+1$ homomorphism and we have
\begin{eqnarray*}
\psi'\varphi'-\psi\varphi&=&\psi'(\varphi'-\varphi)+(\psi'-\psi)\varphi=
          \psi'(d' S + S d)+(d'' S' + S' d')\varphi\\
&=& \psi'd'S+\psi'Sd+d''S'\phi +S'd'\phi\\
&=& d''\psi' S + \psi' S d + d'' S' \varphi + S' \varphi d=d''S''+S''d,
\end{eqnarray*}
which proves that $\psi\varphi$ and $\psi'\varphi'$ are indeed chain homotopic.
\qed

\medskip
We have a covariant functor
$\Ch:\CCR\to\ChCC$ which fixes objects and sends a chain map to its chain homotopy
equivalence class.
Moreover, we say that two chain complexes $(C,d)$ and $(C',d')$ are {\em chain homotopic} if they are isomorphic in $\ChCC$.

Note that $0_{\mathbb{Z}}$, the $\ZZ$-graded zero module,
together with the zero homomorphism as the boundary map, is a chain complex.
We call it the {\em zero chain complex}.
We say that a chain complex $(C,d)$ is {\em homotopically essential} if it is not
chain homotopic to the zero chain complex. Otherwise we say that the chain
complex~$(C,d)$ is {\em homotopically trivial} or {\em homotopically inessential}.
Finally, we call a chain complex~$(C,d)$ {\em\boundaryless\/} if $d=0$.

\begin{prop}
\label{prop:simple-homotopic-isomorphic}
Assume $(C,d)$ and $(C',d')$ are two \boundaryless\  chain complexes.
Then the chain complexes $(C,d)$ and $(C',d')$ are chain homotopic
if and only if $C$ and $C'$ are
isomorphic as $\ZZ$-graded modules.
\end{prop}
\proof
  Suppose that the chain complexes $(C,d)$ and $(C',d')$ are \boundaryless.  
  First assume that $C$ and $C'$ are isomorphic as $\ZZ$-graded modules.
  Let $\phi:C\to C'$ and $\phi':C'\to C$ be  mutually inverse isomorphisms
  of $\ZZ$-graded modules. Since both boundary homomorphisms $d$ and $ d'$
  are zero, we have $\phi d = 0 = d'\phi$ and $\phi' d' = 0 = d\phi'$.
  Hence, $\phi$ and $\phi'$ are chain maps
  and from $\phi'\phi=\id_C$ and $\phi\phi'=\id_{C'}$
we get $[\phi'][\phi]=[\id_C]$ and $[\phi][\phi']=[\id_{C'}]$.
This shows that $[\phi]$ and $[\phi']$
are mutually inverse isomorphisms in $\ChCC$.

  To prove the opposite implication we now assume that $(C,d)$ and $(C',d')$ are chain homotopic.
  Choose $\phi:(C,d)\to (C',d')$ and $\phi':(C',d')\to (C,d)$
such that $[\phi]:(C,d)\to (C',d')$ and $[\phi']:(C',d')\to (C,d)$ are mutually inverse isomorphisms in $\ChCC$.
 Let $S:C\to C$ and $S':C' \to C'$ be the chain homotopies between $\phi'\phi$
 and $\id_C$, and between $\phi\phi'$ and $\id_{C'}$, respectively.
  Then $\id_C-\phi'\phi=Sd+dS = 0$ and $\id_{C'}-\phi\phi'=S'd'+d'S' = 0$,
  which proves that $\phi$ and $\phi'$ are mutually inverse isomorphisms in $\CCR$.
\qed

\begin{cor}
\label{cor:simple-homotopic-isomorphic}
The only chain complex which is both \boundaryless\  and homotopically
trivial is the zero chain complex.
\qed
\end{cor}

We now turn our attention to a discussion of the homology of chain complexes.
The {\em homology module} of a chain complex $(C,d)$ is the $\ZZ$-graded
module $H(C):=(H_n(C,d))_{n\in\ZZ}$, where we define $H_n(C,d):=\ker d_n/\im d_{n+1}$.
In the sequel, we will consider the homology module as a \boundaryless\ chain complex,
that is, as a chain complex with zero boundary homomorphism.

By a {\em homology decomposition} of a chain complex $(C,d)$ we mean a direct
sum decomposition $C=V\oplus H\oplus B$ such that $V$, $H$, $B$ are $\ZZ$-graded
submodules of $C$, $d_{|H}=0$, $d(V)\subset B$ and $d_{|V}:V\to B$ is a module
isomorphism. Note that then also $d_{|B}=0$, since $B = d(V)$ implies $d(B) = d^2(V) = \{0\}$.
Finally, a {\em homology complex} of a chain complex~$(C,d)$ is defined as
a \boundaryless\  chain complex which is chain homotopic to~$(C,d)$. Then
one has the following proposition.

\begin{prop}
\label{prop:hom-decomp-existence}
Assume that~$R$ is a field and that the pair~$(C,d)$ is a chain complex over~$R$.
\begin{itemize}
   \item[(i)] There exists a homology decomposition of $(C,d)$.
   \item[(ii)] If $C=V\oplus H\oplus B$ is a homology decomposition of the
   chain complex~$(C,d)$, then $(H,0)$ is chain homotopic
   to $(C,d)$. Therefore, it is a homology complex of $(C,d)$.
   Moreover, the modules~$H$ and~$H(C)$ are isomorphic as $\ZZ$-graded modules,
   where~$H(C)$ denotes the homology module of~$C$.
\end{itemize}
\end{prop}
\proof
  To prove statement (i), let $Z:=\ker d$
and $B:=\im d$. Since $R$ is a field, we may choose
$\ZZ$-graded submodules $V\subset C$ and $H\subset Z$
such that both $C=V\oplus Z$ and $Z=H\oplus B$ are satisfied, 
thus also $d_{H\oplus B}=0$. Clearly, one has $d(V)\subset\im d=B$.
We will prove that $d_{|V}:V\to B$ is an isomorphism.
Indeed, if $dx=0$ for an $x\in V$, then $x\in Z\cap V=\{0\}$,
hence~$d_{|V}$ is a monomorphism. To see that it is also onto, take
a $y\in B$. Since $B=\im d$, we have $y=dx$ for an $x\in C$. But,
$x=v+z$ for a $v\in V$ and a $z\in Z$. It follows that
$dv=dv+dz=dx=y$, which establishes~$d_{|V}$ as an epimorphism,
and proves (i).

To prove (ii), assume that
$C=V\oplus H\oplus B$ is a homology decomposition of~$(C,d)$.
Let $\iota: H \to C$ denote inclusion, and let $\pi : C \to H$
be the projection map defined via
$\pi(v+h+b) = h$ for $v+h+b \in V\oplus H\oplus B = C$.
One can easily see that $\iota:(H,0)\to (C,d)$
and $\pi:(C,d)\to (H,0)$ are chain maps, and that $\pi\iota=\id_H$.
We will show that $\iota\pi$ is chain homotopic to $\id_C$.
For this, define the degree~$+1$ homomorphism
$S:C\to C$ by $Sx=d^{-1}_{|V}b$, where $x=v+h+b\in C=V\oplus H\oplus B$.
Then we have
\begin{eqnarray*}
  (dS+Sd)x & = & (dS+Sd)(v+h+b)=d d_{|V}^{-1}b+d_{|V}^{-1}dv=b+v\\
  & = & (b+v+h)-h=x-\iota\pi x=(\id_C-\iota\pi)x,
\end{eqnarray*}
which proves that $S$ is a chain homotopy joining $\iota\pi$ and $\id_C$.
Finally, directly from the homology decomposition definition one obtains
$\ker d=H\oplus B$ and $\im d=B$. Therefore,
$H(C)=\ker d/\im d\cong H\oplus B/B\cong H$, where all the isomorphisms are $\ZZ$-graded.
\qed

\medskip
In the case of field coefficients the concepts of homology module and homology complex
are essentially the same as the following theorem shows.

\begin{thm}
\label{thm:existence-and-uniqueness}
Assume that~$R$ is a field. Then the following hold.
\begin{itemize}
   \item[(i)] Every chain complex admits a homology complex.
   \item[(ii)] Two chain complexes are chain homotopic
   if and only if the associated homology complexes are
   isomorphic as $\ZZ$-graded modules.
   \item[(iii)]  The homology module of a chain complex is its homology complex.
   \item[(iv)] Two chain complexes are chain homotopic
   if and only if the associated homology modules are
   isomorphic as $\ZZ$-graded modules.
\end{itemize}
\end{thm}
\proof
Property (i) follows immediately from Proposition~\ref{prop:hom-decomp-existence}.
In order to prove~(ii), observe that two chain complexes are chain homotopic if and only if
the associated homology complexes are chain homotopic. Therefore,
property (ii) follows immediately from Proposition~\ref{prop:simple-homotopic-isomorphic}.
To prove (iii), take a chain complex~$(C,d)$. By (i) we may consider a homology complex~$(H,0)$
of~$(C,d)$. Since $(C,d)$ and $(H,0)$ are chain homotopic, by a standard theorem of homology
theory~\cite[\S 13]{Munkres1984} the homology modules of $(C,d)$ and $(H,0)$ are isomorphic
as $\ZZ$-graded modules. Hence, it follows from Proposition~\ref{prop:hom-decomp-existence}(ii)
that the homology module of $(C,d)$ is its homology complex.
Finally, property (iv) follows immediately from (ii) and (iii).
\qed

\subsection{Lefschetz complexes}
\label{sec:lefschetz}
We end our section on preliminaries by recalling basic definitions
and facts about Lefschetz complexes. The following definition goes back
to S.~Lefschetz, see~\cite[Chapter~III, Section~1, Definition~1.1]{Le1942}.

\begin{defn} \label{def:lefschetz}
{\em
We say that $(X,\kappa)$ is a {\em Lefschetz complex} over a ring $R$
if $X=(X_q)_{q\in{\scriptsize\NN_0}}$ is a finite set with $\NN_0$-gradation,
$\kappa : X \times X \to R$ is a map such that
\begin{equation}
\label{eq:kappa-condition-1}
\kappa(x,y)\neq 0 \quad\implies\quad x\in X_q,\;\;y\in X_{q-1},
\end{equation}
and for any $x,z\in X$ we have
\begin{equation}
\label{eq:kappa-condition-2}
    \sum_{y\in X}\kappa(x,y)\kappa(y,z)=0.
\end{equation}
We refer to the elements of $X$ as {\em cells}, to $\kappa(x,y)$ as the {\em incidence
coefficient} of the cells $x$ and $y$, and to $\kappa$ as the {\em incidence coefficient map}.
We define the {\em dimension} of a cell $x\in X_q$ as~$q$, and denote it by $\dim x$.
Whenever  the incidence coefficient map is clear from context
we often just refer to~$X$ as the Lefschetz complex.
We say that $(X,\kappa)$ is {\em regular} if for any $x,y\in X$
the incidence coefficient
$\kappa(x,y)$ is either zero or it is invertible in $R$.
}
\end{defn}

Let $(X,\kappa)$ be a given Lefschetz complex.
We denote by $C_k(X):=R\spn{X_k}$ the free $R$-module
spanned by the set $X_k$ of cells of dimension~$k$ for $k\in\NN_0$,
and let~$C_k(X)$ denote the zero module for $k<0$.
Then it is clear that the sum $C(X):=\bigoplus_{k\in\ZZ} C_k(X)$
is a free $\ZZ$-graded $R$-module generated by $X$.
Finally, define the module homomorphism $\bdy^\kappa:C(X)\to C(X)$
on generators by
\begin{equation}
\label{eq:kappa-bdy}
   \bdy^\kappa(x) := \sum_{y\in X}\kappa(x,y)y.
\end{equation}

\begin{propdef}
  The pair $(C(X),\bdy^\kappa)$ is a chain complex.
We call it the {\em chain complex of $(X,\kappa)$},
and refer to the homology of this chain complex as the {\em Lefschetz homology} of $(X,\kappa)$.
\end{propdef}
\proof
Condition \eqref{eq:kappa-condition-1} guarantees that $\bdy^\kappa$
is a degree~$-1$ module homomorphism, and condition \eqref{eq:kappa-condition-2}
implies that $(\bdy^\kappa)^2=0$.
\qed

\medskip
Note that every finitely generated free chain complex is the chain complex of a Lefschetz
complex obtained by selecting a basis. More precisely, assume that~$(C,\bdy)$ is a finitely
generated free chain complex over a ring~$R$ and $U\subset C$ is a fixed basis of~$C$.
Suppose further that~$C_k = 0$ for all $k < 0$. Then for every $v\in U$ there are
uniquely determined coefficients $a_{vu}\in R$ such that
$$
   \bdy v=\sum_{u\in U}a_{vu}u.
$$
Let $\kappa_\bdy:U\times U\to R$ be defined by $\kappa_\bdy(v,u)=a_{vu}$.
The following  proposition is straightforward.

\begin{prop}
\label{prop:U-Lefschetz}
 The pair  $(U,\kappa_\bdy)$ is a Lefschetz complex.
 \qed
\end{prop}

The family of cells of a simplicial complex~\cite[Definition 11.8]{KaMiMr2004} and the family of
elementary cubes of a cubical set~\cite[Definition 2.9]{KaMiMr2004} provide simple
but important examples of Lefschetz complexes.
In these two cases the respective formulas for the incident coefficients are explicit and
elementary, see for example~\cite{MB2009}.
Also a general regular cellular complex, or a regular finite CW complex as
considered in~\cite[Section IX.3]{Ma1991},
is an example of a Lefschetz complex. In this case
the incident coefficients may be obtained
from a system of equations as described in~\cite[Section~IX.5]{Ma1991},
and the Lefschetz homology may be computed efficiently as outlined
in~\cite{dlotko:etal:11a}.
Note that a Lefschetz complex over a field is always regular.

Given $x,y\in X$ we say that $y$ is a {\em facet} of~$x$, and we write $y\adhl_{\kappa} x$,
if we have $\kappa(x,y)\neq 0$.
It is easily seen that the reflexive and transitive closure of the relation~$\adhl_{\kappa}$
is a partial order.
We denote this partial order by $\leq_\kappa$, call it the {\em face relation},
and denote the associated strict order by~$<_\kappa$.
As an immediate consequence of~\eqref{eq:kappa-condition-1}
we have the following proposition.

\begin{prop}
\label{prop:dim-monotoonicity}
The map
$
  \dim: (X,\leq_\kappa)\to(\ZZ,\leq)
$
which assigns to a cell $x\in X$ its dimension $\dim x\in\ZZ$ is order preserving.
Moreover, if the inequality $x \adhl_{\kappa} y$ holds, then one has $\dim y=\dim x+1$.
\end{prop}

We say that $y$ is a {\em face} of $x$ if $y\leq_\kappa x$.
The $T_0$ topology defined via the Alexandrov Theorem~\cite{Al1937}
by the partial order $\leq_\kappa$ is called the
{\em Lefschetz topology} of the Lefschetz complex $(X,\kappa)$.
The Lefschetz topology makes the Lefschetz complex a finite topological space. 
Notice that in this topology the closure of any set $A\subset X$
consists of all faces of all cells in~$A$.

We say that the subset $A\subset X$ is a {\em $\kappa$-subcomplex} or  {\em Lefschetz subcomplex}
of~$X$, if $(A,\kappa_{|A\times A})$ with the $\ZZ$-gradation induced from $X$ is a Lefschetz complex
in its own right.
The following proposition provides sufficient conditions for a subset of a Lefschetz complex to be
a Lefschetz subcomplex. For more details see~\cite[Proposition~5.3]{Mr2017}
and~\cite[Theorem~3.1]{MB2009}.
\begin{prop}
\label{prop:Lefschetz-subcomplex}
If $A\subset X$ is locally closed in the Lefschetz topology, then~$A$
is a Lefschetz subcomplex of $(X,\kappa)$.
In particular, every open and every closed subset of $X$ is a Lefschetz subcomplex.
\end{prop}

\begin{figure}
  \includegraphics[width=0.3\textwidth]{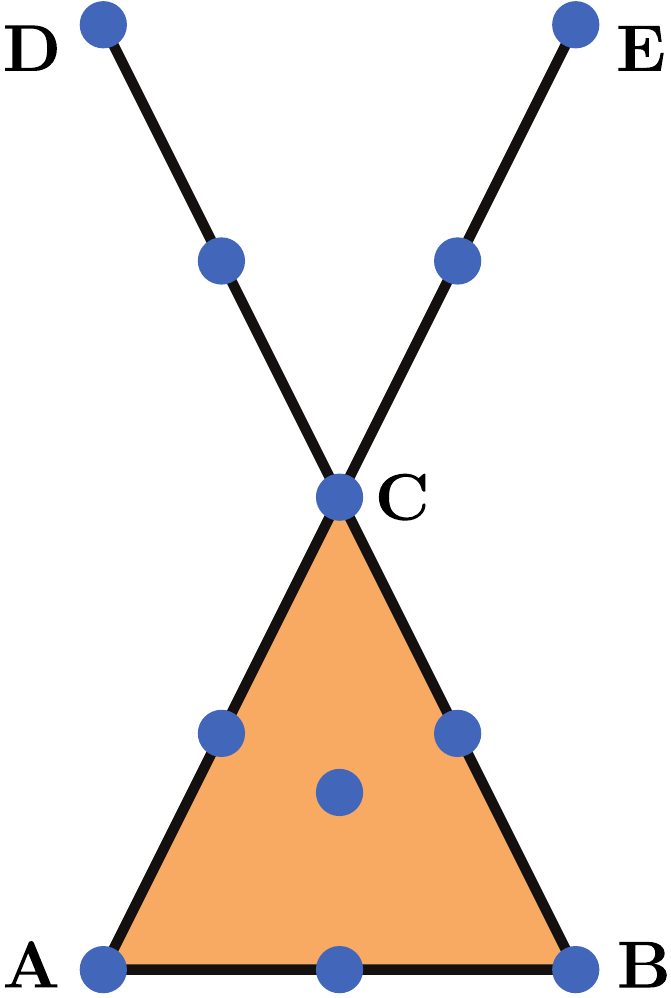}\qquad\qquad\qquad
  \includegraphics[width=0.3\textwidth]{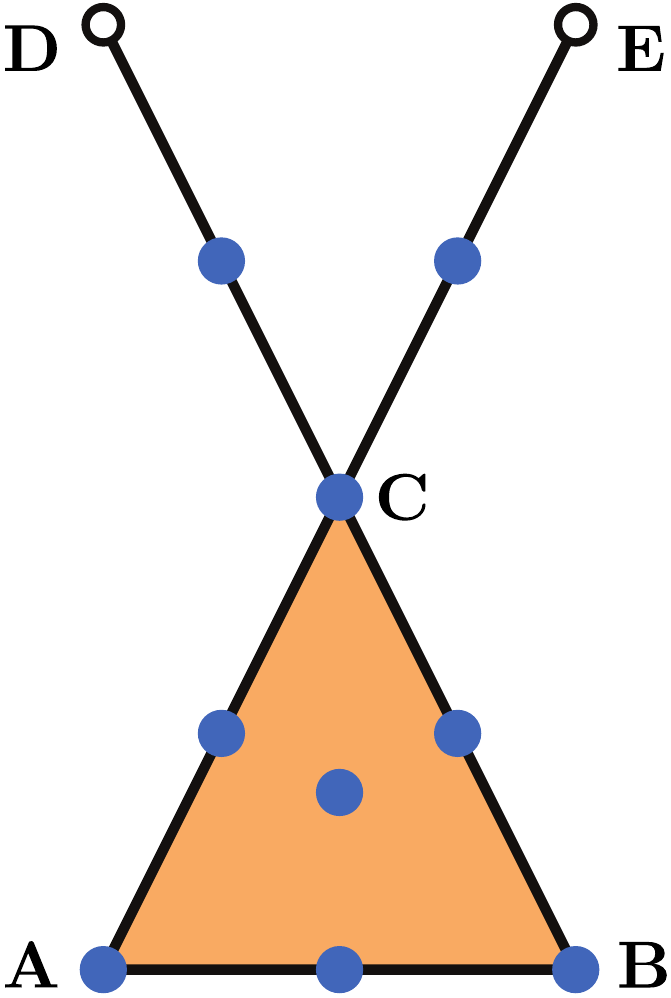}
  \caption{
  The simplicial complex shown on the left and consisting of the triangle $\mathbf{ABC}$, the edges
  $\mathbf{AB}$, $\mathbf{AC}$, $\mathbf{BC}$, $\mathbf{CD}$, $\mathbf{CE}$,
  and the vertices $\mathbf{A}$, $\mathbf{B}$, $\mathbf{C}$, $\mathbf{D}$, and~$\mathbf{E}$
  is an example of a Lefschetz complex. Simplices  are marked with a small circle in
  the center of mass of each simplex. In the right panel, the simplices marked with a blue circle
  constitute a locally closed (convex) subset of the set of all simplices.
  Therefore, it is another example of a Lefschetz complex $X$ which consists
  of all simplices of the simplicial complex except the vertices~$\mathbf{D}$
  and~$\mathbf{E}$, shown as circles with white interiors.
  }
  \label{fig:t2eLef}
\end{figure}

\begin{ex}
\label{ex:lefschetz}
{\em 
  The easiest way to visualize a Lefschetz complex is by presenting it as a $\kappa$-subcomplex
  of a simplicial complex.  A sample Lefschetz complex is presented in Figure~\ref{fig:t2eLef}
  as a locally closed collection of simplices of a simplicial complex. Note, however, that
  not all Lefschetz complexes can be written in this form.
\exend}
\end{ex}

A Lefschetz subcomplex~$A$ of a Lefschetz complex~$X$ has its own Lefschetz topology 
as well as the associated chain complex~$C(A)$. As a topological space~$A$ is always 
a subspace of~$X$, but~$C(A)$ typically is not a chain subcomplex of~$X$ unless~$A$ is closed in~$X$
(see~\cite[Theorem~5.4]{Mr2017}), as the following result shows.

\begin{prop}
\label{prop:closed-subcomplex}
If $A$ is closed in $X$ in the Lefschetz topology, then
\begin{equation}
\label{eq:closed-subcomplex-1}
\bdy^{\kappa_{|A\times A}}=\bdy^\kappa_{|C(A)}
\end{equation}
and
\begin{equation}
\label{eq:closed-subcomplex-2}
\bdy^\kappa(C(A))\subset C(A).
\end{equation}
In particular, the chain complex  $(C(A),\bdy^{\kappa_{|A\times A}})$ is
a chain subcomplex of the chain complex $(C(X),\bdy^\kappa)$.
\end{prop}
\proof
To see \eqref{eq:closed-subcomplex-1} it suffices to verify the equality on basis elements. Thus, take an $x\in A$. Then
\[
\bdy^{\kappa_{|A\times A}}x=\sum_{y\in A} \kappa(x,y)y=\sum_{y\in X} \kappa(x,y)y=\bdy^\kappa x,
\]
because $\kappa(x,y)\neq 0$ implies $y\in\cl x\subset\cl A=A$.
Property \eqref{eq:closed-subcomplex-2} and the remaining assertion are obvious.
\qed

\medskip
Given a closed subset $A\subset X$ in the Lefschetz topology
we define the relative Lefschetz homology $H(X,A)$ as the homology
of the quotient chain complex $(C(X,A),\tilde{\bdop})$,
where $C(X,A):=C(X)/C(A)$ and $\tilde{\bdop}$ stands for the
induced boundary map. In addition, we also have the following
proposition, which follows immediately from~\cite[Theorem~5.4]{Mr2017}
(see also~\cite[Theorem~3.5]{MB2009}), and which uses such a relative
homology to characterize the Lefschetz homology of a locally
closed $A\subset X$.
\begin{prop}
\label{prop:rel-homology}
If $A\subset X$ is locally closed in the Lefschetz topology of $X$, then $H(A)$ is isomorphic 
to the homology $H(\cl A,\mo A)$ of the pair $(\cl A,\mo A)$ of closed subsets of $X$.
\qed
\end{prop}
The following result was established in~\cite[Proposition~5.2]{Mr2017}, and it computes
the homology of some simple Lefschetz complexes.
\begin{prop}
\label{prop:singl-dubl-homology}
For every $x\in X$ the singleton $\{x\}$ is a Lefschetz
subcomplex of $X$ and
\[
     H_q(\{x\})\cong\begin{cases}
                   R & \text{ if $q=\dim x$,}\\
                   0 & \text{ otherwise.}
                \end{cases}
\]
For every $x,y\in X$ such that $x$ is a facet of $y$
the doubleton $\{x,y\}$ is a Lefschetz
subcomplex of $X$ and for all $q \in \ZZ$ one has
\[
     H_q(\{x,y\})\cong 0.
\]\end{prop}


\section{Poset filtered chain complexes}
\label{sec:posetfilteredcc}

In the preliminaries section we recalled basic definitions and properties
of chain complexes and their homology. For the definition of the connection 
matrix, these considerations have to be extended to the case of chain
complexes which are poset filtered. In addition, we need to study maps
between these poset filtered chain complexes, even in the situation where
the filtration posets of the involved chain complexes differ. All of these
concepts will be introduced in the current section.

\subsection{Graded and filtered module homomorphisms}

\begin{figure}
\begin{center}
  \includegraphics[height=0.40\textwidth]{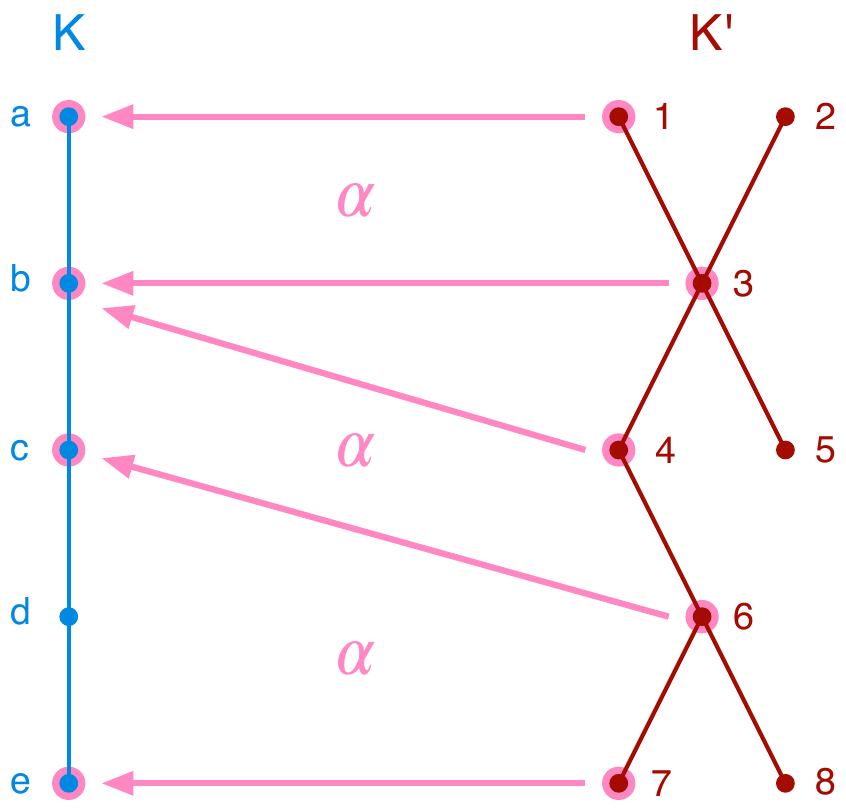}\quad\quad\quad
  \includegraphics[height=0.40\textwidth]{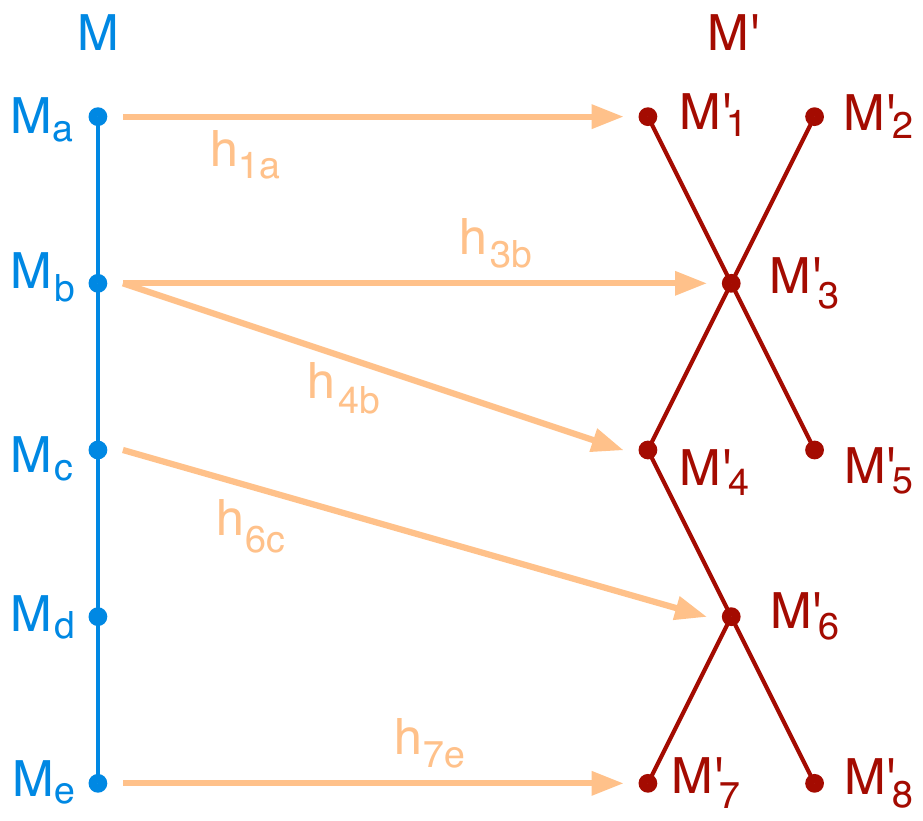}
\end{center}
  \caption{Sample $\alpha$-graded morphism. The left image shows two posets~$K$ and
           $K'$, together with an order preserving partial map~$\alpha : K' \pto K$.
           The panel on the right indicates the potential nonzero module
           homomorphisms~$h_{ij}$ if~$h$ is $\alpha$-graded.}
  \label{fig:hom-graded}
\end{figure}

\begin{figure}
\begin{center}
  \includegraphics[height=0.40\textwidth]{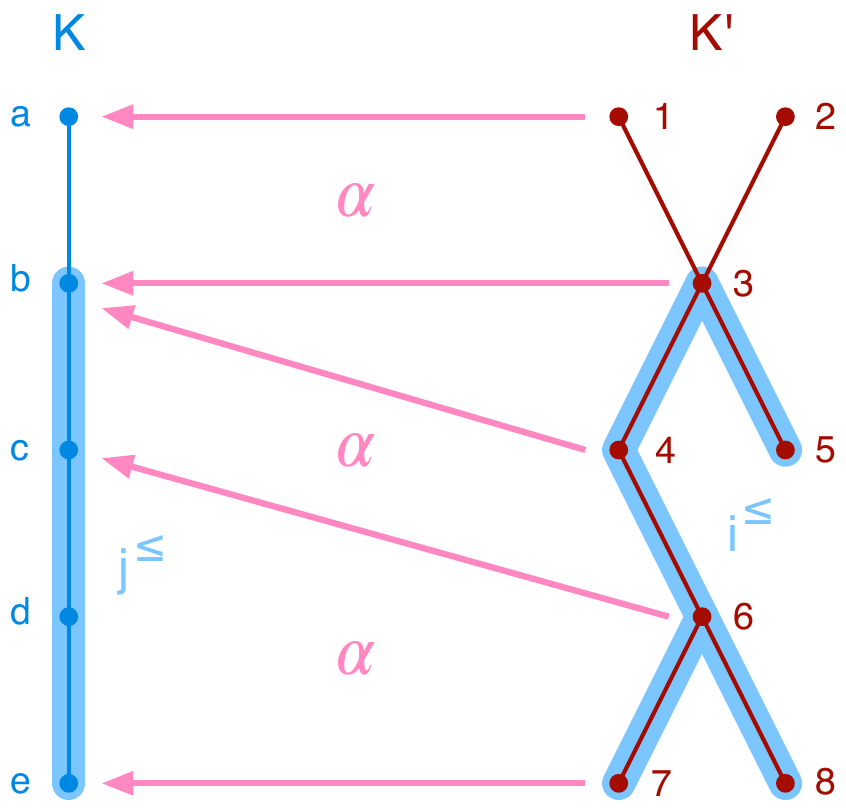}\quad\quad\quad
  \includegraphics[height=0.40\textwidth]{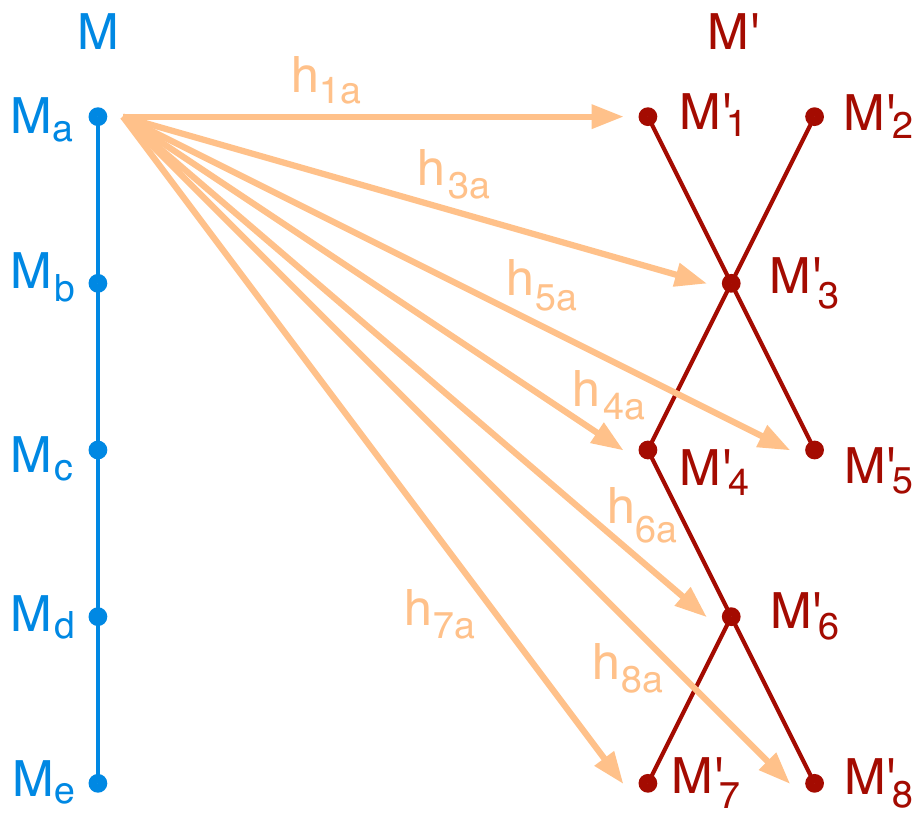}\\[4ex]
  \includegraphics[height=0.40\textwidth]{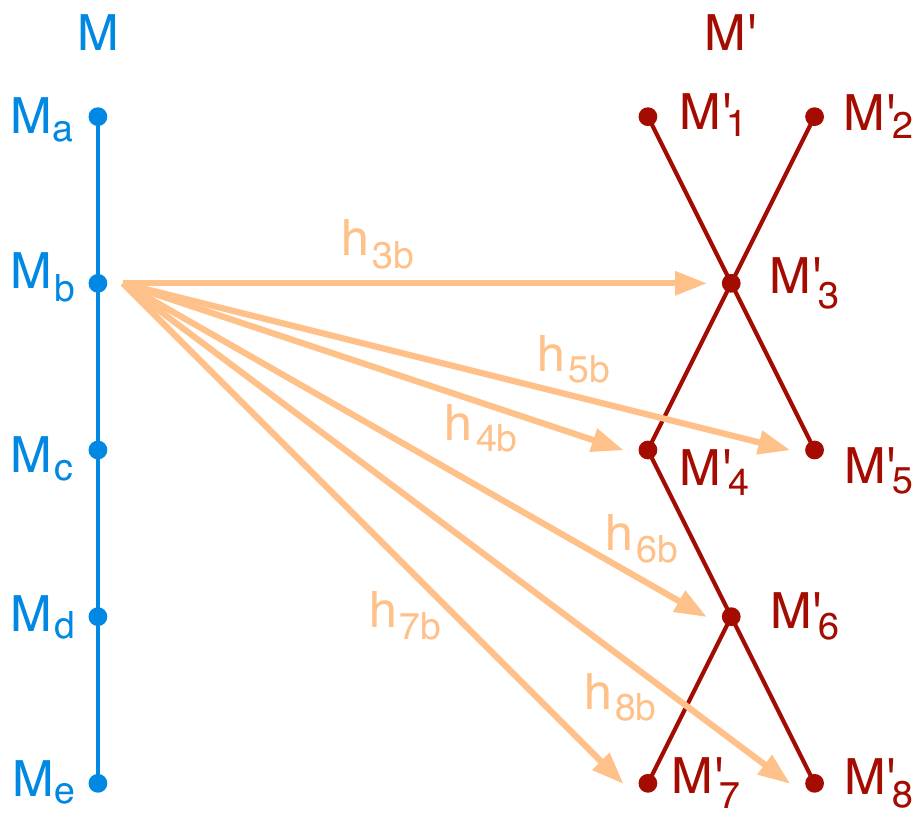}\quad\quad\quad
  \includegraphics[height=0.40\textwidth]{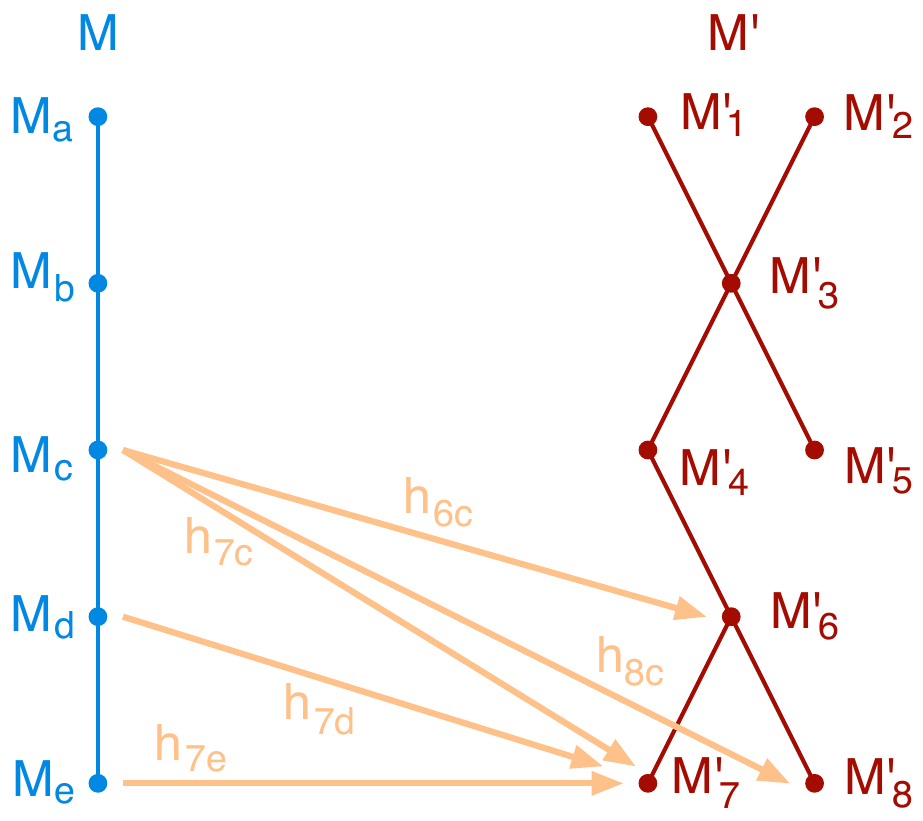}
\end{center}
  \caption{Sample $\alpha$-filtered morphism. The top left image shows two posets~$K$ and
           $K'$, together with an order preserving partial map~$\alpha : K' \pto K$.
           The remaining three panels indicate all potentially nonzero module
           homomorphisms~$h_{ij}$ if~$h$ is $\alpha$-filtered (top right for $j=a$, 
           bottom left for $j=b$ and bottom right for $j=c,d,e$).
           }
  \label{fig:hom-filtered}
\end{figure}

We begin our discussion on the level of general graded modules.
Consider a $K$-graded module~$M$ and a $K'$-graded module~$M'$.
Let $\alpha: K'\pto K$ be a partial map.
We say that a homomorphism
$h:M\to M'$ is $\alpha$-{\em graded} if for every pair of indices
$j\in K$ and $i\in K'$ we have
\[
h_{ij}\neq 0 \quad\implies\quad i\in\alpha^{-1}(j),
\]
which is equivalent to the implication
\begin{equation}
\label{eq:grad-hom}
h_{ij}\neq 0 \quad\implies\quad i\in\dom \alpha
\ \ \ \text{and}\ \ \
\alpha(i)=j.
\end{equation}
If in addition the sets~$K$ and~$K'$ are partially ordered sets
and~$\alpha$ is order preserving, then we say that a homomorphism
$h:M\to M'$ is $\alpha$-{\em filtered} if for every $j\in K$ and $i\in K'$
\begin{equation}
\label{eq:filt-hom}
h_{ij}\neq 0 \quad\implies\quad i\in \alpha^{-1}(j^\leq)^\leq
\end{equation}
which is equivalent to
\[
h_{ij}\neq 0 \quad\implies\quad
\exists_{i'\in \dom\alpha}\;i\leq i'
\ \ \ \text{ and } \ \ \ \alpha(i')\leq j.
\]
These two definitions are illustrated in Figures~\ref{fig:hom-graded}
and~\ref{fig:hom-filtered}. The left image in Figure~\ref{fig:hom-graded}
shows two partially ordered sets together with a partial map~$\alpha$ between
them. In the right image of the same figure arrows indicate which of the
homomorphisms~$h_{ij}$ can be nontrivial if~$h$ is $\alpha$-graded. In 
contrast, Figure~\ref{fig:hom-filtered} depicts all possible nontrivial
homomorphisms~$h_{ij}$ if~$h$ is $\alpha$-filtered. The more involved
implication~(\ref{eq:filt-hom}) is illustrated in the top left image, 
which indicates the sets~$3^\le$ and~$b^\le$ in light blue. The remaining
three panels in the figure depict all possible nontrivial
homomorphisms~$h_{ij}$. Clearly, $\alpha$-filtered homomorphisms are
far less restricted than $\alpha$-graded ones. In fact, an $\alpha$-filtered
homomorphism can have a nontrivial~$h_{ij}$ even if~$i \not\in \dom\alpha$
or if~$j \not\in \im\alpha$.

Conditions~\eqref{eq:grad-hom} and~\eqref{eq:filt-hom} for $\alpha$-graded and
$\alpha$-filtered homomorphisms, respectively, may be simplified if $\dom\alpha = K'$,
i.e., if~$\alpha$ is actually a function. This is the subject of the following
straightforward proposition.

\begin{prop}
\label{prop:filt-grad-hom-full-dom}
Assume $\dom\alpha = K'$. Then $h:M\to M'$ is $\alpha$-graded if and only if 
for every $i\in K'$ and $j\in K$
\begin{equation}
\label{eq:grad-hom-full-dom}
   h_{ij}\neq 0 \quad\implies\quad \alpha(i)=j,
\end{equation}
and $h$ is $\alpha$-filtered if and only if 
\begin{equation}
\label{eq:filt-hom-full-dom}
   h_{ij}\neq 0 \quad\implies\quad \alpha(i)\leq j
\end{equation}
for every $i\in K'$ and $j\in K$. \qed
\end{prop}

Obviously, if $\alpha$ is order preserving and~$h$ is an $\alpha$-graded
homomorphism then~$h$ is an $\alpha$-filtered homomorphism. However, the 
converse is not true in general. This can readily be seen from 
Figures~\ref{fig:hom-graded} and~\ref{fig:hom-filtered}.

In the case when $K=K'$, we say that a homomorphism
$h:M\to M'$ is {\em $K$-filtered} (respectively {\em $K$-graded})
if $h$ is $\id_K$-filtered (respectively $\id_K$-graded).
We simplify the terminology to {\em filtered} (respectively {\em graded}) if $K$ is clear from the context.
As an immediate consequence of Proposition~\ref{prop:filt-grad-hom-full-dom}
we get the following corollary.
\begin{cor}
An endomorphism  $h:M\to M$ on a module~$M$ is graded if and only if
for every $i,j\in K$ we have
\[
h_{ij}\neq 0 \quad\implies\quad i=j,
\]
and $h$ is filtered if and only if one has
\[
h_{ij}\neq 0 \quad\implies\quad i\leq j
\]
for every $i,j\in K$. \qed
\end{cor}

Note that in  the special case when $K=K'=\ZZ$ an $\alpha$-graded homomorphism with 
$\alpha:\ZZ\ni i\mapsto i-k\in \ZZ$
given as the left-shift by~$k$ on the integers coincides with a $\ZZ$-graded homomorphism of degree $k$.
We recall that in such a case, when $k$ is clear from the context, we shorten the notation $h_{j+k,j}$ to $h_j$.

\begin{prop}
\label{prop:filt-hom}
Assume that $P$ and $P'$ are posets, $\alpha:P'\pto P$ is order preserving,
 $M$ is a $P$-graded module, and $M'$ is a $P'$-graded module.
A module homomorphism $h:M\to M'$ is an $\alpha$-filtered homomorphism
if and only if
\begin{equation}
\label{eq:filt-hom-downset}
   h(M_L)\subset M'_{\alpha^{-1}(L)^\leq}
\end{equation}
for every $L\in\Down(P)$.
\end{prop}
\proof
Assume that $h: M \to M'$ is a module homomorphism such
that~\eqref{eq:filt-hom-downset} is satisfied for every $L\in\Down(P)$.
Moreover, assume that $h_{pq}\neq 0$ for an element $q\in P$ and a $p\in P'$.
Let $x\in M_q$ be such that $h_{pq}(x)\neq 0$.
If we consider $L:= q^\leq \in\Down(P)$, then~\eqref{eq:filt-hom-downset}
implies $h(M_L)\subset M'_{\alpha^{-1}(L)^\leq}$.
Since we have the inclusions $x\in M_q\subset M_L$, we get $h(x)\in M'_{\alpha^{-1}(L)^\leq}$.
Thus, the identity $h_{p'q}(x)=\pi_{p'}(h(x))=0$ holds for $p'\not\in\alpha^{-1}(L)^\leq$,
which in turn implies that $p\in \alpha^{-1}(L)^\leq$, in view of $h_{pq}(x)=\pi_p(h(x))\neq 0$.
In consequence, $p\leq\bar{p}$ for a $\bar{p}$ with $\alpha(\bar{p})\in L= q^\leq$.
Therefore, $p\in\alpha^{-1}(q^\leq)^\leq$ which implies \eqref{eq:filt-hom} and proves that $h$ is $\alpha$-filtered.
To prove the opposite implication, assume that
 $h:M\to M'$ is a module homomorphism such that
property \eqref{eq:filt-hom} holds.
Let $L\in\Down(P)$ and consider $x\in M_L$. Without loss of generality we may
assume that $x\in M_q$ for a $q\in L$. Then an application of \eqref{eq:filt-hom} yields
\[
   h(x)=\sum_{p\in P} h_{pq}(x)=\sum_{p\in\alpha^{-1}(q^\leq)^\leq} h_{pq}(x),
\]
which means that $h(x)\in M'_{\alpha^{-1}(L)^\leq}$,
since we have $\alpha^{-1}(q^\leq)^\leq \subset \alpha^{-1}(L)^\leq$.
This establishes the inclusion in~\eqref{eq:filt-hom-downset}.
\qed

\medskip
In the special case when $P=P'$ and $\alpha=\id_P$ we have the following corollary of Proposition~\ref{prop:filt-hom}.

\begin{cor}
\label{cor:filt-hom}
   Assume that $P$ is a poset and that $M$ is a $P$-graded module.
   Then $h:M\to M$ is a filtered homomorphism if and only if
\begin{equation}
\label{eq:id-filt-hom-downset}
   h(M_L)\subset M_L
\end{equation}
for every $L\in\Down(P)$.
\end{cor}
\proof
  By applying our assumptions $P=P'$ and $\alpha=\id_P$ to the setting
  of Proposition~\ref{prop:filt-hom}, one obtains that
  $\alpha^{-1}(L)^\leq=L^\leq=L$. Therefore, the condition
  in~\eqref{eq:filt-hom-downset} reduces to the one
  in~\eqref{eq:id-filt-hom-downset}.
\qed

\begin{defn}
\label{defn:filt-equiv}
{\em
   We say that two $P$-gradations $(M_p)_{p\in P}$ and  $(M'_p)_{p\in P}$  of the same module $M$ are
{\em filtered equivalent}, if for every $L\in\Down(P)$ we have $M_L=M'_L$.
}
\end{defn}

As an immediate consequence of Corollary~\ref{cor:filt-hom} one obtains the following proposition.
\begin{prop}
\label{prop:filt-equiv}
   Let $P$ be a poset. Assume that $M$ and $M'$ are $P$-graded modules
   and that $h:M\to M'$ is a filtered homomorphism.
   Then $h$ is a filtered homomorphism with respect to any
   filtered equivalent gradation of the modules~$M$ and~$M'$.
\qed
\end{prop}

\subsection{The category of graded and filtered moduli}
We define the two categories $\GMod$ of graded moduli and $\FMod$ of filtered moduli as follows.
The objects of $\GMod$ (respectively $\FMod$) are pairs $(P,M)$
where $P$ is an object of $\DSet$ (respectively $\DPSet$) and $M$ is a $P$-graded module;
(see Section~\ref{sec:modules}).
The morphisms in $\GMod$ (respectively $\FMod$) are pairs $(\alpha,h):(P,M)\to (P',M')$
such that $\alpha:P'\pto P$ is a morphism in $\DSet$ (respectively $\DPSet)$
and $h:M\to M'$ is an $\alpha$-graded (respectively $\alpha$-filtered)
module homomorphism. In the following, we will briefly refer to morphisms in~$\GMod$
(respectively~$\FMod$) as graded (respectively filtered) morphisms $(\alpha,h):(P,M)\to (P',M')$.
In the special case when one has the identities $(P,M)=(P',M')$ and $\alpha=\id_P$
we simplify the terminology by referring to~$h$ as a filtered (respectively graded)
morphism, meaning that~$h$ is $\id_P$-filtered (respectively $\id_P$-graded).

One can easily verify that $\id_{(P,M)}:=(\id_P,\id_M)$ is the identity morphism
on~$(P,M)$ in both~$\GMod$ and~$\FMod$. Also, if $(\alpha,h): (P,M) \to (P',M')$
and $(\alpha',h'): (P',M') \to (P'',M'')$ are two graded or filtered morphisms,
then we define their composition
\[
(\alpha'',h''):=(\alpha',h')\circ (\alpha,h): (P,M) \to
               (P'',M'')
\]
by $\alpha'':=\alpha\alpha'$ and $h'':=h'h$. This leads to the following fundamental
results.
\begin{prop}
\label{prop:GMod}
   $\GMod$ is a well-defined category.
\end{prop}
\proof
   The only nontrivial part of the proof is the verification that the composition
   of two graded morphisms is again a graded morphism.
   For this, let $(\alpha'',h''):=(\alpha',h')\circ (\alpha,h)$ be given as above.
Obviously  $h'h$ is a module homomorphism.
We will show that $h''_{pr}\neq 0$ for $p\in P''$ and $r\in P$ implies the
equality $(\alpha\alpha')(p)= r$.
Hence, we assume that $h''_{pr}\neq 0$ for some $p\in P''$ and $r\in P$.
It follows from Proposition~\ref{prop:composition} that
there exists a $q\in P'$ such that both $h'_{pq}\neq 0$ and $h_{qr}\neq 0$
are satisfied. This further shows that one has to have $p\in \dom\alpha'$ and $\alpha'(p)= q$,
as well as $q\in \dom\alpha$ and $\alpha(q)= r$.
It follows that $p\in\dom(\alpha\alpha')$ and
 $(\alpha\alpha')(p)= r$.
This proves that $(\alpha'',h'')$ is a graded morphism.
\qed

\medskip
A given pair of poset graded morphisms $f: (P,M)\to (P,M)$
and $f': (P',M')\to (P',M')$ are called {\em graded-conjugate}, if there exists
a poset graded isomorphism $(\alpha,h):(P,M)\to(P',M')$ such that
\[
  (\alpha,h)\circ f=f'\circ (\alpha,h).
\]
If the morphisms $f$ and $f'$ are graded-conjugate, then we say that their
respective $(P,P)$- and $(P',P')$-matrices are {\em graded-similar}.
It is easily seen that the graded similarity of matrices
is an equivalence relation.

\begin{prop}
\label{prop:FMod}
   $\FMod$ is a well-defined category.
\end{prop}
\proof
   The only nontrivial property to verify is that
the composition of filtered morphisms is again a filtered morphism.
  For this, consider a composition $(\alpha'',h''):=(\alpha',h')\circ (\alpha,h)$
  of filtered morphisms as above.
Obviously, the map $\alpha\alpha'$ is order preserving and $h'h$ is a module homomorphism.
Thus, we only need to verify that
that $h''_{pr}\neq 0$ for $p\in P''$ and $r\in P$ implies the inclusion $p\in(\alpha\alpha')^{-1}(r^\leq)^\leq$.
Hence, assume that $h''_{pr}\neq 0$ for $p\in P''$ and $r\in P$.
It follows from Proposition~\ref{prop:composition} that
there exists a $q\in P'$ such that $h'_{pq}\neq 0$ and $h_{qr}\neq 0$.
Therefore, there exist $\bar{p}\geq p$ and $\bar{q}\geq q$
such that
$\alpha'(\bar{p})\leq q$ and $\alpha(\bar{q})\leq r$.
Since $\alpha$ is order preserving, it follows that $(\alpha\alpha')(\bar{p})\leq \alpha(\bar{q})\leq r$.
Hence, $\bar{p}\in(\alpha\alpha')^{-1}(r^\leq)$
and $p\in (\alpha\alpha')^{-1}(r^\leq)^\leq$.
This proves that $(\alpha',h'')$ is a filtered morphism.
\qed

\begin{prop}
\label{prop:FMod-homo}
Assume that the two maps $(\alpha,h):(P,M)\to (P',M')$ and $(\alpha',h'):(P',M')\to (P'',M'')$
are morphisms in $\FMod$ which satisfy the identities $\dom\alpha=P'$ and
$\dom \alpha'=P''$. Suppose further that~$\alpha$ is injective. 
If $p\in P$ and $p''\in P''$ are such that $\alpha\alpha'(p'')=p$, then
\begin{equation}
\label{eq:FMod-homo-2}
(h'h)_{p''p}=h'_{p''\alpha'(p'')}h_{\alpha'(p'')p}\quad \text{for $p\in P$.}
\end{equation}
\end{prop}
\proof
We get from Proposition~\ref{prop:composition} that
\begin{equation}
\label{eq:FMod-homo-3}
  (h'h)_{p''p}=\sum_{p'\in P'} {h'}_{p''p'} h_{p'p}.
\end{equation}
It follows from \eqref{eq:filt-hom-full-dom}   that the index $p'$ of a non-zero term
on the right-hand side
of \eqref{eq:FMod-homo-3} must satisfy $\alpha'(p'')\leq p'$
and $\alpha(p')\leq p$.
Hence, 
\[
   p=\alpha\alpha'(p'')\leq\alpha(p')\leq p.
\]
It follows  that $\alpha(p') = p$ and, since $\alpha$ is injective, we get $p'=\alpha'(p'')$.
This proves \eqref{eq:FMod-homo-2}.
\qed

\medskip
As an immediate consequence of Proposition~\ref{prop:FMod-homo} one obtains the following
corollary. 
\begin{cor}
\label{cor:FMod-homo}
Assume that the two maps $(\alpha,h):(P,M)\to (P',M')$ and $(\alpha',h'):(P',M')\to (P,M)$
are morphisms in $\FMod$ such that  $\alpha:P'\to P$ and $\alpha':P \to P'$ are mutually inverse
bijections. Then
\begin{equation*}
(h'h)_{pp}=h'_{p\alpha'(p)}h_{\alpha'(p)p}\quad \text{for $p\in P$.}
\end{equation*}
\end{cor}

\begin{lem}
\label{lem:FMod-iso}
A filtered morphism $(\alpha,h):(P,M)\to (P',M')$ is an isomorphism in $\FMod$
if and only if the map $\alpha:P'\to P$ is an order isomorphism and
$h_{p'\alpha(p')}:M_{\alpha(p')}\to M'_{p'}$ is a module isomorphism for every $p'\in P'$.
\end{lem}
\proof
  First assume that $(\alpha,h)$ is an isomorphism in the category~$\FMod$.
Let $(\alpha',h'):(P',M')\to (P,M)$ be the inverse of  $(\alpha,h)$.
Then $\alpha\alpha'=\id_P$ and $\alpha'\alpha=\id_{P'}$.
Since $\alpha$ and $\alpha'$ are order preserving, we see that  $\alpha: P'\to P$ is an order isomorphism.
It follows from Corollary~\ref{cor:FMod-homo} that
\begin{equation}
\label{eq:FMod-iso-1}
   \id_{M_p}=(\id_M)_{pp}=(h'h)_{pp}=h'_{p\alpha'(p)}h_{\alpha'(p)p} \text{ for $p\in P$}
\end{equation}
and
\begin{equation}
\label{eq:FMod-iso-2}
   \id_{M_{p'}}=(\id_M)_{p'p'}=(hh')_{p'p'}=h_{p'\alpha(p')}h'_{\alpha(p')p'} \text{ for $p'\in P'$.}
\end{equation}
But, $\alpha$ and $\alpha'$ are mutually inverse bijections.
Hence, substituting $p':=\alpha'(p)$ in \eqref{eq:FMod-iso-1}
we get $p=\alpha(p')$ and
\[
     \id_{M_{\alpha(p')}}=h'_{\alpha(p')p'}h_{p'\alpha(p')} \text{ for $p'\in P'$.}
\]
It follows that $h'_{\alpha(p')p'}$ is the inverse of $h_{p'\alpha(p')}$.
Therefore, $h_{p'\alpha(p')}$ is a module isomorphism for every $p'\in P'$.

To see the opposite implication, assume that $(\alpha,h):(P,M)\to (P',M')$
is a filtered homomorphism such that $\alpha:P'\to P$ is an order isomorphism and 
$h_{p'\alpha(p')}:M_{\alpha(p')}\to M'_{p'}$ is a module isomorphism for all $p'\in P'$.
Since~$P$ and~$P'$ as objects of~$\DPSet$ are finite sets, we may proceed by
induction on their cardinality $n:=\card P=\card P'$.
If $n=1$, then $P=\{p\}$, $P'=\{p'\}$ and $h_{p'\alpha(p')}=h_{p'p}=h$
is a module isomorphism.
Clearly, its inverse is a filtered homomorphism.
Hence, $h$ is an isomorphism in $\FMod$. Thus, we now assume $n>1$.
Let $\bar{p}'$ be a maximal element in $P'$.
Set $\bar{p}:=\alpha(\bar{p}')$.
Since~$\alpha'$ is an order isomorphism, we see that $\bar{p}$ is a maximal element in $P$.
Let $\bar{P}':=P'\setminus\{\bar{p}'\}$, $\bar{P}:=P\setminus\{\bar{p}\}$,
$\bar{M}':=M_{\bar{P}'}$, $\bar{M}:=M_{\bar{P}}$. Let $Q:=\{\bar{P},\bar{p}\}$
and  $Q':=\{\bar{P}',\bar{p}'\}$ be linearly ordered respectively by $\bar{P}<\bar{p}$ and $\bar{P}'<\bar{p}'$.
Then $(Q,M)$ and $(Q',M')$ are graded modules.
Define $\bar{\alpha}:Q'\to Q$ through the identities $\bar{\alpha}(\bar{P}'):=\bar{P}$ and $\bar{\alpha}(\bar{p}'):=\bar{p}$.
We will prove that
\begin{equation}
\label{eq:FMod-iso-3}
h_{\bar{p}'p}=0 \text{ for $p\in\bar{P}.$}
\end{equation}
Arguing by contradiction, assume that $h_{\bar{p}'p}\neq 0$
for a $p\in\bar{P}$. Since $h$ is an $\alpha$-filtered homomorphism,
we get from \eqref{eq:filt-hom}
 that $h_{\bar{p}'p}\neq 0$ implies
$\bar{p}'\leq \bar{\bar{p}}'$ for some $\bar{\bar{p}}'\in\dom\alpha$
such that $\alpha(\bar{\bar{p}}')\leq p$. Since $\bar{p}'$ is maximal
in $P'$, we obtain both $\bar{p}'=\bar{\bar{p}}'$ and $\bar{p}=\alpha(\bar{p}')\leq p$.
Since $\bar{p}$ is maximal in $P$, this further implies the equality
 $p=\bar{p}$, a contradiction proving \eqref{eq:FMod-iso-3}.
Hence, the identity in~\eqref{eq:matrix-rep} yields $h_{\bar{p}'\bar{P}}=\sum_{p\in\bar{P}}h_{\bar{p}'p}\circ\pi_p=0$.
Therefore, the $(Q,Q')$-matrix of $h$ is
\[
\left[
  \begin{array}{cc}
    h_{\bar{P}'\bar{P}} & h_{\bar{P}'\bar{p}} \\
    0    & h_{\bar{p}'\bar{p}}
  \end{array}
\right].
\]
By induction assumption $(\alpha_{|\bar{P}'}, h_{\bar{P}'\bar{P}}):M_{\bar{P}}\to M_{\bar{P}'}$
is an isomorphism in $\FMod$.
Set $\alpha':=\alpha^{-1}$ and let $h': M'\to M$ be the  homomorphism given by the $(Q',Q)$-matrix
\[
\left[
  \begin{array}{cc}
    h_{\bar{P}'\bar{P}}^{-1} & -h_{\bar{P}'\bar{P}}^{-1}h_{\bar{P}'\bar{p}}h_{\bar{p}'\bar{p}}^{-1} \\
    0    & h_{\bar{p}'\bar{p}}^{-1}
  \end{array}
\right].
\]
One easily verifies that $h'$ is $\alpha'$-filtered and a straightforward computation shows that $(\alpha',h')\circ(\alpha,h)=\id_{(P,M)}$
and $(\alpha,h)\circ(\alpha',h')=\id_{(P',M')}$.
Hence, the morphism $(\alpha,h)$ is indeed an isomorphism in $\FMod$.
\qed

\begin{cor}
\label{cor:FMod-iso}
Assume $(\alpha,h):(P,M)\to(P',M')$ is both a homomorphism in~$\GMod$
and an isomorphism in~$\FMod$. Then it is automatically an isomorphism in~$\GMod$.
\end{cor}
\proof
   Due to Lemma~\ref{lem:FMod-iso} the map $\alpha:P'\to P$ is an order
   isomorphism and $h_{p'\alpha(p')}:M_{\alpha(p')}\to M'_{p'}$
   is a module isomorphism. Let $g:M'\to M$ be the $\alpha^{-1}$-graded
   homomorphism with  $(P',P)$-matrix given by
 \[
  g_{pp'}:=\begin{cases}
             h_{p'\alpha(p')}^{-1} & \text{ if $p=\alpha(p')$,}\\
             0 & \text{ otherwise.}
           \end{cases}
 \]
  It is straightforward to verify that $(\alpha^{-1},g)$ is the inverse of
  $(\alpha,h)$ in $\GMod$. Hence, $(\alpha,h)$ is
  an isomorphism in $\GMod$.
\qed

\subsection{Poset filtered chain complexes}

Let~$P$ be an arbitrary finite set and~$(C,d)$ a chain complex. We call~$(C,d)$
a {\em $P$-graded chain complex}, if~$C$ is a $P$-graded module in which each
$C_p\subset C$ is a $\ZZ$-graded submodule of~$C$. Since the chain complex~$C$
is $P$-graded, its boundary homomorphism $d$  has a $(P,P)$-matrix.
A partial order $\leq$ in $P$ is called $(C,d)$-admissible, or briefly {\em $d$-admissible},
if $d$ is a filtered homomorphism with respect to $(P,\leq)$, i.e., if
\begin{equation}
\label{eq:admissible}
d_{pq}\neq 0\implies p\leq q
\end{equation}
for all $p,q\in P$.
Note that a $(C,d)$-admissible partial order on $P$ may not always exist.
However, we have the following straightforward proposition and definition.
\begin{propdef}
\label{prop:d-admissible}
If a given $P$-graded chain complex $(C,d)$ admits a $(C,d)$-admissible
partial order on $P$, then the intersection of all such $(C,d)$-admissible
partial orders on~$P$ is again a $(C,d)$-admissible partial order on~$P$.
We call it the {\em native} partial order of~$d$.
\qed
\end{propdef}

\begin{defn}
\label{defn:poset-filtered-cc}
{\em
We say that the triple $(P,C,d)$ is a {\em poset filtered chain complex}
if $P$ is a poset, $(C,d)$ is a $P$-graded chain complex and the partial order in $P$
is $d$-admissible.
}
\end{defn}
Note that for a poset filtered chain complex $(P,C,d)$ the module $C$ is not only $P$-graded
but also $\ZZ$-graded where the $n$th summand of the $\ZZ$-gradation is the direct sum
over $p\in P$  of the $n$th summands in the $\ZZ$-gradation of $C_p$.

\begin{ex}
\label{ex:pfcc}
{\em 
  Consider the set of words 
 \[
    X:=\setof{\mathbf{A},\mathbf{B},\mathbf{C},\mathbf{AB},\mathbf{AC},
              \mathbf{BC},\mathbf{CD},\mathbf{CE},\mathbf{ABC}}
 \] 
 and the free module $C:=\ZZ_2\scalprod{X}$ with basis $X$ and
 coefficients in the field~$\ZZ_2$. 
 For a word $x\in X$ define its {\em dimension} as 
one less than the number of characters in $x$
and denote the set of words of dimension~$i$ by~$X_i$.
Setting
\[
  C_i:=\begin{cases}
            \ZZ_2\scalprod{X_i} & \text{ if $X_i\neq\emptyset$,}\\
            0 & \text{ otherwise.}
           \end{cases}
\]
for $i\in\ZZ$ we obtain a $\ZZ$-gradation of $C$ given by
\begin{equation}
\label{ex:Z-gradation}
   C=\bigoplus_{i\in\ZZ} C_i.
\end{equation}
Let $P=\{\xa,\xb,\xc,\xd,\xe,\xf\}$ be a poset with its partial order $\leq$ defined by the Hasse diagram
\begin{equation}
\label{ex:hasse0}
   \begin{diagram}
  \dgARROWLENGTH 1.1em 
    \node{\xd}
    \arrow{se,-}
    \node[2]{\xe}
    \arrow{sw,-}\\
    \node[2]{\xf}
    \arrow{s,-}\\
    \node[2]{\xc}
    \arrow{se,-}
    \arrow{sw,-}\\
    \node{\xa}
    \node[2]{\xb}
  \end{diagram}.
\end{equation}
For $p=\xa,\xb,\xc,\xd,\xe,\xf\in P$ we define $X_p$ respectively as
\[
\{\mathbf{A}\}, \; \{\mathbf{B}\}, \; \{\mathbf{AB}\}, \;
\{\mathbf{CD}\}, \; \{\mathbf{CE}\}, \;
\{\mathbf{C},\mathbf{AC},\mathbf{BC},\mathbf{ABC}\}.
\]
Then $\{X_p\}_{p\in P}$ is a partition of $X$ which induces a gradation 
\begin{equation}
\label{eq:ex0}
C=\bigoplus_{p\in P} C_p
\end{equation}
with $C_p:=\ZZ_2\scalprod{X_p}$.
Consider the homomorphism $d:C\to C$ of degree~$-1$ which is defined on the basis~$X$ by the matrix
\[
\begin{array}{c||c|c|c|cccc|c|c|}
    &  \mathbf{A} &  \mathbf{B} &    \mathbf{AB} &  \mathbf{C} &    \mathbf{AC} &    \mathbf{BC} & \mathbf{ABC}   & \mathbf{CD}    & \mathbf{CE}    \\
  \hline
  \hline
  \mathbf{A} &    &    & \cg 1 &    & \cg 1 &       & \cg   & \cg   & \cg   \\
  \hline                                                                                    
  \mathbf{B} &    &    & \cg 1 &    & \cg   & \cg 1 & \cg   & \cg   & \cg   \\
  \hline
 \mathbf{AB} &    &    &       &    & \cg   & \cg   & \cg 1 & \cg   & \cg   \\
  \hline
  \mathbf{C} &    &    &       &    &     1 &     1 &       & \cg 1 & \cg 1 \\
 \mathbf{AC} &    &    &       &    &       &       &  1    & \cg   & \cg   \\
 \mathbf{BC} &    &    &       &    &       &       &  1    & \cg   & \cg   \\
\mathbf{ABC} &    &    &       &    &       &       &       & \cg   & \cg   \\
  \hline
 \mathbf{CD} &    &    &       &    &       &       &       &       &       \\
  \hline
 \mathbf{CE} &    &    &       &    &       &       &       &       &       \\
  \hline
\end{array}.
\]
One can immediately verify that $d$ is $P$-filtered and $d^2=0$. Therefore,
the triple $(P,C,d)$ with the $\ZZ$-gradation \eqref{ex:Z-gradation}
and the $P$-gradation \eqref{eq:ex0} is a well-defined poset filtered
chain complex.
\exend}
\end{ex}

We would like to point out that if $J\subset P$, then $(C_J,d_{J J})$ does not
need to be a chain complex in general. However, we have the following proposition.
\begin{prop}
\label{prop:M-convex}
Assume that~$J$ is a convex subset of~$P$. Then we have:
\begin{itemize}
   \item[(i)] $(C_J,d_{J J})$ is a chain complex.
   \item[(ii)] $(C_J,d_{J J})$ is chain isomorphic to the quotient
   complex $(C_{J^\leq}/C_{J^<},d')$, where $d'$ denotes the homomorphism
   induced by $d_{J^\leq J^\leq}$. Notice that the quotient complex is
   well-defined since~$J^<$ is a down set due to the convexity of~$J$.
\end{itemize}
\end{prop}
\proof
  In order to prove (i) we need to verify that $d_{J J}^2=0$.
Assume first  that $J$ is a down set.  Let $x\in C_J$.
Then Corollary~\ref{cor:filt-hom} immediately implies that $dx\in C_J$.
Hence, one obtains both $d_{JJ}x=(\iota_J\circ d\circ\pi_J)(x)=dx$ and
$d_{JJ}^2x=d_{J J}dx=d^2x=0$, which yields $d_{JJ}^2=0$.
If $J$ is just convex, we consider the down sets $I:=J^<$ and $K:=J^\leq$,
which clearly satisfy the identity $K=I\cup J$.
Since $d$ is a filtered homomorphism, we have $d_{pq}=0$
for $p\in J$ and $q\in I $. Therefore,
the matrix of $d_{KK}$ takes the form
\[
\left[
  \begin{array}{cc}
   d_{I I } & d_{I J } \\
    0    & d_{J J}
  \end{array}
\right].
\]
Since $K $ is a down set, we already have verified that $d_{KK}^2=0$.
It follows that
\[
  0=d_{K  K }^2=\left[
  \begin{array}{cc}
   d_{I I }^2 & d_{I I }d_{I J }+d_{I J } d_{J J}\\
    0    & d_{J J}^2.
  \end{array}
\right]
\]
Thus, $d_{J J}^2=0$, which proves (i).

To establish (ii), we consider again the down sets $I:=J^<$ and $K:=J^\leq$, as well
as the homomorphism $\kappa: C_J\ni x\mapsto [x]_I\in C_K/C_I$, where $[x]_I$ denotes
the equivalence class of $x$ in the quotient module $ C_K/C_I$.
Let $x\in C_J$ and let~$\bar{d}_{KK}$ denote the boundary homomorphism
induced by $d_{KK}$ on the quotient module $C_K/C_I$. We then have
$$\bar{d}_{KK}\kappa x=[dx]_I=[d_{IJ}x+d_{JJ}x]_I=[d_{JJ}x]_I=\kappa d_{JJ}x,$$
and this implies that the map $\kappa$ is a chain map.
Assume now that $[x]_I=0$ for an $x\in C_J$. Then $x\in C_I\cap C_J=\{0\}$, that is, $x=0$.
This proves that~$\kappa$ is a monomorphism.
Finally, given the class $[y]_I\in C_K/C_I$ generated by $y\in C_K$, we have $y=y_I+y_J$ for a $y_I\in C_I$
and a $y_J\in C_J$. It follows that $[y]_I=[y_J]_I=\kappa y_J$, proving that
$\kappa$ is an epimorphism. Hence, $\kappa$ is an isomorphism.
\qed

\medskip
It is straightforward to observe that the chain complex $(C_J,d_{J J})$ in
Proposition~\ref{prop:M-convex}(i) is in fact $J$-filtered.
Therefore, we have the following proposition.
\begin{prop}
\label{prop:induced-chain-complex}
Given a filtered chain complex $(P,C,d)$ and a convex subset $J\subset P$, the triple $(J,C_J,d_{J J})$
is also a filtered chain complex. We call it the {\em filtered chain complex induced by a convex subset}.
\qed
\end{prop}

Since a singleton $\{p\}\subset P$   is always convex, one immediately obtains
the following corollary.
\begin{cor}
\label{cor:M-convex-sing}
Let~$(P,C,d)$ be an arbitrary poset filtered chain complex. Then the pair $(C_p,d_{pp})$
is a chain complex for every $p\in P$.
\qed
\end{cor}

Given a poset filtered chain complex $(P,C,d)$ we
consider the poset $P$ as an object of $\DPSet$ with the distinguished subset $P_\star$ satisfying
\begin{equation}
\label{eq:P-star-condition}
\setof{p\in P\mid \text{ $C_p$ is homotopically essential}}\subset P_\star.
\end{equation}
The concrete choice of $P_\star$  depends on applications.
Typically, the simplest choice which gives equality in \eqref{eq:P-star-condition}
suffices. In other situations, however, we might take the other extreme case
$P_\star=P$, and in this case we call the poset filtered chain complex
$(P,C,d)$ {\em peeled}. 

In terms of applications to dynamics, a choice of a larger~$P_\star$ can
mean that some Morse sets may only be known via isolating neighborhoods
and potentially may have zero Conley index. Recall that~$C_p$ is
homotopically essential if it is not chain homotopic to the zero
chain complex. 

We note that in Example~\ref{ex:pfcc} the chain complex $(C_p,d_{pp})$
is essential for  $p=\xa,\xb,\xc,\xd,\xe$ and inessential
for $p=\xf$.  Therefore, for the poset $P$ in Example~\ref{ex:pfcc}
there are only two possible choices of $P_\star$, namely $P_\star=\{\xa,\xb,\xc,\xd,\xe\}$ and $P_\star=P$.

\subsection{The category of filtered and graded chain complexes}

As our next step, we study morphisms between poset filtered chain complexes.
For this, consider poset filtered chain complexes $(P,C,d)$ and $(P',C',d')$.
Then we say that the map $(\alpha,h):(P,C,d)\to (P',C',d')$ is a {\em filtered
chain morphism}, if  $h:(C,d)\to (C',d')$ is a chain map and 
$h$ is a $\alpha$-filtered, that is $(\alpha,h):(P,C)\to (P',C')$ is a morphism in $\FMod$.
We say that the map~$(\alpha,h)$ is a
{\em graded chain morphism} if $h:(C,d)\to (C',d')$ is a chain map
and $h$ is a $\alpha$-graded, that is $(\alpha,h):(P,C)\to (P',C')$ is a morphism in $\GMod$.

We define the category $\PfCC$ of poset filtered chain complexes
by taking poset filtered chain complexes as objects and filtered chain morphisms
as morphisms. One easily verifies that this is indeed a category.
We also define the subcategory $\PgCC$ of poset graded chain complexes
by taking the same objects as in $\PfCC$ and graded chain morphisms as morphisms.

\begin{defn}
\label{defn:chain-equivalence}
{\em
Let 
$
(\alpha,\varphi),(\alpha',\varphi'):(P,C,d)\to (P',C',d')
$
denote a pair of filtered chain morphisms, and let $(\gamma,\Gamma):(P,C)\to (P',C')$
be a filtered module morphism. Then $(\gamma,\Gamma)$  is called an
{\em elementary filtered chain homotopy} between $(\alpha,\phi)$ and
$(\alpha',\phi')$ if the following conditions are satisfied:
\begin{itemize}
\item[(i)] $\Gamma$ is a module homomorphism of degree $+1$ with respect to the $\ZZ$-gradation of $C$ and $C'$,
\item[(ii)] $\Gamma$ is a chain homotopy between $\phi$ and $\phi'$, that is, $\phi'-\phi=\Gamma d+d'\Gamma$,
\item[(iii)] we have the equalities $\alpha_{|P'_\star}=\gamma_{|P'_\star}=\alpha'_{|P'_\star}$.
\end{itemize}
We say that two filtered chain morphisms $(\alpha,h)$ and $(\alpha',h')$
are {\em elementarily filtered chain homotopic}, and we write $(\alpha,h)\sim_e(\alpha',h')$,
if  there exists an elementary filtered chain homotopy between $(\alpha,h)$ and $(\alpha',h')$.
We say that filtered chain morphisms $(\alpha,h)$, $(\alpha',h')$
are {\em filtered chain homotopic}, and we write $(\alpha,h)\sim(\alpha',h')$,
if  there exists a sequence
\begin{equation} \label{defn:chain-equivalence1}
  (\alpha,h)=(\alpha_0,h_0) \; \sim_e \; (\alpha_1,h_1) \; \sim_e
  \; \ldots \; \sim_e \; (\alpha_n,h_n)=(\alpha',h')
\end{equation}
of filtered chain morphisms such that successive pairs are elementarily filtered chain homotopic.
}
\end{defn}
The following proposition is straightforward.
\begin{prop}
\label{prop:chain-equivalence}
The relation $\sim$ in the set of morphisms from $(P,C,d)$ to $(P',C',d')$ in $\PfCC$
is an equivalence relation.
\qed
\end{prop}

In certain situations, a filtered chain homotopy automatically is
an elementary filtered chain homotopy. To state this result,
recall that a poset filtered chain complex $(P,C,d)$ is
{\em peeled} if we have $P_\star=P$. Then the following holds.
\begin{prop}
\label{prop:filt-homot-essential}
Assume that the poset filtered chain complex $(P',C',d')$ is peeled.
If~$(\alpha,h), (\alpha',h'):(P,C,d)\to (P',C',d')$ are filtered
chain homotopic, then $\alpha=\alpha'$, the domain of $\alpha=\alpha'$ is~$P'$,
and we have in fact that $(\alpha,h)\sim_e(\alpha',h')$.
\end{prop}
\proof
Choose $n+1$ filtered morphisms $(\alpha_i,h_i):(P,C,d)\to (P',C',d')$
as in~(\ref{defn:chain-equivalence1}) for $i = 0,\ldots,n$.
Let $(\gamma_i,\Gamma_i):(P,C)\to (P',C')$ for $i\in\{1,2,\ldots,n\}$
be an elementary filtered chain homotopy between
$(\alpha_{i-1},h_{i-1})$ and $(\alpha_i,h_i)$.
Since $P'_\star=P'$, we have $\alpha_{i-1}=(\alpha_{i-1})_{|P'_\star}=
(\gamma_i)_{|P'_\star}=\gamma_i$, as well as the identity
$\alpha_{i}=(\alpha_{i})_{|P'_\star}=(\gamma_i)_{|P'_\star}=\gamma_i$
for  $i\in\{1,2,\ldots,n\}$. It follows then that one has $\alpha_i=\alpha_{i-1}$
for all~$i$, and thus $\alpha=\alpha'$.
Moreover, we obtain
\begin{equation}
\label{eq:filt-homot-essential}
   h_{i}-h_{i-1}=d'\Gamma_i+\Gamma_i d
\end{equation}
for  $i\in\{1,2,\ldots,n\}$.
Let $\Gamma:=\sum_{i=1}^n\Gamma_i$.
One easily verifies that $\Gamma$ is an $\alpha$-filtered module homomorphism.
Summing \eqref{eq:filt-homot-essential} for $i = 1,\ldots,n$ we get
\[
   h'-h=h_n-h_0=d'\Gamma+\Gamma d.
\]
Thus, $(\alpha,\Gamma)$  is an elementary filtered chain homotopy between
the filtered morphisms~$(\alpha,h)$ and~$(\alpha',h')$.
\qed

\subsection{Homotopy category of filtered chain complexes}

We refer to the equivalence classes of $\sim$ as the {\em homotopy
equivalence classes}. They enable us to introduce the {\em homotopy category}
of poset filtered chain complexes, denoted by $\ChPfCC$,
by taking poset filtered chain complexes as objects,
homotopy equivalence classes of filtered chain morphisms in~$\PfCC$
as morphisms in~$\ChPfCC$, and using the formula
\begin{equation}
\label{eq:chpfcc-comp}
   [(\beta,g)]_\sim\circ[(\alpha,h)]_\sim:=[(\alpha\circ\beta,g\circ h)]_\sim
\end{equation}
for any two given filtered morphisms $(\alpha,h)\in \PfCC((P,C,d),(P',C',d'))$ and
$(\beta,g)\in \PfCC((P',C',d'),(P'',C'',d''))$
as the definition of composition of morphisms in~$\ChPfCC$. Finally,
equivalence classes of identities in~$\PfCC$ are identities in~$\ChPfCC$.
The following result verifies that these definitions indeed lead to
a category.

\begin{prop}
\label{prop:filt-homotopy-composition}
Assume that the filtered chain morphisms
\[
(\alpha,h),(\alpha',h'):(P,C,d)\to (P',C',d')
\]
and
\[
(\beta,g),(\beta',g'):(P',C',d')\to (P'',C'',d'')
\]
are filtered chain homotopic. Then the compositions
\[
(\beta,g)\circ(\alpha,h),(\beta',g')\circ(\alpha',h'):(P,C,d)\to (P'',C'',d'')
\]
are also filtered chain homotopic.
In particular, the category~$\ChPfCC$ is well-defined.
\end{prop}
\proof
Assume that $(\alpha,h)\sim(\alpha',h')$ and $(\beta,g)\sim(\beta',g')$.
Since $\sim$ is an equivalence relation, it suffices to prove that
$(\beta,g)\circ(\alpha,h)\sim(\beta,g)\circ(\alpha',h')$
and $(\beta,g)\circ(\alpha',h')\sim(\beta',g')\circ(\alpha',h')$.
For the same reason, we may assume that $(\alpha,h)\sim_e(\alpha',h')$ and $(\beta,g)\sim_e(\beta',g')$.
Let $(\gamma,S):(P,C)\to (P',C')$ be an elementary filtered chain homotopy between $(\alpha,h)$ and $(\alpha',h')$.
Consider the filtered morphism $(\eta,T):=(\beta,g)\circ(\gamma,S):(P,C)\to (P'',C'')$.
Then $\eta=\gamma\beta$ and $T=gS$. We will prove that $(\eta,T)$ is an elementary filtered chain homotopy
between $(\beta,g)\circ(\alpha,h)$ and $(\beta,g)\circ(\alpha',h')$.
For this, we need to verify properties (i)-(iii) of Definition~\ref{defn:chain-equivalence}.
Clearly, $T$ is a degree $+1$ $\ZZ$-graded module homomorphism.
Hence, property (i) is satisfied.
To see property (ii), observe that
\begin{equation*}
   gh'-gh=g(h'-h)=g(d'S+Sd)=d''gS+gSd=d''T+Td.
\end{equation*}
Finally, in view of the inclusion $\beta(P''_\star) \subset P'_\star$, which is a consequence of
the definition of~$\DPSet$, one obtains
\begin{equation*}
   \alpha\beta_{|P''_\star}=\alpha_{|P'_\star}\beta_{|P''_\star}=\gamma_{|P'_\star}\beta_{|P''_\star}=
   \eta_{|P''_\star}=\gamma_{|P'_\star}\beta_{|P''_\star}=\alpha'_{|P'_\star}\beta_{|P''_\star}= \alpha'\beta_{|P''_\star}.
\end{equation*}
Altogether, we have verified that $(\beta,g)\circ(\alpha,h) \sim_e
(\beta,g)\circ(\alpha',h')$, which in turn immediately implies that also
$(\beta,g)\circ(\alpha,h)\sim(\beta,g)\circ(\alpha',h')$. The
equivalence $(\beta,g)\circ(\alpha',h')\sim(\beta',g')\circ(\alpha',h')$
can be shown similarly.
\qed

\medskip
After these preparations the following notions are immediate. We say
that a filtered chain morphism $(\alpha,h):(P,C,d)\to (P',C',d')$ is a
{\em filtered chain equivalence}, if its equivalence class~$[(\alpha,h)]$
is an isomorphism in~$\ChPfCC$. In addition, we call two poset filtered
chain complexes~$(P,C,d)$ and~$(P',C',d')$ {\em filtered chain homotopic},
if they are isomorphic in~$\ChPfCC$, that is, if there exist two
filtered chain morphisms $(\alpha,\varphi):(P,C,d)\to (P',C',d')$ and
$(\alpha',\varphi'):(P',C',d')\to (P,C,d)$ such that $(\alpha',\varphi')
\circ(\alpha,\varphi)$ is filtered chain homotopic to~$\id_{(P,C)}$, and
$(\alpha,\varphi)\circ(\alpha',\varphi')$ is filtered chain homotopic
to~$\id_{(P',C')}$. In this situation, the filtered chain
morphisms~$(\alpha,\varphi)$ and~$(\alpha',\varphi')$ are
referred to as {\em mutually inverse} filtered chain equivalences. 
If $(\alpha',\varphi')\circ(\alpha,\varphi)$ is elementarily 
filtered chain homotopic to~$\id_{(P,C)}$, and
$(\alpha,\varphi)\circ(\alpha',\varphi')$ is elementarily filtered chain homotopic
to~$\id_{(P',C')}$, we call morphisms~$(\alpha,\varphi)$ and~$(\alpha',\varphi')$
mutually inverse {\em elementary filtered chain equivalences}.

We finish this section with  the following auxiliary proposition.

\begin{prop}
\label{prop:W-gradation}
Assume that $(P,C,d)$ is a poset filtered chain complex with the $P$-gradation
$(C_p)_{p\in P}$, and let $(W_p)_{p\in P}$ denote another $P$-gradation of~$C$.
If $(C_p)_{p\in P}$ and $(W_p)_{p\in P}$ are filtered equivalent in the sense
of Definition~\ref{defn:filt-equiv}, then the triple $(P,W,d)$ with
$W:=\bigoplus_{p\in P}W_p$ and $P$-gradation given by $(W_p)_{p\in P}$ is also a
poset filtered chain complex. Moreover,
$(P,W,d)$ and $(P,C,d)$ are isomorphic in~$\PfCC$, and therefore also in~$\ChPfCC$.
\end{prop}
\proof
   Note that $W$ and $C$ are in fact the same modules, but the gradations $(W_p)_{p\in P}$ and $(C_p)_{p\in P}$
need not be the same. We know that $d$ is a filtered homomorphism with respect to the  $(C_p)_{p\in P}$
gradation of $C$.
Since the $(W_p)_{p\in P}$ gradation of $W=C$ is filtered equivalent to the~$(C_p)_{p\in P}$ gradation,
we see from Proposition~\ref{prop:filt-equiv} that $d$ is a filtered homomorphism with respect to
the $(W_p)_{p\in P}$ gradation of~$C$ as well. Hence, $(P,W,d)$ is a poset filtered chain complex.
We will prove that
\begin{equation}
\label{eq:W-gradation}
   (\id_P,\id_C):(P,C,d)\to (P,W,d)
\end{equation}
and
\begin{equation}
\label{eq:W-gradation2}
 (\id_P,\id_W):(P,W,d)\to (P,C,d)
\end{equation}
are mutually inverse isomorphisms in $\PfCC$.
Let $p\in P$. Since both~$(C_p)_{p\in P}$ and~$(W_p)_{p\in P}$  are filtered equivalent, using
Proposition~\ref{prop:M-convex}(ii) we get
\[
  C_p\cong C_{p^\leq}/C_{p^<} = W_{p^\leq}/W_{p^<}\cong W_p .
\]
Hence, $C_p$ is essential if and only if $W_p$ is essential, and from this
it follows immediately that~$W$ satisfies~\eqref{eq:P-star-condition}, because~$C$
satisfies~\eqref{eq:P-star-condition}. 

Since $(C_p)_{p\in P}$ and $(W_p)_{p\in P}$  are filtered equivalent, we get from
Corollary~\ref{cor:filt-hom} that $(\id_P,\id_C):(P,C)\to (P,W)$ and $(\id_P,\id_W):(P,W)\to (P,C)$
are filtered morphisms. Clearly, both are chain maps. It follows that the morphisms
\eqref{eq:W-gradation} and~\eqref{eq:W-gradation2} are well-defined. Since they are
also mutually inverse, the conclusion follows.
\qed


\section{Algebraic connection matrices}
\label{sec:algconnmatrix}

After the preparations of the previous section, we are now in a position 
to introduce our notion of connection matrices in a purely algebraic way.
While the definition is modeled on previous work by Robbin and
Salamon~\cite{RoSa1992}, as well as Harker, Mischaikow, and
Spendlove~\cite{HMS2021}, their approach has to be extended to allow
for varying underlying posets. For this, we first introduce the notion
of reduced filtered chain complexes, before discussing Conley complexes
and connection matrices, as well as their existence for arbitrary poset
filtered chain complexes. We close the section with a new equivalence
relation for Conley complexes, which enables us to precisely formulate
the uniqueness question of connection matrices for the first time. Despite
being a completely algebraic criterion based on the notion of essentially
graded morphisms, it will later allow us to detect underlying bifurcations
in combinatorial dynamics.

\subsection{Reduced filtered chain complexes}

Recall that by Corollary~\ref{cor:M-convex-sing}, every poset filtered
chain complex~$(P,C,d)$ gives rise to the induced chain complexes~$(C_p,d_{pp})$
for every $p\in P$. The following definition lies at the heart of the notion
of Conley complexes.
\begin{defn}
\label{def:reduced}
{\em
We say that a poset filtered chain complex $(P,C,d)$ is {\em \reduced\/},
if it is peeled, i.e., if we have $P_\star = P$, and if the chain
complex~$(C_p,d_{pp})$ is \boundaryless\ for all $p\in P$.
}
\end{defn}
\begin{ex}
\label{ex:reduced}
{\em 
   The filtered chain complex in Example~\ref{ex:pfcc} is not reduced,
   because one can check that $d_p$ is not boundaryless for $p=\xf$. Consider
   therefore the different set of words 
 \[
    X':=\setof{\mathbf{A},\mathbf{B},\mathbf{AB},\mathbf{AD},\mathbf{AE}},
 \] 
 a subset $\Psub:=\{\xa,\xb,\xc,\xd,\xe\}$ of $P$
 and a poset $\PPsub=(\Psub,\leqsub)$  with its partial order $\leqsub$ defined by the Hasse diagram
\begin{equation}
\label{ex:hasse1s}
   \begin{diagram}
  \dgARROWLENGTH 1.6em
    \node{\xd}
    \arrow{se,-}
    \node[2]{\xe}
    \arrow{sw,-}\\
    \node[2]{\xc}
    \arrow{se,-}
    \arrow{sw,-}\\
    \node{\xa}
    \node[2]{\xb}
  \end{diagram}.
\end{equation}
Notice that the partial order~$\leqsub$ is just the restriction to~$\Psub$ of
the partial order~$\leq$ in the poset~$P$ defined in Example~\ref{ex:pfcc}.
Proceeding as in Example~\ref{ex:pfcc} we obtain a $\PPsub$-filtered
$\ZZ_2$-module~$(\PPsub,C',d')$ with $C'_p:=\ZZ_2\scalprod{X'_p}$,
where $X'_p$ for $p=\xa,\xb,\xc,\xd,\xe$ in $\Psub$  are defined respectively as
\[
\{\mathbf{A}\},\{\mathbf{B}\},\{\mathbf{AB}\},\{\mathbf{AD}\},\{\mathbf{AE}\}
\]
and the homomorphism $d':C'\to C'$ is defined on the basis $X'$ by the matrix
\begin{equation}
\label{eq:d1m}
\begin{array}{c||c|c|c|c|c|}
d'  &  \mathbf{A} &  \mathbf{B} & \mathbf{AB} & \mathbf{AD} & \mathbf{AE} \\
  \hline
  \hline
  \mathbf{A} &    &    &  1 &  1 &  1 \\
  \hline                                                                                    
  \mathbf{B} &    &    &  1 &    &    \\
  \hline
 \mathbf{AB} &    &    &    &    &    \\
  \hline
 \mathbf{AD} &    &    &    &    &    \\
  \hline
 \mathbf{AE} &    &    &    &    &    \\
  \hline
\end{array}\;.
\end{equation}
One can easily check that the partial order $\leqsub$ is $d'$-admissible.
However, it is not the native partial order of $d'$. As the reader may
verify, the native partial order of~$d'$ is the partial order~$\leq'$
given by the Hasse diagram
\begin{equation}
\label{ex:hasse1}
   \begin{diagram}
  \dgARROWLENGTH 1.6em
    \node{\xd}
    \arrow{se,-}
    \node{\xe}
    \arrow{s,-}
    \node{\xc}
    \arrow{sw,-}
    \arrow{se,-}\\
    \node[2]{\xa}
    \node[2]{\xb}
  \end{diagram},
\end{equation}
which gives a poset~$\PP'=(\Psub,\leq')$. Replacing the poset~$\PPsub$
in~$(\PPsub,C',d')$ by~$\PP'$, we obtain a $\PP'$-filtered $\ZZ_2$-module~$(\PP',C',d')$
which differs from~$(\PPsub,C',d')$ only in the $d'$-admissible partial order.
Nevertheless, it is another object in~$\PfCC$. 
Note that both~$(\PPsub,C',d')$ and~$(\PP',C',d')$ are peeled, because it follows from
Proposition~\ref{prop:singl-dubl-homology} that $C'_p$ is homotopically essential for $p\in\Psub$. 
Therefore, we get from Proposition~\ref{prop:filt-homot-essential}
that~$(\PPsub,C',d')$ and~$(\PP',C',d')$ are neither isomorphic in~$\PfCC$
nor in~$\ChPfCC$, because otherwise~$\PPsub$ and~$\PP'$ would be isomorphic
as posets. We note that both~$(\PPsub,C',d')$ and~$(\PP',C',d')$ are also reduced, 
because $d'_{pp}=0$ for every index $p\in \Psub$.
\exend}
\end{ex}

\begin{prop}
\label{prop:diag-equality}
  Assume that $(P,C,d)$ and $(P',C',d')$ are two reduced poset filtered chain complexes
  and that $(\alpha,h),(\beta,g):(P,C,d)\to (P',C',d')$ are filtered chain homotopic morphisms.
  If $\alpha$ is injective, then for $p\in P$ and $q\in P'$ we have
 \[ 
     \alpha(q)=p\implies h_{qp}=g_{qp}.
\]
\end{prop}
\proof
  Since $(P',C',d')$ is reduced, it is peeled. Thus, we get 
  from Proposition~\ref{prop:filt-homot-essential} that $\alpha=\beta: P'\to P$
  and there exists an $\alpha$-filtered degree $+1$ homomorphism $\Gamma:C\to C'$ such that
  $g-h=d'\Gamma+\Gamma d$. Since $(P,C,d)$ and $(P',C',d')$ are reduced, 
  we get from Proposition~\ref{prop:FMod-homo} the identities
  \begin{eqnarray*}
    g_{qp}-h_{qp} & = & (d'\Gamma)_{qp}+(\Gamma d)_{qp}=
    d'_{q\id(q)}\Gamma_{\id(q)p}+\Gamma_{q\alpha(q)}d_{\alpha(q)p} \\
    & = & 0\Gamma_{qp}+\Gamma_{qp}0=0.
  \end{eqnarray*}
  For this, recall also that~$d$ and~$d'$ are $\id_P$-filtered, and that~$\id_P$ is injective.
  This completes the proof.
\qed

\medskip
If two poset filtered chain complexes $(P,C,d)$ and $(P',C',d')$ are filtered chain
homotopic, then the equivalence class of every filtered chain equivalence
$(\alpha,\phi) : (P,C,d)\to(P',C',d')$ is automatically an isomorphism
in~$\ChPfCC$. However, in the category~$\PfCC$ one has to prove that fact.
This is addressed in the following result.

\begin{thm}
\label{thm:Conley-uniqueness}
Assume that the two poset filtered chain complexes $(P,C,d)$ and $(P',C',d')$ are \reduced.
If they are filtered chain homotopic,
then every filtered chain equivalence $(\alpha,\phi):(P,C,d)\to(P',C',d')$
is an isomorphism in $\PfCC$. 
In particular,  $(P,C,d)$ and $(P',C',d')$ are isomorphic in $\PfCC$.
\end{thm}
\proof
   Assume that both $(P,C,d)$ and $(P',C',d')$ are \reduced.
Then we have $P_\star=P$ and $P'_\star=P'$.
Let   $(\alpha',\varphi'):(P',C',d')\to (P,C,d)$
be a filtered chain morphism such that one has both $(\alpha',\varphi')\circ(\alpha,\varphi)\sim\id_{(P,C,d)}$ 
and $(\alpha,\varphi)\circ(\alpha',\varphi')\sim\id_{(P',C',d')}$.
Clearly, $\id_P$ and $\id_{P'}$ are injective. Hence, it
follows from Proposition~\ref{prop:diag-equality} that
\[
(\varphi'\varphi)_{pp}=(\id_C)_{pp}=\id_{C_p} \quad\text{  and  }\quad
(\varphi\varphi')_{qq}=(\id_{C'})_{qq}=\id_{C'_q}
\]
for arbitrary $p\in P$ and $q\in P'$.
Therefore, it follows from Corollary~\ref{cor:FMod-homo} that
$(\id_{C})_{pp}=\varphi'_{p\alpha'(p)}\varphi_{\alpha'(p)p}$ for $p\in P$.
Similarly we obtain the identity $(\id_{C'})_{p'p'}=\varphi_{p'\alpha(p')}\varphi'_{\alpha(p')p'}$ for $p'\in P'$.
Since $\alpha$, $\alpha'$ are mutually inverse bijections, we may substitute
$\alpha'(p)$ for $p'$ and we get $(\id_{C'})_{\alpha'(p)\alpha'(p)}=\varphi_{\alpha'(p)p}\varphi'_{p\alpha'(p)}$
Hence, we conclude that $\varphi_{\alpha'(p)p}$ and $\varphi'_{p\alpha'(p)}$ are mutually inverse module homomorphisms.
Thus, it follows from Lemma~\ref{lem:FMod-iso} that $(\alpha,\varphi)$
is an isomorphism in $\FMod$.
Since $\varphi$ is a chain map, by Proposition~\ref{prop:inverse-chain-map}
its inverse is also a chain map.
Hence, $(\alpha,\varphi)$ is an isomorphism in $\PfCC$, which proves that $(P,C,d)$ and $(P',C',d')$
are isomorphic in $\PfCC$.
\qed

\subsection{Conley complexes and connection matrices}
We now turn our attention the the definition of the Conley complex
and the associated connection matrix. For this, we need two additional
concepts.

\begin{defn}\label{def:representation:pfcc}
{\em  
 Let $(P,C,d)$ be a given poset filtered chain complex.
 By a {\em representation} of~$(P,C,d)$ we mean a triple  $((P',C',d'),(\alpha,\phi),(\beta,\psi))$
 such that $(P',C',d')$ is an object in~$\PfCC$ and $(\alpha,\phi):(P,C,d)\to (P',C',d')$, 
$(\beta,\psi) : (P',C',d')\to (P,C,d)$ are mutually inverse 
elementary chain equivalences such that $\alpha$ and $\beta$ are strict, 
that is, $\dom\alpha = P'_\star$ and $\dom\beta=P_\star$.
To simplify the terminology, in the sequel we refer to the object $(P',C',d')$ as a representation of $(P,C,d)$,
assuming that the associated mutually inverse elementary chain equivalences $(\alpha,\phi)$ and $(\beta,\psi)$ are implicitly given.
}
\end{defn}
\begin{prop}
\label{prop:alpha-beta}
Given a representation $((P',C',d'),(\alpha,\phi),(\beta,\psi))$ of~$(P,C,d)$, the maps 
$\alpha:P'_\star\to P_\star$ and $\beta:P_\star\to P'_\star$ are mutually inverse
order preserving bijections.
\end{prop}
\proof Since $\alpha$ and $\beta$ are strict, we get from the definition of elementary chain homotopy
\begin{align*}
  \id_{P_\star}&=(\alpha\beta)_{P_\star}=\alpha_{P'_\star}\beta_{P_\star}=\alpha\beta,\\
  \id_{P'_\star}&=(\beta\alpha)_{P'_\star}=\beta_{P_\star}\alpha_{P'_\star}=\beta\alpha.
\end{align*}
\qed

\medskip
Let $((P',C',d'),(\alpha,\phi),(\beta,\psi))$ and
$((P'',C'',d''),(\alpha',\phi'),(\beta',\psi'))$ be two representations 
of  a filtered chain complex  $(P,C,d)$.
\begin{defn}
{\em
We define the {\em transfer  morphism} from the filtered chain complex~$(P',C',d')$ to 
the filtered chain complex~$(P'',C'',d'')$
as the filtered chain morphism $(\beta\alpha',\phi'\psi)$.
}
\end{defn}
Note that the pair of transfer morphisms~$(\beta\alpha',\phi'\psi)$ from~$(P',C',d')$
to~$(P'',C'',d'')$ and $(\beta'\alpha,\phi\psi')$ from~$(P'',C'',d'')$
to~$(P',C',d')$ are mutually inverse chain equivalences.

After these preparations we can now define the central concept of this paper.
For this, let~$(P,C,d)$ denote an arbitrary poset filtered chain complex.
Then the following definition builds upon Definitions~\ref{def:reduced}
and~\ref{def:representation:pfcc} and introduces both the notions of Conley
complex and of connection matrix.

\begin{defn}
\label{defn:Conley-complex}
{\em
By a {\em Conley complex} of $(P,C,d)$ we mean every reduced representation~$(\bar{P},\bar{C},\bar{d})$
of~$(P,C,d)$. The $(\bar{P},\bar{P})$-matrix of the boundary homomorphism~$\bar{d}$
is then called a {\em connection matrix} of the poset filtered chain complex~$(P,C,d)$.
}
\end{defn}

Since a Conley complex $(\bar{P},\bar{C},\bar{d})$ of a filtered chain complex $(P,C,d)$ is
reduced, we have $\bar{P}_\star=\bar{P}$ and, in view of Proposition~\ref{prop:alpha-beta},
we see that~$\bar{P}$ is order isomorphic to~$P_\star$. 
Therefore, in a Conley complex~$(\bar{P},\bar{C},\bar{d})$ of~$(P,C,d)$ we can
identify the poset~$\bar{P}$ with~$P_\star$, the map~$\alpha$ with the inclusion
map $\iota_P:P_\star \hookrightarrow P$,
and the map~$\beta$ with $\iota_P^{-1}:P\pto P_\star$.
We do so in the sequel. 
Clearly, the Conley complex $(\bar{P},\bar{C},\bar{d})$ is chain homotopic 
to $(P,C,d)$. Therefore, we have the following straightforward proposition. 
\begin{prop}
\label{prop:hom-conley-complex}
If $(P_\star,\bar{C},\bar{d})$ is a Conley complex of a filtered chain complex $(P,C,d)$,
then $H(\bar{C},\bar{d})$ is isomorphic to $H(C,d)$.
\qed
\end{prop}

As an immediate consequence of Theorem~\ref{thm:Conley-uniqueness},
we obtain the following two corollaries, which show that these notions
are well-defined.

\begin{cor}
\label{cor:Conley-uniqueness}
The Conley complex of a poset filtered chain complex is unique up to an isomorphism in~$\PfCC$.
In particular, the transfer homomorphism between two Conley complexes of a given poset filtered
chain complex is an isomorphism in~$\PfCC$.
\qed
\end{cor}
\begin{cor}
\label{cor:Conley-homotopy-equivalence}
If two poset filtered chain complexes are filtered chain homotopic,
that is, isomorphic in $\ChPfCC$, then their Conley complexes are isomorphic in~$\PfCC$.
\qed
\end{cor}

\begin{ex}
\label{ex:pfcc-conley1}
{\em 
   One can verify that the reduced filtered chain complex $(\PPsub,C',d')$ in Example~\ref{ex:reduced} 
   with poset $\PPsub=(\Psub,\leqsub)$ and $\leqsub$ given by the Hasse diagram~\eqref{ex:hasse1s}
   is a Conley complex of the filtered chain complex in Example~\ref{ex:pfcc}
   with the associated connection matrix~\eqref{eq:d1m}. 
   More precisely, consider the inclusion map $\epsilon:\Psub\ni x\mapsto x\in P$ 
   and homomorphisms $h':C\to C'$ and $g':C'\to C$ given by the matrices 
\[
\begin{array}{c||c|c|c|cccc|c|c|}
  h' &  \mathbf{A} &  \mathbf{B} &    \mathbf{AB} &  \mathbf{C} &    \mathbf{AC} &    \mathbf{BC} & \mathbf{ABC}   & \mathbf{CD}    & \mathbf{CE}    \\
  \hline
  \hline
  \mathbf{A} &  1 &    &       &  1 &       &       &       &       &       \\
  \hline                                                                                    
  \mathbf{B} &    &  1 &       &    &       &       &       &       &       \\
  \hline
 \mathbf{AB} &    &    &    1  &    &       &     1 &       &       &       \\
  \hline
 \mathbf{AD} &    &    &       &    &       &       &       &  1    &       \\
  \hline
 \mathbf{AE} &    &    &       &    &       &       &       &       &  1    \\
  \hline
\end{array}
\]
and
\[
\begin{array}{c||c|c|c|c|c|}
  g' &  \mathbf{A} &  \mathbf{B} & \mathbf{AB} & \mathbf{AD} & \mathbf{AE} \\
  \hline
  \hline
  \mathbf{A} &  1 &    &    &    &    \\
  \hline                                                                                    
  \mathbf{B} &    &  1 &    &    &    \\
  \hline
 \mathbf{AB} &    &    &  1 &    &    \\
  \hline
  \mathbf{C} &    &    &    &    &    \\
 \mathbf{AC} &    &    &    &  1 &  1 \\
 \mathbf{BC} &    &    &    &    &    \\
\mathbf{ABC} &    &    &    &    &    \\
  \hline
 \mathbf{CD} &    &    &    &  1 &    \\
  \hline
 \mathbf{CE} &    &    &    &    &  1 \\
  \hline
\end{array}.
\]
One can verify that 
\[
  (\epsilon,h'):(P,C,d)\to (\PPsub,C',d') \text{ and }(\epsilon^{-1},g'):(\PPsub,C',d')\to (P,C,d),
\]
with the partial map $\epsilon^{-1}:P\pto \Psub$ defined as the inverse relation of $\epsilon$,
are mutually inverse elementary filtered chain equivalences. This implies that
a Conley complex of the filtered chain complex $(P,C,d)$ in Example~\ref{ex:pfcc}
is given by the representation $((\PPsub,C',d'),(\epsilon,h'),(\epsilon^{-1},g'))$.
\exend}
\end{ex}

\subsection{Existence of Conley complexes}
We now turn our attention to the question of existence, i.e., does
every poset filtered chain complex have an associated Conley complex
and connection matrix? Although in our setting the poset in a poset
filtered chain complex is not fixed as in~\cite{HMS2021, RoSa1992},
the existence proof of a Conley complex for a poset filtered chain
complex in our sense can be adapted from the argument
in~\cite[Theorem 8.1, Corollary 8.2]{RoSa1992}. For the sake
of completeness, we present the details. We begin with a technical
lemma, which is a counterpart to~\cite[Theorem~8.1]{RoSa1992}.
\begin{lem}
\label{lem:splitting}
Assume that $P\neq\emptyset$ and that~$(P,C,d)$ denotes an arbitrary
poset filtered chain complex with field coefficients. Then there exist
four families
$\{W_p\}_{p\in P}$,
$\{V_p\}_{p\in P}$,
$\{B_p\}_{p\in P}$,
$\{H_p\}_{p\in P}$
of $\ZZ$-graded submodules of~$C$ such that the following statements
are satisfied:
\begin{itemize}
   \item[(i)]   The family~$\{W_p\}_{p\in P}$ is a $P$-gradation of~$C$ which is
                filtered equivalent to the gradation $\{C_p\}_{p\in P}$.
                In particular, $C=\bigoplus_{p\in P}W_p$.
   \item[(ii)]  For every $p\in P$ we have $W_p=V_p\oplus H_p\oplus B_p$.
   \item[(iii)] The inclusion $d(V_p)\subset B_p$ holds, and $d_{|V_p}:
                V_p\to B_p$ is a module isomorphism, for every $p\in P$.
   \item[(iv)]  The identity $d_{pp}(H_p)=0$ is satisfied for all $p\in P$.
   \item[(v)]   For $H:=\bigoplus_{p\in P}H_p$ we have
                $d(H)\subset H$ and $(P,H,d_{|H})$ is a poset filtered
                chain complex.
\end{itemize}
\end{lem}
\proof
Assume that $(P,C,d)$ is a poset filtered chain complex.
We proceed by induction on $n:=\card P$. First assume that $n=1$.
Let $p_*$ be the unique element of $P$. Then $C=C_{p_*}$. We set $W_{p_*}:=C$.
It follows that the $P$-gradations $(C_p)_{p\in P}$ and $(W_p)_{p\in P}$ are
identical, and therefore they are filtered equivalent, i.e., (i) is satisfied.
The existence of $V_{p_*}$, $B_{p_*}$, $H_{p_*}$
satisfying properties (ii)-(v) then follows from
Proposition~\ref{prop:hom-decomp-existence}.

Now assume that $n>1$. Let $r\in P$ be a maximal element in $P$
and consider the down set $P':=P\setminus \{r\}\in\Down(P)$.
Let $C':=\bigoplus_{p\in P'}C_p$.
Since~$d$ is a filtered homomorphism, we have  $d(C')\subset C'$.
Thus, $(P',C',d_{|C'})$ is a $P'$-filtered chain complex.
Since $\card P'=n-1$, by our induction hypothesis, there exist families
$\{W_p\}_{p\in P'}$,
$\{V_p\}_{p\in P'}$,
$\{B_p\}_{p\in P'}$,
$\{H_p\}_{p\in P'}$
of submodules such that properties (i)-(iv) hold for $(P',C',d_{|C'})$.
In order to have respective families for $(P,C,d)$
we will extend the families over $P'$ to families over $P$
by constructing  in turn the modules $V_r$, $B_r$, $H_r$ and $W_r$.

To begin with, for a family $\{M_p\}_{p\in P'}$ of submodules of $C$ which
satisfies the identity $M_p\cap M_q=\{0\}$ for all $p\neq q$,
and an $L\in\Down(P')$, we recall the notation
\[
    M_L:=\bigoplus_{p\in L}M_p.
\]
Set $C':=C_{P'}$, $W':=W_{P'}$, $V':=V_{P'}$, $B':=B_{P'}$, and $H':=H_{P'}$.
Then from~(i) and~(ii) applied to $(P',C',d_{|C'})$ we get
\begin{equation}
\label{eq:L-VHB}
  C'_L=C_L=W_L=V_L\oplus H_L\oplus B_L
  \quad\text{ for every }\quad
  L\in\Down(P').
\end{equation}
We will show that
\begin{equation}
\label{eq:VpHp}
   V'\cap (H'+d(C))=0.
\end{equation}
Indeed, if $x\in V'$ and $x=x_1+dx_2$ for an $x_1\in H'$ and an $x_2\in C$, then
one obtains $dx=dx_1$. By (v) of the induction assumption, $(P',H',d_{|H'})$ is a
poset filtered chain complex. Hence, $dx=dx_1=d_{|H'}x_1\in H'$.
Also, by~(iii) of the induction assumption we have $dx\in B'$, and therefore
$dx\in H'\cap B'=0$. Since $x\in V'$ one then obtains $x=0$.
This proves \eqref{eq:VpHp}.

Since $d^{-1}(C_{r^<})\cap C_{r^\leq}$ is a $\ZZ$-graded submodule of $C_{r^\leq}$, we can find a $\ZZ$-graded
submodule $V_r$ of $C_{r^\leq}$ such that
\begin{equation}
\label{eq:Vr}
   C_{r^\leq}=(d^{-1}(C_{r^<})\cap C_{r^\leq})\oplus V_r,
\end{equation}
where we also use the fact that the modules have field coefficients.
We will prove that
\begin{equation}
\label{eq:cC-geq}
   d^{-1}(C_{r^<}) \cap C_{r^\leq} =
   C_{r^<} + (d^{-1}(H_{r^<}) \cap C_{r^\leq}).
\end{equation}
To see that the right-hand side of \eqref{eq:cC-geq} is contained in the left-hand side
observe that the right-hand side is obviously contained in $C_{r^\leq}$.
Since $C$ is a filtered chain complex and $r^<$ is a down set, we have $d(C_{r^<})\subset C_{r^<}$ and, in consequence,
$C_{r^<}\subset d^{-1}(C_{r^<})$.
We also have $d^{-1}(H_{r^<})\subset d^{-1}(C_{r^<})$, because  $H_{r^<}\subset W_{r^<}$ by (ii)
and $W_{r^<}=C_{r^<}$ by (i).
To prove the opposite inclusion take an $x\in C_{r^\leq}\cap d^{-1}(C_{r^<})$.
Then $dx\in C_{r^<}$. Hence, by \eqref{eq:L-VHB}, we can find $x_V\in V_{r^<}$,
$x_H\in H_{r^<}$ and $x_B\in B_{r^<}$ such that
$dx=x_V+x_H+x_B$. From (iii) we get $x_B=dy_V$ for some $y_V\in V_{r^<}$.
It follows that $x_V=d(x-y_V)-x_H$.
Hence, $x_V\in V_{r^<}\cap (d(C)+H_{r^<})\subset V'\cap (d(C)+H')$.
Thus, from \eqref{eq:VpHp} we get $x_V=0$ and $d(x-y_V)=x_H\in H_{r^<}$ which means
$x-y_V\in d^{-1}(H_{r^<})$.
We also have $x\in C_{r^\leq}$ and
$y_V\in V_{r^<}\subset C_{r^<}\subset C_{r^\leq}$.
Therefore, one finally obtains $x=y_V+(x-y_V)\in C_{r^<}+(C_{r^\leq}\cap d^{-1}(H_{r^<}))$.
This completes the proof of \eqref{eq:cC-geq}.

Now set $B_r:=d(V_r)$.
We will prove that
\begin{equation}
\label{eq:VrCrd-1Hr}
B_r\cap  C_{r^<}=0.
\end{equation}
Let $x\in B_r\cap  C_{r^<}$.
Then $x=dy$ for some $y\in V_r\subset C_{r^{\leq}}$.
Since $x\in C_{r^<}$, we get $y\in d^{-1}(C_{r^<})$.
Therefore, $y\in V_r\cap C_{r^{\leq}}\cap d^{-1}(C_{r^<})$.
It follows from \eqref{eq:Vr} that $y=0$.
Hence, $x=d0=0$, which proves \eqref{eq:VrCrd-1Hr}.

Since $(P,C,d)$ is a poset filtered chain complex,
the boundary homomorphism $d$ is a filtered homomorphism.
Therefore, we have $C_{r^<}\subset d^{-1}(C_{r^<})$.
Obviously, $C_{r^<}\subset C_{r^{\leq}}$
and $B_r=d(V_r)\subset d(C_{r^\leq})\cap d^{-1}(0)\subset  C_{r^{\leq}}\cap d^{-1}(H_{r^<})$.
Thus, by \eqref{eq:VrCrd-1Hr},
we have a direct sum of $\ZZ$-graded submodules
\begin{equation}
\label{eq:CrBr}
(C_{r^<}\cap d^{-1}(H_{r^<}))\oplus B_r\subset C_{r^{\leq}}\cap d^{-1}(H_{r^<}).
\end{equation}
Hence, we can choose a $\ZZ$-graded submodule $H_r$ such that
\begin{equation}
\label{eq:CrBrHr}
   C_{r^{\leq}}\cap d^{-1}(H_{r^<})=(C_{r^<}\cap d^{-1}(H_{r^<}))\oplus B_r\oplus H_r.
\end{equation}
Thus, it follows from \eqref{eq:cC-geq} and Proposition~\ref{prop:XAB} that
\begin{equation}
\label{eq:CrBrHr2}
   C_{r^{\leq}}\cap d^{-1}(C_{r^<})=C_{r^<}\oplus B_r\oplus H_r.
\end{equation}
Therefore, setting $W_r:= V_r\oplus B_r\oplus H_r$ we get from \eqref{eq:Vr} that
\begin{equation}
\label{eq:CrBrHr3}
    C_{r^\leq}= C_{r^<}\oplus W_r.
\end{equation}
We now have well-defined families
$\{W_p\}_{p\in P}$,
$\{V_p\}_{p\in P}$,
$\{B_p\}_{p\in P}$,
$\{H_p\}_{p\in P}$ of submodules of~$C$.
We will prove that they indeed satisfy properties (i)-(v) for~$(P,C,d)$.
To prove (i) observe that by \eqref{eq:CrBrHr3}
\[
C=C_{P'\cup r^\leq}=C'+C_{r^\leq}=C'+C_{r^<}+W_r=C'+W_r.
\]
We claim that $C=C'\oplus W_r$. Indeed,
by \eqref{eq:CrBrHr3} we have $W_r\subset C_{r^\leq}$.
Therefore, $W_r\cap C'\subset W_r\cap C_{r^\leq} \cap C'=W_r\cap  C_{r^<}=0$.
This together with the induction assumption shows that
\[
   C=\bigoplus_{p\in P}W_p.
\]
To show that $(W_p)_{p\in P}$ is filtered equivalent to $(C_p)_{p\in P}$,
one needs to verify that $C_L=W_L$ for all $L\in\Down(P)$. Note that
by the induction assumption
\begin{equation}
\label{eq:CLWL}
   C_L=W_L \quad\text{ for }\quad L\in\Down(P').
\end{equation}
Thus, we only need to consider the case when $r\in L$.
Let $L':=L\setminus\{r\}$.
Then $L=L'\cup r^\leq$, and~\eqref{eq:CLWL} and \eqref{eq:CrBrHr3} yield
\[
C_L=C_{L'}+C_{r^\leq}=W_{L'}+C_{r^<}+W_r=W_{L'}+W_{r^<}+W_r=W_{L'}+W_{r^\leq}=W_L.
\]
This proves property (i).
By the induction assumption properties (ii)-(v) need to be verified only
for $p=r$.
Property (ii) for $p=r$ follows from the definition of $W_r$.
To see (iii) for $p=r$ take an $x\in V_r$ such that $dx=0\in C_{r^<}$.
Since the inclusion $V_r\subset C_{r^\leq}$ holds, it follows that $x\in V_r\cap C_{r^\leq}\cap d^{-1}(C_{r^<})$
and we get from \eqref{eq:Vr} that $x=0$. Thus, $d_{|V_r}$ is a monomorphism.
By the definition of $B_r$ it is an epimorphism, which proves (iii).
Finally, by \eqref{eq:CrBrHr} we have $H_r\subset d^{-1}(H_{r^<})$ which implies
\begin{equation}
\label{eq:dHr}
   d(H_{r})\subset H_{r^<}.
\end{equation}
Therefore,  $d_{rr}(H_r)=0$ which proves (iv) for $p=r$.
To show that $(P,H,d_{|H})$ is a poset filtered chain complex,
we will prove that
\begin{equation}
\label{eq:CrBrHr4}
\text{$d(H_L)\subset H_L$ for every $L\in\Down(P)$.}
\end{equation}
Property \eqref{eq:CrBrHr4} holds by our induction assumption if
one has $r\not\in L$. Thus, assume $r\in L$ and set $L':=L\setminus\{r\}$.
Then we have $H_L=H_{L'}+H_{r^<}+H_r$, and by \eqref{eq:dHr} and the
induction assumption one further obtains
\[
   d(H_L)\subset d(H_{L'})+d(H_{r^<})+d(H_r)\subset H_{L'}+H_{r^<}\subset H_L,
\]
which proves \eqref{eq:CrBrHr4}.
Since $P\in\Down(P)$ and $H_P=H$, we get from \eqref{eq:dHr} that
the restriction $d_{|H}:H\to H$ is well-defined.
Since $d^2=0$ and $d(H)\subset H$, we get $d_{|H}^2=0$.
This proves that $(H,d_{|H})$ is a chain complex.
From \eqref{eq:CrBrHr4} and Corollary~\ref{cor:filt-hom} one finally can
conclude that~$d_{|H}$ is a filtered homomorphism.
This proves~(v) and completes the proof of the lemma.
\qed

\medskip
After these preparations we can now prove the main result of this
section, which guarantees the existence of a Conley complex and 
associated connection matrix for every poset filtered chain complex.

\begin{thm}
\label{thm:Conley-complex-existence}
Every poset filtered chain complex admits a Conley complex and a
connection matrix.
\end{thm}
\proof
Assume that~$(P,C,d)$ is an arbitrary poset filtered chain complex.
Let
$\{W_p\}_{p\in P}$,
$\{V_p\}_{p\in P}$,
$\{B_p\}_{p\in P}$,
$\{H_p\}_{p\in P}$
be four families of submodules of~$C$ satisfying properties (i)-(v)
of Lemma~\ref{lem:splitting}. Then, by Lemma~\ref{lem:splitting}(i),
the collection~$\{W_p\}_{p\in P}$ is a $P$-gradation of the module
$W:=\bigoplus_{p\in P}W_p$ which coincides with the module~$C$ and, by
Proposition~\ref{prop:W-gradation}, the triple~$(P,W,d)$ is a poset filtered
chain complex which is filtered chain isomorphic to $(P,C,d)$.
Thus, it suffices to prove that $(P,W,d)$ admits a Conley complex
and a connection matrix, since the composition of a filtered chain
isomorphism and an elementary filtered chain equivalence gives again
an elementary filtered chain equivalence.

For this, set $V:=\bigoplus_{p\in P}V_p$, $B:=\bigoplus_{p\in P}B_p$, 
and $H:=\bigoplus_{p\in P}H_p$. Then one has $W=V\oplus H\oplus B$.
Moreover, $(P,V)$, $(P,B)$, and~$(P,H)$ are all objects of~$\GMod$,
and by Lemma~\ref{lem:splitting}(iv) the triple~$(P,H,d_{|H})$ is a
poset filtered chain complex. Clearly, it is a filtered chain subcomplex
of~$(P,W,d)$.

Now take $Q:=Q_\star:=P_\star$.
Let $\alpha:Q\to P$ denote the inclusion map.
Clearly, $\alpha$ is order preserving.
Moreover, $\dom\alpha=Q=Q_\star$, therefore $\alpha$ is strict.
Since $\alpha$ is injective, the inverse relation
$\beta:=\alpha^{-1}:P\pto Q$ is a well-defined partial map
which is also order preserving and strict, because $\dom\beta=\im\alpha=P_\star$.
Hence, 
$\alpha:(Q,Q_\star)\to (P,P_\star)$ and $\beta:(P,P_\star)\to (Q,Q_\star)$
are well-defined strict morphisms in $\DPSet$.

Recall that by Proposition~\ref{prop:M-convex}(i)
for every $p\in P$ we have
a chain complex~$(W_p,d_{pp})$. By Lemma~\ref{lem:splitting} it has the homology
decomposition $W_p=V_p\oplus H_p\oplus B_p$.
Thus, by Proposition~\ref{prop:hom-decomp-existence}(ii), the chain complexes $(W_p,d_{pp})$
and $(H_p,0)$ are chain homotopic. It follows that~$(W_p,d_{pp})$
is essential if and only if~$(H_p,0)$ is essential. In addition,
$(H_p,0)$ is essential if and only if $H_p\neq 0$ in view of
Corollary~\ref{cor:simple-homotopic-isomorphic}. This proves that 
\begin{equation}
\label{eq:Q-star}
\setof{p\in P\mid H_p\neq 0}\subset P_\star.
\end{equation}

We claim that $(Q,H,d_{|H})$ is also a poset filtered chain complex. 
Since $Q=P_\star$, by \eqref{eq:Q-star} we have $H=\bigoplus_{p\in Q}H_p$. 
Let $L\in\Down(Q)$ and let $L':=\setof{p\in P\mid\exists_{q\in L}\; p\leq q}$.
Then $L'\in\Down(P)$ and $H_p=0$ for $p\in L'\setminus L$. Therefore, $H_{L'}=H_L$.
Thus, since~$d_{|H}$ is $P$-filtered, it follows from Corollary~\ref{cor:filt-hom}
that~$d_{|H}$ is also $Q$-filtered, and~$(Q,H,d_{|H})$ is indeed a poset filtered
chain complex. Moreover, we get from $Q_\star=Q=P_\star$ and \eqref{eq:Q-star} 
that $(Q,H,d_{|H})$ satisfies \eqref{eq:P-star-condition}.
From  Lemma~\ref{lem:splitting}(iv)
we get that~$(Q,H,d_{|H})$ is boundaryless. Since $Q=Q_\star$ by definition, 
we conclude that~$(Q,H,d_{|H})$ is reduced.

Next, we will prove that~$(Q,H,d_{|H})$ is filtered chain homotopic
to the poset filtered chain complex~$(P,W,d)$. Let $\iota:H\to W$ and
$\pi:W\to H$ denote the inclusion and projection homomorphisms, respectively.
It follows from Lemma~\ref{lem:splitting}(v) that for $x\in H$ we have
$d\iota x=d x =\iota d x$, that is, $\iota$ is a chain map.
Also, for $x\in W$ we have $x=x_V+x_H+x_B \in V \oplus H \oplus B$.
Hence, one obtains $dx=dx_V+dx_H+dx_B=dx_V+dx_H\in B\oplus H$ in view
of Lemma~\ref{lem:splitting}(iii). This in turn implies
$\pi d x =d x_H= d\pi x$ and proves that~$\pi$ is a chain map.

The inclusion $\iota:H\to W$ is $\beta$-graded and thus also
$\beta$-filtered. Similarly, the projection $\pi:W\to H$ is $\alpha$-filtered.
Therefore, we have two well-defined filtered morphisms
$(\beta,\iota):(Q,H,d_{|H})\to (P,W,d)$
and~$(\alpha,\pi):(P,W,d)\to (Q,H,d_{|H})$.

Obviously, one has $(\alpha,\pi)\circ(\beta,\iota)=(\beta\alpha,\pi\iota)
=(\id_Q,\id_H)=\id_{(Q,H)}$ which means that $(\alpha,\pi)\circ(\beta,\iota)$
is, trivially, elementarily filtered chain homotopic to $\id_{(Q,H)}$.
We will now show that $(\beta,\iota)\circ(\alpha,\pi)
=(\alpha\beta,\iota\pi)$ is elementarily filtered chain homotopic to~$\id_{(P,W)}$.
Let $\mu:W\ni x=x_V+x_H+x_B\mapsto x_B\in B$ be the projection map and let
$\nu: V\ni x\mapsto x\in W$ be the inclusion map.
Clearly, $\mu$ and $\nu$ are graded and, in consequence, filtered homomorphisms.
By Lemma~\ref{lem:splitting}(iii) we have a well-defined $P$-graded degree~$-1$
isomorphism $d_{|V}:V\ni x\mapsto dx\in B$ with a $P$-graded inverse, which is a
degree~$+1$ isomorphism $d_{|V}^{-1}:B\to V$. Then the definition
$\Gamma:=\nu\circ d_{|V}^{-1}\circ \mu: W\to W$ gives a degree~$1$ filtered
module homomorphism. We claim that
\begin{equation}
\label{eq:Conley-complex-existence}
   \id_C-\iota\pi=\Gamma d + d\Gamma.
\end{equation}
To see \eqref{eq:Conley-complex-existence}, take an arbitrary $x\in W=C$.
Then $x=x_V+x_H+x_B$, where we have $x_V\in V$, $x_H\in H$, and $x_B\in B$.
Hence, $(\id_C-\iota\pi)(x)=x_V+x_B$, as well as
$dx=d x_V+ d x_H\in B+H$, $\Gamma d x=x_V$,
$\Gamma x=d_{|V}^{-1}(x_B)$, and also $d \Gamma x=x_B$.
It follows that $(\Gamma d + d\Gamma)(x)=x_V+x_B$ which proves the
identity~\eqref{eq:Conley-complex-existence}. Clearly, the identity
$(\beta\alpha)_{|P_\star}=\id_{P_\star}=\id_{P|P_\star}$ holds.
Hence, $(\id_P,\Gamma)$ is an elementary filtered chain homotopy between
$(\alpha\beta,\iota\pi)$ and $\id_{(P,W)}$. Thus, the poset filtered
chain complexes~$(Q,H,d_{|H})$ and~$(P,W,d)$ are filtered chain homotopic.
Since we have already seen that~$(Q,H,d_{|H})$ is \reduced, it follows
that the poset filtered chain complex~$(Q,H,d_{|H})$ is in fact
a Conley complex for~$(P,W,d)$ --- and that the $(Q,Q)$-matrix of~$d_{|H}$
is a connection matrix of the poset filtered chain complex~$(P,W,d)$.
\qed

\subsection{Conley complexes of subcomplexes}

As we stated in Proposition~\ref{prop:induced-chain-complex}, 
given a filtered chain complex $(P,C,d)$, a convex subset $J\subset P$
induces a filtered chain complex $(J,C_J,d_{J J})$. Hence, a natural question arises
whether one can obtain a Conley complex of $(J,C_J,d_{J J})$ from a Conley complex of $(P,C,d)$.
A positive answer to this question is given by the following theorem.

\begin{thmdef}
\label{thm:Conley-subcomplex}
Assume $(P_\star,\bar{C},\bar{d})$ is a Conley complex of a given poset filtered chain complex $(P,C,d)$
and $(\id_{P_\star},\phi):(P,C,d)\to (P_\star,\bar{C},\bar{d})$, $(\id_{P_\star},\psi):(P_\star,\bar{C},\bar{d})\to (P,C,d)$ are
the associated mutually inverse elementary chain equivalences. 
Let $J\subset P_\star$ be a convex subset and let $\hat{J}:=\conv_P(J)$ denote the convex hull of $J$ in $P$. 
Consider restrictions
\[
  \phir:=\phi_{J\hat{J}}, \quad \psir:=\psi_{\hat{J}J}.
\]
Then $(J,\bar{C}_{J},\bar{d}_{JJ})$ is a Conley complex of $(\hat{J},C_{\hat{J}},d_{\hat{J} \hat{J}})$
with $(\id_J,\phir)$ and  $(\id_J,\psir)$
the associated mutually inverse elementary chain equivalences.
We will call it the {\em restriction} of the Conley complex $(P_\star,\bar{C},\bar{d})$ to $J$.
\end{thmdef}
\proof
   We will first show that $(\id_J,\psir)$ is a morphism in $\PfCC$.
To prove that $\psir$ is $\id_J$-filtered, assume that $(\psir)_{pq}\neq 0$ for some $p\in\hat{J}$ and $q\in J$.
Since $(\psir)_{pq}=\psi_{pq}$ and $\psi$ is $\id_{P_\star}$-filtered, we see that $p\leq p_1$ for some $p_1\in P_\star$
such that $p_1\leq q$. 
Hence, we get $p_1\in\hat{J}$ by \eqref{eq:conv-hull}, because $p\leq p_1\leq q$ and $p,q\in J$.

   To complete the proof that $(\id_J,\psir)$ is a morphism in $\PfCC$, we still need to verify that 
$\psir:(\bar{C}_{J},\bar{d}_{JJ})\to (C_{\hat{J}},d_{\hat{J}\hat{J}})$ is a chain map. 
Observe that since $\psi:(\bar{C},\bar{d})\to (C,d)$, as a chain map,  satisfies $\psi \bar{d}=d\psi$, 
in order to prove that $\psir$ is a chain map
it suffices to verify the following two equalities
\begin{eqnarray}
\label{eq:psi-d}
(\psi \bar{d})_{\hat{J}J}&=&\psi_{\hat{J}J}\bar{d}_{JJ}\\
\label{eq:d-psi}
(d\psi)_{\hat{J}J}&=&d_{\hat{J}\hat{J}}\psi_{\hat{J}J}.
\end{eqnarray}
Since for $p\in P$, $q\in P_\star$ we have
\[
(\psi \bar{d})_{pq}=\sum_{r\in P_\star}\psi_{pr}\bar{d}_{rq},
\]
in order to prove \eqref{eq:psi-d} it suffices to verify that 
\begin{equation}
\label{eq:psi-d-1}
\psi_{pr}\bar{d}_{rq}\neq 0,\; p\in\hat{J}, \; q\in J \implies r\in J. 
\end{equation}
Hence, assume that $\psi_{pr}\bar{d}_{rq}\neq 0$ for some $ p\in\hat{J}$ and $q\in J$.
Then $\psi_{pr}\neq 0$, $\bar{d}_{rq}\neq 0$ and, using the fact that $\psi$ is $\id_{P_\star}$-filtered, 
we can find a $p_1\in P_\star$ such that $p\leq p_1$ and $p_1\leq r$.
Since $p\in\hat{J}=\conv_p(J)$, by \eqref{eq:conv-hull} we can find a $p_-\in J\subset P_\star$ 
such that $p_-\leq p\leq p_1$ and, in consequence, $p_-\leq p_1\leq r$.
Since  $\bar{d}_{rq}\neq 0$ and $\bar{d}$ is filtered, we also have $r\leq q$. 
But, $q\in J$ and also $p_-\in J$.
Now, again by \eqref{eq:conv-hull}, we get from $p_-\leq r\leq q$ that $r\in J$. 
This proves \eqref{eq:psi-d-1} and \eqref{eq:psi-d}.

Similarly, to prove \eqref{eq:d-psi} it suffices to verify that 
\begin{equation}
\label{eq:d-psi-1}
d_{pr}\psi_{rq}\neq 0,\; p\in\hat{J}, \; q\in J \implies r\in \hat{J}. 
\end{equation}
Hence, assume that $d_{pr}\psi_{rq}\neq 0$ for some $ p\in \hat{J}$ and $q\in J$.
Then $\psi_{rq}\neq 0$, $d_{pr}\neq 0$ and, using the fact that $\psi$ is $\id_{P_\star}$-filtered, 
we can find an $r_1\in P_\star$ such that $r\leq r_1$ and $r_1\leq q$.
From $d_{pr}\neq 0$ and the fact that $d$ is filtered, we get $p\leq r$. 
Since $p\in\hat{J}=\conv_P(J)$, we can find a $p_-\in J$ such that $p_-\leq p\leq r$.
Hence, we get $p_-\leq r\leq r_1\leq q$. Since $p_-,q\in J$, we get $r\in\hat{J}$, 
which verifies \eqref{eq:d-psi-1} and \eqref{eq:d-psi}.
This completes the proof that $\psi_{\hat{J}J}$ is a chain map and a morphism in $\PfCC$. 

By an analogous argument, we will prove that 
\[
\phir:(C_{\hat{J}},d_{\hat{J}\hat{J}})\to (\bar{C}_{J},\bar{d}_{JJ})
\]
is a morphism in $\PfCC$. 
First, assume that $(\phir)_{pq}\neq 0$ for some $p\in J$ and $q\in\hat{J}$.
Then $p\leq q$,  which  proves that $\phir$ is $\id_J$-filtered.
To see that $\phir$ is a chain map it suffices to verify the following two equalities
\begin{eqnarray}
\label{eq:phi-d}
(\phi d)_{J\hat{J}}&=&\phi_{J\hat{J}}\bar{d}_{\hat{J}\hat{J}}\\
\label{eq:d-phi}
(\bar{d}\phi)_{J\hat{J}}&=&\bar{d}_{JJ}\phi_{J\hat{J}}, 
\end{eqnarray}
which reduces to verifying that 
\begin{eqnarray}
\label{eq:phi-d-1}
&&\phi_{pr}d_{rq}\neq 0,\; p\in J, \; q\in \hat{J} \implies r\in \hat{J} \\
\label{eq:d-phi-1}
&&\bar{d}_{pr}\phi_{rq}\neq 0,\; p\in J, \; q\in \hat{J} \implies r\in J. 
\end{eqnarray}
To see \eqref{eq:phi-d-1} assume that  $\phi_{pr}d_{rq}\neq 0$ for some $p\in J$ and $q\in \hat{J}$.
Then $r\leq q$ and also $p\leq r$, because $p\in J\subset P_\star$. 
Hence, we get from \eqref{eq:conv-hull} that $r\in\hat{J}$.
This proves \eqref{eq:phi-d-1} and \eqref{eq:phi-d}.
To see \eqref{eq:d-phi-1} assume that  $\bar{d}_{pr}\phi_{rq}\neq 0$ for some $p\in J$ and $q\in \hat{J}$.
Then $p\leq r$ and also $r\leq q$, because $r\in P_\star$. 
Since $q\in\hat{J}$, by \eqref{eq:conv-hull} we can choose a $q^+\in J$ such that $q\leq q^+$.
Then $r\leq q^+\in J$.
We also have $p\leq r$, because $\bar{d}$ is filtered.
Hence, since $J$ is convex, we get $r\in J$, which proves \eqref{eq:d-phi-1} and \eqref{eq:d-phi}.
This completes the proof that also $\phi_{J\hat{J}}$ is a chain map and a morphism in $\PfCC$.

We still need to verify that $(\id_J,\phir)$ and  $(\id_J,\psir)$ are mutually inverse elementary chain equivalences. 
Since $(\id_{P_\star},\phi)$ and  $(\id_{P_\star},\psi)$ are mutually inverse elementary chain equivalences, we can find
filtered, degree $+1$ homomorphisms $\Gamma: C\to C$ and $\Gamma': \bar{C}\to \bar{C}$ such that
\begin{eqnarray*}
  \psi\phi&=&\id_C+\Gamma d+ d\Gamma,\\
  \phi\psi&=&\id_{\bar{C}}+\Gamma' \bar{d}+ \bar{d}\Gamma'.
\end{eqnarray*}
Then also
\begin{eqnarray}
  (\psi\phi)_{\hat{J}\hat{J}}&=&\id_{C_{\hat{J}}}+(\Gamma d)_{\hat{J}\hat{J}}+ (d\Gamma)_{\hat{J}\hat{J}},\label{eq:psi-phi-G}\\
  (\phi\psi)_{JJ}&=&\id_{\bar{C}_{J}}+(\Gamma' \bar{d})_{JJ}+ (\bar{d}\Gamma')_{JJ}.\label{eq:phi-psi-G}
\end{eqnarray}
Analogously to (\ref{eq:psi-d}-\ref{eq:d-psi}) and (\ref{eq:phi-d}-\ref{eq:d-phi}) one can verify that
\begin{eqnarray*}
  (\psi\phi)_{\hat{J}\hat{J}}&=&\psi_{\hat{J}J}\phi_{J\hat{J}},\\ 
  (\Gamma d)_{\hat{J}\hat{J}}&=&\Gamma_{\hat{J}\hat{J}}d_{\hat{J}\hat{J}},\\ 
  (d\Gamma)_{\hat{J}\hat{J}}&=&d_{\hat{J}\hat{J}}\Gamma_{\hat{J}\hat{J}},\\ 
  (\phi\psi)_{JJ}&=&\phi_{J\hat{J}}\psi_{\hat{J}J},\\ 
  (\Gamma' \bar{d})_{JJ}&=&\Gamma'_{JJ}\bar{d}_{JJ},\\ 
  (d\Gamma)_{JJ}&=&\bar{d}_{JJ}\Gamma'_{JJ}. 
\end{eqnarray*}
Applying these formulas to (\ref{eq:psi-phi-G}-\ref{eq:phi-psi-G}) we obtain
\begin{eqnarray*}
  \psir\phir&=&\id_{C_{\hat{J}}}+\Gamma_{\hat{J}\hat{J}} d_{\hat{J}\hat{J}}+ d_{\hat{J}\hat{J}}\Gamma_{\hat{J}\hat{J}},\\
  \phir\psir&=&\id_{\bar{C}_{J}}+\Gamma'_{JJ} \bar{d}_{JJ}+ \bar{d}_{JJ}\Gamma'_{JJ}.
\end{eqnarray*}
which completes the proof.

As an immediate consequence of Theorem~\ref{thm:Conley-subcomplex} and Proposition~\ref{prop:hom-conley-complex}
we get the following corollary. 
\begin{cor}
\label{cor:hom-conley-subcomplex}
Under the assumptions of Theorem~\ref{thm:Conley-subcomplex} for every convex $J\subset P_\star$ the homology 
$H(C_{\hat{J}},d_{\hat{J}\hat{J}})$ is isomorphic to $H(\bar{C}_J,\bar{d}_{JJ})$.
\qed
\end{cor}

\subsection{Equivalence of Conley complexes}

The last three subsections demonstrated that every poset filtered
chain complex~$(P,C,d)$ does have an associated Conley complex which
is unique up to an isomorphism in~$\PfCC$. Moreover, while different
representations of the Conley complex might lead to different connection
matrices, any two of them are similar to each other via an isomorphism
in the category~$\PfCC$. However, in the application of the classical
connection matrix theory to dynamical systems it was already pointed
out by Franzosa~\cite{Fr1989} and Reineck~\cite{reineck:90a, reineck:95a}
that the connection matrix of a given Morse decomposition may be not unique
in dynamic terms, which reflects the underlying bifurcations in the dynamical system.

This dynamics-related lack of uniqueness of connection matrices can also be described in
the algebraic setting of the present section, if one considers a
stronger equivalence relation between Conley complexes. For this,
we first have to discuss the notion of essentially graded morphisms,
which can be defined as follows.

\begin{defn}
\label{defn:ess-graded}
{\em
Consider two poset filtered chain complexes~$(P,C,d)$ and~$(P',C',d')$,
and let $(\alpha,h):(P,C,d)\to (P',C',d')$ be a filtered chain morphism.
Then~$(\alpha,h)$ is called {\em essentially graded}, if there exists a
graded chain morphism $(\beta,g):(P,C,d)\to (P',C',d')$ which is filtered
chain homotopic to $(\alpha,h)$.
}
\end{defn}

\begin{propdef}
\label{propdef:graded-representation}
Assume  $(\alpha,h):(P,C,d)\to (P',C',d')$ is an essentially graded chain equivalence and 
that both poset filtered chain complexes $(P,C,d)$ and $(P',C',d')$ are reduced. Then
there exists a unique graded chain morphism 
$(\beta,g):(P,C,d)\to (P',C',d')$ in the homotopy equivalence class $[(\alpha,h)]$.
We call it the {\em graded representation} of $(\alpha,h)$.
\end{propdef}
\proof
  The existence of $(\beta,g)$ follows directly from the definition of an essentially graded morphism. 
  To prove uniqueness, we assume that $(\beta,g)$ and $(\beta',g')$ are two graded chain morphisms
  in the homotopy equivalence class $[(\alpha,h)]$. Since $(\alpha,h)$ is a chain equivalence,
  we can select a filtered chain morphism $(\alpha',h'):(P',C',d')\to (P,C,d)$
  such that $(\alpha\alpha',h'h)\sim\id_{(P,C,d)}$ and $(\alpha'\alpha,hh')\sim\id_{(P',C',d')}$.
  Hence, it follows from Proposition~\ref{prop:filt-homot-essential} that $\alpha\alpha'=\id_P$
  and $\alpha'\alpha=\id_{P'}$. In particular, $\alpha$ is injective. 
  Since $(\beta,g)\sim(\alpha,h)$ and $(\beta',g')\sim(\alpha,h)$, we see that $(\beta,g)\sim(\beta',g')$
  and from Proposition~\ref{prop:filt-homot-essential} we get $\beta=\alpha=\beta'$.
  Hence, $\beta$ is injective. Consider $p\in P$ and $q\in P'$. 
  If we have $p=\alpha(q)$, then one obtains $g_{qp}=g'_{qp}$ from Proposition~\ref{prop:diag-equality}.
  If $p\neq\alpha(q)$, we get $g_{qp}=0=g'_{qp}$ from \eqref{eq:grad-hom-full-dom},
  because both $g$ and $g'$ are $\alpha$-graded.
\qed

\medskip
The following proposition is an immediate consequence of
Proposition~\ref{prop:filt-homotopy-composition}  and the fact that the composition
of graded chain morphisms is again a graded chain morphism. The later statement in
turn follows from the fact that the composition of graded morphisms is a graded
morphism, and that the composition of chain maps is again a chain map. Thus, we
have the following result.

\begin{prop}
\label{prop:composition-essentially-graded}
  Assume that both
  $(\alpha,h):(P,C,d)\to (P',C',d')$ and $(\alpha',h'):(P',C',d')\to (P'',C'',d'')$
  are filtered chain morphisms which are essentially graded, and which have the
  respective graded representations $(\beta,g):(P,C,d)\to (P',C',d')$ and
  $(\beta',g'):(P',C',d')\to (P'',C'',d'')$.
  Then the composition $(\alpha',h')\circ(\alpha,h)$ is also essentially graded
  and its graded representation is given by~$(\beta',g')\circ(\beta,g)$. \qed
\end{prop}

It is straightforward to observe that every identity morphism in~$\PfCC$ is essentially
graded, and that it is its own graded representation. 
Thus, in view of Proposition~\ref{prop:composition-essentially-graded} 
we have a well-defined wide subcategory $\EgPfCC$ of $\PfCC$ whose morphisms are
essentially graded morphisms in $\PfCC$.

\begin{prop}
\label{prop:inverse-essentially-graded}
  Assume that both $(P,C,d)$ and $(P',C',d')$ are reduced poset filtered chain complexes, 
  and that both $(\alpha,h):(P,C,d)\to (P',C',d')$ and $(\alpha',h'):(P',C',d')\to (P,C,d)$
  are mutually inverse chain equivalences. Suppose further that $(\alpha,h)$ is essentially
  graded and has the graded representation~$(\bar{\alpha},\bar{h})$. Then $(\alpha',h')$ is
  essentially graded as well,  $(\bar{\alpha},\bar{h})$ is an isomorphism
  in $\PgCC$, and $(\bar{\alpha},\bar{h})^{-1}$ is the graded representation of $(\alpha',h')$.
\end{prop}
\proof
By Definition~\ref{defn:ess-graded} we have  $(\alpha,h)\sim(\bar{\alpha},\bar{h})$.
Hence, we get from Proposition~\ref{prop:filt-homotopy-composition} the equivalences
\begin{gather}
    (\bar{\alpha},\bar{h})\circ (\alpha',h')\sim (\alpha,h)\circ (\alpha',h')\sim\id_{(P',C',d')},\label{eq:ieg-1}\\
    (\alpha',h')\circ(\bar{\alpha},\bar{h}) \sim  (\alpha',h')\circ(\alpha,h)\sim\id_{(P,C,d)},\nonumber
\end{gather}
which imply that $(\bar{\alpha},\bar{h})$ is a filtered chain equivalence. 
Now Theorem~\ref{thm:Conley-uniqueness} shows that $(\bar{\alpha},\bar{h})$ is an isomorphism
in $\PfCC$, and it follows that $(\bar{\alpha},\bar{h})$  is an isomorphism in $\FMod$.
Thus, we get from Proposition~\ref{prop:filt-homot-essential} that $\bar{\alpha}=\alpha$
and from Lemma~\ref{lem:FMod-iso} that $\alpha$ is an order isomorphism.
In consequence, the inverse of $(\bar{\alpha},\bar{h})=(\alpha,\bar{h})$ in $\PfCC$
takes the form $(\alpha^{-1},\bar{h}^{-1})$ and, clearly, the inverse $\bar{h}^{-1}$ is $\alpha^{-1}$-graded.
Moreover, from \eqref{eq:ieg-1}  we get
$$
   (\alpha',h')=(\alpha^{-1},\bar{h}^{-1})\circ(\alpha,\bar{h})\circ (\alpha',h')\sim (\alpha^{-1},\bar{h}^{-1})
$$
which proves that $(\alpha',h')$ is essentially graded with $(\alpha,\bar{h})^{-1}=(\bar{\alpha},\bar{h})^{-1}$ 
as its graded representation.
\qed

\medskip
After this algebraic intermission we return to our study of Conley complexes.
As mentioned earlier, dynamical considerations require a stronger notion of
equivalence of Conley complexes than the one provided by isomorphisms
in~$\PfCC$. This equivalence is based on essentially graded morphisms,
and can be defined as follows.

\begin{defn} \label{def-equivalent-conley-complexes}
{\em 
  We say that two Conley complexes of a poset filtered chain complex
  are {\em equivalent} if the associated transfer homomorphism is
  essentially graded.
}
\end{defn}

Note that the transfer homomorphism is always an isomorphism in~$\PfCC$. For equivalence,
one needs in addition that it is essentially graded, i.e., that it is filtered chain homotopic
to a graded chain isomorphism.

It follows from Proposition~\ref{prop:composition-essentially-graded} that the equivalence of Conley 
complexes of a given poset filtered chain complex is indeed an equivalence relation. 
Two equivalent Conley complexes, as well as the associated connection matrices,
are basically the same. This justifies the following definition. 

\begin{defn}
{\em 
  A poset filtered chain complex has a {\em uniquely determined Conley complex
  and connection matrix} if any two of its Conley complexes are equivalent. 
}
\end{defn}

Note that two non-equivalent Conley complexes of a filtered chain complex may be graded-conjugate 
and the associated connection matrices may be graded similar. Thus, graded similarity of connection 
matrices is not sufficient for their uniqueness. It is indeed necessary to consider the specific
transfer homomorphisms. This is illustrated in the following
example.

\begin{ex}
\label{ex:pfcc-conley2}
{\em 
   The filtered chain complex in Example~\ref{ex:pfcc} does not have a uniquely determined Conley complex. 
    To see this consider the set of words 
 \[
    X'':=\setof{\mathbf{A},\mathbf{B},\mathbf{AB},\mathbf{BD},\mathbf{BE}}.
 \] 
  Proceeding as in Example~\ref{ex:pfcc} we obtain a reduced, $\PPsub$-filtered $\ZZ_2$-module~$(\PPsub,C',d')$ with $C''_p:=\ZZ_2\scalprod{X''_p}$
 where $X''_p$ for $p=\xa,\xb,\xc,\xd,\xe\in \Psub$  are defined respectively as
 \[
 \{\mathbf{A}\},\{\mathbf{B}\},\{\mathbf{AB}\},\{\mathbf{BD}\},\{\mathbf{BE}\}
 \]
 and the homomorphism $d'':C''\to C''$ is defined on the basis $X''$ by the matrix
\begin{equation}
\label{eq:d2m}
\begin{array}{c||c|c|c|c|c|}
d'' &  \mathbf{A} &  \mathbf{B} & \mathbf{AB} & \mathbf{BD} & \mathbf{BE} \\
  \hline
  \hline
  \mathbf{A} &    &    &  1 &    &    \\
  \hline                                                                                    
  \mathbf{B} &    &    &  1 &  1 &  1 \\
  \hline
 \mathbf{AB} &    &    &    &    &    \\
  \hline
 \mathbf{BD} &    &    &    &    &    \\
  \hline
 \mathbf{BE} &    &    &    &    &    \\
  \hline
\end{array}.
\end{equation}
As in Example~\ref{ex:reduced} one can readily check that the partial order~$\leqsub$ in~$\PPsub$
is $d''$-admissible, but not the native partial order of~$d''$. The native partial order of~$d''$
is~$\leq''$ defined by the Hasse diagram
 \begin{equation}
 \label{ex:hasse2}
    \begin{diagram}
   \dgARROWLENGTH 1.13em
     \node[2]{\xc}
     \arrow{se,-}
     \arrow{sw,-}
     \node{\xd}
     \arrow{s,-}
     \node{\xe}
     \arrow{sw,-}\\
     \node{\xa}
     \node[2]{\xb}
   \end{diagram},
 \end{equation}
and gives the poset~$\PP''=(\Psub,\leq'')$. Then $(\PP'',C'',d'')$  
is a well-defined $\PP''$-filtered, reduced chain complex which differs from $(\PPsub,C'',d'')$ only in the partial order.
However, our interest in $(\PPsub,C'',d'')$ comes from the fact that it provides another Conley complex of the filtered chain complex in Example~\ref{ex:pfcc}.
For this, consider again the inclusion map $\epsilon:\Psub\ni x\mapsto x\in P$ 
and homomorphisms $h'':C\to C''$ and $g'':C''\to C$ given by the matrices 
\[
\begin{array}{c||c|c|c|cccc|c|c|}
h'' &  \mathbf{A} &  \mathbf{B} &    \mathbf{AB} &  \mathbf{C} &    \mathbf{AC} &    \mathbf{BC} & \mathbf{ABC}   & \mathbf{CD}    & \mathbf{CE}    \\
  \hline
  \hline
  \mathbf{A} &  1 &    &       &    &       &       &       &       &       \\
  \hline                                                                                    
  \mathbf{B} &    &  1 &       &  1 &       &       &       &       &       \\
  \hline
 \mathbf{AB} &    &    &    1  &    &     1 &       &       &       &       \\
  \hline
 \mathbf{BD} &    &    &       &    &       &       &       &  1    &       \\
  \hline
 \mathbf{BE} &    &    &       &    &       &       &       &       &  1    \\
  \hline
\end{array}
\]
and
\[
\begin{array}{c||c|c|c|c|c|}
g'' &  \mathbf{A} &  \mathbf{B} & \mathbf{AB} & \mathbf{BD} & \mathbf{BE} \\
  \hline
  \hline
  \mathbf{A} &  1 &    &    &    &    \\
  \hline                                                                                    
  \mathbf{B} &    &  1 &    &    &    \\
  \hline
 \mathbf{AB} &    &    &  1 &    &    \\
  \hline
  \mathbf{C} &    &    &    &    &    \\
 \mathbf{AC} &    &    &    &    &    \\
 \mathbf{BC} &    &    &    &  1 &  1 \\
\mathbf{ABC} &    &    &    &    &    \\
  \hline
 \mathbf{CD} &    &    &    &  1 &    \\
  \hline
 \mathbf{CE} &    &    &    &    &  1 \\
  \hline
\end{array}.
\]
Then one can verify that $((\PPsub,C'',d''),(\epsilon,h''),(\epsilon^{-1},g''))$
is another reduced representation, that is, a Conley complex
of the filtered chain complex~$(P,C,d)$ presented in Example~\ref{ex:pfcc}.
We claim that this Conley complex is not equivalent to the Conley complex $(\PPsub,C',d')$ 
presented in Example~\ref{ex:pfcc-conley1}.
To see this, assume to the contrary that the transfer homomorphism
from~$( \PPsub,C',d')$ to~$( \PPsub,C'',d'')$
is essentially graded. 
Observe that this transfer homomorphism 
is $(\epsilon^{-1}\epsilon,h''g')=(\id_{\Psub},h''g')$
and $h''g'$ has the matrix
\[
\begin{array}{c||c|c|c|c|c|}
h''g'  &  \mathbf{A} &  \mathbf{B} & \mathbf{AB} & \mathbf{AD} & \mathbf{AE} \\
  \hline
  \hline
  \mathbf{A} &  1 &    &    &    &    \\
  \hline                                                                                    
  \mathbf{B} &    &  1 &    &    &    \\
  \hline
 \mathbf{AB} &    &    &  1 &  1 &  1 \\
  \hline
 \mathbf{BD} &    &    &    &  1 &    \\
  \hline
 \mathbf{BE} &    &    &    &    &  1 \\
  \hline
\end{array}.
\]
Now let $(\gamma,f)$ be a graded representation of $(\id_{\Psub},h''g')$.
Then we clearly have $(\gamma,f)\sim (\id_{\Psub},h''g')$ and from
Proposition~\ref{prop:filt-homot-essential} one obtains
$\gamma=\id_{\Psub}$, as well as 
\begin{equation}
\label{ex:pfcc-conley2-1}
h''g'-f=d''\Gamma+\Gamma d'
\end{equation}
for some $\id_{\Psub}$-filtered degree $+1$
homomorphism $\Gamma: C'\to C''$.  
Evaluating both sides of \eqref{ex:pfcc-conley2-1} at $\mathbf{AD}$ we get
\begin{equation}
\label{ex:pfcc-conley2-2}
   \mathbf{AB}+\mathbf{BD}-f(\mathbf{AD})=
   (d''\Gamma)(\mathbf{AD})+\Gamma(\mathbf{A})=\Gamma(\mathbf{A}),
\end{equation}
because $C''_2=0$.  
Since $f$ is $\id_{\Psub}$-graded, one obtains
$\scalprod{f(\mathbf{AD}),\mathbf{AB}}=0$. Hence, from \eqref {ex:pfcc-conley2-2}
we can deduce the identity
$$
   1=\scalprod{\mathbf{AB}+\mathbf{BD}-f(\mathbf{AD}),\mathbf{AB}}=
   \scalprod{\Gamma(\mathbf{A}),\mathbf{AB}}.
$$
It follows that $\Gamma_{\xc,\xa}\neq 0$ and, since $\Gamma$ is $\id_{\Psub}$-filtered, we immediately obtain
the inequality $\xc\leq \xa$. However, we see from the Hasse
diagram~\eqref{ex:hasse1s}
of the poset $ \PPsub$ that $\xa<\xc$. This contradiction proves that the two Conley complexes
$( \PPsub,C',d')$  and $( \PPsub,C'',d'')$  of $(P,C,d)$ are not equivalent.
The associated connection matrices are the matrices of~$d'$ and~$d''$ which are shown
in~\eqref{eq:d1m} and~\eqref{eq:d2m}, respectively. 
Note that despite the fact that these matrices are not equivalent
as connection matrices of $(P,C,d)$, they are graded similar. To see this
define $\chi: \PPsub\to \PPsub$ by
\[
  \chi(x):= \begin{cases}
     \xb & \text{ if $x=\xa$,}\\
     \xa & \text{ if $x=\xb$,}\\
     x & \text{otherwise.}
   \end{cases}
\]
and $h:C'\to C''$ by the matrix
\[
\begin{array}{c||c|c|c|c|c|}
h &  \mathbf{A} &  \mathbf{B} & \mathbf{AB} & \mathbf{AD} & \mathbf{AE} \\
  \hline
  \hline
  \mathbf{A} &   & 1 &    &    &    \\
  \hline                                                                                    
  \mathbf{B} &  1 &    &    &    &    \\
  \hline
 \mathbf{AB} &    &    &   1 &    &    \\
  \hline
 \mathbf{BD} &    &    &    &  1 &    \\
  \hline
 \mathbf{BE} &    &    &    &    &  1 \\
  \hline
\end{array}.
\]
Then one easily verifies that $\chi$ is an order preserving bijection, 
 that~$(\chi,h)$ is a graded chain isomorphism, and that we have
\[
  (\id_{ \Psub},d'')\circ(\chi,h)=(\chi,h)\circ(\id_{ \Psub},d').
\]
This readily establishes their graded similarity.
\exend}
\end{ex}


\section{Connection matrices in Lefschetz complexes}
\label{sec:lefconnmatrix}

With this section we turn our attention to a more specialized situation.
Rather than continuing to study connection matrices for general poset
filtered chain complexes, we now consider the setting of Lefschetz
complexes. Their definition and basic properties have already been
recalled in Section~\ref{sec:lefschetz}, and therefore this section
focuses on the introduction of acyclic partitions in Lefschetz complexes,
as well as their associated connection matrices. We then concentrate
on the fundamental properties of a special case, the so-called
singleton partition, and end the section with a discussion of
refinements.

\subsection{Connection matrices of acyclic partitions}
\label{sec:c-matr-acyc-part}

Let~$(X,\kappa)$ be a Lefschetz complex as in Definition~\ref{def:lefschetz},
equipped with the Lefschetz topology induced by the face relation partial
order~$\leq_\kappa$. Furthermore, let~$\cE$ denote a partition of~$X$ into
locally closed sets in this Lefschetz topology. Finally, we consider the
relation~$\preceq_\cE$ in~$\cE$ defined by
\begin{equation}
\label{eq:preceq}
   E \preceq_\cE E' \quad\text{ if and only if }\quad E\cap \cl E'\neq\emptyset.
\end{equation}
We say that $\cE$ is an {\em acyclic partition} of $X$
if $~\preceq_\cE$ may be extended to  a partial order on $\cE$.
Every such extension is called $\cE$-admissible. Note that the
smallest $\cE$-admissible partial order on an acyclic partition $\cE$
is the transitive closure of $\preceq_\cE$. We call it the {\em inherent}
partial order of $\cE$ and denote it by the symbol~$\leq_\cE$.

We recall (see Section~\ref{sec:sets-and-maps}) that given a subfamily $\cD\subset\cE$ we 
use the compact notation $|\cD|:=\bigcup\cD\subset X$ to denote the subset of~$X$ associated with~$\cD$.
The following proposition characterizes some topological features of~$|\cD|$ in terms of poset features of~$\cD$.
\begin{prop}
\label{prop:convex-lcl}
   Assume that~$\cE$ is an acyclic partition of a Lefschetz complex~$X$,
   that~$\leq$ is an $\cE$-admissible partial order, and that $\cD\subset\cE$.
   Then the following statements hold:
\begin{itemize}
   \item[(i)] If $\cD$ is a down set with respect to $\leq$, then $|\cD|$ is closed.
   \item[(ii)] If $\cD$  is convex with respect to $\leq$, then $|\cD|$ is locally closed.
\end{itemize}
\end{prop}
\proof
   In order to prove (i), we assume that $\cD\subset\cE$ is a down set and
   let $x\in\cl |\cD|$. Then there exists an $E'\in\cD$ such that $x\in \cl E'$.
   Since $\cE$ is a partition of $X$, there exists
   an $E\in\cE$ such that $x\in E$. It follows that $E\cap \cl E'\neq\emptyset$
   and, since $\leq$ is an $\cE$-admissible partial order, we get
    $E\leq E'$. Since $\cD$ is a down set, this in turn implies $E\in \cD$.
   Thus, $x\in E\subset |\cD|$, and therefore~$|\cD|$ is closed. This completes
   the proof of (i).
   To prove (ii), we now assume that $\cD$ is convex.
   It follows from Proposition~\ref{prop:interval} that both~$\cD^\leq$ and~$\cD^<$
   are down sets, and due to~(i) the sets $|\cD^\leq|$ and $|\cD^<|$ are closed.
   Moreover, since $\cD=\cD^\leq\setminus \cD^<$ and $\cE$ is a partition,
   we have $|\cD|=|\cD^\leq|\setminus |\cD^<|$. Hence, $|\cD|$ is locally closed
   by Proposition~\ref{prop:lcl}.
\qed

\medskip
With any acyclic partition of a Lefschetz complex one can associate
a connection matrix. In order to make use of the theory developed in
the previous section, we first need to recognize the Lefschetz complex
as a poset filtered chain complex. Thus, assume that~$\cE$ is an acyclic
partition of a Lefschetz complex~$X$. Then the module~$C(X)$ spanned
by the cells of~$X$ admits the $\cE$-gradation
\begin{equation}
\label{eq:cE-gradation}
    C(X)=\bigoplus_{E\in\cE} C(E).
\end{equation}
We have the following proposition.
\begin{prop}
\label{prop:Lefschetz-filtered-chain-complex}
The triple $(\cE,C(X),\bdy^\kappa)$, with an $\cE$-admissible partial
order~$\leq$ on~$\cE$, and the $\cE$-gradation~\eqref{eq:cE-gradation},
is a poset filtered chain complex.
\end{prop}
\proof
We need to prove that $\bdy^\kappa$ is a filtered module homomorphism,
that is, $\bdy^\kappa$ is an $\alpha$-filtered module homomorphism with $\alpha=\id_{\cE}$.
For this, we will verify \eqref{eq:id-filt-hom-downset} of Corollary~\ref{cor:filt-hom}
for $M=C(X)$ and $L=\cD\in\Down(\cE)$.
Then one immediately obtains $M_L=C(X)_\cD=C(|\cD|)$,
and we need to verify that the inclusion $\bdy^\kappa(C(|\cD|))\subset C(|\cD|)$ holds.
But, this follows from Proposition~\ref{prop:closed-subcomplex}, because,
by Proposition~\ref{prop:convex-lcl}(i), the set $|\cD|$ is closed.
Thus, by Corollary~\ref{cor:filt-hom}, we see that the boundary homomorphism
is $\id_{\cE}$-filtered. Therefore, $(\cE,C(X),\bdy^\kappa)$ is a poset filtered chain complex.
\qed

\medskip
Proposition~\ref{prop:Lefschetz-filtered-chain-complex} lets us define the
Conley complex and connection matrix of an acyclic partition of a Lefschetz
complex. For this we consider poset $\cE$ as an object of $\DPSet$ with the distinguished subset $\cE_\star$
satisfying \eqref{eq:P-star-condition}, that is we assume that $C(E)$ is homotopically inessential for $E\in\cE\setminus\cE_\star$.
\begin{defn}
{\em
By a {\em filtration} of an acyclic partition~$\cE$ of a Lefschetz
complex~$X$ we mean a poset filtered chain complex~$(\cE,C(X),\bdy^\kappa)$
with~$\cE$ ordered by an $\cE$-admissible partial order.
The {\em Conley complex} and associated {\em connection matrix}
of~$\cE$ is defined as the Conley complex and associated connection
matrix of the filtration of~$\cE$, that is, the Conley complex and
respective connection matrix of~$(\cE,C(X),\bdy^\kappa)$.
}
\end{defn}

\subsection{The singleton partition}
\label{subsec-singletonpart}

While a given Lefschetz complex can have many different acyclic
partitions, one of the most important one from a technical point
of view is its finest one, which is given by the partition into
singletons. For this special case, we have the following proposition
and definition.
\begin{propdef}
\label{prop:singleton-partition}
Assume that~$(X,\kappa)$ is a Lefschetz complex.
Then the following hold:
\begin{enumerate}
  \item[(i)]   The family $\cX:=\setof{\{x\}\mid x\in X}$ of 
               all singletons is a partition of~$X$ into locally
               closed sets.
  \item[(ii)]  The map $\sing: X\ni x\mapsto\{x\}\in\cX$ is a
               bijection which preserves in both directions the
               face relation~$\leq_\kappa$ on~$X$ and the
               relation~$\preceq_\cX$ on~$\cX$ given
               by~\eqref{eq:preceq}.
  \item[(iii)] The relation~$\preceq_\cX$ is a partial order
               on~$\cX$ which coincides with the inherent partial 
               order~$\leq_\cX$ on~$\cX$.
  \item[(iv)]  The family~$\cX$ is an acyclic partition of~$X$.
\end{enumerate}
We call~$\cX$ the {\em singleton partition} of~$X$.
\end{propdef}
\proof
Clearly $\cX$ is a partition of $X$ and every singleton $\{x\}$ is locally closed,
because we have $\{x\}=x^{\leq_\kappa}\setminus x^{<_\kappa}$ and the sets
$x^{\leq_\kappa}$ and $x^{<_\kappa}$ are both down sets, hence closed. This proves (i).
Now let $x,y\in X$. Then one easily verifies the sequence of equivalences
\begin{eqnarray*}
x\leq_\kappa y &\iff & x\in\cl y\\
                         &\iff & \{x\}\cap\cl\{y\}\neq\emptyset\\
                         &\iff & \{x\}\preceq_\cX\{y\},
\end{eqnarray*}
which proves (ii).
Hence, $\preceq_\cX$ is a partial order on $\cX$.
In particular, $\preceq_\cX$ is transitively closed, and this
implies that $\preceq_\cX$ coincides with the inherent partial order on $\cX$.
This proves (iii).
 It follows that $\cX$ is an acyclic partition of $X$, which finally establishes (iv).
\qed

\medskip
Proposition~\ref{prop:singleton-partition} shows that we have a well-defined
poset filtered chain complex~$(\cX,C(X),\bdy^\kappa)$. Using the order
isomorphism $\sing:X\to\cX$ we can then identify it with~$(X,C(X),\bdy^\kappa)$.
Note that via this identification closed sets in~$X$ are in one-to-one
correspondence with down sets in~$\cX$.

\begin{prop}
\label{prop:native-porder}
Let~$(X,\kappa)$ be a Lefschetz complex. Then the native partial order
of~$\bdy^\kappa$ is precisely the face relation~$\leq_\kappa$.
\end{prop}
\proof
  According to Definition~\ref{prop:d-admissible}
  we have to prove that for every admissible partial order
  $\leq$ on $X$ and $x,y\in X$ the inequality $x\leq_\kappa y$
  implies the inequality $x\leq y$. Clearly, it suffices to prove
  that $x\prec_\kappa y$ implies $x\leq y$, because $\leq_\kappa$
  is the transitive closure of $\prec_\kappa$. But,
  $x\prec_\kappa y$ by definition means that $\kappa(y,x)\neq 0$
  which implies $\bdop^\kappa_{yx} \neq 0$, and the admissibility
  of~$\leq$ finally gives the inequality~$x\leq y$.
\qed

\medskip
Note that the bijection $\sing$ lets us identify partial orders in $X$
with partial orders in the associated singleton partition $\cX$ of $X$.
As an immediate consequence of Proposition~\ref{prop:singleton-partition}
and  Proposition~\ref{prop:native-porder} we get the following corollary.
\begin{cor}
Let $(X,\kappa)$ be a  Lefschetz complex and let $\cX$ denote the associated singleton partition.
A partial order is $\bdop^\kappa$-admissible if and only if it is $\cX$-admissible. 
\end{cor}

We say that a $\bdy^\kappa$-admissible partial order $\leq$ in $X$ is {\em natural}  if
\begin{equation} \label{def-natural-order}
  x\leq y \;\;\;\text{ and }\;\; \dim x=\dim y
  \quad\implies\quad x=y.
\end{equation}
Note that the native  partial order of $\bdy^\kappa$ is natural,
but that a $\bdy^\kappa$-admissible partial order $\leq$ in $X$ does not
need to be natural in general.

\begin{defn}
\label{defn:filtration-Lefschetz}
{\em
A {\em filtration of a Lefschetz complex $(X,\kappa)$}
is defined as a filtration of the singleton partition of~$X$,
that is, the poset filtered chain complex $(X,C(X),\bdy^\kappa)$
with $X$ ordered by a $\bdy^\kappa$-admissible partial order.
 When we consider $X$ as a poset ordered by a natural partial
 order in $X$ we refer to  the filtration $(X,C(X),\bdy^\kappa)$  as a {\em natural filtration of $X$}.
 When we consider $X$ as a poset ordered by the native partial
 order of $X$ we refer to the filtration $(X,C(X),\bdy^\kappa)$  as the {\em native filtration of $X$}.
}
\end{defn}

\begin{thm}
\label{thm:singleton-partition}
Let~$(X,\kappa)$ be a Lefschetz complex. Then the following hold:
\begin{itemize}
   \item[(i)]  Every filtration of the Lefschetz complex~$(X,\kappa)$  is \reduced,
   that is, it is a Conley complex of itself. In particular, $X_\star=X$.
   \item[(ii)] Every natural filtration of a Lefschetz complex~$(X,\kappa)$ has a uniquely
   determined Conley complex and  connection matrix.
   The connection matrix coincides with the $(X,X)$-matrix of the boundary homomorphism
    $\bdy^\kappa$.
\end{itemize}
\end{thm}
\proof
Consider an $x\in X$. We have from \eqref{eq:kappa-condition-1} that
$\kappa_{|\{x\}\times\{x\}}=0$. Therefore, $\bdy^\kappa_{|C(\{x\})}=0$.
Moreover, since $C_{\dim x}(\{x\})=Rx\neq 0$, we see 
that~$X_\star=X$ which means that~$(X, C(X),\bdy^\kappa)$
is \reduced. In consequence, it is a Conley complex of itself. This proves~(i).   

In order to prove (ii), it suffices to verify that every transfer morphism from
$(X, C(X),\bdop^\kappa)$ to a Conley complex $(P,C,d)$ of $(X, C(X),\bdop^\kappa)$
is essentially graded. Actually, we will prove the stronger fact that every such transfer
morphism is graded. Thus, assume that $(X, C(X),\bdop^\kappa)$ is a natural filtration
of~$X$, that~$(P,C,d)$ is another Conley complex of $X$, and that the map $(\alpha,\phi):
(X,C(X),\bdy^\kappa)\to(P,C,d)$ is a filtered chain isomorphism.

Since in view of~(i) the filtration $(X,C(X),\bdy^\kappa)$ is a Conley complex of itself,
the transfer morphism from $(X,C(X),\bdy^\kappa)$ to $(P,C,d)$ is just $(\alpha,\phi)$. 
Since both $(X, C(X),\bdop^\kappa)$ and $(P,C,d)$, as Conley complexes, are reduced, we see from 
Proposition~\ref{prop:inverse-essentially-graded} that
$(\alpha,\phi)$ is also an isomorphism in $\FMod$. 
Thus, it follows from Lemma~\ref{lem:FMod-iso}
that $\alpha:P\to X$ is an order isomorphism and $\phi_{p\alpha(p)}:C(\{\alpha(p)\}) \to C_p$
is an isomorphism of $\ZZ$-graded moduli for every $p\in P$. 
But, Proposition~\ref{prop:singl-dubl-homology} provides the structure of $C(\{\alpha(p)\})$.
Hence, we can choose a non-zero $c_p\in C$ such that
\[
   (C_p)_n=\begin{cases}
                   Rc_p & \text{ if $n=\dim\alpha(p),$}\\
                   0 & \text{ otherwise.}
            \end{cases}
\]
In other words, $C_p$ is a one-dimensional $\ZZ$-graded module generated by $c_p$ such that
\begin{equation}
\label{eq:lefc-as-pfcc}
\dim c_p=\dim\alpha(p).
\end{equation}
In order to  prove that $(\alpha,\phi)$ is an isomorphism in $\GMod$
we will first show that $\phi$ is $\alpha$-graded by verifying \eqref{eq:grad-hom}.
Assume therefore that we have $\varphi_{py}\neq 0$ for $p\in P$ and $y\in X$.
Furthermore, let $x = \alpha(p)$. Since $(\alpha,\phi)$ is a filtered
module homomorphism, we get from \eqref{eq:filt-hom} that
$p\in\alpha^{-1}(y^\leq)^\leq$. This implies that there exists a $p' \in P$
such that $p \le p'$ and $\alpha(p') \le y$. Since $\alpha$ is an order isomorphism,
one further obtains $x = \alpha(p) \le \alpha(p') \le y$.

Note that $C_p$ is a one-dimensional $\ZZ$-graded module which is non-zero only
in dimension $n=\dim\alpha(p)=\dim x$. Since $\varphi_{py}:C(\{y\}) \to C_p$ is
$\ZZ$-graded homomorphism of degree zero, both $C(\{y\})$ and $C_p$ are one-dimensional
$\ZZ$-graded modules, and we have $\varphi_{py}\neq 0$, one immediately obtains
that~$\varphi_{py}$ maps~$Ry$ isomorphically onto~$R c_p$. It follows that
$\dim c_p=\dim y$. In view of $x = \alpha(p)$ and~\eqref{eq:lefc-as-pfcc}
we therefore have $\dim x=\dim y$. Since $x\leq y$ and the partial order is
natural, this implies $x=y$, as well as $\alpha(p)=x=y$. This in turn proves
that $\phi$ is $\alpha$-graded, that is, $(\alpha,\phi)$ is a module homomorphism
in $\GMod$. Since it is an isomorphism in $\FMod$, one obtains from
Corollary~\ref{cor:FMod-iso} that it is also an isomorphism in $\GMod$.
Since it is a chain map, Proposition~\ref{prop:inverse-chain-map} shows
that it is an isomorphism in~$\PgCC$. This completes the proof of (ii).
\qed

\subsection{Refinements of acyclic partitions}

Consider two objects of~$\DPSet$, that is, posets~$P$ and~$Q$ with
distinguished subsets $P_\star$ and $Q_\star$,
an order preserving surjection $\mu: Q\to P$ such that 
\begin{equation}
\label{eq:mu-Q-star}
\mu(Q_\star)=P_\star,
\end{equation}
and a module $M$ which is both $P$-filtered and $Q$-filtered, 
that is, we have objects $(P,M)$ and $(Q,M)$ in $\FMod$.

\begin{defn}
\label{defn:refinement}
{\em 
We say that $(Q,M)$ is a {\em $\mu$-refinement} of $(P,M)$
if for every element $p\in P$ we have
\begin{equation}
\label{eq:refinement}
   M_p=\bigoplus_{q\in\mu^{-1}(p)}M_q.
\end{equation}
Given a chain complex $(C,d)$ which is both $P$-filtered and $Q$-filtered, 
we say that the object~$(Q,C,d)$ in~$\PfCC$ is a {\em $\mu$-refinement}
of~$(P,C,d)$ in~$\PfCC$ if~$(Q,C)$ is a $\mu$-refinement of~$(P,C)$
as a module.
}
\end{defn}

\begin{prop}
\label{prop:homom-refinement}
Assume that the object~$(Q,M)$ of~$\FMod$ is a {\em $\mu$-refinement}
of the object~$(P,M)$, and that $(\id_Q,h):(Q,M)\to (Q,M)$ is a morphism
in~$\FMod$. Then also $(\id_P,h):(P,M)\to (P,M)$ is a morphism in~$\FMod$.
\end{prop}
\proof
We need to verify that
   $h_{p'p}\neq 0$ for $p,p'\in P$ implies $p'\leq p$. 
   Hence, assume that $h_{p'p}\neq 0$ holds for some $p,p'\in P$. Then there exist
   elements $q\in\mu^{-1}(p)$ and $q'\in\mu^{-1}(p')$ such that $h_{q'q}\neq 0$.
   Since $(\id_Q,h)$ is a morphism in~$\FMod$, the homomorphism~$h$ is $\id_Q$-filtered.
   Therefore, we get $q'\leq q$ and since $\mu$ is order preserving, 
   we obtain $p'=\mu(q')\leq\mu(q)=p$. 
   This proves that also $(\id_P,h)$ is a morphism in $\FMod$.
\qed

\begin{prop}
\label{prop:C-mu}
Assume that $(Q,C,d)$ is a {\em $\mu$-refinement} of $(P,C,d)$. Then, for every $A\subset P$ we have
\begin{equation}
\label{eq:C-mu}
   C_A=C_{\mu^{-1}(A)}.
\end{equation}
\end{prop}
\proof
  From  the definition \eqref{eq:M-I} of $C_A$ and \eqref{eq:refinement} we get
\[
  C_A=\bigoplus_{p\in A}C_p=\bigoplus_{p\in A}\left(\bigoplus_{q\in\mu^{-1}(p)}C_q\right)=\bigoplus_{q\in\mu^{-1}(A)}C_q = C_{\mu^{-1}(A)}.
\]
This completes the proof of the proposition.
\qed

\medskip
Consider now a poset filtered chain complex~$(P,C,d)$, as well as a Conley
complex $(Q_\star,\bar{C},\bar{d})$ of a $\mu$-refinement~$(Q,C,d)$ of~$(P,C,d)$.
In addition, define the map $\bar{\mu}:=\mu\circ \iota_Q : Q_\star \to P_\star$,
where $\iota_Q:Q_\star\hookrightarrow Q$ stands for the inclusion map. Note that,
in view of \eqref{eq:mu-Q-star}, the map~$\bar{\mu}$ is a surjection. 
For $p\in P_\star$ set 
\begin{displaymath}
  \bar{C}_p:=\bigoplus_{q\in \bar{\mu}^{-1}(p)}\bar{C}_q.
\end{displaymath}
Since $\bar{\mu}$ is a surjection, we clearly have
\begin{equation}
\label{eq:bar-C-p}
\bar{C}=\bigoplus_{p\in P_\star}\bar{C}_p,
\end{equation}
which makes $\bar{C}$ a $P_\star$-graded module.

\begin{prop}
\label{prop:P-star-chain-complex}
The triple $(P_\star,\bar{C},\bar{d})$ is a poset filtered chain complex. 
\end{prop}
\proof
   Clearly, $(Q_\star,\bar{C})$  is a $\bar{\mu}$-refinement of $(P_\star,\bar{C})$.
   Since $(Q_\star,\bar{C},\bar{d})$ is a poset filtered chain complex, the boundary map~$\bar{d}$ is $\id_{Q_\star}$-filtered.
   By Proposition~\ref{prop:homom-refinement} it is also $\id_{P_\star}$-filtered.
   Hence, $(P_\star,\bar{C},\bar{d})$ is indeed a poset filtered chain complex.
\qed

\medskip
The following proposition shows that, under suitable assumptions, the grouping \eqref{eq:bar-C-p} 
of components in the Conley complex $(Q_\star,\bar{C},\bar{d})$ of $(Q,C,d)$ provides a Conley complex of $(P,C,d)$.
\begin{prop}
\label{prop:refinement}
In the situation described above, assume further that in the poset filtered
chain complex $(P_\star,\bar{C},\bar{d})$ the homomorphisms 
$\bar{d}_{pp}$ are boundaryless for $p\in P_\star$.
Then $(P_\star,\bar{C},\bar{d})$ is a Conley complex of $(P,C,d)$.
\end{prop}
\proof
   Clearly, $(P_\star,\bar{C},\bar{d})$ is peeled, therefore, it is also reduced,
   because we assume that $\bar{d}_{pp}$ is boundaryless for $p\in P_\star$.
   Thus, we only need to verify that~$(P,C,d)$ is
   elementarily filtered equivalent to~$(P_\star,\bar{C},\bar{d})$.
   Since~$(Q_\star,\bar{C},\bar{d})$ is a Conley complex of $(Q,C,d)$,
   we have the associated mutually inverse elementary chain equivalences 
   $(\iota_Q,\varphi):(Q,C,d)\to (Q_\star,\bar{C},\bar{d})$
   and
   $(\iota_Q^{-1},\psi):(Q_\star,\bar{C},\bar{d})\to (Q,C,d)$.
   We claim that also 
   $(\iota_P,\varphi):(P,C,d)\to (P_\star,\bar{C},\bar{d})$
   and
   $(\iota_P^{-1},\psi):(P_\star,\bar{C},\bar{d})\to (P,C,d)$,
   with $\iota_P: P_\star \hookrightarrow P$ denoting the inclusion,
   are mutually inverse elementary chain equivalences. 
   We will first prove that $(\iota_P,\varphi)$ is a morphism in $\PfCC$.
   For this it suffices to prove that it is $\iota_P$-filtered, that is,
   $(\iota_P,\varphi)$ is a morphism in $\FMod$,
   because, since $(\iota_Q,\varphi)$ is a morphism in $\PfCC$, the
   homomorphism~$\varphi$ is already a chain map. 
   For this, let $L\in\Down(P)$ be a down set.
   By Proposition~\ref{prop:filt-hom} it suffices to verify that
   the inclusion
\begin{equation}
\label{eq:C-L-0}
   \varphi(C_L)\subset \bar{C}_{\iota_P^{-1}(L)^\leq}
\end{equation}
holds.
   We get from Proposition~\ref{prop:C-mu} and from Proposition~\ref{prop:filt-hom} applied to $\varphi$
   as an $\iota_Q$-filtered homomorphism that
\begin{equation}
\label{eq:C-L-1}
   \varphi(C_L)=\varphi(C_{\mu^{-1}(L)})\subset \bar{C}_{\iota_Q^{-1}(\mu^{-1}(L))^\leq}= \bar{C}_{\bar{\mu}^{-1}(L)^\leq}.
\end{equation}
   Note that $\bar{\mu}^{-1}(L)=\bar{\mu}^{-1}(L\cap P_\star)$, because $\dom \bar{\mu}=Q_\star$ and, by \eqref{eq:mu-Q-star}, 
   we have the equality $\mu(Q_\star)=P_\star$. Therefore, 
\begin{equation}
\label{eq:C-L-2}
   \bar{C}_{\bar{\mu}^{-1}(L)^\leq}=\bar{C}_{\bar{\mu}^{-1}(L\cap P_\star)^\leq}.
\end{equation}
   Observe that $\bar{\mu}:Q_\star\to P_\star$ is a surjection and $(Q_\star,\bar{C},\bar{d})$ 
   is a $\bar{\mu}$-refinement of~$(P_\star,\bar{C},\bar{d})$ by the very definition~\eqref{eq:bar-C-p} of  $(P_\star,\bar{C},\bar{d})$.
   Hence, applying Proposition~\ref{prop:C-mu} to $\bar{\mu}$, we get
\begin{equation}
\label{eq:C-L-3}
   \bar{C}_{\bar{\mu}^{-1}(L\cap P_\star)^\leq}=\bar{C}_{(L\cap P_\star)^\leq}=\bar{C}_{\iota_P^{-1}(L)^\leq}.
\end{equation}
   Combining \eqref{eq:C-L-1}, \eqref{eq:C-L-2}, and \eqref{eq:C-L-3} we finally get \eqref{eq:C-L-0}.
   This completes the proof that~$(\iota_P,\varphi)$ is a morphism in~$\PfCC$.

   Next we will prove that $(\iota_P^{-1},\psi)$ is a morphism in $\PfCC$.
   Again, it suffices to verify that $\psi$ is $\iota_P^{-1}$-filtered.
   Hence, assume now that $L\in\Down(P_\star)$ is a down set.
   We first observe that  
\begin{equation}
\label{eq:mu-L-leq}
   \mu^{-1}(L)^\leq\subset\mu^{-1}(L^\leq).
\end{equation}
   Indeed, if we assume $q\in \mu^{-1}(L)^\leq$, then $q\leq q'$ for some $q'\in\mu^{-1}(L)$, 
   which gives $\mu(q)\leq\mu(q')\in L$ and proves that $q\in \mu^{-1}(L^\leq)$.
   Hence, applying Proposition~\ref{prop:C-mu} and using the fact that $\psi$
   is $\iota_Q^{-1}$-filtered, in combination with
   the observation that $(\iota_Q^{-1})^{-1}(A)=A$ for $A\subset Q_\star$, we get
\[
   \psi(\bar{C}_L)=\psi(\bar{C}_{\bar{\mu}^{-1}(L)})\subset C_{(\iota_Q^{-1})^{-1}(\bar{\mu}^{-1}(L))^\leq}
   =C_{\bar{\mu}^{-1}(L)^\leq}=C_{\iota_Q^{-1}(\mu^{-1}(L))^\leq}.
\]
Clearly, 
\begin{equation}
\label{eq:iota-Q}
\iota_Q^{-1}(\mu^{-1}(L))=Q_\star\cap\mu^{-1}(L)\subset\mu^{-1}(L).
\end{equation}
Therefore, applying \eqref{eq:iota-Q}, \eqref{eq:mu-L-leq}, and again Proposition~\ref{prop:C-mu}, we further obtain
\[
   C_{\iota_Q^{-1}(\mu^{-1}(L))^\leq}\subset C_{\mu^{-1}(L)^\leq}
   \subset C_{\mu^{-1}(L^\leq)}=C_{L^\leq}=C_{(\iota_P^{-1})^{-1}(L^\leq)}.
\]
This proves that $\psi$ is $\iota_P^{-1}$-filtered and $(\iota_P^{-1},\psi)$ is a morphism in $\PfCC$.

Finally, consider an elementary filtered chain homotopy $T$ 
between 
\[
(\iota_Q^{-1},\psi)\circ(\iota_Q,\varphi)=(\id_{Q_\star},\psi\varphi):(Q,C,d)\to(Q,C,d)
\]
and 
\[
\id_{(Q,C,d)}=(\id_Q,\id_C):(Q,C,d)\to(Q,C,d).
\]
By Proposition~\ref{prop:homom-refinement}, the homomorphism $T$ is $\id_P$-filtered.
Therefore, it is straightforward to see that $T$ is also an elementary filtered chain homotopy
between 
\[
   (\iota_P^{-1},\psi)\circ(\iota_P,\varphi)=(\id_{P_\star},\psi\varphi):(P,C,d)\to(P,C,d)
\]
and 
\[
   \id_{(P,C,d)}=(\id_P,\id_C):(P,C,d)\to(P,C,d). 
\]
Similarly, we find an elementary chain homotopy between
\[
   (\iota_P,\varphi)\circ(\iota_P^{-1},\psi)=(\id_{P_\star},\varphi\psi):(P_\star,\bar{C},\bar{d})\to(P_\star,\bar{C},\bar{d})
\]
and 
\[
   \id_{(P_\star,\bar{C},\bar{d})}=(\id_{P_\star},\id_{\bar{C}}):(P_\star,\bar{C},\bar{d})\to(P_\star,\bar{C},\bar{d}).
\]
This proves that $(P_\star,\bar{C},\bar{d})$ is a Conley complex of $(P,C,d)$.
\qed

\medskip
Proposition~\ref{prop:refinement} facilitates the computation of Conley complexes of
acyclic partitions of Lefschetz complexes. For this,
consider two acyclic partitions~$\cF$ and~$\cE$ of a Lefschetz complex~$X$.
\begin{defn}
\label{defn:part-refinement}
{\em 
We say that $\cF$ is a {\em refinement} of $\cE$ 
if for every $F\in\cF$ there is an $E\in\cE$ such that $F\subset E$. 
}
\end{defn}
Assume that the acyclic partition~$\cF$ of~$X$ is a refinement of the acyclic partition~$\cE$ of~$X$.
Then, for every $F\in\cF$ there is exactly one $E\in\cE$ such that $F\subset E$.
Therefore, we have a well-defined map
\[
\mu=\mu_{\cF,\cE}:\cF\ni F\mapsto E\in\cE.
\]
Since $\cE$ and $\cF$ are partitions, $\mu$ is clearly a surjection.
It is straightforward to verify that $\mu$ preserves the inherent partial
orders of~$\cF$ and~$\cE$. Therefore, we have the following proposition.

\begin{prop}\label{prop:aprefinement}
Assume that the acyclic partitions~$\cF$ and~$\cE$ of the Lefschetz complex~$X$
are objects of~$\DPSet$ with the distinguished subsets~$\cE_\star$
and~$\cF_\star$ satisfying~\eqref{eq:P-star-condition}.
If~$\cF$  is a refinement of~$\cE$ and $\mu_{\cF,\cE}(\cF_\star)=\cE_\star$,
then the filtration $(\cF,C(X),\bdy^\kappa)$ of~$X$ is 
a $\mu_{\cF,\cE}$-refinement of the filtration~$(\cE,C(X),\bdy^\kappa)$.
\qed
\end{prop}


\section{Dynamics of combinatorial multivector fields}
\label{sec:dcmvf}

Multivector fields are a natural generalization of Forman's combinatorial
vector fields. They provide greater flexibility in the description of
a variety of dynamical phenomena, and they can be used to combinatorialize
even such concepts as chaotic behavior or multiflows. In the present 
section, we review basic notions from combinatorial multivector fields
on Lefschetz complexes, and show that they provide a natural framework
for our theory of connection matrices. This is achieved through a
straightforward construction of an associated acyclic partition of
the underlying Lefschetz complex.

\subsection{Combinatorial multivector fields}
\label{sec:cmvf}
The concept of a combinatorial multivector field was proposed
in~\cite[Definition 5.10]{Mr2017} as a dynamically oriented extension
of the notion of combinatorial vector field in the sense of
Forman~\cite{Fo98b,Fo98a}. The definition of combinatorial multivector
field considered here is based on~\cite{LKMW2020}, yet restricted to the
special setting of Lefschetz complexes. Thus, let~$X$ be an arbitrary
Lefschetz complex. By a {\em combinatorial multivector} in~$X$ we mean
any non-empty subset of~$X$ which is locally closed with respect to the
Lefschetz topology. Then a {\em combinatorial multivector field} is a
partition~$\cV$ of~$X$ into combinatorial multivectors. Since in this
paper we do not use any other concepts of vector fields, in the sequel
we simplify the terminology by dropping the adjective combinatorial in
combinatorial multivector and combinatorial multivector field.

We say that a multivector~$V$ is {\em critical}, if the relative Lefschetz
homology~$H(\cl V,\mouth V)$ is non-zero. A multivector~$V$ which is not
critical is called {\em regular}. For each $x\in X$ we denote by~$\vclass{x}$
the unique multivector in~$\cV$ which contains~$x$. If the multivector
field~$\cV$ is clear from context, we abbreviate the notation by writing
$[x]:=\vclass{x}$. We say that $x\in X$ is {\em critical} (respectively
{\em regular}) with respect to the multivector field~$\cV$, if the
multivector~$\vclass{x}$ is critical (respectively regular).
We say that a subset $A\subset X$ is $\cV$-{\em compatible} if
for each $x\in X$ either $\vclass{x}\cap A=\emptyset$ or
$\vclass{x}\subset A$.

We associate with every multivector field a multivalued map $\Pi_\cV: X\mto X$
given by the definition
\begin{equation} \label{def-multivalflow}
   \Pi_\cV(x):=\cl x\cup \vclass{x}.
\end{equation}
The multivalued map~$\Pi_\cV$ may be interpreted as a digraph with
vertices in~$X$ and an arrow from $x\in X$ to $y\in X$ if $y\in\Pi_\cV(x)$.
Clearly, every multivector $V\in\cV$ forms a clique in this digraph.
By collapsing all vertices in a multivector to a point we obtain an
induced digraph with vertices in $\cV$ and an arrow from $V\in\cV$ to $W\in\cW$
if there exist an $x\in V$ and a $y\in W$ such that $y\in\Pi_\cV(x)$.
We refer to this digraph as the $\cV$-digraph.

Since~$\cV$ is a partition of~$X$ into locally closed subsets,
we can consider the relation $\preceq_\cV$ in~$\cV$ introduced
at the beginning of Section~\ref{sec:c-matr-acyc-part}. The following
proposition shows that the $\cV$-digraph and the relation $\preceq_\cV$
are the same concepts.

\begin{prop}
\label{prop:cV-digraph}
There is an arrow from $V$ to $W$ in the $\cV$-digraph of $\cV$
if and only if $W\preceq_\cV V$, that is, if and only if $W\cap\cl V\neq\emptyset$.
\end{prop}
\proof
   Assume there is an arrow from $V$ to $W$ in the $\cV$-digraph of $\cV$, and
   let $x,y\in X$ be such that $x\in V$, $y\in W$, and
   $y\in\Pi_\cV(x)=[x]_\cV\cup\cl x$.
   If $y\in[x]$, then $W=[y]=[x]=V$ which yields $W\cap\cl V=V\neq\emptyset$.
   If $y\in\cl[x]$, then $y\in W\cap\cl V$ proving that $W\cap\cl V\neq\emptyset$.
   Vice versa, if $W\cap\cl V\neq\emptyset$,
   then we can take a $y\in W\cap\cl V$ and an $x\in V$ such that $y\in \cl x$.
   It follows that $y\in\Pi_\cV(x)$.
   Hence, there is an arrow in the $\cV$-digraph from $V$ to $W$.
\qed

\begin{figure}
  \includegraphics[width=0.25\textwidth]{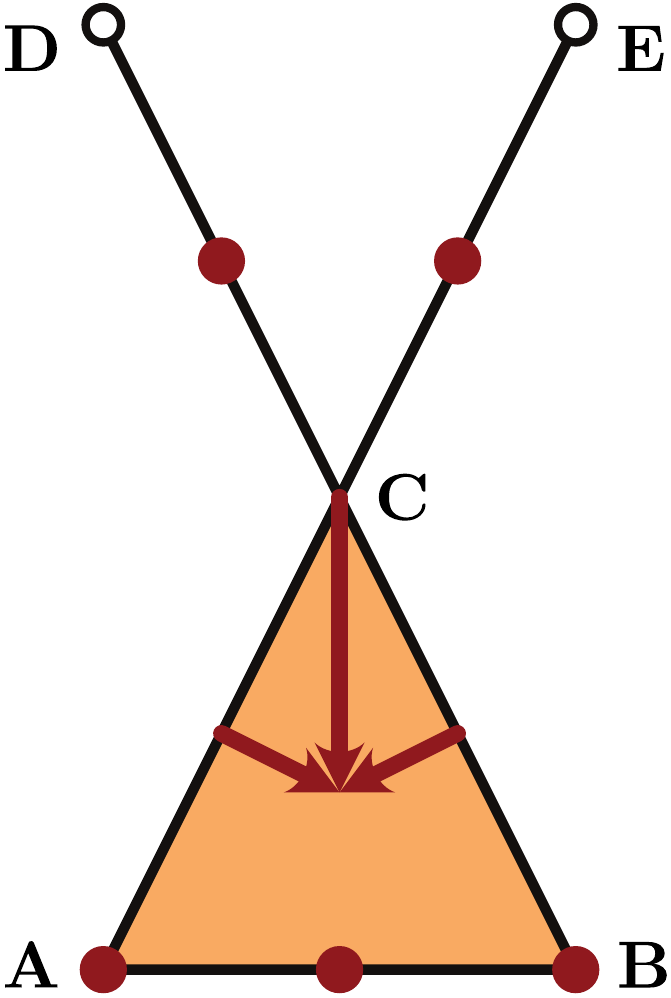}\qquad\quad
  \includegraphics[width=0.25\textwidth]{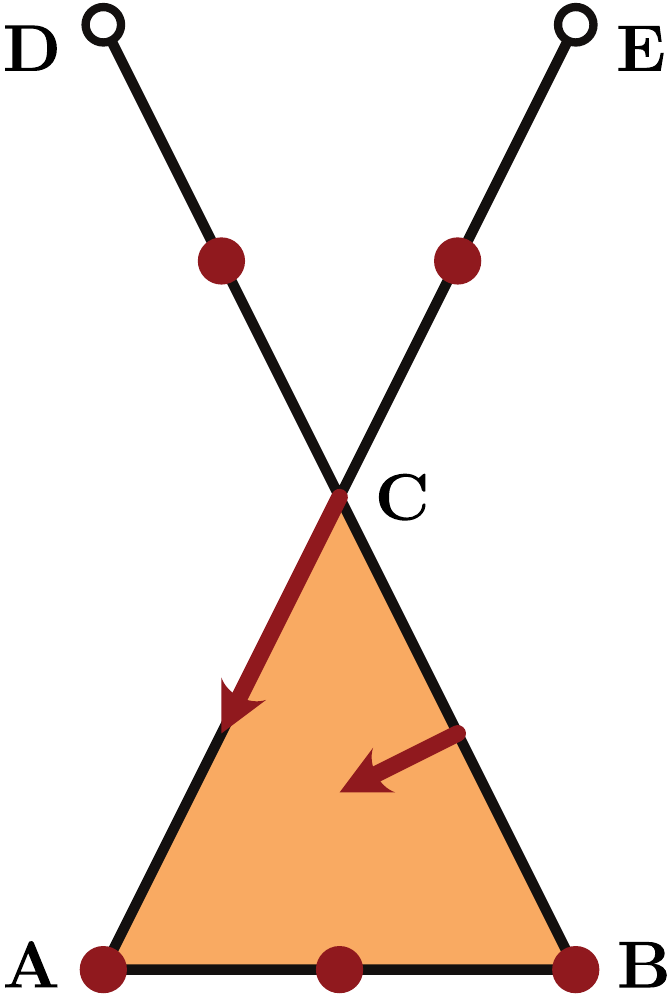}\qquad\quad
  \includegraphics[width=0.25\textwidth]{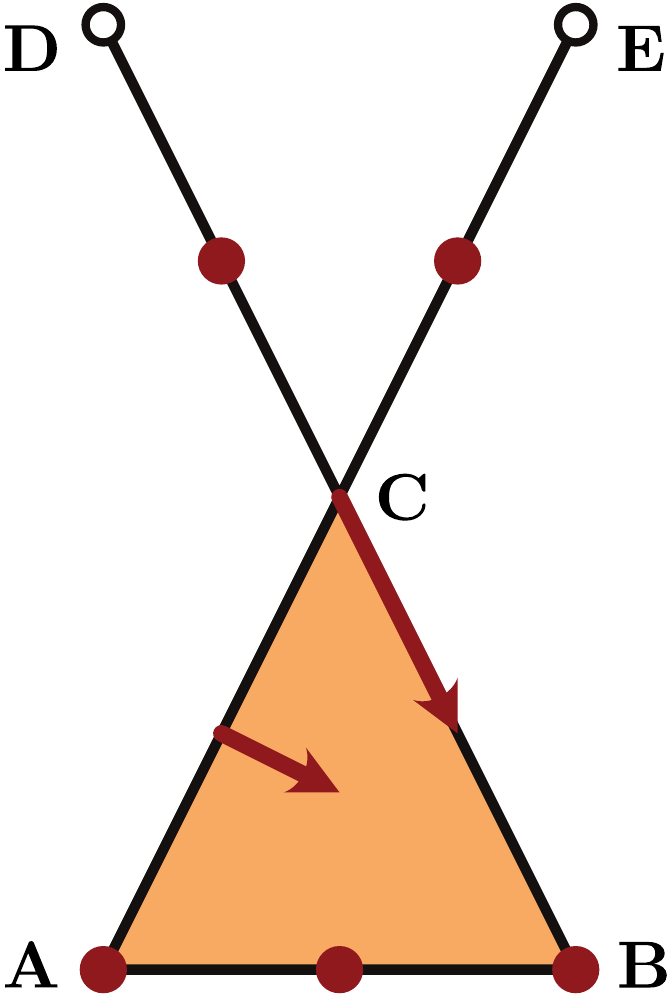}
  \caption{
  Three different multivector fields~$\cV_0$ (left), $\cV_1$ (middle),
  and~$\cV_2$ (right) on the Lefschetz complex $X$ introduced in
  Example~\ref{ex:lefschetz}. Critical cells are marked with a fat
  circle in the center of mass of a simplex. Vectors and multivectors
  are indicated by an arrow from simplex~$x$ to simplex~$y$ whenever
  $y\in \Pi_\cV(x)$ and $y\not\in\cl x$, and as long as~$y$ is a
  top-dimensional coface in the multivector containing~$x$. The
  remaining cases of arrows $x \to y$ are not marked in order to
  keep the image readable.
  }
  \label{fig:t2eLefMVF}
\end{figure}

\begin{ex}
\label{ex:lefschetzMVF}
{\em 
  Figure~\ref{fig:t2eLefMVF} presents three different combinatorial multivector
  fields~$\cV_0$, $\cV_1$, and~$\cV_2$ on the Lefschetz complex~$X$ introduced
  in Example~\ref{ex:lefschetz}. The $\cV_0$-digraph coincides with the Hasse
  diagram~\eqref{ex:hasse0}. Notice that by using the Hasse diagram representation
  of the digraph, we implicitly assume that arrows always point downwards, i.e.,
  we represent them without arrow heads. Moreover, to keep the diagrams as simple
  as possible, we do not indicate the loops which are present at every node.
  Similarly, the $\cV_1$-digraph is
  \begin{equation}
  \label{ex:hasseV1}
     \begin{diagram}
    \dgARROWLENGTH 2mm
      \node{\{\mathbf{CD}\}}
      \arrow{se,-}
      \node{\{\mathbf{CE}\}}
      \arrow{s,-} 
      \node[2]{\{\mathbf{BC},\mathbf{ABC}\}}
      \arrow{wsw,-}
      \arrow{s,-}\\
      \node[2]{\{\mathbf{C},\mathbf{AC}\}}
      \arrow{se,-}
      \node[2]{\{\mathbf{AB}\}}
      \arrow{se,-}
      \arrow{sw,-}\\
      \node[3]{\{\mathbf{A}\}}
      \node[2]{\{\mathbf{B}\}}
    \end{diagram}~~
  \end{equation}
and the  $\cV_2$-digraph is given by
\begin{equation}
\label{ex:hasseV2}
   \begin{diagram}
  \dgARROWLENGTH 2mm
    \node[2]{\{\mathbf{AC},\mathbf{ABC}\}}
    \arrow{s,-}    
    \arrow{ese,-}    
    \node[2]{\{\mathbf{CD}\}}
    \arrow{s,-}
    \node{\{\mathbf{CE}\}}
    \arrow{sw,-}\\
    \node[2]{\{\mathbf{AB}\}}
    \arrow{se,-}
    \arrow{sw,-}
    \node[2]{\{\mathbf{C},\mathbf{BC}\}}
    \arrow{sw,-}\\
    \node{\{\mathbf{A}\}}
    \node[2]{\{\mathbf{B}\}}
  \end{diagram}~~
\end{equation}
We will return to these three multivector fields later in more detail.
\exend}
\end{ex}

We call a subset $A\subset \ZZ$ {\em left bounded}
(respectively {\em right bounded})
if it has a minimum (respectively maximum).
Otherwise we call it {\em left unbounded} or {\em left infinite}
(respectively {\em right unbounded} or {\em right infinite}).
We call $A\subset \ZZ$ {\em bounded} if $A$ has both a minimum and a maximum.
Finally, we call $A\subset \ZZ$ a $\ZZ$-interval if $A=\ZZ\cap I$
where $I$ is an interval in $\RR$.

A {\em solution} of a multivector field $\cV$ in $A\subset X$ is a partial
map $\rho:\ZZ\pto A$ whose {\em domain}, denoted by $\dom \rho$,  is a
$\ZZ$-interval and for any $i,i+1\in\dom\rho$ the inclusion
$\rho(i+1)\in \Pi_\cV(\rho(i))$ holds.
The solution {\em passes} through $x\in X$ if $0\in\dom \rho$ and $x=\rho(0)$.
The solution $\rho$ is {\em full} if $\dom \rho=\ZZ$.
It is {\em periodic} if there is a $k>0$ such that $\rho(i+k)=\rho(i)$ for all $i\in\ZZ$.
It is a {\em partial solution} or simply a {\em path} if $\dom \rho$ is bounded.
We refer to the cardinality of the domain of a path as the {\em length} of the path.
If the maximum of $\dom\rho$ exists, we call the value of~$\rho$ at this maximum
the {\em right endpoint} of $\rho$.
If the minimum of $\dom\rho$ exists, we call the value of~$\rho$ at this minimum
the {\em left endpoint} of $\rho$.
We denote the left and right endpoints of $\rho$ by $\lep{\rho}$ and $\rep{\rho}$,
respectively. Given a full solution~$\rho$ through $x\in X$,
we denote by $\rho^+:=\rho_{|\scriptsize{\ZZ}_0^+}$ the {\em forward solution} through $x$,
and by $\rho^-:=\rho_{|\scriptsize{\ZZ}_0^-}$ the {\em backward solution} through $x$.

By a {\em shift} of a solution $\rho$ we mean the composition
$\rho\circ\tau_n$, where the translation map is defined as
$\tau_n:\ZZ\ni m\mapsto m+n\in \ZZ$.
Given two solutions~$\varphi$ and~$\psi$ such that~$\lep{\psi}$
and~$\rep{\varphi}$ exist and $\lep{\psi}\in\Pi_\cV(\rep{\varphi})$,
there is a unique shift~$\tau_n$ such that $\varphi\cup\psi\circ\tau_n$ is a solution.
We call this union of paths the {\em concatenation} of~$\varphi$ and~$\psi$
and denote it $\varphi\cdot\psi$. We also identify each $x\in X$
with the path of length one whose image is~$\{x\}$.

The following proposition is straightforward.
\begin{prop}
\label{prop:X-path}
Let $x,y\in X$. Then $y\in\Pi_\cV(x)$ if and only if there is an arrow
from~$[x]$ to~$[y]$ in the $\cV$-digraph. Consequently, the multivector
field~$\cV$  admits a path from~$x$ to~$y$ if and only if
there is a path from~$[x]$ to~$[y]$ in the $\cV$-digraph.
\qed
\end{prop}
Finally, and for later reference, we define the {\em backward
and forward ultimate image} of a full solution $\rho:\ZZ\to X$,
respectively, as the sets
\begin{eqnarray*}
   \uim^-(\rho)&:=&\bigcap_{t\in\ZZ^-}\rho((-\infty,t]),\\
   \uim^+(\rho)&:=&\bigcap_{t\in\ZZ^+}\rho([t,\infty)).
\end{eqnarray*}
One can easily see that both of these sets are necessarily nonempty.
Moreover, the following proposition is straightforward.
\begin{prop}
\label{prop:omega-periodic}
If $\rho:\ZZ\to X$ is a periodic solution, then its backward and forward 
ultimate images satisfy $\uim^-(\rho)=\im\rho=\uim^+(\rho)$.
\qed
\end{prop}
While the above notion of solution is a natural generalization of the
classical case, the specific definition of the multivalued map~$\Pi_\cV$
causes some slight complication. Since we clearly have $x \in \Pi_\cV(x)$
for every $x \in X$, every point in the Lefschetz complex lies on a constant
full solution through itself. Thus, if we are interested in extending
notions such as invariance from the classical situation, our concept of
solution implies that every subset of~$X$ would be an invariant set, which
clearly is not appropriate. This can be remedied by strengthening our
notion of solutions to the concept of essential solutions.

A full solution $\rho:\ZZ\to X$ is called {\em left-essential} (or,
respectively, {\em right-essential}), if for every regular $x\in \im\rho$
the set $\setof{t\in\ZZ\mid\rho(t)\not\in\vclass{x}}$ is {\em left-infinite}
(or, respectively, {\em right-infinite}). We say that~$\rho$ is an
{\em essential solution} if it is both left- and right-essential.
In other words, essential solutions have to leave every regular multivector
in forward and in backward time before they can re-enter it.

\subsection{Conley index and Morse decompositions}

We say that $S\subset X$ is {\em $\cV$-invariant}, or briefly {\em invariant},
if for every $x\in S$ there exists an essential solution through~$x$ in~$S$.
A closed set $N\subset X$ is an {\em isolating set} for a $\cV$-invariant
subset $S\subset N$ if $\Pi_\cV(S)\subset N$ and any path in~$N$ with
endpoints in~$S$ is a path in~$S$. Finally, we say that~$S$ is an
{\em isolated invariant set} if~$S$ admits an isolating set. One can
show that isolated invariant sets have an easy direct characterization. This
is the subject of the following proposition (see also~\cite[Propositions~4.10,
4.12, and~4.13]{LKMW2020}).
\begin{prop}
\label{prop:isolation-criterion}
A subset $S\subset X$ is an isolated invariant set for~$\cV$ if and
only if it is an invariant set which is both $\cV$-compatible and
locally closed.
\end{prop}

Consider an isolated invariant set $S$ of a combinatorial multivector field $\cV$
on a Lefschetz complex $X$. The Conley index of $S$ 
is defined in~\cite[Section~5.2]{LKMW2020} as the homology of any index pair of $S$, 
in particular the homology $H(\cl S,\mo S)$ of a special index pair $(\cl S,\mo S)$
(see~\cite[Theorem~5.3]{LKMW2020}).
Since, in view of Proposition~\ref{prop:rel-homology} and Proposition~\ref{prop:isolation-criterion},
$H(\cl S,\mo S)$ is isomorphic to $H(S)$, for the purposes of this paper it suffices
to assume that the Conley index of an isolated invariant set $S$ is 
\begin{equation*}
\Con(S):=H(S),
\end{equation*}
that is, the Lefschetz homology of~$S$.

Observe that given an isolated invariant set $S$ for a multivector field  $\cV$
on a Lefschetz complex $X$, by Proposition~\ref{prop:isolation-criterion} the family
\[
\cV_S:=\setof{V\in\cV\mid V\subset S}
\]
is a partition of $S$ and, in consequence, a multivector field on $S$ considered as a Lefschetz subcomplex of $X$.
We call it the multivector field {\em induced} by $\cV$ on $S$.

Multivector fields on Lefschetz complexes can exhibit a variety
of different dynamical behaviors, ranging from critical cells which
correspond to equilibrium solutions, to both transient and more 
complicated recurrent behavior. In this context, an essential, and
therefore, full solution $\rho:\ZZ\to X$ is called {\em recurrent}
if for every $x\in \im\rho$ the set $\rho^{-1}(x)$ is both left-
and right-infinite. Examples of recurrent solutions include constant
solutions whose image is a single critical multivector, or nonconstant
periodic solutions which loop through a sequence of regular multivectors.

In classical dynamics, one has the notion of a gradient vector field,
which rules out the existence of nonconstant recurrent solutions --- and in 
fact, Forman's original combinatorial vector fields were motivated by
precisely this situation. While our definition of combinatorial multivector
fields is considerably more general and allows for recurrence, we will
see that assuming an additional gradient structure can provide more
detailed information about the associated connection matrices.

Thus, we call a multivector field~$\cV$ {\em gradient-like} if for
every recurrent solution $\rho:\ZZ\to X$ there exists a multivector
$V\in\cV$ such that $\im\rho\subset V$. Moreover, a multivector field~$\cV$
is called a {\em gradient multivector field} if the only  recurrent
solutions are constant solutions with its value in a singleton of~$\cV$.
We would like to point out that this rules out multivectors of size at
least two which have nontrivial Conley index. Moreover,
recall that a singleton in~$\cV$ is always critical, see
Proposition~\ref{prop:singl-dubl-homology}.

Gradient-like multivector fields on Lefschetz complexes are of significant
importance on our theory of connection matrices, since they correspond
precisely to acyclic partitions. This is the subject of the following
result.
\begin{prop}
\label{prop:grdient-like-mvf}
Let~$\cV$ be a combinatorial multivector field on a Lefschetz complex~$X$.
The following properties are equivalent.
\begin{itemize}
  \item[(i)]   The multivector field~$\cV$ is gradient-like.
  \item[(ii)]  After discarding self-loops, the $\cV$-digraph is acyclic.
  \item[(iii)] The collection~$\cV$ of multivectors is an acyclic partition of~$X$.
\end{itemize}
\end{prop}
\proof
The equivalence (i)$\iff$(ii) follows from Proposition~\ref{prop:X-path},
and the equivalence (ii)$\iff$(iii) follows  from Proposition~\ref{prop:cV-digraph}.
\qed

\medskip
In order to define connection matrices for  multivector
fields on a Lefschetz complex~$X$ we need the concept of a
Morse decomposition of an isolated invariant set. 
Let $S$ be an isolated invariant set of a multivector field $\cV$ on a Lefschetz complex $X$.
By {\em Morse decomposition of~$S$} we mean in this paper a collection~$\cM$ of mutually
disjoint, non-empty isolated invariant subsets of~$S$, referred to as {\em Morse sets},
together with a partial order~$\leq$ on~$\cM$ such that for every essential
solution $\rho:\ZZ\to S$ one either has $\im \rho\subset M$ for a Morse set
$M\in \cM$, or there exist two Morse sets $M^\pm\in\cM$ with $M^+ < M^-$
and for which $\uim^-(\rho)\subset M^-$ and $\uim^+(\rho)\subset M^+$ are satisfied.
By a {\em Morse decomposition of $\cV$} we mean a Morse decomposition of the maximal invariant
subset of $\cV$ in $X$.

We would like to point out that this definition of Morse decomposition differs in two aspects from
the one given in~\cite{LKMW2020}.  The definition in~\cite{LKMW2020},
instead of the condition with $\uim^-(\rho)$ and $\uim^+(\rho)$, uses the condition with  $\alpha$-
and $\omega$-limit sets of $\rho$. The reader familiar with~\cite{LKMW2020}
will, however, immediately notice that these conditions are equivalent,
in view of Proposition~\ref{prop:isolation-criterion}.
Also, paper~\cite{LKMW2020} defines a {\em global Morse decomposition}, that is a Morse decomposition of the whole space $X$
instead of an isolated invariant set $S$ of $\cV$ in $X$. Since the definition in~\cite{LKMW2020}
is given under the assumption that $X$ is invariant, 
one can easily regain the definition of Morse decomposition of an isolated invariant set $S$
from the definition in~\cite{LKMW2020} by replacing $X$ with $S$ and $\cV$ with the multivector
field $\cV_S$ induced by $\cV$ on $S$.
Global Morse decompositions of a combinatorial multivector field can easily be
determined via the strongly connected components of the associated
$\cV$-digraph. In fact, in contrast to the classical case, there always
exists a finest Morse decomposition for~$\cV$. 

In our definition of Morse decomposition we assume that all Morse sets are non-empty.
Although this is in the spirit of the classical definition by Conley~\cite{Conley1978},
in terms of computational applications the assumption may be inconvenient, because 
there may be areas where it is impossible to decide whether they contain a non-empty 
invariant set due to limited data or computational power.
As a remedy to this situation we define connection matrices for block decompositions, 
a concept slightly more general than Morse decompositions. 

\begin{defn}
  {\em Let $\cV$ be a combinatorial multivector field on a Lefschetz complex $X$ and let $S\subset X$ be
  an isolated invariant set.
  A {\em block} in $S$ is a locally closed and $\cV$-compatible subset of  $S$.
  A {\em block decomposition} of $S$ is a family $\cM$
  of mutually disjoint, non-empty blocks in $S$, together with a partial order $\leq$ on $\cM$
  and such that:
\begin{itemize}
   \item[(BD1)] If $\tau$ is a path in~$S$ with $\lep{\tau}\in M$ and $\rep{\tau}\in M'$
                for some $M,M'\in \cM$, then $M\geq M'$, and additionally $M=M'$ implies $\im\tau\subset M=M'$.
   \item[(BD2)] For every essential solution $\rho\in S$ the set $\rho^{-1}\left(\bigcup\cM\right)$
                is left and right infinite.
\end{itemize}
  }
\end{defn}

It follows from Proposition~\ref{prop:isolation-criterion} that every isolated invariant set is a block. 
A block does not have to be an isolated invariant set, because it may fail to be invariant. 
But, if $B$ is a block, then one can easily verify that
\[
   \Inv B:=\setof{x\in B\mid \text{there is an essential solution through $x$ in $B$}}
\]
is invariant, locally closed, and $\cV$-compatible, hence an isolated invariant set, again by 
Proposition~\ref{prop:isolation-criterion}.
We leave the proof of the following simple proposition to the reader.
\begin{prop}
If $\cM$ is a Morse decomposition of $S$ then it is also a block decomposition of $S$. 
Conversely, if $\cM$ is a block decomposition of $S$, then 
\begin{equation}
\label{eq:cM-bullet}
\cM^\bullet:=\setof{\Inv M\mid M\in\cM,\, \Inv M\neq\emptyset}
\end{equation}
is a Morse decomposition of $S$.
\qed
\end{prop}

The problem with potentially empty invariant sets in a Morse decomposition
is also addressed in~\cite{KMV+2022}, and our proposed concept of block
decomposition fits well with the solution of this reference. The reader familiar
with the terminology of~\cite{KMV+2022} will easily notice that given a
block decomposition~$\cN$, the family $\cM:=\cN^\bullet$ becomes a Morse
representation in the sense of~\cite{KMV+2022}, and the order embedding
$\cM\to\cN$ which sends a nonempty invariant set $M\in\cM$ to the unique
block $N\in\cN$ containing~$M$ is then a Morse decomposition in their
sense.

It is often convenient to consider a Morse or block decomposition as a
family $\cM=\{M_p\mid p\in P\}$ indexed by a poset $P$ with the partial
order in~$\cM$ given by
\begin{equation}
\label{eq:Mr-Ms}
M_r \leq M_s\iff r\leq s.
\end{equation}
A block decomposition~$\cM$ shares with a Morse decomposition its crucial property, 
as the following proposition shows.

\begin{prop}
\label{prop:block-decomp}
Assume that $\cM=\{M_p\}_{p\in P}$ is a block decomposition of an isolated invariant set $S$.
Then for every essential solution $\rho$ in $S$ there exist $p^-,p^+\in P$ 
such that $p^-\geq p^+$ and
\begin{equation}
\label{eq:block-decomp-1}
\uim^-\rho\subset M_{p^-} \quad\mbox{and}\quad \uim^+\rho\subset M_{p^+}.
\end{equation}
Moreover, if $p^- = p^+$, then one also has
\begin{equation}
\label{eq:block-decomp-2}
\im\rho\subset M_{p^-}.
\end{equation}
\end{prop}
\proof
It follows from (BD1) that 
\[
P_\rho:=\setof{p\in P\mid M_p\cap\im\rho\neq\emptyset}
\]
is a linearly ordered subset of $P$.
Set $p^-:=\max P_\rho$, $p^+:=\min P_\rho$ and choose $k^-,k^+\in\ZZ$ such that 
$\rho(k^-)\in M_{p^-}$,  $\rho(k^+)\in M_{p^+}$.
We claim that
\begin{equation}
\label{eq:k-plus}
\rho([k^+,\infty))\subset M_{p^+}.
\end{equation}
To see this, assume to the contrary that $\rho(n)\not\in M_{p^+}$ for some $n>k^+$.
Since, by (BD2), the preimage $\rho^{-1}\left(\bigcup\cM\right)$  is right infinite, 
there exists an index $m>n$ such that $\rho(m)\in  \bigcup\cM$.
Then $\rho(m)\in M_{p^+}$, because otherwise we contradict the definition of $p^+$.
Consider now the path $\tau:=\rho_{|[k^+,m]}$. We have $\lep{\tau}=\rho(k^+)\in M_{p^+}$
and $\rep{\tau}=\rho(m)\in M_{p^+}$. Hence, we get from (BD1) that $\rho(n)=\tau(n)\in\im\tau\subset M_{p^+}$, 
a contradiction proving \eqref{eq:k-plus}.
Similarly we prove that 
\begin{equation}
\label{eq:k-minus}
\rho((-\infty,k^-])\subset M_{p^-}.
\end{equation}
Consider now the case $k^-\geq k^+$. Applying (BD1) to $\sigma:=\rho_{|[k^+,k^-]}$
we obtain the inequality $p^+\geq p^-$, which in turn forces the equality
$p^+= p^-$, because we have $p^+=\min\ P_\rho\leq \max\ P_\rho=p^-$.
We also get from \eqref{eq:k-plus} and  \eqref{eq:k-minus} that
\[
\im\rho\subset\rho((-\infty,k^-])\cup\rho([k^+,\infty))\subset M_{p^-},
\]
which shows that $\rho$ satisfies \eqref{eq:block-decomp-2}.
Hence, this solution also satisfies \eqref{eq:block-decomp-1}, because
$\uim^-\rho\subset\im\rho$ and $\uim^+\rho\subset\im\rho$.
This completes the proof in the case $k^-\geq k^+$.

Consider finally the case  $k^-< k^+$. Then, again by (BD1), one has the
inequality $p^-\geq p^+$, together with 
$\uim^-\rho\subset\rho((-\infty,k^-])\subset M_{p^-}$ by \eqref{eq:k-minus}, as well as
$\uim^+\rho\subset\rho([k^+,\infty))\subset M_{p^+}$ by \eqref{eq:k-plus}.
Hence,  $\rho$ satisfies \eqref{eq:block-decomp-1}. Moreover, 
if $p^-=p^+$, then
\[
\im\rho\subset \rho((-\infty,k^-])\cup\rho([k^-,k^+])\cup\rho([k^+,\infty))\subset M_{p^+},
\]
which shows that~$\rho$ satisfies \eqref{eq:block-decomp-1} when $p^-=p^+$.
\qed

\begin{cor}
\label{cor:block-decomp}
Assume $\cM=\{M_p\}_{p\in P}$ is a block decomposition of an isolated invariant set $S$.
Then every multivector $V\not\in\bigcup\cM$ is regular.
\end{cor}
\proof
Assume, to the contrary, that $V\not\in\bigcup\cM$ is critical. 
Let $\rho$ be a constant, full solution through an element  $x\in V$.
Then $\rho$ is essential and $\uim^-\rho=\uim^+\rho=\{x\}=\im\rho$.
Hence, by Proposition~\ref{prop:block-decomp}
there is a $p\in P$ such that $\im\rho\subset M_p$. It follows that $V\subset M_p$,
a contradiction.
\qed

\medskip
Let~$\cM$ be a Morse or block decomposition of $\cV$.
Then the family
\[
   \cM\cup\setof{V\in\cV\mid V\cap\bigcup\cM=\emptyset}
\]
is clearly a partition of~$X$. We call it the {\em partition induced by
the decomposition $\cM$} and denote it $\cE_{\cM,\cV}$, or just $\cE_\cM $
when $\cV$ is clear from the context.
Each $E\in\cE_\cM$ is either a block,
or it is a multivector which is not contained in any of the blocks.
For two sets $E,E'\in\cE_\cM$ we write~$E\preccurlyeq E'$ if and only if
there is a path~$\sigma$ such that~$\lep{\sigma}\in E'$ and~$\rep{\sigma}\in E$.
Then the following statement holds.
\begin{prop}
\label{prop:curly-poset}
  The relation $\preccurlyeq$ is a partial order in $ \cE_\cM$.
\end{prop}
\proof
  Clearly, the relation is reflexive and transitive.
  To see that it is antisymmetric, take $E,E'\in \cE_\cM$ such that $E\preccurlyeq E'$
  and $E'\preccurlyeq E$. We have to prove that  $E= E'$.
  Let $\sigma$ be a path such that $\lep{\sigma}\in E'$ and $\rep{\sigma}\in E$.
  In addition, let $\tau$ denote a path which satisfies both $\lep{\tau}\in E$
  and $\rep{\tau}\in E'$.

  Consider first the case when $E$ and $E'$ are multivectors disjoint from~$\bigcup\cM$.
  Then $[\lep{\sigma}]=E'=[\rep{\tau}]$ and $[\rep{\sigma}]=E=[\lep{\tau}]$.
  It follows that the paths~$\sigma$ and~$\tau$ can be concatenated to a full, periodic
  solution~$\rho$. Assume now that we have $E\neq E'$. Then the solution~$\rho$ is essential.
  Hence, by Proposition~\ref{prop:block-decomp} we can find $p^-,p^+\in P$
  such that $p^-\geq p^+$, $\uim^-\rho\subset M_{p^-}$ and $\uim^+\rho\subset M_{p^+}$.
  However, it follows from Proposition~\ref{prop:omega-periodic}
  that $\uim^-(\rho)=\im\rho=\uim^+(\rho)$. Therefore, we must have $p^-= p^+$.
  In consequence, $\lep{\sigma}\in E'\cap M_{p^-}$, which contradicts our assumption
  that both~$E$ and~$E'$ are multivectors disjoint from $\bigcup\cM$. Thus, we
  have $E=E'$ in this case.

  Consider now the case  when both $E$ and $E'$ are in~$\cM$.
  Then $E=M_p$ and~$E'=M_{p'}$ for some $p,p'\in P$.
  Hence, applying (BD1) to $\sigma$ we obtain the inequality $p'\geq p$,
  and after applying (BD1) to $\tau$ we further get $p\geq p'$.
  It follows that $p=p'$ and $E=M_p=M_{p'}=E'$.

  Finally, consider the case when one of the sets $E$, $E'$
  is in $\cM$ and the other is a multivector disjoint from $\bigcup\cM$.
  Without loss of generality we may assume that $E=M_p$ for some $p\in P$ 
  and $E'\cap \bigcup\cM=\emptyset$. 
  Then we have $[\lep{\sigma}]=E'=[\rep{\tau}]$.
  It follows that the paths~$\sigma$ and~$\tau$ can be concatenated  
  to a path $\pi$ such that $\lep{\pi}=\lep{\tau}\in M_p$ and
  $\rep{\pi}=\rep{\sigma}\in M_p$.
  Applying (BD1) to~$\pi$ we obtain the inclusion $\im\pi\subset M_p=E$. 
  Since $\lep{\sigma}\in\im\pi$ and $\lep{\sigma}\in E'$, 
  we get $E\cap E'\neq\emptyset$ which implies $E=E'$, 
  because $\cE_\cM$ is a partition. 
\qed

\medskip
Since~$\cE_\cM$ is a partition of the Lefschetz complex, one can also
consider the relation~$\preceq_{\cE_\cM}$ defined in~\eqref{eq:preceq}.
As the following result demonstrates, its transitive closure is the
partial order $\preccurlyeq$ from above.
\begin{prop}
\label{prop:Morse-to-acyclic-decomp}
The relation~$\preccurlyeq$ is the transitive closure of the
relation~$\preceq_{\cE_\cM}$ defined in~\eqref{eq:preceq}. Thus,
the two relations~$\preccurlyeq$ and~$\leq_{\cE_\cM}$ coincide, and
in particular, the family~$\cE_\cM$ is an acyclic partition of~$X$.
\end{prop}
\proof
We need to prove that for all $E,E'\in \cE_\cM$ one has the equivalence
\begin{equation}
\label{eq:Morse-to-acyclic-decomp-1}
  E\leq_{\cE_\cM}E' \quad\iff\quad  E\preccurlyeq E'.
\end{equation}
First, we will prove that
\begin{equation}
\label{eq:Morse-to-acyclic-decomp-2}
    E\cap \cl E'\neq\emptyset \quad\implies\quad E\preccurlyeq E'.
\end{equation}
Let $x\in E\cap \cl E'$ be arbitrary. Then $x\in\cl y$ for a $y\in E'$.
Thus, $\sigma:=y\cdot x$ is a path. Moreover,   $\lep{\sigma}=y\in E'$ and $\rep{\sigma}=x\in E$,
proving  \eqref{eq:Morse-to-acyclic-decomp-2}.
Since, by Proposition~\ref{prop:curly-poset},
the relation $\preccurlyeq$ is transitive,
it follows from property \eqref{eq:Morse-to-acyclic-decomp-2} that the left-hand side of
\eqref{eq:Morse-to-acyclic-decomp-1}  implies the right-hand side of
\eqref{eq:Morse-to-acyclic-decomp-1}.

For the reverse implication assume that $E\preccurlyeq E'$.
Let $\rho=x_0\cdot x_1\cdot\ldots\cdot x_n$ be a path from $E'$ to $E$.
Let $V_i:=[x_i]$. Then $V_0\subset E'$ and $V_n\subset E$.
Moreover, since $x_{i+1}\in \cl x_i\cup[x_i]$, we see that 
$V_{i+1}\cap\cl V_i\neq\emptyset$ or $V_{i+1}=V_i$.
It follows that $V_{i+1}\preceq_{\cE_\cM} V_i$, and therefore $E\leq_{\cE_\cM} E'$.
Thus, \eqref{eq:Morse-to-acyclic-decomp-1} is proved. Hence, it follows from
Proposition~\ref{prop:curly-poset} that $\leq_{\cE_\cM}$ is a partial order
on~$\cE_\cM$, and this finally implies that~$\cE_\cM$ is an acyclic partition.
\qed

\medskip
Proposition~\ref{prop:Morse-to-acyclic-decomp} lets us consider $(\cE_\cM,\preccurlyeq)$
as a poset. We also consider it as an object of $\DPSet$ with the distinguished subset 
given by
\[ 
  (\cE_\cM)_\star:=\cM.
\]
Note that such a definition of $(\cE_\cM)_\star$ satisfies \eqref{eq:P-star-condition},
because it follows from Corollary~\ref{cor:block-decomp} that 
every $V\in\cV\setminus\bigcup\cM$ must be regular which implies that~$C(V)$ is inessential. 
%
\subsection{Connection matrices and heteroclinics}

After these preparations, Proposition~\ref{prop:Morse-to-acyclic-decomp}
finally enables us to define the Conley complex and associated connection
matrix of a given Morse or block decomposition.
\begin{defn}
{\em
\label{defn:conn-matrix-Morese-decomp}
The {\em Conley complex} and the associated {\em connection matrix} of a
Morse or block decomposition~$\cM$ of a combinatorial multivector field~$\cV$
on a Lefschetz complex~$X$ is defined as the Conley complex and the associated
connection matrix of the acyclic partition~$\cE_\cM$, i.e., the Conley
complex and respective connection matrix of the poset filtered chain
complex~$(\cE_\cM,C(X),\bdy^\kappa)$.
}
\end{defn}

\begin{ex}
\label{ex:lefschetzMVF-connMatr}
{\em 
  All three combinatorial multivector fields $\cV_0$, $\cV_1$, $\cV_2$ in
  Figure~\ref{fig:t2eLefMVF} have a common Morse decomposition given as the family
\[
\{\{\mathbf{A}\},\{\mathbf{B}\},\{\mathbf{AB}\},
  \{\mathbf{CD}\},\{\mathbf{CE}\}\},
\]
  which consists of singletons.
  However, the acyclic partition~$\cE_i:=\cE_{\cM,\cV_i}$ is different for each
  index $i\in\{0,1,2\}$. In these simple examples, $\cE_i$ coincides
  with~$\cV_i$. The filtered  chain complex $(\cE_0,C(X),\bdop^\kappa)$
  coincides with the filtered chain complex considered in Example~\ref{ex:pfcc}.
  Therefore, the connection matrices given by \eqref{eq:d1m} and \eqref{eq:d2m}
  are the connection matrices of~$\cM$ for $\cV_0$. We know that these matrices are not
  equivalent.
\exend}
\end{ex}

For a subfamily $\cL\subset\cM$ we define $M(\cL)$ as the collection of cells $x\in X$ for
which there exists an essential solution $\rho$ through $x$ and $M_-,M_+\in\cL$ such that
$\uim^-(\rho)\subset M^-$ and $\uim^+(\rho)\subset M^+$. 

\begin{prop}
\label{prop:M-cL}
The set $M(\cL)\subset X$ is an isolated invariant set for $\cV$. Moreover, $\cL$ is a Morse
decomposition of $M(\cL)$.
\end{prop}
\proof The first statement is the contents of~\cite[Theorem~7.4]{LKMW2020}.
The second statement is straightforward.
\qed

\medskip
As an immediate consequence of Theorem~\ref{thm:Conley-subcomplex} and Corollary~\ref{cor:hom-conley-subcomplex}
we obtain the following proposition. 
\begin{prop}
\label{prop:Conley-index-subcomplex}
Assume $(\cM,\bar{C},\bar{d})$ is a Conley complex of a Morse decomposition $\cM$ of an isolated invariant set $S$.
If $\cL\subset\cM$ is convex with respect to the partial order $\leq$ in $\cM$, then 
$(\cL,\bar{C}_\cL,\bar{d}_{\cL\cL})$ is a Conley complex of $\cL$ considered as a Morse decomposition of $M(\cL)$.
Moreover, the Conley index of $M(\cL)$ is isomorphic to $H(\bar{C}_\cL,\bar{d}_{\cL\cL})$.
\qed
\end{prop}

Clearly, every singleton $\{M\}\subset\cM$, consisting of just one Morse set $M$, is convex in $\cM$.
Since, by the definition of Conley complex, $\bar{d}_{MM}=0$, we get 
the following corollary from Proposition~\ref{prop:M-cL}.

\begin{cor}
\label{cor:Conley-index-M}
Given a Conley complex $(\cM,\bar{C},\bar{d})$ of a Morse decomposition $\cM$ of an isolated invariant set
$S$ of $\cV$, the Conley index $\Con(M)$ is isomorphic to $\bar{C}_M$ for every $M\in\cM$.
\qed
\end{cor}

In terms of applications, the most important result of classical connection matrix theory
states that if $\{M_p\}_{p\in P}$ is a Morse decomposition indexed by a poset $P$
and $p,q\in P$ are such that $q$ covers $p$ and the entry in the $q$th row and $p$th column
of the connection matrix is non-zero, then there must be a heteroclinic connection running from $M_q$ to $M_p$.
We have an analogous result for combinatorial multivector fields. 

Let $(\cM,\bar{C},\bar{d})$ be a Conley complex of a Morse decomposition $\cM=\{M_p\}_{p\in P}$.
Denote by~$A$ the associated connection matrix. 
Then~$A$ may be viewed as a block matrix $[A_{rs}]_{r,s\in P}$ by identifying~$P$ with~$\cM$. 

\begin{thm}
\label{thm:h-conn}
Consider $p,q\in P$ such that $p<q$ and $\{p,q\}$ is convex in~$P$. If $A_{pq}\neq 0$, then there is a heteroclinic connection
from $M_q$ to $M_p$, that is an essential solution $\rho$ satisfying
\begin{equation}
\label{eq:h-conn}
\uim^-(\rho)\subset M_q \quad\text{ and }\quad \uim^+(\rho)\subset M_p.
\end{equation}
\end{thm}
\proof
Assume to the contrary that there is no heteroclinic connection from $M_q$ to $M_p$.
Set $M_{pq}:=M(\{M_p,M_q\})$. It follows from Proposition~\ref{prop:M-cL}
that $M_{pq}$ is an isolated invariant set with Morse decomposition $\{M_p,M_q\}$.
We will prove that 
\begin{equation}
\label{eq:Sol-Mpq-1}
M_{pq}=M_p\cup M_q.
\end{equation}
For this, consider an arbitrary $x\in M_{pq}$.
By the definition of $M_{pq}$, there are $M_-,M_+\in \{M_p,M_q\}$ 
and an essential solution $\rho$ through $x$ in $M_{pq}$ such that 
$\uim^-(\rho)\subset M^-$ and $\uim^+(\rho)\subset M^+$.
We cannot have $M_-=M_p$ and $M_+=M_q$, because then $M_p>M_q$, and by \eqref{eq:Mr-Ms} we get $p>q$, which contradicts 
assumption $p<q$. We also cannot have $M_-=M_q$ and $M_+=M_p$, because otherwise $\rho$ would be a heteroclinic
connection from $M_q$ to $M_p$. Hence, either $M_-=M_+=M_p$ or $M_-=M_+=M_q$. Since $\{M_p,M_q\}$ is a Morse decomposition of $M_{pq}$,
in the fist case we get $x\in\im\gamma\subset M_p$ and in the second case we get $x\in\im\gamma\subset M_q$.
This proves that  $M_{pq}\subset M_p\cup M_q$.
Since the opposite inclusion is obvious, the proof of \eqref{eq:Sol-Mpq-1} is completed.  
We also have
\begin{align}
&M_q\cap\cl M_p=\emptyset,\label{eq:Sol-Mpq-2}\\
&M_p\cap \cl M_q=\emptyset.\label{eq:Sol-Mpq-3}
\end{align}
Indeed, to prove~\eqref{eq:Sol-Mpq-2} assume to the contrary that $M_q\cap\cl M_p\neq\emptyset$.
Then $\bdop^\kappa_{M_qM_p}\neq 0$, which gives $M_q<M_p$, because $\bdop^\kappa$ is filtered. 
But then $q<p$, which contradicts assumption $p<q$.  
To see~\eqref{eq:Sol-Mpq-3} assume to the contrary that there is an $x\in M_p\cap\cl M_q$. 
Then $x\in\cl y$ for some $y\in M_q$ and since $M_p\neq\emptyset\neq M_q$, we can find an
essential solution $\gamma_q$ through $y$ in $M_q$ and  $\gamma_p$ through $x$ in $M_p$.
Then $\gamma_q^-\cdot\gamma_p^+$ is a heteroclinic connection from $M_q$ to $M_p$, a contradiction. 

It follows from (\ref{eq:Sol-Mpq-1}-\ref{eq:Sol-Mpq-3}) that $M_{pq}$, $M_p$, $M_q$ satisfy the assumptions
of \cite[Theorem 5.19]{LKMW2020}. This lets us conclude that 
\[
\Con(M_{pq})=\Con(M_p)\oplus\Con(M_q)
\]
and, in consequence, that 
\begin{equation}
\label{eq:Con-Mpq}
n_{pq}=n_p+n_q
\end{equation}
where $n_{pq}$, $n_p$ and $n_q$ denote respectively the dimension of $\Con(M_{pq})$, $\Con(M_{p})$ and $\Con(M_{q})$.
However, from Proposition~\ref{prop:Conley-index-subcomplex} we know that the block matrix
\[
\bar{A}_{pq}:=\left[
              \begin{array}{cc}
                0 & A_{pq} \\
                0 & 0
              \end{array}
              \right]
\]
is a connection matrix of the Conley complex for the Morse decomposition  
$\{M_p,M_q\}$ of $M_{pq}$ and $\Con(M_{pq})$ is isomorphic to the homology of the Conley complex 
of the Morse decomposition $\cM$ restricted to $\{M_p,M_q\}$.
By Corollary~\ref{cor:Conley-index-M} we know that the $2 \times 2$ block matrix $\bar{A}_{pq}$
has exactly $n_p$ columns in the first block and $n_q$ columns in the second block. 
Since $A_{pq}\neq 0$, we see that the kernel of $\bar{A}_{pq}$ is of dimension less than $n_p+n_q$.
Therefore, the dimension of the homology of the restricted Conley complex and, in consequence, 
also the dimension of $\Con(M_{pq})$ is less than $n_p+n_q$, which contradicts \eqref{eq:Con-Mpq}
and proves the theorem. 
\qed


\section{Connection matrices for Forman's gradient vector fields}
\label{sec:cm-gcvf}

This last theoretical section of the paper is devoted to the study of the special case
of Forman's gradient vector fields on regular Lefschetz complexes. In this
situation, one can show that the associated connection matrix is uniquely
determined, and in fact can be determined in a direct way. This is accomplished
in a number of steps. We begin by recalling basic properties of combinatorial
vectors in the sense of Forman, before we present the definition and basic
properties of his concept of combinatorial flow. After discussing the
long-term limit of the latter, we can finally explain how it can be used 
to obtain the unique connection matrix in this setting. The uniqueness
part of the assertion will rely heavily on the singleton partition
which has already been discussed in Section~\ref{subsec-singletonpart}.

\subsection{Forman's combinatorial flow}

In the following, we assume that~$\cV$ is a gradient combinatorial vector
field in the sense of Forman on a Lefschetz complex~$X$. In terms of its
definition, we recall that a {\em combinatorial vector} is just a combinatorial
multivector~$W$ with $\card W\leq 2$. A combinatorial multivector field
whose multivectors are just vectors is then called a {\em combinatorial
vector field}, a concept originally introduced by Forman~\cite{Fo98b, Fo98a}
in a slightly different but equivalent form.
In the case of a combinatorial vector field the critical vectors are precisely
the singletons, and the regular vectors are precisely the doubletons, i.e.,
vectors of cardinality two (see Proposition~\ref{prop:singl-dubl-homology}). Therefore, every gradient-like combinatorial
vector field is a gradient combinatorial vector field.

Let $\cC\subset\cV$ denote the collection of critical vectors in~$\cV$.
Since we assume that~$\cV$ is a combinatorial vector field, it follows
from Proposition~\ref{prop:singl-dubl-homology} that~$\cC$ is exactly
the collection of singletons in~$\cV$. Moreover, the following proposition
is straightforward.
\begin{prop}
\label{prop:crit-M-decomp}
Assume that~$\cV$ is a gradient combinatorial vector field.
Then the collection~$\cC$ is a Morse decomposition of~$\cV$ and
the induced acyclic partition of the Lefschetz complex~$X$ is
given by~$\cE_{\cC}=\cV$.
\qed
\end{prop}
It follows that the Conley complex of the Morse decomposition~$\cC$ of~$X$
is the Conley complex of the filtered chain complex~$(\cV,C(X),\bdy^\kappa)$.
The aim of this section is to prove that the Morse decomposition $\cC$ has
precisely one connection matrix, and that this connection matrix coincides
with the matrix of the boundary operator of the associated Morse complex,
see also~\cite[Section~7]{Fo98a}. In order to make this statement precise,
we first need to recall some concepts, with the current subsection focusing
on properties of combinatorial vectors.

Unlike a general multivector, a combinatorial vector~$W\in\cV$ contains a
unique minimal element and a unique maximal element in~$W$ with respect
to the face relation~$\leq_\kappa$. We denote them by~$W^-$ and~$W^+$,
respectively, and we extend this notation to cells by writing $x^-:=[x]^-$
and $x^+:=[x]^+$. Note that a combinatorial vector is given by
$W=\{W^-,W^+\}$, and~$W$ is critical if and only if we have $W^+=W^-$,
which in turn is equivalent to assuming that~$W$ is a singleton. As a
consequence of Proposition~\ref{prop:dim-monotoonicity} one obtains for
a vector~$W$ the inequality
\begin{equation}
\label{eq:dim-W+W-}
0 \le \dim W^+-\dim W^-\leq 1,
\end{equation}
and for all $x\in W$ one has
\begin{equation}
\label{eq:dim-W-xW+}
\dim W^-\leq\dim x\leq\dim W^+.
\end{equation}
We say that a cell~$x$ is a {\em tail} if $x=x^-\neq x^+$ and a
{\em head} if $x=x^+\neq x^-$. Clearly, if a cell is neither a tail
nor a head, then it is critical. We denote the subsets of critical
cells, tails, and heads by~$X^c(\cV)$, $X^-(\cV)$, and~$X^+(\cV)$,
respectively, and shorten this notation to~$X^c$, $X^-$, and~$X^+$
whenever~$\cV$ is clear from context. Note that the collection~$\cC$
of critical vectors of~$\cV$ is precisely $\setof{\{x\}\mid x\in X^c}$.

For a combinatorial vector $V\in\cV$ we set $\dim V:=\dim V^+$. Moreover,
recall that~$\cV$, as a gradient vector field, is an acyclic partition
of~$X$ according to Proposition~\ref{prop:grdient-like-mvf}. In particular,
the collection~$\cV$ is a poset with the partial order~$\leq_\cV$. Then
one has the following proposition.
\begin{prop}
\label{prop:declining-dim}
Assume that $\cV$ is a gradient vector field on a Lefschetz complex $X$.
Then the following hold:
\begin{itemize}
  \item[(i)] The map $\dim:(\cV,\leq_\cV)\to(\NN_0,\leq)$ is order preserving.
  \item[(ii)] If $V<_\cV W$ and~$V$ is critical, then $\dim V<\dim W$.
  \item[(iii)] The choice $\cV_\star=\cC$ is possible, since the set of $V\in\cV$ such
     that~$C(V)$ is homotopically essential coincides with the set of
     critical vectors~$\cC$.
\end{itemize}
\end{prop}
\proof
  Since, by definition, the relation $\leq_\cV$ is the transitive closure of
  the relation~$\preceq_\cV$, to prove (i)
  it suffices to verify that $V\cap\cl W\neq\emptyset$ for $V,W\in\cV$ implies $\dim V\leq\dim W$.
  Thus, assume that  $V\cap\cl W\neq\emptyset$.
  Let $x\in V\cap\cl W$.  If $x\in W$, then $V\cap W\neq\emptyset$
  which implies $V=W$, because $\cV$ is a partition.
  In particular, this gives $\dim V=\dim W$.
  Thus, consider now case $x\not\in W$. Since $W^-\in\cl W^+$,
  we have $\cl W=\cl\{W^-,W^+\} = \cl W^-\cup\cl W^+=\cl W^+$.
  It follows that  $x\in\cl W^+\setminus\{W^+\}$, and therefore $\dim x<\dim W^+$.
  Since we assumed $x\in V$, we get from \eqref{eq:dim-W+W-} and \eqref{eq:dim-W-xW+}
  that $\dim x\geq\dim V^-\geq\dim V^+-1=\dim V-1$.
  Thus, $\dim V \le \dim x+1\leq\dim W^+=\dim W$.

  To see (ii), without loss of generality, we may assume that $V\prec_\cV W$,
  that is, $V\neq W$ and $V\cap\cl W\neq\emptyset$.
  Observe that since $V$ is critical, both $V=\{x\}$ and $\dim V=\dim x$ are satisfied.
  Since $x\not\in W$ and $x\in\cl W=\cl W^+$, we have $\dim x<\dim W^+$.
  Therefore, $\dim V=\dim x<\dim W^+=\dim W$.

  To see (iii) observe that, by Proposition~\ref{prop:singl-dubl-homology}, the
  Lefschetz homology of~$V$ is zero if and only if~$V$ is not critical.
  Therefore, property (iii) follows immediately from
  Theorem~\ref{thm:existence-and-uniqueness}(iv).
\qed

\medskip
The next piece of the puzzle is Forman's combinatorial flow, which
he introduced in~\cite[Definition~6.2]{Fo98a} and which we use in a 
slightly more general context. More precisely, we now consider a
regular Lefschetz complex~$X$. Then with a combinatorial gradient vector
field~$\cV$ on~$X$ we can associate the degree~$+1$ map
$\Gamma_\cV: C(X)\to C(X)$, which for $x\in X$ is defined by
\begin{equation}
\label{eq:Gamma-def}
     \Gamma_\cV x:= \left\{
     \begin{array}{cl}
         0 & \text{if $x^-=x^+$,}\\[1ex]
         -\kappa(x^+,x^-)^{-1}\scalprod{x,x^-}x^+ & \text{otherwise.}
     \end{array}\right.
\end{equation}
In the sequel, we drop the subscript in~$\Gamma_\cV$ whenever~$\cV$ is
clear from context. Then the following propositions are straightforward.
\begin{prop}
\label{prop:Gamma-supp}
If~$X$ is a regular Lefschetz complex and~$\cV$ is a combinatorial
gradient vector field on~$X$, then for every $c\in C(X)$ we have
$|\Gamma c|\subset X^+$.
\end{prop}
\begin{prop}
\label{prop:Gamma-on-x}
Suppose that~$X$ is a regular Lefschetz complex and~$\cV$ is a
combinatorial gradient vector field on~$X$. Then for all $x,u\in X$ we have
both the implication
\begin{equation}
\label{eq:Gamma-on-x-1}
  \scalprod{u,\Gamma\bdy x}\neq 0 \;\;\implies\;\;
  u\in X^+\text{ and }\; \scalprod{u,\Gamma\bdy x}=-\kappa(x,u^-)\kappa(u^+,u^-)^{-1}
\end{equation}
and the implication
\begin{equation}
\label{eq:Gamma-on-x-2}
  \scalprod{u,\bdy\Gamma x}\neq 0 \;\;\implies\;\;
  x\in X^- \text{ and }\;\scalprod{u,\bdy\Gamma x}=-\kappa(x^+,x^-)^{-1}\scalprod{u,\bdy x^+}.
\end{equation}
\end{prop}
\proof
The definition of the boundary operator of a Lefschetz complex given in~\eqref{eq:kappa-bdy}
implies
\begin{eqnarray*}
  \Gamma\bdy x&=&\Gamma\left(\sum_{y\in X}\kappa(x,y)y\right)
  =\sum_{y\in X}\kappa(x,y)\Gamma y \nonumber \\[1ex]
  &=&-\sum_{y\in X, \; \Gamma y\neq 0}\kappa(x,y)\kappa(y^+,y^-)^{-1}\scalprod{y,y^-}y^+.
\end{eqnarray*}
Hence, we have
\begin{displaymath}
   \scalprod{u,\Gamma\bdy x}=
  -\sum_{y\in X, \;\Gamma y\neq 0}\kappa(x,y)\kappa(y^+,y^-)^{-1}\scalprod{y,y^-}\scalprod{u,y^+}.
\end{displaymath}
Thus, if $\scalprod{u,\Gamma\bdy x}$ is non-zero, then there is a $y\in X$
such that $y^-\neq y^+$ and that
all factors in $\kappa(x,y)\scalprod{y,y^-}\scalprod{u,y^+}$ are non-zero.
In particular, the inequality $\scalprod{y,y^-}\neq 0$ implies
$y=y^-$. Hence, $y\in X^-$ and from $\scalprod{u,y^+}\neq 0$ we get $u=y^+\in X^+$
and $y=y^-=[y^+]^-=[u]^-=u^-$. It follows that such a $y$ is unique and
$\scalprod{u,\Gamma\bdy x}=-\kappa(x,y)\kappa(y^+,y^-)^{-1}=-\kappa(x,u^-)\kappa(u^+,u^-)^{-1}$.
This proves \eqref{eq:Gamma-on-x-1}.

To see  \eqref{eq:Gamma-on-x-2} observe that if $\scalprod{u,\bdy\Gamma x}\neq 0$,
then $\Gamma x\neq 0$. Hence, $x^-\neq x^+$ and one further obtains
$$
  \scalprod{u,\bdy\Gamma x}=\scalprod{u,-\kappa(x^+,x^-)^{-1}\scalprod{x,x^-}\bdy x^+}=
  -\kappa(x^+,x^-)^{-1}\scalprod{x,x^-}\scalprod{u,\bdy x^+}.
$$
Therefore, if $\scalprod{u,\bdy\Gamma x}\neq 0$, then
$\scalprod{x,x^-}\neq 0\neq \scalprod{u,\bdy x^+}$.
In particular, we get $x=x^-$, which in turn means that $x\in X^-$ and
$$
  \scalprod{u,\bdy\Gamma x}= -\kappa(x^+,x^-)^{-1}\scalprod{u,\bdy x^+}.
$$
This identity finally proves~\eqref{eq:Gamma-on-x-2}.
\qed

\medskip
For a combinatorial gradient vector field~$\cV$, we
follow~\cite[Theorem~6.4]{Fo98a} and define the associated
{\em combinatorial flow\/} $\Phi:=\Phi_\cV: C(X)\to C(X)$
on generators $x\in X$ through the formula
\begin{equation}
\label{eq:Phi}
   \Phi_\cV x:=x+\bdy \Gamma_\cV x + \Gamma_\cV\bdy x.
\end{equation}
In the sequel, whenever~$\cV$ is clear from the context, we drop the
subscript in~$\Phi_\cV$. Note that $\Phi$ is a degree zero module
homomorphism which satisfies the identity
$\Phi\bdy=\bdy+\bdy\Gamma\bdy=\bdy\Phi$.
Hence, $\Phi$ is a chain map. The definition of the combinatorial
flow is illustrated in the following example.
\begin{figure}
  \begin{center}
    \includegraphics[width=0.45\textwidth]{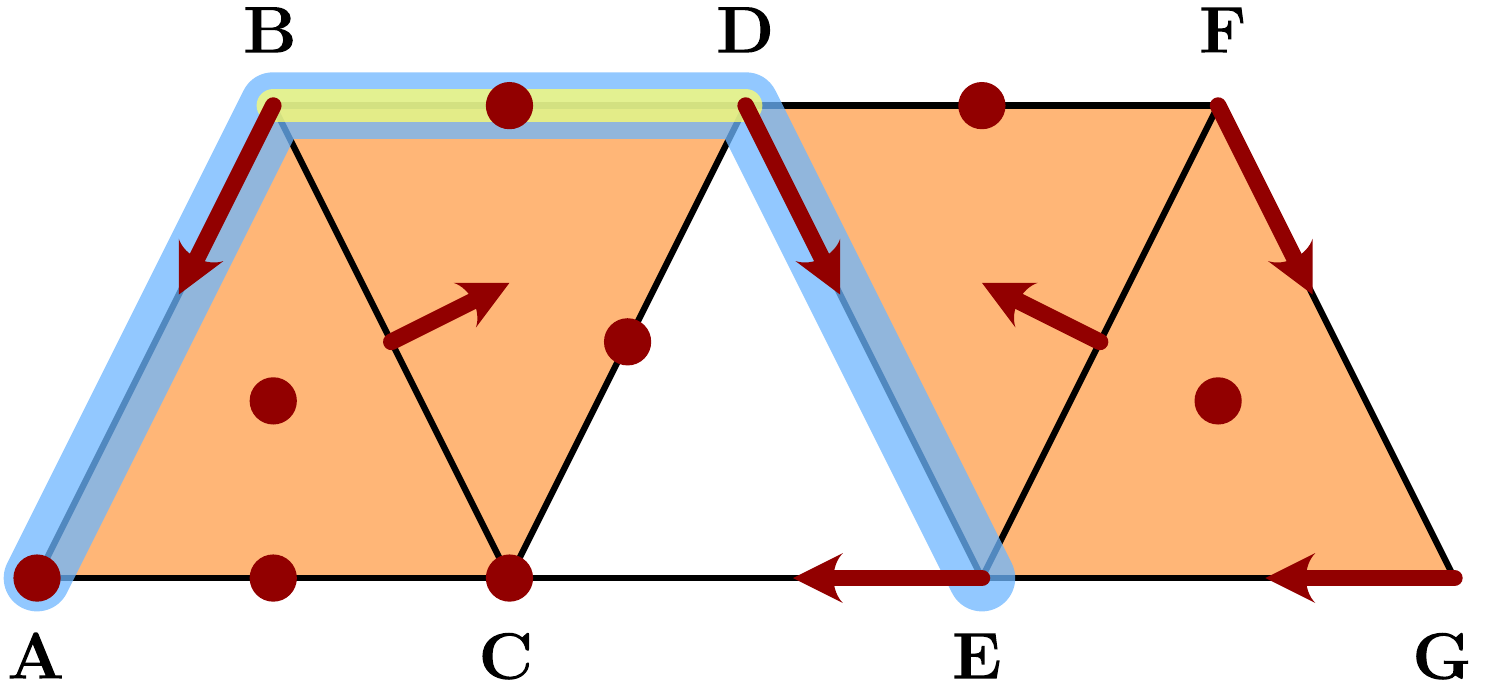}\quad\quad
    \includegraphics[width=0.45\textwidth]{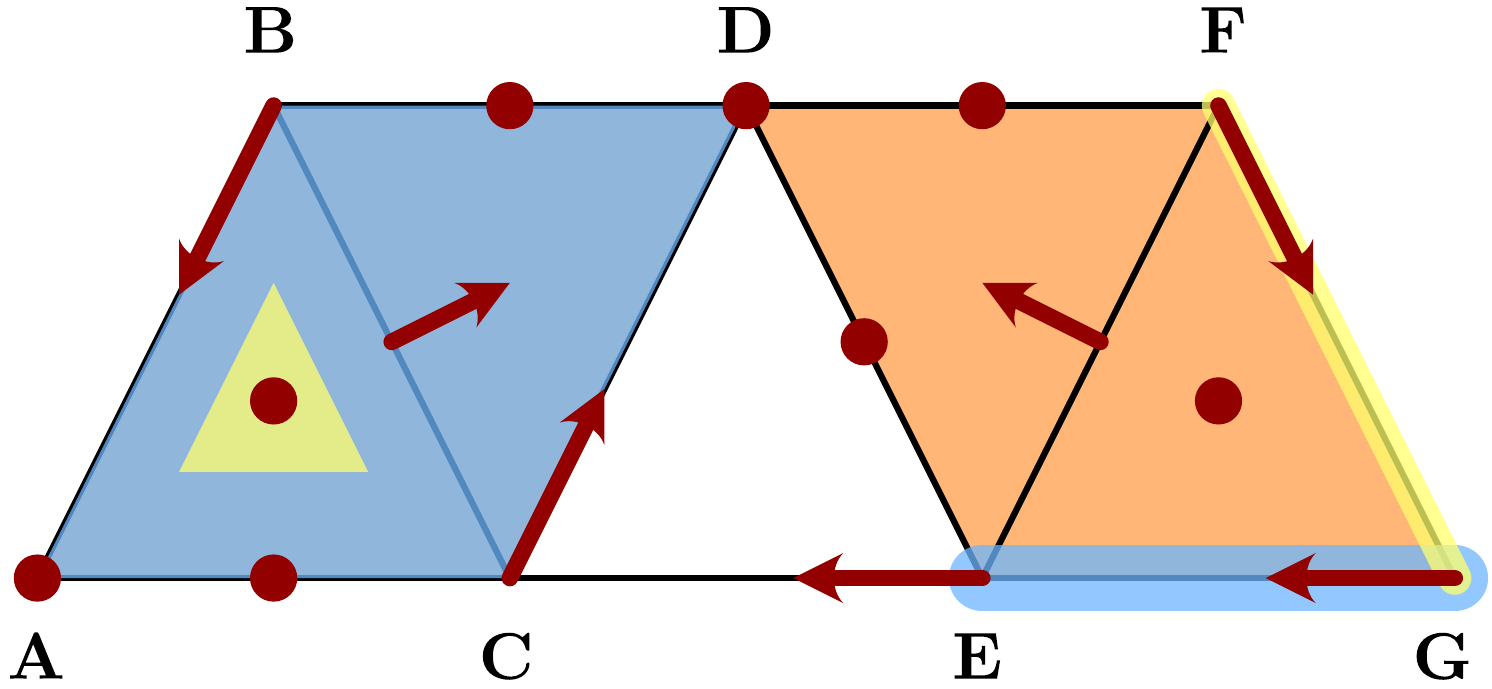}
  \end{center}
  \caption{{\em The combinatorial flow~$\Phi$}.
           The above two panels illustrate the action of the
           combinatorial flows associated with the two
           combinatorial vector fields introduced in the top two
           rows of Figure~\ref{fig:periodicex123}. In the left panel,
           both the chain~$\bB\bD$ and its image~$\Phi_1(\bB\bD)$
           are shown in yellow and blue, respectively. In the
           panel on the right one can find~$\Phi_2(\bA\bB\bC)$
           and~$\Phi_2(\bF\bG)$, where again the arguments are
           depicted in yellow and the images under the
           combinatorial flow in blue.
           }
  \label{fig:phiex}
\end{figure}

\begin{ex}[{\em Three Forman gradient vector fields}]
\label{ex:formangradient-2}
{\em
We return to the setting of Example~\ref{ex:formangradient-1},
where we introduced the three combinatorial vector fields
shown in Figure~\ref{fig:periodicex123}. For this, consider the 
underlying simplicial complex~$X$ as a Lefschetz complex with
$\ZZ_2$-coefficients, that is, the simplices making up the
complex~$X$ are not oriented.

In order to describe the action of the combinatorial flows~$\Phi$
associated with the two gradient combinatorial vector fields
shown in Figure~\ref{fig:phiex}, one first needs to understand
the chain map~$\Gamma$ defined in~\eqref{eq:Gamma-def}. It can
easily be seen that due to the used $\ZZ_2$-coefficients, the
image~$\Gamma(x)$ is nonzero if and only if $x = x^- \neq x^+$,
and in this case we have~$\Gamma(x^-) = x^+$.

With this, we can now turn our attention to a few sample
combinatorial flow computations. We begin by considering the
combinatorial vector field~$\cV_1$ shown in the left panel of
Figure~\ref{fig:phiex}. For the edge~$\bB\bD$ one obtains
\begin{eqnarray*}
  \Phi_1(\bB\bD) & = & \bB\bD + \partial(\Gamma_1(\bB\bD))
    + \Gamma_1(\partial \bB\bD) \\
  & = & \bB\bD + \partial (0) + \Gamma_1(\bB + \bD)
    \;\; = \;\; \bB\bD + \bA\bB + \bD\bE ,
\end{eqnarray*}
and this is illustrated in the left image of Figure~\ref{fig:phiex}.
The edge~$\bB\bD$ is shown in yellow, while the chain~$\Phi_1(\bB\bD)$
is drawn in blue. Notice that the action of the combinatorial
flow~$\Phi_1$ encodes potential solution paths of the combinatorial
vector field~$\cV_1$ starting at~$\bB\bD$.

For the gradient combinatorial vector field~$\cV_2$ shown in the
right panel of Figure~\ref{fig:phiex} we consider two examples.
First, the edge~$\bF\bG$ has the image
\begin{displaymath}
  \Phi_2(\bF\bG) =
  \bF\bG + \partial(0) + \Gamma(\bF + \bG) =
  \bF\bG + \bF\bG + \bE\bG =
  \bE\bG,
\end{displaymath}
since both $\Gamma(\bF\bG) = 0$ and $\bF\bG + \bF\bG = 0$ are satisfied.
In addition, for the two-dimensional cell~$\bA\bB\bC$ one obtains
the identity
\begin{eqnarray*}
  \Phi_2(\bA\bB\bC) & = &
    \bA\bB\bC + \partial(0) + \Gamma_2(\bA\bB + \bA\bC + \bB\bC) \\
  & = & 
  \bA\bB\bC + 0 + 0 + \bB\bC\bD \;\; = \;\;
  \bA\bB\bC + \bB\bC\bD,
\end{eqnarray*}
in view of $\Gamma_2(\bA\bB\bC) = 0$ and $\Gamma_2(\bA\bB) =
\Gamma_2(\bA\bC) = 0$, as well as using our above discussion to
show that $\Gamma_2(\bB\bC) = \bB\bC\bD$. Both of these images
under the combinatorial flow~$\Phi_2$ are depicted in blue in
the right panel of Figure~\ref{fig:phiex}, while the considered
arguments are shown in yellow.
\exend
}
\end{ex}

The above sample computations demonstrate that the combinatorial
flow in some sense encodes the potential paths that solutions of
the combinatorial vector field can take. For later use, we now
present four propositions, which provide more rigorous insight
into the action of the combinatorial flow~$\Phi$.
\begin{prop}
\label{prop:Phi-on-x}
Suppose that~$X$ is a regular Lefschetz complex and~$\cV$ is a
combinatorial gradient vector field on~$X$. Then for all $x\in X$ we have
\begin{equation*}
  \scalprod{x,\Phi x}=
  \begin{cases}
     1 & \text{if $x\in X^c$,}\\
     0 & \text{otherwise.}
  \end{cases}
\end{equation*}
\end{prop}
\proof
We have $\scalprod{x,\Phi x}=\scalprod{x,x}+\scalprod{x,\bdy\Gamma x}+\scalprod{x,\Gamma \bdy x}$.
In the case $x\in X^c$
we get from Proposition~\ref{prop:Gamma-on-x} that
$\scalprod{x,\bdy\Gamma x}=0$ and $\scalprod{x,\Gamma \bdy x}=0$,
which implies $\scalprod{x,\Phi x}=\scalprod{x,x}=1$.

Consider now the case  $x\in X^-$.
Then one has $x=x^-\not \in X^+$ and one obtains from~\eqref{eq:Gamma-on-x-1} that
$\scalprod{x,\Gamma\bdy x}=0$.
Therefore, the definition~\eqref{eq:Gamma-def} yields
\begin{multline*}
\scalprod{x,\Phi x}=\scalprod{x,x}+\scalprod{x,\bdy\Gamma  x}=\\
1-\kappa(x^+,x^-)^{-1}\scalprod{x,x^-}\scalprod{x,\bdy x^+}=1-\kappa(x^+,x^-)^{-1}\kappa(x^+,x^-)=0.
\end{multline*}
Finally, we consider the remaining case $x\in X^+$.
Then $x\not \in X^-$, and~\eqref{eq:Gamma-on-x-2} shows that
$\scalprod{x,\bdy \Gamma  x}=0$ .
Therefore, since $\bdy x=\sum_{y\in X}\kappa(x,y)y$, one has
\begin{multline*}
\scalprod{x,\Phi x}=\scalprod{x,x}+\scalprod{x,\Gamma \bdy x}=
1+\scalprod{x,\sum_{y\in X}\kappa(x,y)\Gamma y}\\
= 1-\scalprod{x,\sum_{y\in X,\;\Gamma y\neq 0}\kappa(x,y)\kappa(y^+,y^-)^{-1}\scalprod{y,y^-} y^+}\\
= 1-\sum_{y\in X,\;\Gamma y\neq 0}\kappa(x,y)\kappa(y^+,y^-)^{-1}\scalprod{y,y^-}\scalprod{x,y^+}\\
= 1-\kappa(x^+,x^-)\kappa(x^+,x^-)^{-1}=0,
\end{multline*}
because the only non-zero term in the last sum occurs for the point $y\in X$
which satisfies both $y^+=x=x^+$ and $y=y^-=x^-$.
\qed

\begin{prop}
\label{prop:Phi-on-X+}
Suppose that~$X$ is a regular Lefschetz complex and~$\cV$ is a combinatorial
gradient vector field on~$X$. Then we have
$$
   \Phi(C(X^+))\subset C(X^+).
$$
\end{prop}
\proof
 Since $\Phi$ is linear, it suffices to show that the inclusion
 $x\in X^+$ implies $|\Phi x|\subset X^+$.
 Thus, take $x\in X^+$ and a $u\in |\Phi x|$. Then $\scalprod{u,\Phi x}\neq 0$.
 From  \eqref{eq:Gamma-on-x-2} we get  $\scalprod{u,\bdy\Gamma x}=0$.
 Therefore, $0\neq\scalprod{u,\Phi x}=\scalprod{u,x}+\scalprod{u,\Gamma \bdy x}$,
 which implies $\scalprod{u,x}\neq 0$ or $\scalprod{u,\Gamma\bdy x}\neq 0$.
 If $\scalprod{u,x}\neq 0$, then $u=x\in X^+$.
 If the inequality $\scalprod{u,\Gamma \bdy x}\neq 0$ holds, then we get from
 \eqref{eq:Gamma-on-x-1}
 that $u\in X^+$. Hence, $|\Phi x|\subset X^+$.
\qed

\begin{prop}
\label{prop:Phi-on-down sets}
Suppose that~$X$ is a regular Lefschetz complex and that~$\cV$ is a given
combinatorial gradient vector field on~$X$. Assume further that the inclusion
$\cA\in\Down(\cV,\leq_\cV)$ holds. Then $A:=|\cA|$ is a $\cV$-compatible and
closed subset of~$X$, and we have the inclusion $\Phi(C(A))\subset C(A)$.
Therefore, the restriction $\Phi_{|C(A)}:(C(A),\bdy^\kappa_{|C(A)})\to
(C(A),\bdy^\kappa_{|C(A)})$ is again a well-defined chain map.
\end{prop}
\proof
Clearly, $A$ is $\cV$-compatible, and it follows from Proposition~\ref{prop:convex-lcl}(i)
that~$A$ is closed.
Take an $x\in A$. We will prove that $\Phi x\in C(A)$. Since $A$ is closed,
we have  $|\bdy x|\subset A$.
Since $A$ is $\cV$-compatible, we have $|\Gamma_\cV x|\subset A$.
Thus, again by closedness and $\cV$-compatibility of $A$,
we get both $|\bdy \Gamma_\cV x|\subset A$ and  $|\Gamma_\cV \bdy x|\subset A$. Therefore,
$$|\Phi x|=|x+\bdy \Gamma_\cV x+\Gamma_\cV\bdy  x|\subset|x|\cup |\bdy \Gamma_\cV x|
\cup|\Gamma_\cV\bdy  x|\subset A.$$
This finally implies $\Phi x\in C(A)$ and completes the proof.
\qed

\begin{prop}
\label{prop:Phi-is-declining}
Suppose that~$X$ is a regular Lefschetz complex and~$\cV$ is a combinatorial
gradient vector field on~$X$. Assume further that $x,y\in X$. Then the following hold:
\begin{itemize}
  \item[(i)] If $y\in |\Phi x|$, then $[y]\leq_\cV[x]$, and there exists
             a path from~$x$ to~$y$ with respect to the multivalued
             flow map~$\Pi_\cV$ defined in~\eqref{def-multivalflow}.
  \item[(ii)] If $y\in |\Phi x|$ and $[y]=[x]$, then  $y=x$ .
  \item[(iii)] If $x\in |\Phi x|$, then  $\scalprod{x,\Phi x}=1$.
\end{itemize}
In addition, for any chain $c \in C(X)$ the following holds:
\begin{itemize}
  \item[(iv)] If $c = \Phi(c)$ and $c \neq 0$, then~$|c|$ has to contain a critical cell.
\end{itemize}
\end{prop}
\proof
  In order to see (i), observe that the assumption $y\in |\Phi x|$ implies the inequality
  $0\neq\scalprod{y,\Phi x}=\scalprod{y,x}+\scalprod{y,\bdy\Gamma x}+\scalprod{y,\Gamma\bdy x}$.
  Thus, either we have $\scalprod{y,x}\neq 0$, or $\scalprod{y,\bdy\Gamma x}\neq 0$,
  or $\scalprod{y,\Gamma\bdy x}\neq 0$.
  In the first case one immediately obtains $x=y$ and $[x]=[y]$.
  Consider now the second case $\scalprod{y,\bdy\Gamma x}\neq 0$. This inequality
  yields both $\Gamma x\neq 0$ and $y\in |\bdy\Gamma x|$.
  Hence, we have $x=x^-$ and $\Gamma x=\lambda x^+$
  with $\lambda=-\kappa(x^+,x^-)^{-1}\neq 0$.
  This gives $|\bdy\Gamma x|=|\bdy x^+|$, as well as $y\in|\bdy x^+|\subset\cl x^+$.
  It follows that $y\in[y]\cap\cl[x^+]=[y]\cap\cl[x^-]$, and therefore $[y]\leq_\cV[x]$.
  Finally, we consider the third inequality $\scalprod{y,\Gamma\bdy x}\neq 0$.
  Since $\bdy x=\sum_{u\in X}\kappa(x,u)u$, we have
 \[
     0\neq\scalprod{y,\Gamma\bdy x}=\sum_{u\in X}\kappa(x,u)\scalprod{y,\Gamma u},
 \]
 which implies $\scalprod{y,\Gamma u}\neq 0$ for some $u\in|\bdy x|$.
 But this implies that $u=u^-$ and  $\Gamma u=\lambda u^+$
 with $\lambda=-\kappa(u^+,u^-)^{-1}\neq 0$.
 Hence, the identities $y=u^+$ and $[y]=[u^+]=[u^-]$ are satisfied,
 as well as the inclusion $u^-\in|\bdy x|\subset\cl x$.
 Thus, $[y]\cap\cl[x]\neq\emptyset$, which implies $[y]\leq_\cV[x]$.
 Furthermore, in all of these three cases one can easily construct a path
 from~$x$ to~$y$ with respect to the multivalued map~$\Pi_\cV$. This finally
 completes the proof of~(i).
 
 In order to prove~(ii), we observe that $y\in |\Phi x|$ implies $\dim x=\dim y$,
 since~$\Phi$ is a degree zero module homomorphism. Since~$\cV$ is a combinatorial
 vector field, we further have both $[x]=\{x^-,x^+\}$ and $[y]=\{y^-,y^+\}$. Thus,
 the validity of $[x]=[y]$ and $y\in |\Phi x|$ implies the identity $x=y$,
 which in turn establishes~(ii).
 
 We now turn our attention to~(iii). If  $x\in |\Phi x|$, then $\scalprod{x,\Phi x}\neq 0$.
 Thus, we get from Proposition~\ref{prop:Phi-on-x} that $\scalprod{x,\Phi x}=1$.
 
 Finally, consider~(iv). Let~$y \in |c|$ be such that~$[y]$ is a maximal element
 with respect to~$\leq_\cV$. In other words, there is no $z \in |c|$ for which
 $[y] <_\cV [z]$. Since
 $$y \in |c| = |\Phi(c)| \subset \bigcup_{x \in |c|} |\Phi(x)|,$$
 there exists an $x \in |c|$ with $y \in |\Phi(x)|$. Then 
 part~(i)
 implies $[y] \leq_\cV [x]$, and in view of our choice
 of~$y$ this yields $[y] = [x]$. But then~(ii) furnishes $y = x$, and an
 application of~(iii) gives $\scalprod{x,\Phi x}=1$. Together with
 Proposition~\ref{prop:Phi-on-x} one finally obtains $y = x \in X^c$ as
 claimed.
\qed

\subsection{The stabilized combinatorial flow}

As mentioned earlier, the combinatorial flow describes on the level of
chains the possible future evolutions of solutions of the underlying
combinatorial vector field. Thus, in the case of combinatorial 
gradient vector fields iterations of~$\Phi$ should eventually
stabilize and provide intuition on potential connecting orbits between
critical cells. In other words, such a stabilization of the combinatorial
flow should encode both the Conley complex and the associated connection
matrix.

However, before we can show that such a stabilization exists, we first
have to study the behavior of iterates of the combinatorial flow. We begin
with a result which is an easy application of Proposition~\ref{prop:Phi-is-declining}
and shows that cells in~$|\Phi^n x|$ have to belong to vectors below~$[x]$
in the partial order~$\leq_\cV$, regardless of the iteration number $n \in \NN$.
\begin{prop}
\label{prop:Phi-infty-is-declining}
Assume that $\cV$ is a combinatorial gradient vector field on a regular
Lefschetz complex $X$. If $y\in |\Phi^n x|$ for some $n \in \NN$,
then the inequality $[y]\leq_\cV[x]$ is satisfied.
\end{prop}
\proof
We proceed by induction on~$n$. For $n=1$ the conclusion is the
statement of Proposition~\ref{prop:Phi-is-declining}. Thus, fix an
$n\in\NN$ and assume that the conclusion holds for  all
$k\in\{1,2,\ldots n\}$. Let $y\in |\Phi^{n+1} x|$ and let
$c:=\Phi^n x$. Then $c=\sum_{w\in|c|}\scalprod{c,w}w$
and $\Phi c=\sum_{w\in|c|}\scalprod{c,w}\Phi w$. Since
$$
y\in |\Phi^{n+1} x|= |\Phi c|\subset\bigcup_{w\in|c|}|\Phi w|,
$$
we get $y\in|\Phi w|$ for some $w\in |c|=|\Phi^n x|$. Thus, by our
induction assumption, we obtain $[y]\leq_\cV[w]\leq_\cV[x]$.
\qed

\medskip
Next, we study in more detail the specific form of the image
chain~$\Phi^n x$ if the cell~$x$ is a critical cell and $n \in \NN$.
This  result will prove to be essential for identifying the Conley
complex associated with~$\cV$.
\begin{prop}
\label{prop:bar-x}
Assume that~$\cV$ is a combinatorial gradient vector field on a regular
Lefschetz complex~$X$. Then for every critical cell $x\in X^c$ and 
every positive integer $n \in \NN$ we have the representation
\begin{equation} \label{eq:barx-form-1}
  \Phi^n x=x+r_{x,n}
\end{equation}
where~$r_{x,n}$ is a chain satisfying
\begin{equation} \label{eq:barx-form-2}
  |r_{x,n}| \subset X^+\cap |\,[x]^{<_\cV}\,|.
\end{equation}
\end{prop}
\proof
We proceed by induction on $n$. Assume first that one has $n=1$ and
set $r_{x,1}:=\Gamma\bdy x$. From Proposition~\ref{prop:Gamma-supp} we see
that $|r_{x,1}|\subset X^+$ and, since we assumed $x\in X^c$, one obtains
further $\scalprod{x,r_{x,1}}=0$. The inclusion $x\in X^c$ also shows that
$\Gamma x=0$, and therefore the identity $\Phi x=x+\Gamma\bdy x=x+r_{x,1}$
holds. To show that $|r_{x,1}| \subset |\,[x]^{<_\cV}\,|$, take an arbitrary
$y\in |r_{x,1}|$. Then we have $\scalprod{y,r_{x,1}}\neq 0$, and since
$\scalprod{x,r_{x,1}}=0$, we get $x\neq y$, that is, $\scalprod{y,x}=0$.
It follows that $\scalprod{y,\Phi x}=\scalprod{y,x}+\scalprod{y,r_{x,1}}=
\scalprod{y,r_{x,1}}\neq 0$. Hence, $y\in |\Phi x|$. Thus,
Proposition~\ref{prop:Phi-is-declining}(i) gives $[y]\leq_\cV[x]$ and
$[y]\in[x]^{\leq_\cV}$. Since $y\neq x$, and in view of the fact
that~$\Phi$ is a degree~$0$ map, one further has $\dim y=\dim x$,
and we also get $[y]\neq[x]$ and $[y]\in[x]^{<_\cV}$. This in turn yields
$y\in  |\,[x]^{<_\cV}\,|$. Therefore, $|r_{x,1}|\subset |\,[x]^{<_\cV}\,|$,
and the proof of \eqref{eq:barx-form-1} for $n=1$ is complete.

Next, fix a $k\in\NN$ and assume that \eqref{eq:barx-form-1} holds for all
$n \le k$ with~$r_{x,n}$ satisfying~\eqref{eq:barx-form-2}.
Then $\Phi^{k+1}x=\Phi (\Phi^k x)=\Phi x+\Phi r_{x,k}=x+r_{x,1}+\Phi r_{x,k}$,
and we set $r_{x,k+1}:=r_{x,1}+\Phi r_{x,k}$. It follows from the
induction assumption and Proposition~\ref{prop:Phi-on-X+} that
$|r_{x,k+1}|\subset |r_{x,1}|\cup |\Phi r_{x,k}|\subset X^+$.
Since $|r_{x,1}|\subset |\,[x]^{<_\cV}\,|$, it suffices to prove
that $|\Phi r_{x,k}|\subset |\,[x]^{<_\cV}\,|$ in order to see that
$|r_{x,k+1}|\subset  |\,[x]^{<_\cV}\,|$. For this, take an arbitrary
element $y\in |r_{x,k}|$. Then the induction assumption gives
$y\in  |\,[x]^{<_\cV}\,|$. Since $ |\,[x]^{<_\cV}\,|$ is $\cV$-compatible,
it follows that $[y]\in [x]^{<_\cV}$. Since $[x]^{<_\cV}$ is a down set,
we get $[y]^{\leq_\cV}\subset [x]^{<_\cV}$ and
$|\,[y]^{\leq_\cV}\,|\subset |\,[x]^{<_\cV}\,|$. Thus, we proved that
\begin{equation}
\label{eq:Phi-n-x-2}
    y\in |r_{x,k}| \quad\implies\quad |\,[y]^{\leq_\cV}\,|\subset |\,[x]^{<_\cV}\,|.
\end{equation}
We have $r_{x,k}=\sum_{y\in |r_{x,k}|}\scalprod{r_{x,k},y}y$
and $\Phi r_{x,k}=\sum_{y\in |r_{x,k}|}\scalprod{r_{x,k},y}\Phi y$.
In view of Proposition~\ref{prop:Phi-is-declining}(i) and
property~\eqref{eq:Phi-n-x-2} one then obtains the inclusion
$$
|\Phi r_{x,k}|\subset\bigcup_{y\in |r_{x,k}|}|\Phi y|\subset
\bigcup_{y\in |r_{x,k}|}|\,[y]^{\leq_\cV}|\,\subset |\,[x]^{<_\cV}\,|,
$$
which completes the induction argument, and therefore the proof.
\qed

\medskip
The statement of the above proposition can be illustrated using the
sample computations in Example~\ref{ex:formangradient-2}, two of
which have arguments which are critical cells. As shown in the
left panel of Figure~\ref{fig:phiex}, for the gradient combinatorial
vector field~$\cV_1$ one has the identity $\Phi_1(\bB\bD) = \bB\bD +
\bA\bB + \bD\bE$, i.e., in this situation $r_{\bB\bD,1} = \bA\bB +
\bD\bE$. Similarly, for the vector field~$\cV_2$ in the right panel
one has $\Phi_2(\bA\bB\bC) = \bA\bB\bC + \bB\bC\bD$, which leads
to $r_{\bA\bB\bC,1} = \bB\bC\bD$. In both cases, the
chain~$r_{x,1}$ clearly satisfies~\eqref{eq:barx-form-2}.

After these preparations we can finally show that iterations of
the combinatorial flow have to stabilize. In the context of the
Morse complex this was already shown by Forman in~\cite{Fo98a}.
The following result is modeled after~\cite[Theorem~7.2]{Fo98a},
yet adapted to our situation.
\begin{prop}
\label{prop:Phi-infty-exists}
Assume again that~$\cV$ denotes a combinatorial gradient vector field
on a regular Lefschetz complex~$X$, and consider the combinatorial
flow $\Phi = \Phi_\cV : C(X) \to C(X)$ defined in~\eqref{eq:Phi}.
In this case there exists an integer~$N \in \NN$ such that
\begin{displaymath}
  \Phi^n = \Phi^N
  \quad\mbox{ for all }\quad
  n \ge N.
\end{displaymath}
We then define $\Phi^\infty := \Phi^N$ and call it the
{\em stabilized combinatorial flow}.
\end{prop}
\proof
We show that for every $x \in X$ there exists an integer~$n_x \in \NN$
such that $\Phi^{n_x+1}x = \Phi^{n_x}x$ holds. Then it is straightforward
to see that choosing~$N$ as the maximum of~$\{ n_x \mid x \in X \}$ 
satisfies the statement of the proposition.

In order to prove the above claim, we proceed by induction over the
number of vectors which are strictly below~$[x]$ with respect
to~$<_\cV$, that is, by induction over the number
\begin{displaymath}
  L(x) := \#\{ V \in \cV \, \mid \, V <_\cV [x] \}.
\end{displaymath}
Suppose first that $x \in X$ is such that $L(x) = 0$. If we have $\Phi x = 0$,
then clearly $\Phi^2 x = \Phi x = 0$ and the statement follows. Otherwise there
exists a cell $y \in X$ with $y \in |\Phi x|$. But then
Proposition~\ref{prop:Phi-is-declining}(i) implies $[y] \leq_\cV [x]$, and the
identity $L(x) = 0$ further yields $[y] = [x]$. This in turn gives both $y = x$
and $x \in |\Phi x|$, as well as $\scalprod{x,\Phi x} = 1$, in view of
Proposition~\ref{prop:Phi-is-declining}(ii),(iii). In other words, we have
to have $\Phi x = x$, and the claim follows again.

Now let $k \in \NN$ be arbitrary and assume that we have verified the claim
for all $z \in X$ with $L(z) < k$. Moreover, let $x \in X$ be a cell with
$L(x) = k$. We now distinguish the cases of~$x$ being critical or not.

Assume first that~$x$ is not a critical cell. Then for all $y \in |\Phi x|$
one has to have $[y] <_\cV [x]$. To show this, suppose otherwise that there
is a $y \in |\Phi x|$ with $[y] = [x]$. Then
Proposition~\ref{prop:Phi-is-declining}(ii),(iii) implies $y = x$ and
$\scalprod{x,\Phi x} = 1$, and in view of Proposition~\ref{prop:Phi-on-x}
this gives $x \in X^c$, a contradiction. Together with
Proposition~\ref{prop:Phi-is-declining}(i), we therefore have $[y] <_\cV [x]$
for all $y \in |\Phi x|$ as claimed. Thus, every term in the representation
\begin{displaymath}
  \Phi x = \sum_{y \in |\Phi x|} \scalprod{y,\Phi x} y
\end{displaymath}
satisfies $L(y) < k$, and according to the inductive hypothesis there exists,
after possibly taking the maximum of the individual values, an $m \in \NN$
such that $\Phi^{m+1} y = \Phi^m y$ for all $y \in |\Phi x|$. This finally
furnishes
\begin{displaymath}
  \Phi^{m+2} x = \sum_{y \in |\Phi x|} \scalprod{y,\Phi x} \Phi^{m+1} y =
  \sum_{y \in |\Phi x|} \scalprod{y,\Phi x} \Phi^m y = \Phi^{m+1} x,
\end{displaymath}
and the claim from the beginning of the proof follows for non-critical~$x$.

Finally, suppose that~$x \in X^c$. According to Proposition~\ref{prop:bar-x}
there exists a chain~$r_{x,1}$ satisfying~\eqref{eq:barx-form-2} as well
as $\Phi x = x + r_{x,1}$. A straightforward induction argument then implies
\begin{displaymath}
  \Phi^\ell x =
  x + r_{x,1} + \Phi r_{x,1} + \ldots + \Phi^{\ell-1} r_{x,1}
  \quad\mbox{ for all }\quad
  \ell \in \NN .
\end{displaymath}
Thus, in order to establish that one has $\Phi^{m+1}x = \Phi^m x$ for some
$m \in \NN$ it suffices to show that $\Phi^m r_{x,1} = 0$. Notice that
according to~\eqref{eq:barx-form-2} we can apply the inductive hypothesis
to~$r_{x,1}$, and therefore there exists an integer~$m \in \NN$ such that
$\Phi^{m+1} r_{x,1} = \Phi^m r_{x,1}$. This immediately shows
that~$\Phi^m r_{x,1}$ is a fixed chain under the combinatorial
flow~$\Phi$. Furthermore, the inclusion in~\eqref{eq:barx-form-2} and
Proposition~\ref{prop:Phi-on-X+} yield $|\Phi^m r_{x,1}| \subset X^+$.
But then Proposition~\ref{prop:Phi-is-declining}(iv) implies that we have
to have $\Phi^m r_{x,1} = 0$, since the support of~$\Phi^m r_{x,1}$
contains no critical cells. This finally establishes the induction
step also for a critical cell~$x$, and the proof is complete.
\qed

\medskip
We would like to point out that in view of this result, the stabilized
combinatorial flow~$\Phi^\infty$ also satisfies the statements in
Propositions~\ref{prop:Phi-infty-is-declining} and~\ref{prop:bar-x},
as long as we replace~$n$ by~$\infty$ in their formulations.
\begin{figure}
  \begin{center}
    \includegraphics[width=0.45\textwidth]{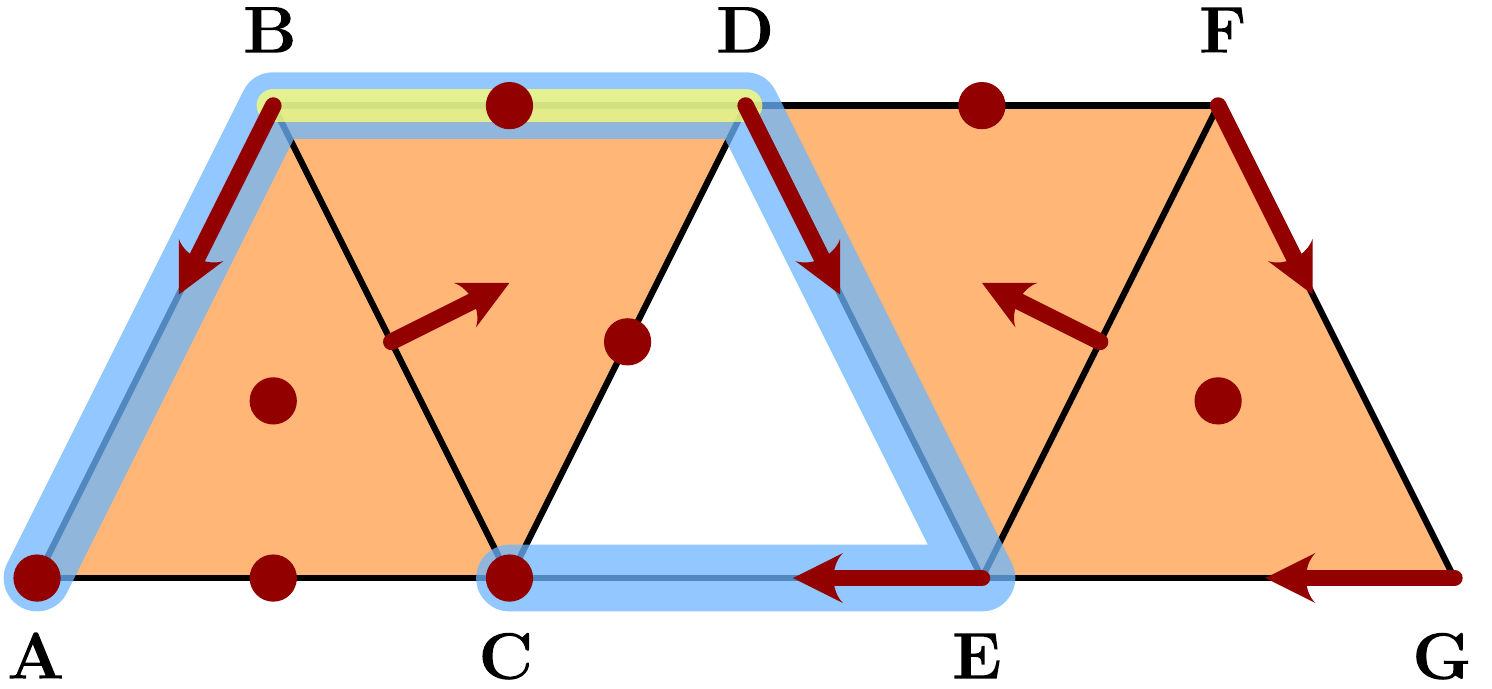}\quad\quad
    \includegraphics[width=0.45\textwidth]{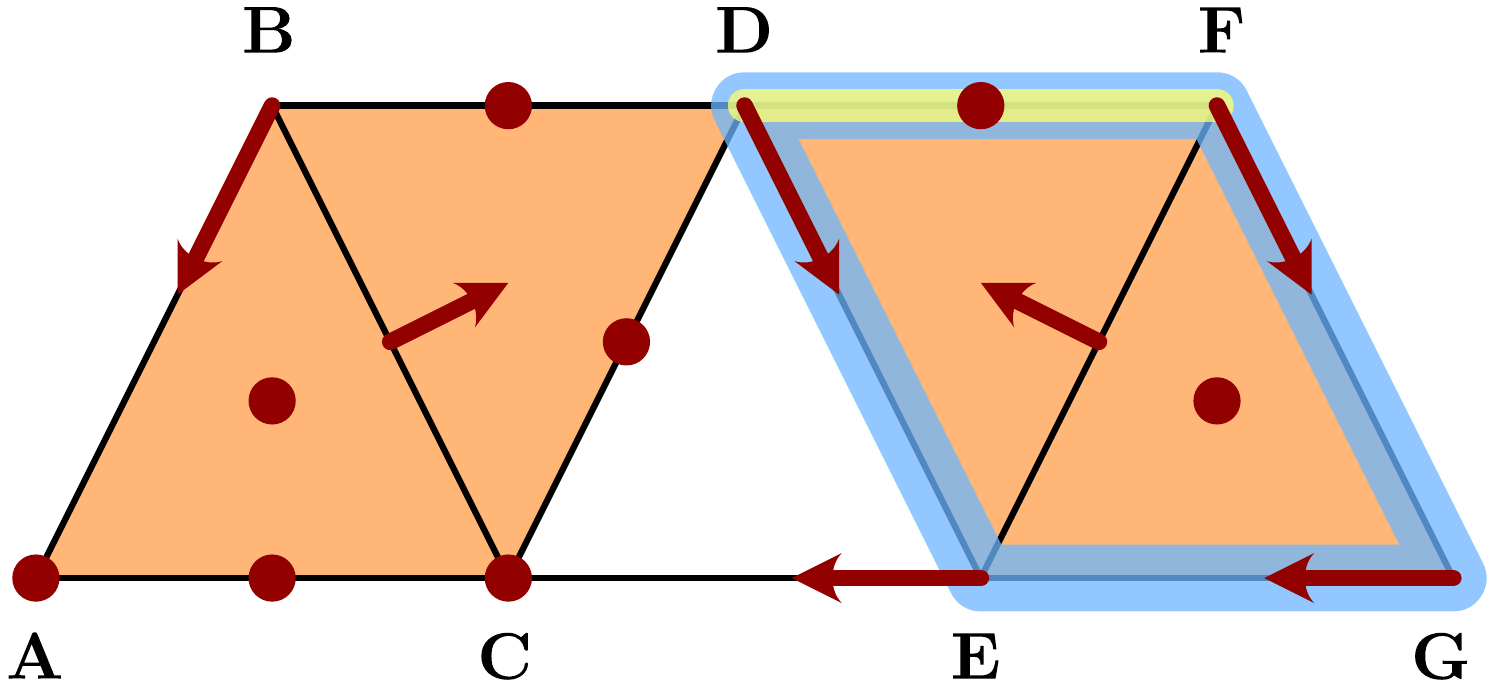} \\[3ex]
    \includegraphics[width=0.45\textwidth]{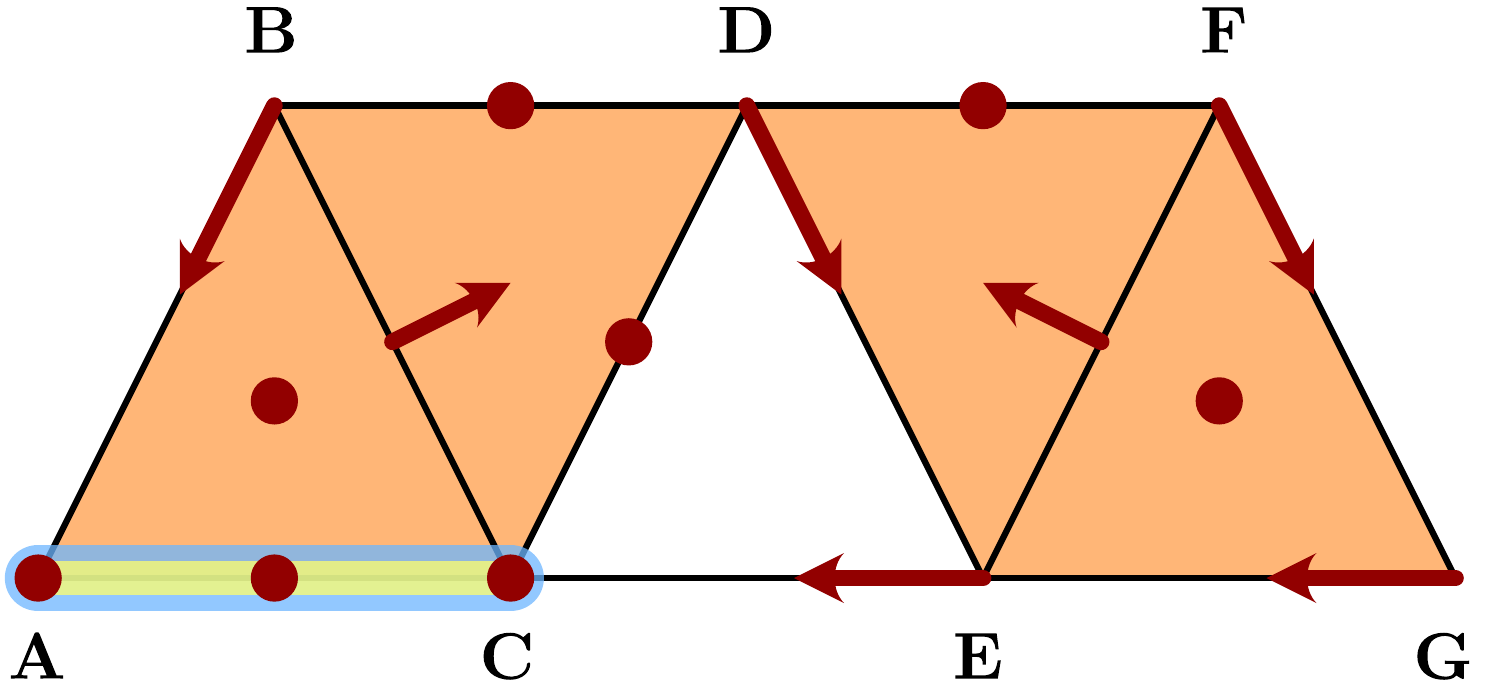}\quad\quad
    \includegraphics[width=0.45\textwidth]{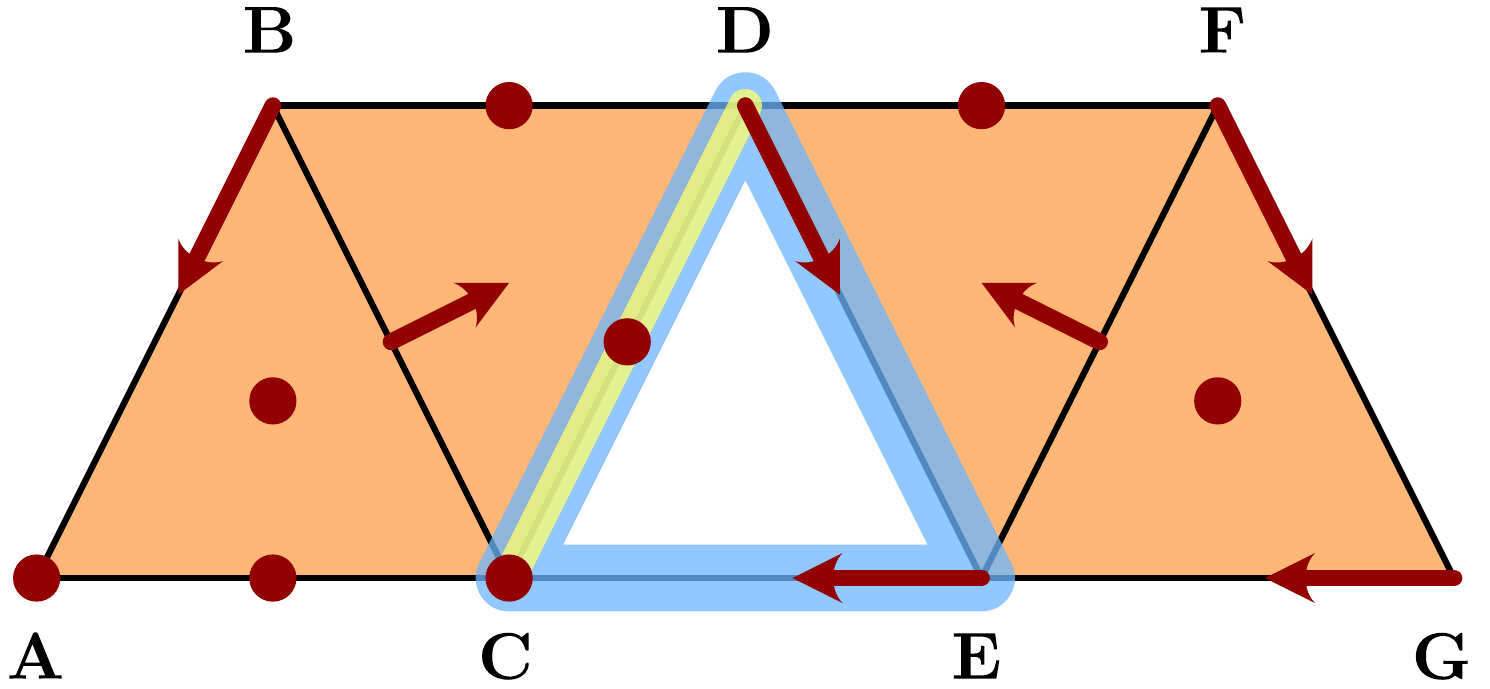}
  \end{center}
  \caption{{\em The stabilized combinatorial flow}.
           For the gradient combinatorial vector field~$\cV_1$
           introduced in the first row of Figure~\ref{fig:periodicex123},
           the four panels show the images~$\Phi_1^\infty(x)$ for the
           one-dimensional critical cells~$x$. From top left to bottom
           right the panels depict the chains~$\Phi_1^\infty(\bB\bD)$,
           $\Phi_1^\infty(\bD\bF)$, $\Phi_1^\infty(\bA\bC)$,
           and~$\Phi_1^\infty(\bC\bD)$, respectively, in blue.
           }
  \label{fig:phiinfex1}
\end{figure}
\begin{figure}
  \begin{center}
    \includegraphics[width=0.45\textwidth]{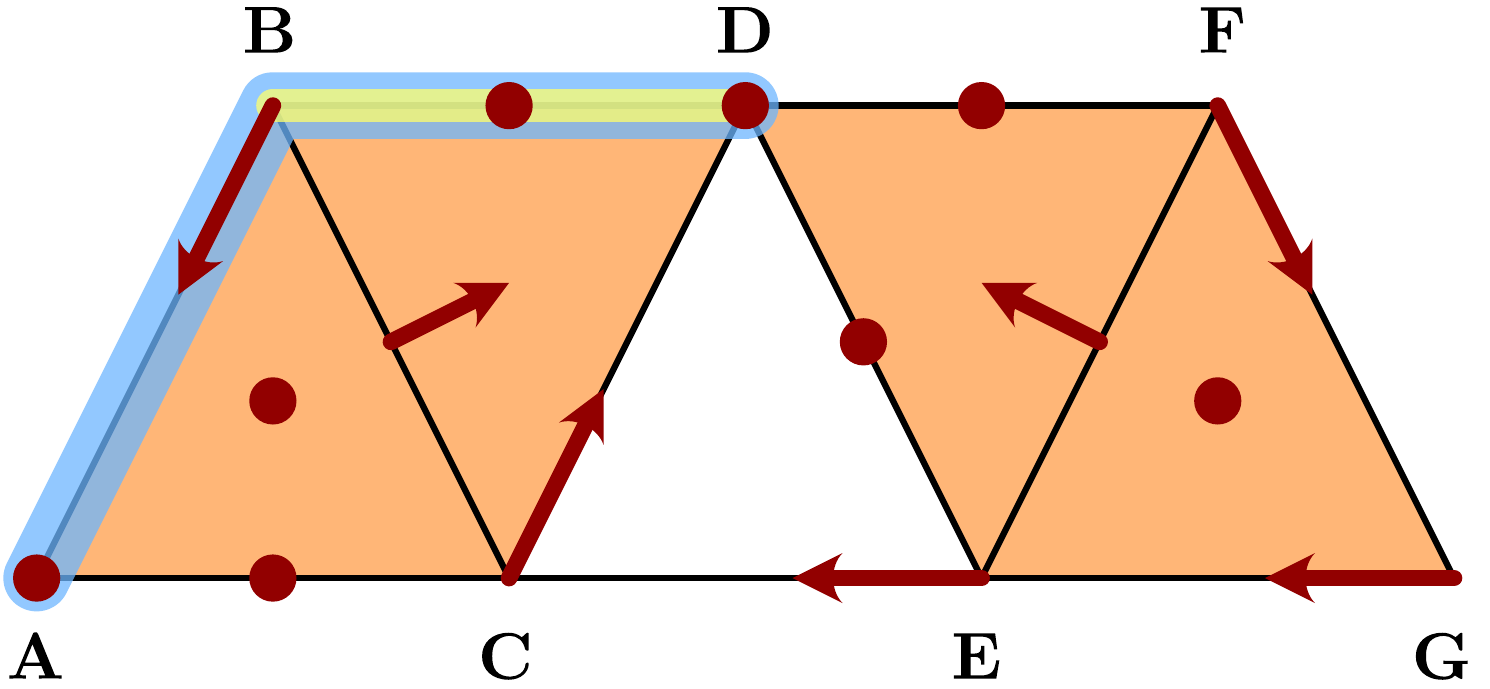}\quad\quad
    \includegraphics[width=0.45\textwidth]{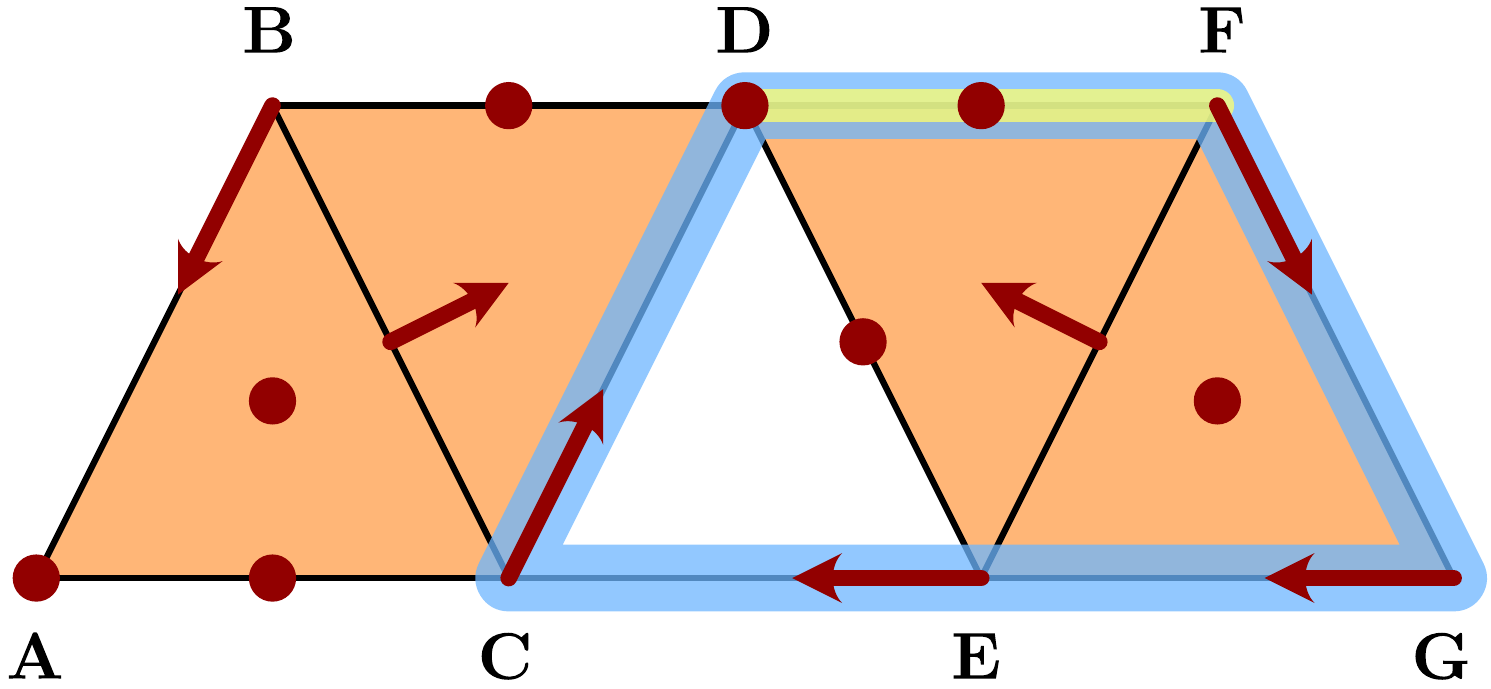} \\[3ex]
    \includegraphics[width=0.45\textwidth]{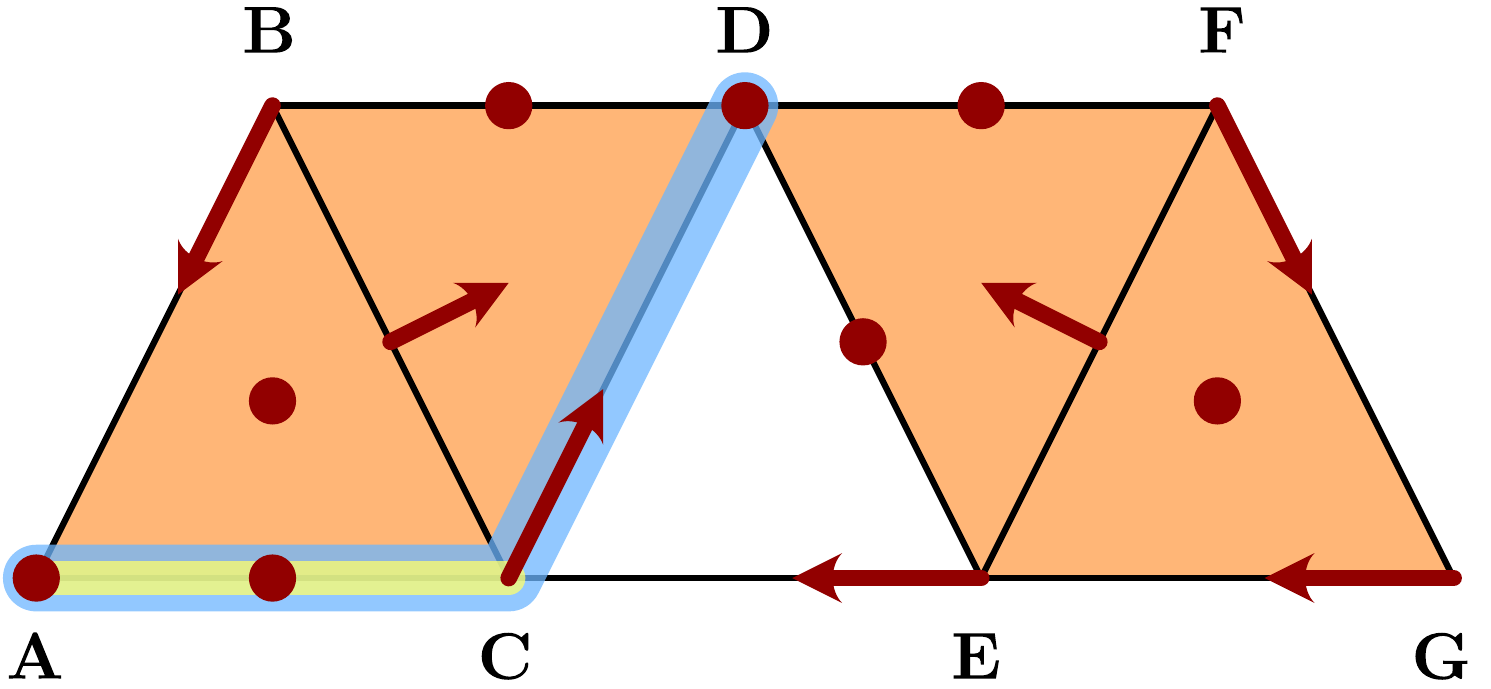}\quad\quad
    \includegraphics[width=0.45\textwidth]{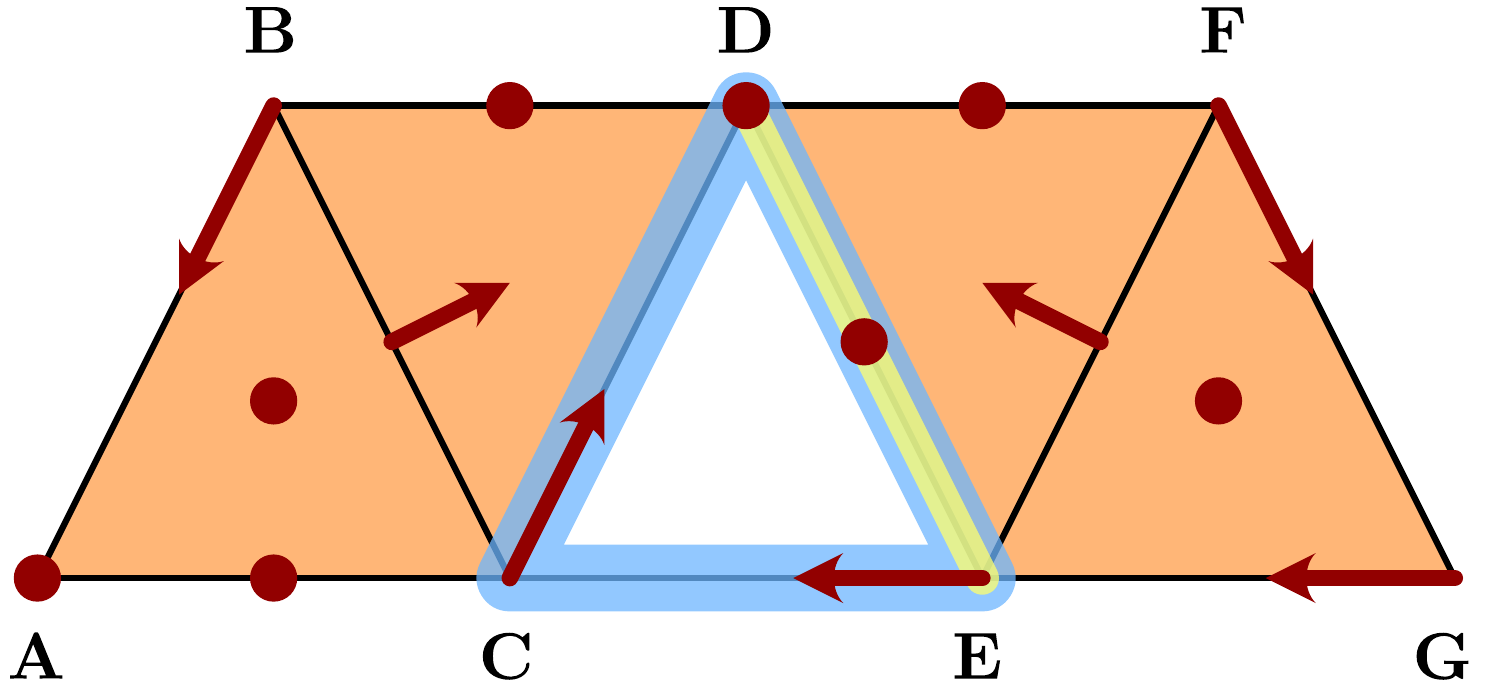}
  \end{center}
  \caption{{\em The stabilized combinatorial flow}.
           For the gradient vector field~$\cV_2$ introduced in the
           second row of Figure~\ref{fig:periodicex123}, the above
           four panels show the images~$\Phi_2^\infty(x)$ for the
           one-dimensional critical cells~$x$. From top left to bottom
           right they depict the chains~$\Phi_2^\infty(\bB\bD)$,
           $\Phi_2^\infty(\bD\bF)$, $\Phi_2^\infty(\bA\bC)$,
           and~$\Phi_2^\infty(\bD\bE)$, respectively, in blue.
           }
  \label{fig:phiinfex2}
\end{figure}
\begin{ex}[{\em Three Forman gradient vector fields}]
\label{ex:formangradient-3}
{\em
To illustrate the action of the stabilized flow, we return
to Examples~\ref{ex:formangradient-1} and~\ref{ex:formangradient-2}.
For the two gradient combinatorial vector fields~$\cV_1$ and~$\cV_2$
shown in the first two rows of Figure~\ref{fig:periodicex123},
and in Figure~\ref{fig:phiex}, what are the images of the critical
cells under the stabilized combinatorial flow?

We begin by considering the zero-dimensional critical cells in the
vector fields. Clearly, both their boundaries and images under~$\Gamma$
are trivial, and therefore we have $\Phi_k^\infty(x) = x$ for all of
these cells. In addition, for the two-dimensional cells one easily
obtains that~$\Phi_k^\infty(\bA\bB\bC) = \bA\bB\bC + \bB\bC\bD$,
as well as~$\Phi_k^\infty(\bE\bF\bG) = \bE\bF\bG + \bD\bE\bF$.
This follows easily from the computation in 
Example~\ref{ex:formangradient-2}, together with the fact that
$\Phi_k(\bB\bC\bD) = \Phi_k(\bD\bE\bF) = 0$, which in turn
implies that the second iterate of~$\Phi_k$ coincides with the
first iterate on critical cells of dimension two.

This leaves us with the one-dimensional critical cells. As it turns
out, their images under the stabilized flow can be different with
respect to~$\cV_1$ or~$\cV_2$. The resulting images~$\Phi_k^\infty(x)$
are shown in Figure~\ref{fig:phiinfex1} for $k = 1$, and in
Figure~\ref{fig:phiinfex2} for $k = 2$. Notice that in almost
all of these cases, the image of a one-dimensional cell~$x$ under
the stabilized flow is given by the support of all connecting orbits
between the cell~$x$ and zero-dimensional critical cells. The only 
exception is the image~$\Phi_1^\infty(\bD\bF)$, which lacks the
edge~$\bC\bE$. Yet, since this edge is contained in both connecting
orbits between~$\bD\bF$ and~$\bC$, the sum of the associated two
chains with respect to $\ZZ_2$-coefficients leads to its cancellation.
\exend
}
\end{ex}

The stabilized combinatorial flow is our main tool for the
construction of the Conley complex and its associated connection
matrix for a combinatorial gradient vector field~$\cV$ on a regular
Lefschetz complex~$X$. It follows immediately from
Proposition~\ref{prop:Phi-infty-exists} that for every chain
$c \in C(X)$ we have
\begin{displaymath}
  \Phi(\Phi^\infty(c)) = \Phi(\Phi^n(c)) = \Phi^{n+1}(c) =
  \Phi^\infty(c),
\end{displaymath}
that is, the image~$\Phi^\infty(c)$ is a fixed point of the combinatorial
flow~$\Phi$. We therefore consider the set
\begin{equation} \label{def:FixPhi}
    \Fix\Phi:=\setof{c\in C(X)\mid \Phi c = c}
\end{equation}
consisting of chains fixed by $\Phi$. Clearly, if $c\in\Fix\Phi$, then
$\bdy c\in\Fix\Phi$, because $\Phi\bdy c =\bdy\Phi c=\bdy c$. It follows
that $(\Fix\Phi,\bdy_{|\Fix\Phi})$ is a chain subcomplex of $(C(X),\bdy)$.
Furthermore, since~$\Phi$ is a chain map, Proposition~\ref{prop:Phi-infty-exists}
further implies that $\Phi^\infty: (C(X),\bdy)\to (\Fix\Phi,\bdy_{|\Fix\Phi})$
is a chain map as well. Yet even more is true, as the following result 
demonstrates.
\begin{prop}
\label{prop:Phi-infty-Xc-Fix}
Assume that~$\cV$ is a combinatorial gradient vector field on a regular
Lefschetz complex~$X$, and let~$\Phi^\infty$ denote the associated
stabilized combinatorial flow. Then the restriction $\Phi^\infty_{|C(X^c)}:
C(X^c) \to \Fix\Phi$ is an isomorphism. Moreover, if we define the
projection $\Pi:\Fix\Phi \to C(X^c)$ via
\begin{equation} \label{prop:Phi-infty-Xc-Fix-1}
  \Pi c := \sum_{x\in X^c} \scalprod{c,x}x
  \quad\mbox{ for every }\quad
  c \in \Fix\Phi,
\end{equation}
then we have $\Phi^\infty \Pi c = c$ for all $c \in \Fix\Phi$.
\end{prop}
\proof
We begin by verifying that $\Phi^\infty_{|C(X^c)}$ is a monomorphism.
For this, suppose that $c\in C(X^c)$ is a chain with~$\Phi^\infty c=0$.
Then $c=\sum_{x\in X^c}a_xx$, and the representation in~\eqref{eq:barx-form-1}
further yields
\[
  0=\Phi^\infty c=\sum_{x\in X^c}a_x \Phi^\infty x=
  c+\sum_{x\in X^c}a_xr_x\in C(X^c)\oplus C(X^+).
\]
Hence, it follows that $c=0$, and therefore~$\Phi^\infty_{|C(X^c)}$
is one-to-one.

In order to prove that $\Phi^\infty_{|C(X^c)}: C(X^c) \to \Fix\Phi$ is
onto, let $c \in \Fix\Phi$ be arbitrary and consider the chain~$\Pi c
\in C(X^c)$ defined in~\eqref{prop:Phi-infty-Xc-Fix-1}. Then we have
\begin{displaymath}
  \Phi^\infty \Pi c = \sum_{x \in X^c} \scalprod{c,x} \Phi^\infty x,
\end{displaymath}
and for every critical cell $y \in X^c$ the representation of~$\Phi^\infty x$
in Proposition~\ref{prop:bar-x} further implies
$\scalprod{\Phi^\infty x,y} = \scalprod{x + r_x,y} = \scalprod{x,y}$,
see also~\eqref{eq:barx-form-1}. Thus one obtains
\begin{displaymath}
  \scalprod{\Phi^\infty \Pi c,y} =
  \sum_{x \in X^c} \scalprod{c,x} \scalprod{x,y} =
  \scalprod{c,y}
  \quad\Rightarrow\quad
  \scalprod{c - \Phi^\infty \Pi c,y} = 0
\end{displaymath}
for all $y \in X^c$. In other words, the support of the chain
$c - \Phi^\infty \Pi c$ does not contain a critical cell. In view
of Proposition~\ref{prop:Phi-is-declining}(iv) and $c - \Phi^\infty \Pi c
\in \Fix\Phi$ this implies $c = \Phi^\infty \Pi c$, and thereby completes
the proof of the result.
\qed

\medskip
The above result shows that the fixed chains under the combinatorial
flow~$\Phi$ are in one-to-one correspondence with chains of critical
cells in the Lefschetz complex~$X$. As the following result demonstrates,
this fact can be used to both find a suitable basis of~$\Fix\Phi$
which is indexed by the critical cells of~$X$, as well as equip this
basis with the structure of a Lefschetz complex.
\begin{prop}
\label{prop:barX}
 Assume that~$\cV$ is a combinatorial gradient vector field on a regular
 Lefschetz complex~$X$, and recall that $X^c\subset X$ denotes the set of 
 all critical cells of~$\cV$. For every $x\in X^c$ set
 \begin{displaymath}
   \bar{x}:=\Phi^\infty x \in\Fix\Phi.
 \end{displaymath}
Then the following statements hold:
 \begin{itemize}
 \item[(i)]   The set $\bar{X}:=\setof{\bar{x}\mid x\in X^c}$
              is a basis of $\Fix\Phi$.
 \item[(ii)]  For every $x\in X^c$ we have the representation
              \begin{equation} \label{eq:barX}
                \bdy\bar{x}=\sum_{z\in X^c}a_{xz}\bar{z},
              \end{equation}
              for uniquely determined coefficients $a_{xz}\in R$.
 \item[(iii)] The pair~$(\bar{X},\bar{\kappa})$, where the
              $\ZZ$-gradation on~$\bar{X}$ is induced by the map
              $\dim: \bar{X}\ni \bar{x}\mapsto\dim x\in\ZZ$, and
              where $\bar{\kappa}(\bar{x},\bar{z})$ is the
              coefficient~$a_{xz}$ in~\eqref{eq:barX}, is a
              Lefschetz complex with $C(\bar{X})=\Fix\Phi$.
 \end{itemize}
\end{prop}
\proof
In order to prove~(i), consider the projection $\Pi:\Fix\Phi\to C(X^c)$
defined in Proposition~\ref{prop:Phi-infty-Xc-Fix}, which has been shown
to satisfy $\Phi^\infty \Pi=\id_{\Fix\Phi}$. In addition, let $c\in \Fix\Phi$
be arbitrary. Then one has
\[
   c=\Phi^\infty\Pi c=
   \Phi^\infty \sum_{x\in X^c}\scalprod{c,x}x =
   \sum_{x\in X^c}\scalprod{c,x}\Phi^\infty x=
   \sum_{x\in X^c}\scalprod{c,x}\bar{x},
\]
which proves that~$\bar{X}$ generates~$\Fix\Phi$. To see that~$\bar{X}$ is
linearly independent, assume that
\[
    \sum_{x\in X^c} \alpha_x \bar{x}=0
\]
for some coefficients $\alpha_x\in R$. Then $0=\sum_{x\in X^c} \alpha_x
\Phi^\infty x = \Phi^\infty \sum_{x\in X^c} \alpha_x  x$. Since
$\sum_{x\in X^c} \alpha_x  x\in C(X^c)$, and $\Phi_{|C(X^c)}:C(X^c) \to
\Fix\Phi$ is an isomorphism due to Proposition~\ref{prop:Phi-infty-Xc-Fix},
we conclude that $\sum_{x\in X^c} \alpha_x  x=0$. This in turn implies
$\alpha_x=0$ for all $x\in X^c$. Thus, the set~$\bar{X}$ is indeed a basis
for~$\Fix\Phi$.

In order to see (ii), we observe that $\bdy \bar{x}\in \Fix\Phi$, since~$\Phi$
is a chain map. Thus, the representation in~\eqref{eq:barX} is an immediate
consequence of~(i). Finally, a direct application of
Proposition~\ref{prop:U-Lefschetz} shows that the pair~$(\bar{X},\bar{\kappa})$
is a Lefschetz complex, and clearly $C(\bar{X})=\Fix\Phi$, which proves~(iii).
\qed

\medskip
The Lefschetz complex~$(\bar{X},\bar{\kappa})$ will take the role of the
Conley complex associated with the combinatorial gradient vector field~$\cV$.
For this, and as a last step, it is necessary to identify a natural partial
order on its cells. This can be accomplished as follows.

Recall that~$\cC$ denotes the collection of all critical cells of~$\cV$.
Then it follows easily from Proposition~\ref{prop:bar-x} that the map
$X^c\ni x\mapsto \bar{x}\in\bar{X}$ is a bijection. Hence, also the map
$\cC\ni \{x\}\mapsto \bar{x}\in\bar{X}$ is a bijection, and this enables 
us to carry over the partial order~$\leq_\cV$ from $\cC \subset \cV$
to~$\bar{X}$. Thus, for arbitrary elements $\bar{x},\bar{y}\in\bar{X}$
we write $\bar{x}\leq_\cV\bar{y}$ if $\{x\}\leq_\cV\{y\}$. Then we have
the following result.
\begin{prop} \label{prop:P-eq-cV-on-barX}
Assume that~$\cV$ is a combinatorial gradient vector field on a regular
Lefschetz complex~$X$. Then the above-defined partial order~$\leq_\cV$
in~$\bar{X}$ is a natural partial order on the Lefschetz
complex~$(\bar{X},\bar{\kappa})$.
\end{prop}
\proof
We begin by proving that the order $\leq_\cV$ in~$\bar{X}$ is admissible,
that is, the validity of $\bar{y}\leq_{\bar{\kappa}}\bar{x}$ for $x,y\in X^c$
implies $\bar{y}\leq_\cV\bar{x}$. Since the partial order~$\leq_{\bar{\kappa}}$
is the transitive closure of $\adhl_{\bar{\kappa}}$, it clearly suffices to
prove that $\bar{y}\adhl_{\bar{\kappa}}\bar{x}$ implies $\bar{y}\leq_\cV\bar{x}$.
Thus, assume that $x,y\in X^c$ and $\bar{y}\adhl_{\bar{\kappa}}\bar{x}$.
Then, $\bar{\kappa}(\bar{x},\bar{y})\neq 0$. By Proposition~\ref{prop:barX}
we have $\bdy\bar{x}=\sum_{z\in X^c}a_{xz}\bar{z}$, with $a_{xz} =
\bar{\kappa}(\bar{x},\bar{z})$. In addition, Propositions~\ref{prop:bar-x}
and~\ref{prop:Phi-infty-exists} imply $\bar{z}=\Phi^\infty z=z+r_{z,\infty}$,
together with the inclusion $|r_{z,\infty}|\subset X^+\cap |\,[z]^{<_\cV}\,|$.
Since $y\in X^c\subset X\setminus X^+$, one further has
\begin{multline*}
    \scalprod{\bdy\bar{x},y}=\sum_{z\in X^c}a_{xz}\scalprod{\bar{z},y}=
    \sum_{z\in X^c}a_{xz}\left(\scalprod{z,y}+\scalprod{r_{z,\infty},y}\right)\\
    = \sum_{z\in X^c}a_{xz}\scalprod{z,y}=a_{xy}=\bar{\kappa}(\bar{x},\bar{y})\neq 0.
\end{multline*}
It follows that
$$
y\in|\bdy\bar{x}|=|\bdy\Phi^\infty x|=|\Phi^\infty\bdy x|=
\left|\Phi^\infty \sum_{w\in |\bdy x|}\kappa(x,w) w \right|\subset
\bigcup_{w\in |\bdy x|}|\Phi^\infty w|.
$$
Thus, $y\in |\Phi^\infty u|$ for some $u\in X$ such that $\kappa(x,u)\neq 0$.
From Proposition~\ref{prop:Phi-infty-is-declining} we get $[y]\leq_\cV[u]$.
Since $\kappa(x,u)\neq 0$ we further have $u\in\cl x$, which in turn implies
$u\in[u]\cap\cl[x]$. In consequence $[u]\leq_\cV[x]$ and $\{y\}=[y]\leq_\cV[x]=\{x\}$.
Thus, by the definition of $\leq_\cV$ in $\bar{X}$ we get $\bar{y}\leq_\cV\bar{x}$.

To see that $\leq_\cV$ is natural, we need to verify the implication
in~\eqref{def-natural-order}. Thus, let $\bar{x},\bar{y}\in\bar{X}$ satisfy
both the inequality $\bar{x}\leq_\cV\bar{y}$ and $\dim\bar{x}=\dim\bar{y}$.
Then both~$\{x\}\leq_\cV\{y\}$ and $\dim x=\dim y$ hold as well. Yet this
implies that one cannot have $\{x\}<_\cV\{y\}$, because in this case, in
view of Proposition~\ref{prop:declining-dim}(ii), one deduces the strict
inequality $\dim x=\dim\{x\}<\dim\{y\}=\dim y$.  It follows that $\{x\}=\{y\}$
and $\bar{x}=\bar{y}$, which proves  that the order is natural.
\qed

\medskip
In the remaining last subsection of the paper, we consider~$\bar{X}$ as a
poset ordered by the natural partial order~$\leq_\cV$ in $\cV$, which has
been transferred to~$\bar{X}$ via the map $\{ x \} \mapsto \bar{x}$. In 
addition, the triple $(\bar{X},C(\bar{X}),\bdy^{\bar{\kappa}})$ is considered
as a natural filtration of~$\bar{X}$.

\subsection{Conley complex and unique connection matrix}

Over the previous four subsections, we have considered combinatorial
gradient vector fields~$\cV$ on a regular Lefschetz complex~$X$.
Through the use of the associated combinatorial flow~$\Phi$, this
study culminated in the creation of a new Lefschetz complex~$\bar{X}$,
together with its associated singleton partition. In fact, we have
introduced~$(\bar{X},C(\bar{X}),\bdy^{\bar{\kappa}})$ as a natural
filtration of~$\bar{X}$. We will now show that this filtration
is the Conley complex associated with~$\cV$, and that it has a 
uniquely determined connection matrix. First, however, we need an
auxiliary result, which allows us to recognize the combinatorial
flow~$\Phi$ and its iterates as filtered morphisms which are
filtered chain homotopic to the identity.
\begin{prop}
\label{prop:Phi-id-filt-homotopy}
Assume that $\cV$ is a gradient vector field on a regular Lefschetz
complex~$X$. Then we have a well-defined filtered morphism
\[
   (\id_\cV,\Phi):( \cV,C(X),\bdy^\kappa)\to( \cV,C(X),\bdy^\kappa).
\]
Moreover, $(\id_\cV,\Phi)^n=(\id_\cV,\Phi^n)$ is filtered chain
homotopic to the identity morphism $\id_{(\cV,C(X))}$ for every $n\in\NN$.
In particular, $(\id_\cV,\Phi^\infty)$ is filtered chain homotopic
to~$\id_{(\cV,C(X))}$.
\end{prop}
\proof
Clearly, $\Phi$ is a chain map and $\id_\cV:(\cV,\cV_\star)\to (\cV,\cV_\star)$
is a morphism in $\DPSet$. Hence, to prove that $(\id_\cV,\Phi)$ is a
well-defined filtered morphism we only have to show that~$\Phi$ is
$\id_\cV$-filtered. We will do so by checking
property \eqref{eq:id-filt-hom-downset} of Corollary~\ref{cor:filt-hom}.
Consider the down set $\cL\in\Down(\cV)$ and set $L:=|\cL|$.
Then $C(X)_\cL=\bigoplus_{V\in\cL}C(V)=C(L)$, and we need to verify
that $\Phi(C(L))\subset C(L)$. But this follows from
Proposition~\ref{prop:Phi-on-down sets}, because~$L$ is
immediately seen to be $\cV$-compatible and closed by
Proposition~\ref{prop:convex-lcl}(i).

Since by Proposition~\ref{prop:filt-homotopy-composition}
filtered chain homotopy between morphisms is preserved by
composition, in order to prove that $(\id_\cV,\Phi)^n=(\id_\cV,\Phi^n)$
is filtered homotopic to the identity morphism it suffices to prove
that~$(\id_\cV,\Phi)$ is filtered homotopic to the identity morphism.
By the definition of $\Phi$ we have $\Phi-\id_{C(X)}=\bdy\Gamma+\Gamma\bdy$.
Thus, we only need to check that $\Gamma$ is $\id_\cV$-filtered.
Hence, assume that $\Gamma_{VW}\neq 0$ for some $V,W\in \cV$.
Then, there exists a $w\in W$ such that $\pi_V(\Gamma w)\neq 0$.
In particular, one has $\Gamma w\neq 0$, which implies both $w=w^-\neq w^+$
and $\Gamma w=-\kappa(w^+,w^-)^{-1}w^+$. It follows that $w^+\in V$.
Thus $V=W=\id_\cV(W)$, which proves that $\Gamma$ is $\id_\cV$-filtered,
in fact, even graded. Finally, since $\Phi^\infty=\Phi^n$ for large $n\in\NN$,
it follows that also $(\id_\cV,\Phi^\infty)$ is filtered chain homotopic
to $\id_{(\cV,C(X))}$.
\qed

\medskip
The next two results show that~$(\bar{X},C(\bar{X}),\bdy^{\bar{\kappa}})$ is
indeed a Conley complex of the combinatorial gradient vector field~$\cV$. First,
we construct a filtered morphism from~$( \cV,C(X),\bdy^\kappa)$ to the Conley
complex.
\begin{prop}
\label{prop:phi-morphism}
In the situation of the last proposition, we have a well-defined filtered morphism
\begin{equation*}
    (\alpha,\phi):( \cV,C(X),\bdy^\kappa)\to(\bar{X},C(\bar{X}),\bdy^{\bar{\kappa}})
\end{equation*}
with $\alpha:\bar{X}\to\cV$ defined for $x\in X^c$ by $\alpha(\bar{x}):=\{x\}$
and $\phi:C(X) \to C(\bar{X})$ defined for $c\in C(X)$ by $\phi(c):=\Phi^\infty(c)$.
\end{prop}
\proof
First observe that by Proposition~\ref{prop:declining-dim}(iii) we
can choose $\cV_\star=\cC$, and by Theorem~\ref{thm:singleton-partition}(i)
we then get $\bar{X}_\star=\bar{X}$. Thus, the poset filtered chain complex
$(\bar{X},C(\bar{X}),\bdy^{\bar{\kappa}})$ is peeled.
Furthermore, $\alpha$ is a morphism
in~$\DSet$, and since we consider~$\bar{X}$ as ordered by the natural
order~$\leq_\cV$, the map~$\alpha$ is trivially order preserving,
hence also a morphism in~$\DPSet$. Obviously, the map~$\phi$ is a
chain map. We will show that~$\phi$ is $\alpha$-filtered by checking
property \eqref{eq:filt-hom-downset} of Proposition~\ref{prop:filt-hom}.
Let $\cL\in\Down(\cV)$ and let $L:=|\cL|$. Then we have
$C(X)_\cL=\bigoplus_{V\in\cL}C(V)=C(L)$. Observe that~$\alpha^{-1}(\cL)$
is a down set in~$\bar{X}$. Indeed, if $\bar{x}\in\alpha^{-1}(\cL)$ and
$\bar{y}\leq_\cV \bar{x}$, then $\{y\}\leq_\cV\{x\}=\alpha(\bar{x})\in\cL$,
which implies $\alpha(\bar{y})=\{y\}\in\cL$ and $\bar{y}\in\alpha^{-1}(\cL)$.
Hence, $\alpha^{-1}(\cL)^{\leq}=\alpha^{-1}(\cL)$. It follows that
$C(\bar{X})_{\alpha^{-1}(\cL)^\leq}=C(\bar{X})_{\alpha^{-1}(\cL)}=
\bigoplus_{\{u\}\in\cL}R\bar{u}$. Thus, in our case, we have to verify
condition~\eqref{eq:filt-hom-downset} in the form
\begin{equation}
\label{eq:phi-morphism-1}
   \Phi^\infty(C(L))\subset \bigoplus_{\{u\}\in\cL}R\bar{u}.
\end{equation}
In order to verify \eqref{eq:phi-morphism-1} take a $c\in C(L)$,
and let $\bar{c} := \Phi(c) = \Phi^\infty(c)$. In view of
Proposition~\ref{prop:bar-x} one has
\begin{equation}
\label{eq:phi-morphism-2}
  \bar{c} = \sum_{u\in X^c}a_u\bar{u} = 
  \sum_{u\in X^c}a_u u+\sum_{u\in X^c}a_u r_{u,\infty}
\end{equation}
for some coefficients $a_u\in R$ and chains~$r_{u,\infty}$
satisfying $|r_{u,\infty}| \cap X^c=\emptyset$, as well as
the inclusion
\begin{equation}
\label{eq:phi-morphism-3}
|r_{u,\infty}|\subset|[u]^{<_\cV}|.
\end{equation}
Consider a $u\in X^c$ such that $\{u\}\not\in\cL$, i.e., one
has~$u\not\in L$.
Since by Proposition~\ref{prop:Phi-on-down sets} we have
$\phi(C(L))=\Phi^\infty(C(L))\subset C(L)$, it then follows
that the identity $\scalprod{\bar{c},u}=0$ holds. But now
equation~\eqref{eq:phi-morphism-2} gives $\scalprod{\bar{c},u}=a_u$,
because~\eqref{eq:phi-morphism-3} implies $\scalprod{r_{u,\infty},u}=0$.
Therefore, we have $a_u=0$ for $u\not\in L$, which subsequently proves
\[
    \bar{c}=\sum_{u\in L}a_u\bar{u}\in \bigoplus_{u\in L}R\bar{u}.
\]
This verifies~\eqref{eq:phi-morphism-1} and completes the proof.
\qed

\medskip
In our next result, we introduce a candidate for the inverse morphism
in~$\PfCC$ to the morphism~$(\alpha,\phi)$ from the last proposition.
\begin{prop}
\label{prop:psi-morphism}
In the situation of the last proposition, we have a well-defined filtered morphism
\begin{equation*}
    (\beta,\psi):(\bar{X},C(\bar{X}),\bdy^{\bar{\kappa}})\to( \cV,C(X),\bdy^\kappa)
\end{equation*}
with $\beta:\cV\pto\bar{X}$ defined for $x\in X^c$ by $\beta(\{x\}):=\bar{x}$
and $\psi:C(\bar{X}) \to C(X)$ defined for $c\in C(X)$ by $\psi(c):=c$.
\end{prop}
\proof
As we already mentioned, Proposition~\ref{prop:declining-dim} enables us to
choose the set $\cV_\star=\cC$, and Theorem~\ref{thm:singleton-partition}(i)
then gives $\bar{X}_\star=\bar{X}$. Therefore, $\beta$ is a morphism in~$\DSet$.
Moreover, since we consider~$\bar{X}$ as ordered by the natural order~$\leq_\cV$,
the map~$\beta$ is trivially order preserving, hence also a morphism in~$\DPSet$.
Obviously, the map~$\psi$ is a chain map. We will show that~$\psi$ is
$\beta$-filtered by checking again property~\eqref{eq:filt-hom-downset}
of Proposition~\ref{prop:filt-hom}. Consider a down set $A\in\Down(\bar{X})$.
Clearly, we have the identity $C(\bar{X})_A=C(A)$. If we define $\cL:=\beta^{-1}(A)$,
then in general~$\cL$ is not a down set in~$\cV$, but $\cL^{\leq_\cV}$ is.
We therefore consider $L:=|\cL^{\leq_\cV}|$, and get
$C(X)_{\beta^{-1}(A)^{\leq_\cV}}=C(X)_{\cL^{\leq_\cV}}=C(L)$.
Thus, in our case, condition \eqref{eq:filt-hom-downset} has to
be verified in the form
\begin{equation}
\label{eq:psi-morphism-1}
   \psi(C(A))=C(A)\subset C(L).
\end{equation}
Consider first a chain $\bar{x}\in A$ for a critical $x\in X^c$. Then
$\beta(\{x\})=\bar{x}\in  A$, which means that $\{x\}\in\beta^{-1}(A)=\cL$
and $|x|=\{x\}=[x]\in \cL\subset\cL^{\leq_\cV}$. It follows from
Proposition~\ref{prop:bar-x} that we have $\bar{x}=\Phi^\infty x=x+r_{x,\infty}$,
where the inclusion $|r_{x,\infty}|\subset X^+\cap |\,[x]^{<_\cV}\,|$ holds.
Hence, one obtains that
\begin{displaymath}
  |\bar{x}| \subset |x|\cup|r_{x,\infty}| \subset
  [x]\cup|[x]^{<_\cV}| = |\,[x]^{\leq_\cV}\,| \subset
  |\cL^{\leq_\cV}|=L
\end{displaymath}
is satisfied, because
$[x]\in \cL \subset \cL^{\leq_\cV}$ and~$\cL^{\leq_\cV}$ is a down set.
This immediately shows that $\bar{x}\in C(L)$. Since every chain in~$C(A)$
is a linear combination of chains~$\bar{x}$ in~$A$ with~$x\in X^c$, the
inclusion~\eqref{eq:psi-morphism-1} follows.
\qed

\medskip
The next result combines Propositions~\ref{prop:phi-morphism}
and~\ref{prop:psi-morphism} to show that the filtered chain complexes
$(\cV,C(X),\bdy^\kappa)$ and $(\bar{X},C(\bar{X}),\bdy^{\bar{\kappa}})$
are indeed elementary filtered chain homotopic.
\begin{thm}
\label{thm:forman-gradient-filtered-equivalence}
  Assume that $\cV$ is a gradient vector field on a regular Lefschetz complex $X$.
  Then the filtered morphisms $(\alpha,\phi)$ and $(\beta,\psi)$ defined in
  the last two propositions are mutually inverse elementary filtered chain equivalences.
  In particular, the poset filtered chain complexes~$( \cV,C(X),\bdy^\kappa)$
  and~$(\bar{X},C(\bar{X}),\bdy^{\bar{\kappa}})$ are elementary filtered chain homotopic,
  and the latter is a representation of the former.
\end{thm}
\proof
Clearly, one has the identity $\beta\alpha=\id_{\bar{X}}$. Moreover, if
$c\in \Fix\Phi$, then we obtain $\Phi^\infty(c)=c$. Therefore, the equality
$\phi\psi=\id_{\Fix\Phi}=\id_{C(\bar{X})}$ holds, and one further has
\[
   (\alpha,\phi)\circ (\beta,\psi)=(\beta\alpha,\phi\psi)=(\id_{\bar{X}},\id_{C(\bar{X})})=\id_{(\bar{X},C(\Bar{X}))}.
\]
Thus, $(\alpha,\phi)\circ (\beta,\psi)$ is filtered homotopic
to $\id_{(\cC,C(\Bar{X}))}$ via the zero filtered homotopy. In the
opposite direction, we have
\begin{equation}
\label{eq:thm:forman-gradient-filtered-equivalence-1}
  (\beta,\psi) \circ (\alpha,\phi)=
  (\alpha\beta,\psi\phi)=
  (\id_{\cV|\cC},\Phi^\infty).
\end{equation}
Since $\cV_\star=\cC$, we see that $(\id_{\cV|\cC},\Phi^\infty)
\sim_e(\id_\cV,\Phi^\infty)$. Hence, it follows from
Proposition~\ref{prop:Phi-id-filt-homotopy}
and~\eqref{eq:thm:forman-gradient-filtered-equivalence-1} that
$(\alpha,\phi)\circ (\beta,\psi)$ is filtered chain homotopic
to the identity~$\id_{(\cV,C(X))}$, and also elementary
filtered chain homotopic --- and this completes the
proof of the theorem.
\qed

\medskip
After these preparations, we have finally reached the main result
of this section, which establishes the uniqueness of the connection
matrix for gradient combinatorial vector fields on regular
Lefschetz complexes.
\begin{thm}
\label{thm:cm-for-forman-gradient}
  Let~$\cV$ be a gradient combinatorial vector field on a regular
  Lefschetz complex~$X$. Then~$\cC$, that is, the collection of
  critical vectors of~$\cV$, is a Morse decomposition of~$\cV$.
  It has exactly one connection matrix which coincides with
  the $(\bar{X},\bar{X})$-matrix of~$\bdy^\kappa_{|\Fix\Phi}$,
  up to a graded similarity.
\end{thm}
\proof
It follows from Proposition~\ref{prop:crit-M-decomp} that~$\cC$
is a Morse decomposition of~$\cV$, and that the partition induced
on the complex~$X$ by this Morse decomposition is~$\cV$. Thus, by
Definition~\ref{defn:conn-matrix-Morese-decomp}, the Conley complex
of~$\cC$ is the Conley complex of~$(\cV,C(X),\bdy^\kappa)$. From
Corollary~\ref{cor:Conley-homotopy-equivalence} and
Theorem~\ref{thm:forman-gradient-filtered-equivalence} we obtain that
the Conley complex of $(\cV,C(X),\bdy^\kappa)$ is isomorphic in~$\PfCC$
to the Conley complex of $(\bar{X},C(\bar{X}),\bdy^{\bar{\kappa}})$,
with~$\bar{X}$ ordered by a natural partial order. Since, in view of
Proposition~\ref{prop:barX} and Definition~\ref{defn:filtration-Lefschetz}
the latter is the Conley complex of a natural filtration of a Lefschetz
complex, the conclusion now follows directly from
Theorem~\ref{thm:singleton-partition}(ii).
\qed

\medskip
In addition, the above result allows us to easily detect the
existence of connecting orbits in gradient combinatorial vector
fields.
\begin{thm}
\label{thm:co-for-forman-gradient}
  Let~$\cV$ be a gradient combinatorial vector field on a regular
  Lefschetz complex~$X$. Furthermore, suppose that for a pair of critical
  cells~$x,y \in X^c$ the associated connection matrix guaranteed by
  Theorem~\ref{thm:cm-for-forman-gradient} has a nonzero entry in the
  $(\bar{x},\bar{y})$-position. Then there exists a path with respect 
  to the multivalued map~$\Pi_{\cV}$ from~$y$ to~$x$, i.e., a
  connecting orbit.
\end{thm}
\proof
In view of Theorem~\ref{thm:cm-for-forman-gradient}, the
chain~$\partial\bar{y}$ has to contain the basis element~$\bar{x}$
in its representation with respect to the basis of~$\Fix\Phi$ 
guaranteed by Proposition~\ref{prop:barX}(i). Moreover, due
to the representation of~$\bar{x} = \Phi^\infty(x)$ given in
Proposition~\ref{prop:bar-x}, one has in fact the inclusion
$x \in |\partial\bar{y}| = |\partial\Phi^\infty(y)|$.
Using a simple inductive argument and
Proposition~\ref{prop:Phi-is-declining}(i) then immediately
establishes the existence of a path from~$y$ to~$x$. This completes
the proof of the theorem.
\qed
\begin{figure}
  \begin{center}
    \includegraphics[width=0.45\textwidth]{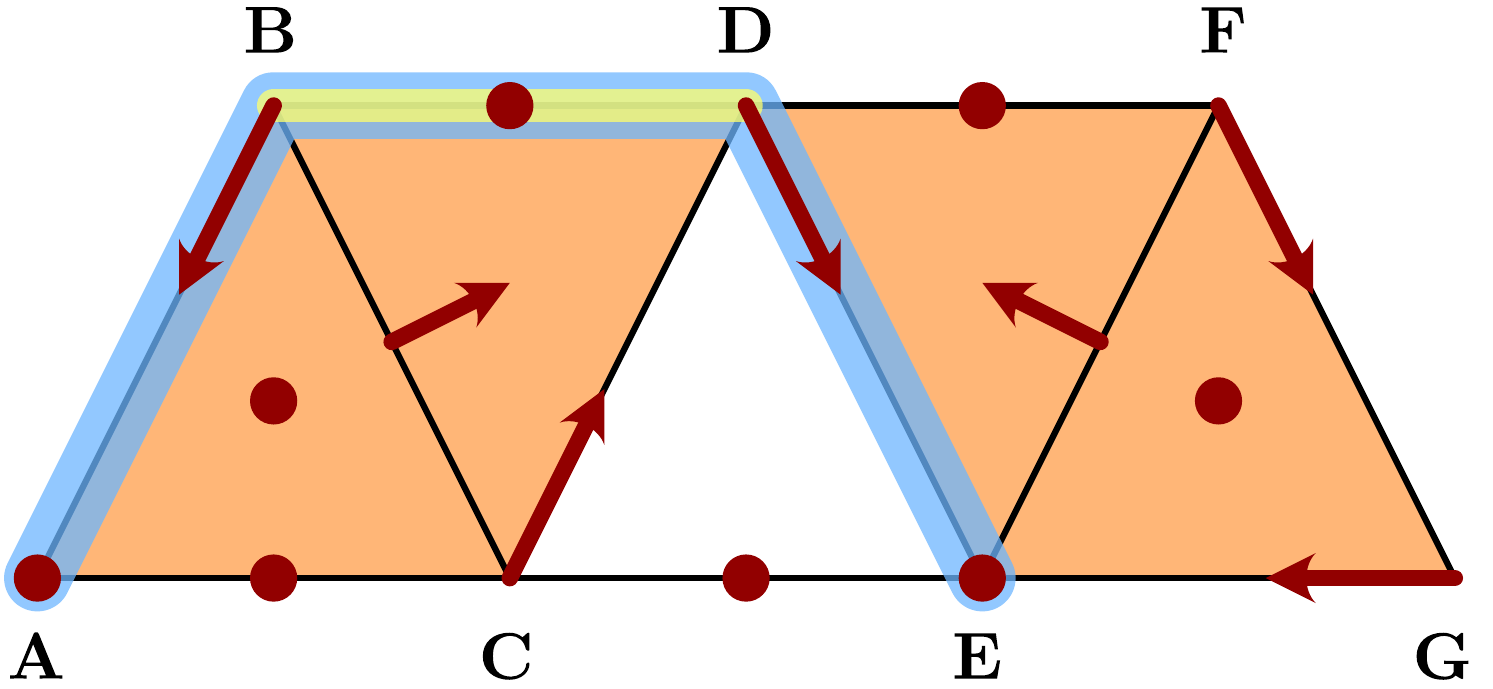}\quad\quad
    \includegraphics[width=0.45\textwidth]{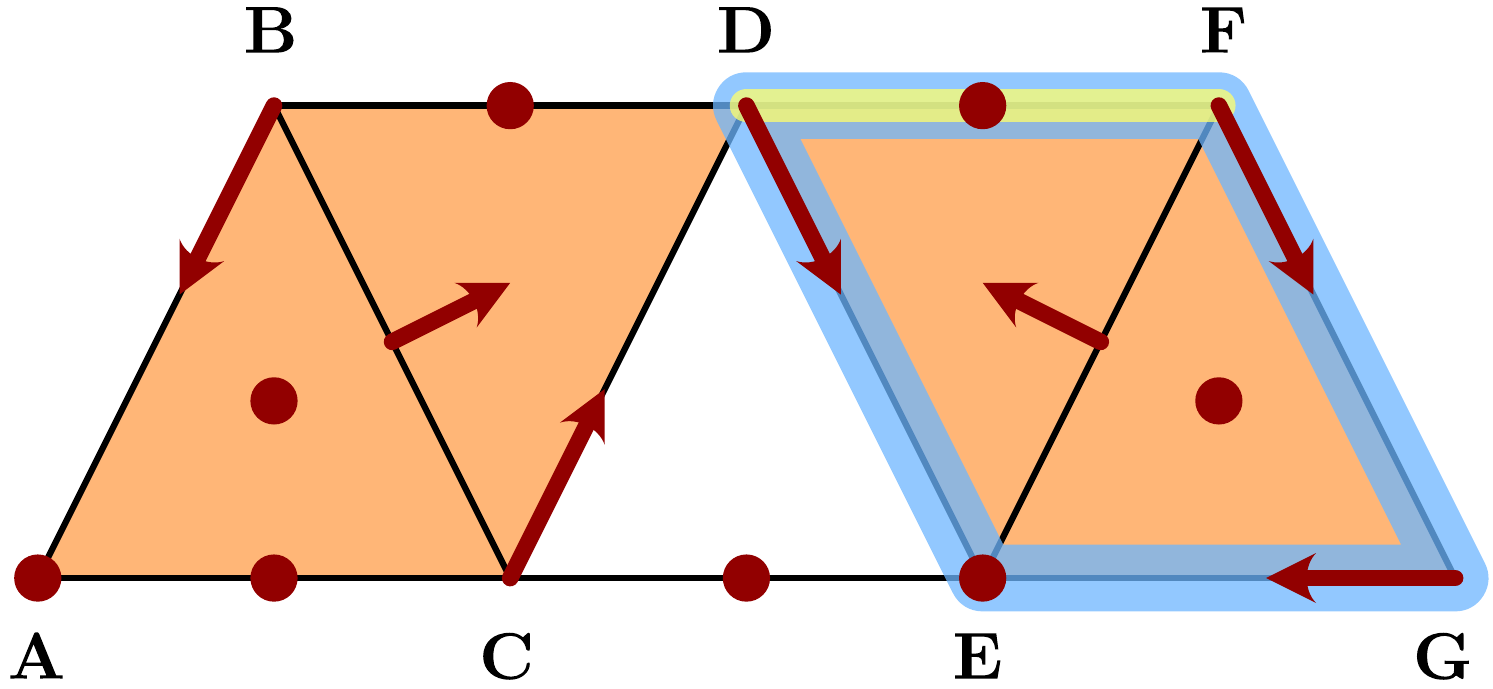} \\[3ex]
    \includegraphics[width=0.45\textwidth]{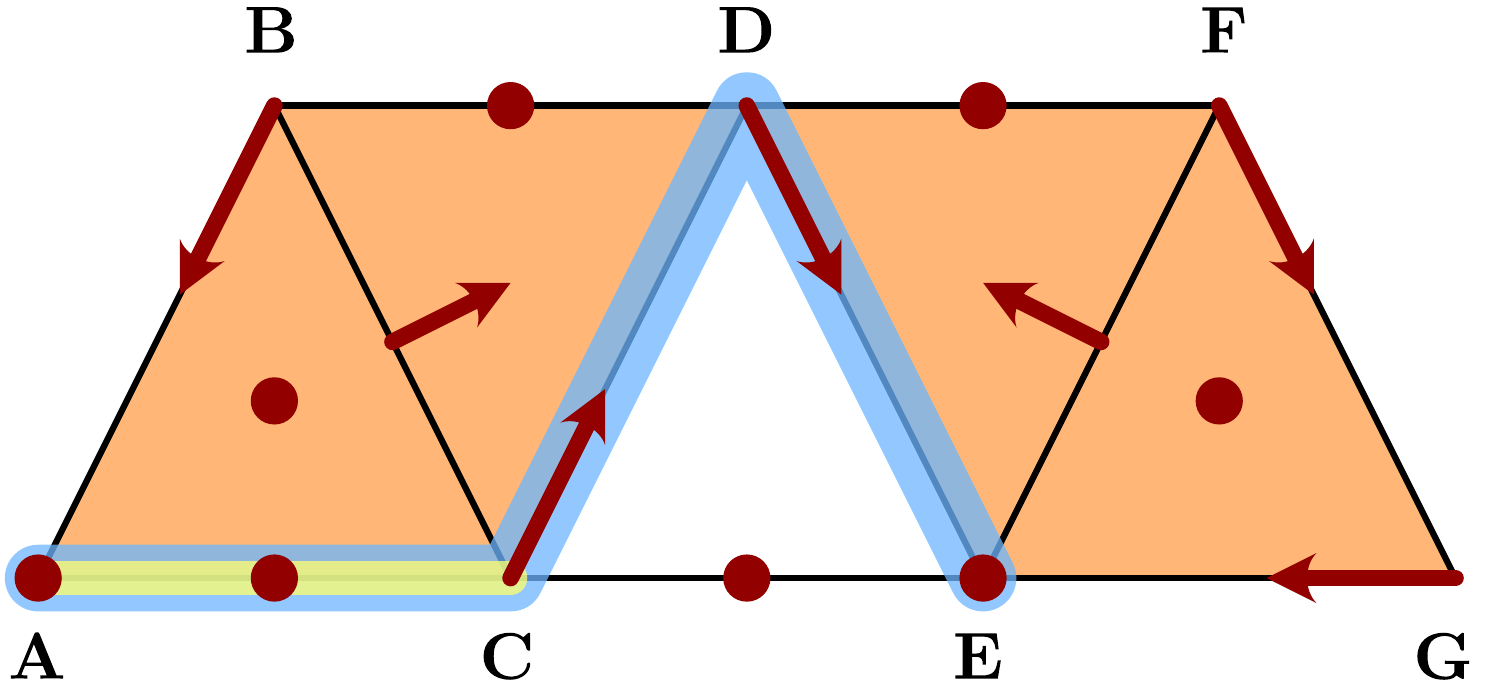}\quad\quad
    \includegraphics[width=0.45\textwidth]{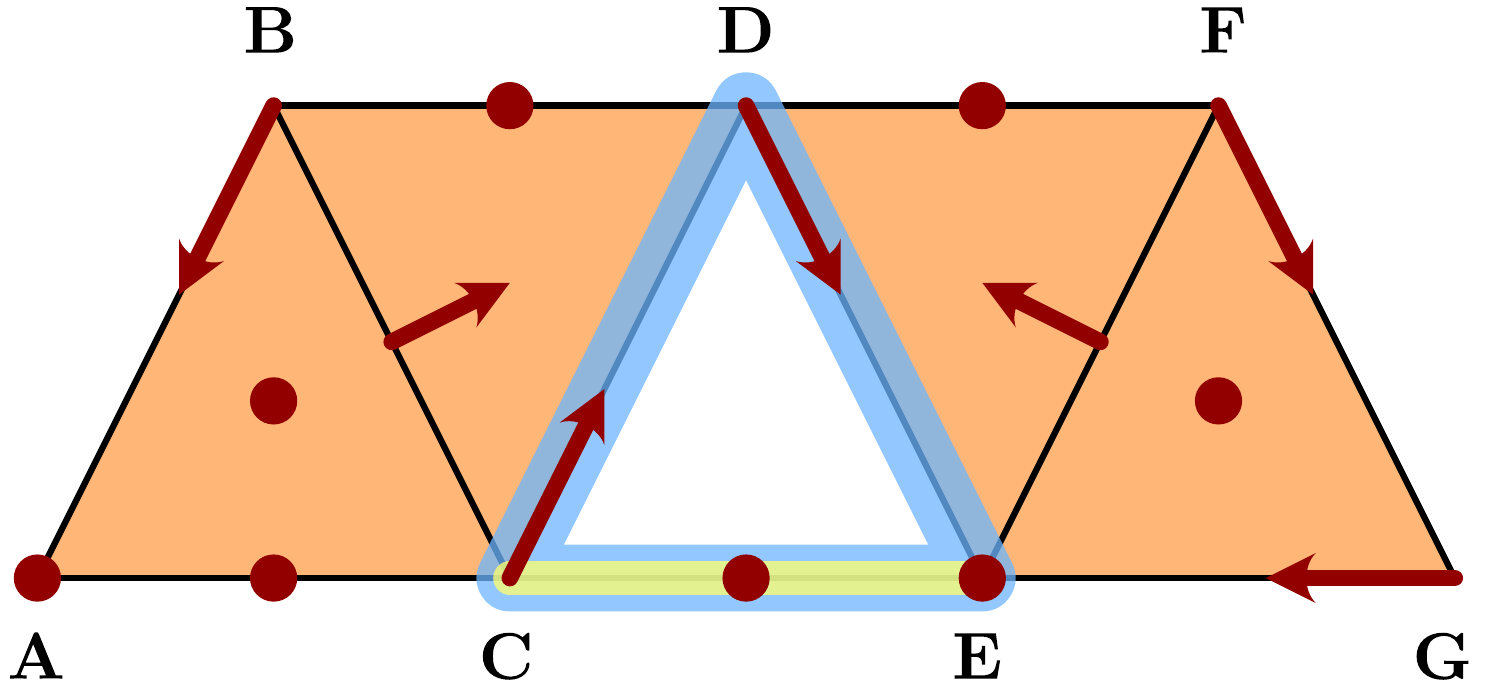}
  \end{center}
  \caption{{\em The stabilized combinatorial flow}.
           For the gradient combinatorial vector field~$\cV_3$
           introduced in the third row of Figure~\ref{fig:periodicex123},
           the above four panels show the images~$\Phi_3^\infty(x)$
           for the one-dimensional critical cells~$x$. From top left
           to bottom right they depict the chains~$\Phi_3^\infty(\bB\bD)$,
           $\Phi_3^\infty(\bD\bF)$, $\Phi_3^\infty(\bA\bC)$,
           and~$\Phi_3^\infty(\bC\bE)$, respectively, in blue.
           }
  \label{fig:phiinfex3}
\end{figure}

\medskip
We close this section with the following three examples. In the
first of these, it will be demonstrated how the connection matrices
in Table~\ref{table:periodicex123} can be determined using
Theorem~\ref{thm:cm-for-forman-gradient}.
\begin{ex}[{\em Three Forman gradient vector fields}]
\label{ex:formangradient-4}
{\em
Consider the gradient combinatorial Forman vector
field~$\cV_3$ which was introduced in the third row of
Figure~\ref{fig:periodicex123}. By following the arguments in
Example~\ref{ex:formangradient-3}, one can easily show that
$\Phi_3^\infty(x) = x$ for all zero-dimensional critical cells~$x$.
Moreover, we have both~$\Phi_3^\infty(\bA\bB\bC) = \bA\bB\bC + \bB\bC\bD$
and~$\Phi_3^\infty(\bE\bF\bG) = \bE\bF\bG + \bD\bE\bF$. For the four
one-dimensional critical cells~$y$, their images~$\Phi_3^\infty(y)$
are illustrated in Figure~\ref{fig:phiinfex3}.

The images under the stabilized flow~$\Phi_3^\infty$ form the
basis with respect to which one can determine the connection matrix.
As before, we abbreviate these by~$\bar{c} = \Phi_3^\infty(c)$ for
all critical cells~$c \in X^c$. Using this notation, for example
in order to determine the last column of the third connection
matrix in Table~\ref{table:periodicex123}, we have to express the
boundary~$\partial\overline{\bE\bF\bG}$ in terms of the chains~$\bar{y}$
for the one-dimensional critical cells~$y$, which leads to the identity
\begin{displaymath}
  \partial\overline{\bE\bF\bG} =
  \bD\bF + \bF\bG + \bE\bG + \bD\bE = 
  \overline{\bD\bF},
\end{displaymath}
and this accounts for the single entry of~$1$ in the corresponding
column given in Table~\ref{table:periodicex123}. Similarly, since
we have
\begin{displaymath}
  \partial\overline{\bA\bB\bC} =
  \bA\bB + \bB\bD + \bC\bD + \bA\bC = 
  \overline{\bA\bC} + \overline{\bB\bD},
\end{displaymath}
since the edge~$\bD\bE$ which is contained in both of the
last two chains cancels upon addition with respect to
$\ZZ_2$-coefficients. This leads to the two nonzero entries
in the second-to-last column of the last connection matrix
in Table~\ref{table:periodicex123}. The remaining columns
can be determined analogously.
\exend
}
\end{ex}
Next, we briefly sketch how the above example can be used to
determine connection matrices in the case of a combinatorial 
vector field which is not gradient. For this, we return to
the setting of Figure~\ref{fig:periodicex0}.
\begin{ex}[{\em A Forman vector field with periodic orbit}]
\label{ex:formanperiodic-3}
{\em
Consider again the combinatorial vector field~$\cV_0$ from
Figure~\ref{fig:periodicex0}, which has already been discussed in
Examples~\ref{ex:formanperiodic-1} and~\ref{ex:formanperiodic-2}.
This vector field has a periodic orbit which traverses the vertices~$\bC$,
$\bD$, and~$\bE$ --- and by breaking this periodic orbit into the union
of two critical cells of dimensions~$0$ and~$1$, as well as the associated
connecting orbits, we arrived at the three gradient systems~$\cV_k$
for $k = 1,2,3$ shown in Figure~\ref{fig:periodicex123}.

We now claim that each of the connection matrices from 
Table~\ref{table:periodicex123} are in fact connection matrices
for~$\cV_0$. Note that this vector field has a Morse decomposition~$\cM_0$
in which one Morse set, say~$M_p$, consists of the periodic orbit. Since 
the associated Conley index has the Conley polynomial~$t + 1$, we now
further subdivide this Morse set into a critical cell~$M_{\tilde{p}_1}$
of index~$1$, one critical cell~$M_{\tilde{p}_0}$ of index~$0$, as well
as two Forman vectors, as shown in Figure~\ref{fig:periodicex123}.
This results in one of the combinatorial vector fields~$\cV_k$,
which in turn leads to a refinement of the associated acyclic
partitions used in Definition~\ref{defn:conn-matrix-Morese-decomp}.
But then a direct application of Propositions~\ref{prop:refinement}
and~\ref{prop:aprefinement} shows that the Conley complex associated
with~$\cV_k$ is also a Conley complex for~$\cV_0$. (Note in particular 
that the above introduction of the critical cells~$M_{\tilde{p}_1}$
and~$M_{\tilde{p}_0}$ implies that the Conley complex has a trivial
boundary operator on the coarsened Morse set which includes the
periodic orbit in~$\cM_0$.) In other words, the connection matrices
shown in Table~\ref{table:periodicex123} are also connection matrices
for~$\cV_0$.

In fact, even more is true. Not only are all three connection
matrices from Table~\ref{table:periodicex123} connection matrices 
for~$\cV_0$, they are also pairwise nonequivalent in the sense of
Definition~\ref{def-equivalent-conley-complexes}, that is, the
vector field~$\cV_0$ does not have a unique connection matrix.
In the remainder of this example, we briefly sketch the necessary
proof of nonequivalence.

For this, we only consider the vector fields~$\cV_1$ and~$\cV_2$. In view
of Propositions~\ref{prop:phi-morphism} and~\ref{prop:psi-morphism}
the transfer morphism~$T : \cC' \to \cC''$ from the Conley
complex~$(\cM_0,\cC',d')$ associated with~$\cV_1$ to the Conley
complex~$(\cM_0,\cC'',d'')$ associated with~$\cV_2$ is defined
via the formula
\begin{displaymath}
  T(\bar{c}) =
  \Phi_2^\infty \left( \Phi_1^\infty(c) \right)
  \quad\mbox{ for }\quad
  \bar{c} = \Phi_1^\infty(c) \; .
\end{displaymath}
Furthermore, the matrix representations of the boundary
operators~$d'$ and~$d''$ can be found in the first two entries
of Table~\ref{table:periodicex123}.

As before we suppose that there exists a graded morphism~$f$ and
a degree~$+1$ homomorphism $\Gamma : C' \to C''$ with
\begin{displaymath}
  T - f = d'' \Gamma + \Gamma d' \; .
\end{displaymath}
We now evaluate both sides of this equation at the
chain~$\overline{\bD\bF}$. The first matrix in
Table~\ref{table:periodicex123} shows that
$d'(\overline{\bD\bF}) = 0$. Furthermore, quick glances at
Figures~\ref{fig:phiinfex1} and~\ref{fig:phiinfex2} imply that
$T(\overline{\bD\bF}) = \overline{\bD\bF} + \overline{\bD\bE}$.
Finally, since~$f$ is a graded morphism we have to have
$f(\overline{\bD\bF}) = \alpha \overline{\bD\bF}$ for some
$\alpha \in \ZZ_2$. Substituting these into the above identity
relating~$T$ and~$f$ one then obtains
\begin{displaymath}
  \overline{\bD\bF} + \overline{\bD\bE} - \alpha \overline{\bD\bF} =
  d''(\Gamma(\overline{\bD\bF})) \; .
\end{displaymath}
This equation, however, is impossible for $\alpha = 1$ due to the
matrix form of~$d''$ given in Table~\ref{table:periodicex123}, since
$(0,0,1,0,0,0)^t$ (which corresponds to the chain $\overline{\bD\bE}$)
is not contained in the column space of the matrix of~$d''$. On the
other hand, if $\alpha = 0$, then the form of~$d''$ immediately implies
$\Gamma(\overline{\bD\bF}) = \overline{\bE\bF\bG}$, and since~$\Gamma$
is filtered this furnishes $\overline{\bE\bF\bG} \le \overline{\bD\bF}$, 
a contradiction. In other words, the connection matrices arising from
the vector fields~$\cV_1$ and~$\cV_2$ are not equivalent in the sense
of Definition~\ref{def-equivalent-conley-complexes}.

The nonequivalence proofs for~$\cV_1 \not\equiv \cV_3$, as well as~$\cV_2
\not\equiv \cV_3$, are similar and thus left to the reader. We would like
to point out, however, that the first of these two cases requires a
slightly more involved argument.
\exend
}
\end{ex}
\begin{figure}
  \includegraphics[width=0.3\textwidth]{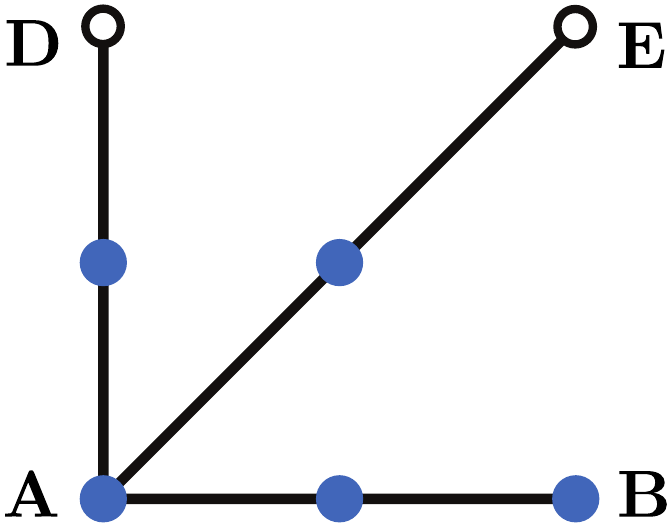}\hspace*{20mm}
  \includegraphics[width=0.3\textwidth]{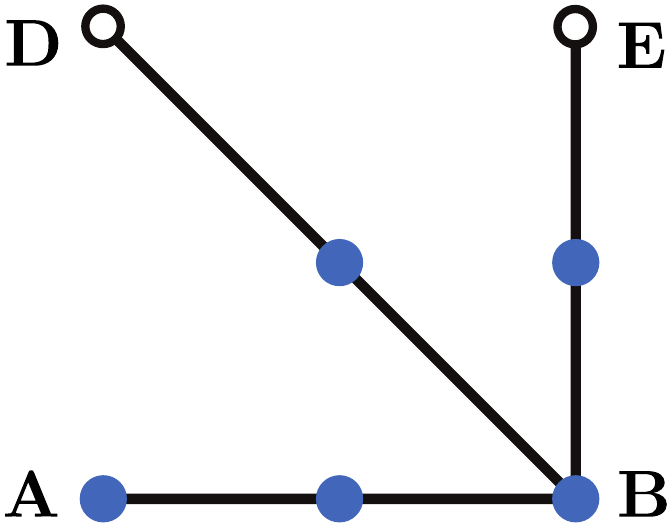}
  \caption{
  Conley complexes $(\bar{X},\bdop^{\bar{\kappa}})$ for the Morse
  decompositions~$\cM_1$ (left) and~$\cM_2$ (right) of the combinatorial
  gradient vector fields~$\cV_1$ and~$\cV_2$ in Figure~\ref{fig:t2eLefMVF},
  visualized as $\kappa$-subcomplexes of simplicial complexes.}
  \label{fig:3eLefA-B}
\end{figure}

As our last example, we return to the discussions begun in
Examples~\ref{ex:pfcc-conley1} and~\ref{ex:pfcc-conley2}, which 
center on the combinatorial flows in Figure~\ref{fig:t2eLefMVF}.
\begin{ex}
{\em
Consider the partition $\cE_1:=\cE_{\cM_1}$ associated with 
the Morse decomposition~$\cM_1$ of the vector field~$\cV_1$ in the
middle of Figure~\ref{fig:t2eLefMVF}. In this simple case it
coincides with the vector field~$\cV_1$. Hence, the Hasse diagram
of the partial order~$\leq_{\cE_1}$ coincides with
diagram~\eqref{ex:hasseV1}. It follows that~$\cV_1$ is a gradient
vector field and by Theorem~\ref{thm:cm-for-forman-gradient} the
Morse decomposition~$\cM_1$ has exactly one connection matrix. 
One can verify that the morphisms~$(\epsilon,h')$ and~$(\epsilon^{-1},g')$ 
constructed as in Example~\ref{ex:pfcc-conley1} are also filtered with
respect to the filtration~$(\cE_1,C(X),\bdop^\kappa)$. Therefore, arguing
as in Example~\ref{ex:pfcc-conley1} we conclude that the Conley complex
of~$\cM_1$ is the filtered chain complex in Example~\ref{ex:pfcc} with
the associated connection matrix \eqref{eq:d1m}. In fact, the
homomorphism~$h'$ coincides with the homomorphism~$\Phi_{\cV_1}$,
and we have $\Phi_{\cV_1}=\Phi_{\cV_1}^\infty$ in this case. An analogous
argument shows that the Conley complex of~$\cM_2$ is the filtered
chain complex in Example~\ref{ex:pfcc-conley2} with the associated
connection matrix~\eqref{eq:d2m}. Both Conley complexes are visualized
in Figure~\ref{fig:3eLefA-B} as $\kappa$-subcomplexes of suitable
simplicial complexes.
\exend}
\end{ex}
%


\section{Future work and open problems}
\label{sec:openproblems}

The main goal of the present paper was to carry over the classical theory
of connection matrices to the case of multivector fields on Lefschetz complexes.
Yet, in the course of this we achieved even more. By allowing for the change of
poset in the underlying Morse decomposition, we now have a mechanism at hand to
compare and ultimately classify connection matrices. In contrast to existing
results, this makes it possible to give a precise meaning to the phenomenon of
multiple connection matrices. As a welcome side effect, we could also shorten
the connection matrix pipeline described in the introduction. Combined with
the recent development of concrete algorithms for the computation of connection
matrices, this should considerably broaden their applicability in concrete examples.

Nevertheless, the results presented in this paper are just a first step.
There are a number of important open questions that have to be addressed
in the future, and which will be briefly outlined in the following.
\begin{itemize}
\item While our definitions give a precise definition that two connection
matrices are the same, verifying this equivalence is not immediately
straightforward. For this, it is necessary to recognize a morphism
as essentially graded. Are there easy ways to do this?
\item Related to the previous point, are there easy ways to see that
two connection matrices are really different, i.e., that one is in
the situation of a multivector field with multiple connection matrices?
Are there easily computable invariants that can be used to detect this?
\item What is the deeper meaning of the occurrence of multiple connection
matrices? In the classical case, this is always seen as an indication
of dynamical objects such as saddle-saddle connections. Is the same true
in the combinatorial multivector case?
\item In the classical case, transition matrices encode global bifurcations
that can be detected using a change in connection matrices. Can this notion
be extended to the combinatorial multivector case? In fact, we suppose that
this is strongly related to our notion of transfer morphisms. In what 
way precisely? In order to address this question, one also will have to
establish a precise notion of continuation of connection matrices.
\item Can every connection matrix be computed using the algorithms
mentioned earlier? We have indicated through examples that in certain
cases it is possible to break up recurrent Morse sets in such a way that
the resulting multivector field is actually gradient. Are the connection
matrices obtained in this way all possible connection matrices? If not,
is it possible to obtain the remaining ones algorithmically?
\item The earlier mentioned persistence based connection matrix algorithm 
indicates the possibility of the interpretation of connection matrix theory
in the language of persistence. Such a possibility might lead to applications
of connection matrices which go beyond dynamics.
\end{itemize}
Progress on any of these topics would only serve to extend the range
of potential applications of this theory.


\newpage

\bibliographystyle{abbrv}
\bibliography{references-ds,references-top,references-misc,references-thw}

\end{document}